%% file: main.tex
\documentclass[review,hidelinks,onefignum,onetabnum]{siamart250211}

\ifpdf
\hypersetup{
  pdftitle={Supervised and Unsupervised Neural Network Solver for First Order Hyperbolic Nonlinear PDEs},
  pdfauthor={Z. BABA, A. BAYEN, A. CANESSE, M.L. DELLE MONACHE, M. DRIEUX, Z. FU, N. LICHTLÉ, Z. LIU, H. NICK ZINAT MATIN and B. PICCOLI}
}
\fi

\headers{Supervised \& Unsupervised Neural Architecture for Hyperbolic PDEs}{BABA ET AL.}

\title{Supervised and Unsupervised Neural Network Solver for First Order Hyperbolic Nonlinear PDEs\thanks{Submitted to the editors on XXX; Authors are listed in alphabetical order
}}

\author{
Zakaria Baba$^{\dagger,1}$
\and Alexandre M. Bayen$^{2}$
\and Alexi Canesse$^{\dagger,3}$
\and Maria Laura Delle Monache$^{2}$
\and Martin Drieux$^{\dagger,1}$
\and Zhe Fu$^{\dagger,\S,2}$
\and Nathan Lichtl\'e$^{\dagger,2}$
\and Zihe Liu$^{2}$
\and Hossein Nick Zinat Matin$^{\dagger,2}$
\and Benedetto Piccoli$^{4}$
}

\DeclareMathOperator{\dt}{d\!}

\definecolor{darkmagenta}{rgb}{.38,0,0.38}
\usepackage{tikzscale}
\usepackage{amssymb}
\usepackage{adjustbox}
\usepackage{tabularx}

\usepackage{afterpage}
\usepackage{adjustbox} 
\usepackage{dsfont}
\usepackage{placeins}
\usepackage{cleveref}
\usepackage{multirow}
\usepackage{makecell}
\usepackage{stackengine}
\usepackage{tikz}
\usepackage{float}
\usepackage{pgfplots}
\usepackage{mathtools}
\usepackage{lscape}
\usepackage{calc}
\usetikzlibrary{positioning, matrix, shapes.geometric, arrows.meta, fit, backgrounds, calc, intersections, fillbetween}
\usepgfplotslibrary{fillbetween}
\pgfplotsset{compat=1.18}
\usetikzlibrary{external}
\usepackage{tabu}
\usepackage{makecell}
\usepackage{array}
\newcolumntype{C}{>{\centering\arraybackslash}X}
\newcommand{\sizestd}{6}
\newcommand{\sizestdinter}{7}
\newcommand{\std}[1]{\color{black!70!white}\fontsize{\sizestd}{\sizestdinter}\selectfont ±#1}

\usepackage{stmaryrd} 

\newcommand{\shortminus}{\mkern-2mu{\scriptscriptstyle-}\mkern-2mu}

\DeclareMathOperator{\mean}{mean}
\usepackage{nicefrac}
\usepackage{lipsum}
\usepackage{graphicx}
\usepackage{epstopdf}
\usepackage{caption}
\usepackage{subcaption}
\usepackage{mathtools}
\usepackage{relsize}
\usepackage{booktabs}
\usepackage{bm}

\usepackage{amsopn}

\newcommand{\set}[1]{\left\{ #1 \right \}}
\newcommand{\Def}{\overset{\textbf{def}}{=}}
\newcommand{\RR}{\mathbb{R}}
\newcommand{\R}{\mathbb{R}}
\newcommand{\ZZ}{\mathbb{Z}}
\newcommand{\NN}{\mathbb{N}}

\newcommand{\eps}{\varepsilon}

\newcommand{\rhomax}{\rho_{\text{max}}}

\newcommand{\abs}[1]{\lvert #1 \rvert}
\newcommand{\norm}[1]{\left\| #1 \right\|}

\usepackage{amsopn}

\DeclareMathOperator{\loc}{loc}

\newcommand{\nnorm}[1]{{\left\vert\kern-0.25ex\left\vert\kern-0.25ex\left\vert #1 
    \right\vert\kern-0.25ex\right\vert\kern-0.25ex\right\vert}}

\newcommand{\mC}{\mathcal C}
\newcommand{\mW}{\mathcal W}

\usepackage{duckuments} 

\usepackage{hyperref}
\usepackage{lipsum}
\usepackage{amsfonts}
\usepackage{amsmath}
\usepackage{amssymb}
\usepackage{graphicx}
\usepackage{epstopdf}
\usepackage{ifthen}
\usepackage{algorithmic}
\ifpdf
  \DeclareGraphicsExtensions{.eps,.pdf,.png,.jpg}
\else
  \DeclareGraphicsExtensions{.eps}
\fi

\newsiamremark{remark}{Remark}
\newsiamremark{hypothesis}{Hypothesis}
\crefname{hypothesis}{Hypothesis}{Hypotheses}
\newsiamthm{claim}{Claim}
\newsiamremark{fact}{Fact}
\crefname{fact}{Fact}{Facts}
\newsiamthm{case}{Case}

\begin{document}

\maketitle

\footnotetext[1]{\'{E}cole Polytechnique, Palaiseau, France, visiting at UC Berkeley, USA.}
\footnotetext[2]{University of California, Berkeley, USA.}
\footnotetext[3]{\'{E}cole Normale Sup\'{e}rieure de Lyon, France, visiting at UC Berkeley, USA.}
\footnotetext[4]{Rutgers University--Camden.}

\begingroup
\renewcommand\thefootnote{}%
\footnotetext{$\dagger$ Equal first author contribution.}%
\footnotetext{$\S$ Corresponding author (zhefu@berkeley.edu).}%
\endgroup

\begin{abstract}
We present a neural network-based method for learning scalar hyperbolic conservation laws. Our method replaces the traditional numerical flux in finite volume schemes with a trainable neural network while preserving the conservative structure of the scheme. The model can be trained both in a supervised setting with efficiently generated synthetic data or in an unsupervised manner, leveraging the weak formulation of the partial differential equation. We provide theoretical results that our model can perform arbitrarily well, and provide associated upper bounds on neural network size. Extensive experiments demonstrate that our method often outperforms efficient schemes such as Godunov's scheme, WENO, and Discontinuous Galerkin for comparable computational budgets. Finally, we demonstrate the effectiveness of our method on a traffic prediction task, leveraging field experimental highway data from the Berkeley DeepDrive drone dataset.
\end{abstract}



\section{Introduction and Related Work}

Conservation laws arise in a wide range of physical and engineering applications, including wave propagation, fluid dynamics, and traffic flow modeling~\cite{LWR}. The corresponding hyperbolic partial differential equations (PDEs) are often non-linear. A key characteristic of such systems is that they develop discontinuities in finite time, even when starting from smooth initial data. This behavior poses significant challenges for numerical methods, as schemes designed for smooth solutions often break down near discontinuities~\cite[Chapter~1]{leveque2002finite}.

In parallel, once discontinuities form, the PDEs can no longer be interpreted in the classical differentiable sense. Instead, solutions are sought in the weak sense, where the equations are satisfied in an integrated form against smooth test functions~\cite{holden2015front}. However, weak solutions are generally not unique, and additional admissibility criteria, such as entropy conditions, must be imposed to ensure physical relevance and restore uniqueness~\cite{holden2015front}.

While explicit solutions to Cauchy problems for conservation laws can be constructed in certain simple cases, such as Riemann problems, where the initial data consists of a single discontinuity~\cite{leveque2002finite}, or for piecewise constant initial conditions using the Lax–Hopf formula~\cite{ClaudelBayen2010}, numerical methods are generally required to approximate solutions in more realistic or complex settings.

Classical approaches to solving conservation laws rely on discretizing the spatial domain and updating solution values in a way that preserves conservation of the total quantity. Finite volume methods achieve this by evolving cell averages through numerical fluxes at cell interfaces. Notable examples include Godunov’s method (first-order, based on exact or approximate Riemann solvers) \cite[Chapter~4]{leveque2002finite} and high-order extensions such as ENO and WENO schemes~\cite{ENO, WENO}, which maintain non-oscillatory behavior near discontinuities. In parallel, conservative finite element variants such as the discontinuous Galerkin (DG)~\cite{first_paper_DG} method have been developed and demonstrated to be effective. These methods are well understood in terms of error estimates, consistency and stability. They offer interpretability and can achieve higher resolution with increasingly refined meshes~\cite{kuznetsov1976accuracy, shu1997essentially, harten1987uniformly,amat2025convergence}.

Several classes of methods leveraging the progress in deep learning have emerged as powerful ways to compute solutions to general PDEs over the recent years. Physics-Informed Neural Networks (PINNs) can directly approximate the solution of the PDE using a neural network. The model is trained using the residual of the PDE. Theoretical error estimates have been established for elliptic and parabolic PDEs~\cite{shin2020convergence,shin2020error,mishra2023generalization,mishra2022inverse,deryck2021kolmogorov,deryck2022navier,deryck2022generic}, however they typically rely on smoothness assumptions and fail for hyperbolic PDEs~\cite{deryck2022generic,mishra2023generalization}. In such cases, minimizing the residual of a smooth approximation of the PDE becomes difficult and may lead to convergence issues 
~\cite{deryck2022generic,mishra2023generalization}. 
The weak formulation provides an alternative which is more robust and adapted to a proper formulation of the problem: by testing the equation against smooth functions and applying integration by parts, the derivatives are transferred to the test functions, which allows a better treatment of non-smooth features for PINN-inspired networks~\cite{deryck2022wpinnsweakphysicsinformed,de2024wpinns}.

Another fundamental issue with PINNs is that they require training for each initial condition, which can take significant computational time. To overcome this, neural operators provide a general framework for learning mappings between infinite-dimensional function spaces, such as from initial conditions to solution spaces of PDEs~\cite{first_Deep0Net,first_FNO,nelsen2021random, patel2018operator, tripura2023wavelet, raonic2023iclr, raonic2023neurips}. Notable examples include Deep Operator Networks (DeepONets)~\cite{first_Deep0Net} and Fourier Neural Operators (FNOs)~\cite{first_FNO}. These models are typically trained in a supervised setting using synthetic data generated from known PDE solvers. Once trained, they allow for fast evaluations on new initial conditions, making them appealing for surrogate modeling tasks.

Universal approximation theorems have been established for mappings between spaces of scalar-valued functions~\cite{chen1995universal, lanthaler2022error, kovachki2021fourier, bhattacharya2021model, castro2023kolmogorov}. However, these results often assume smoothness, linearity or holomorphicity of the target operator \cite{herrmann2022neural, dehoop2023convergence, boulle2023learning}, which fails in the case of hyperbolic conservation laws.  The solution operator for such PDEs is nonlinear and typically non-smooth due to the formation of shocks. Numerical evidence further confirms that neural operator models tend to struggle when facing discontinuities~\cite{thodi2023learning, kim2025approximating}. In the presence of non-smooth dynamics, accurate approximation has been proven to require prohibitively large networks~\cite{lanthaler2023parametric}. While recent architectural modifications have improved approximation in specific cases~\cite{lanthaler2022nonlinear,wang2023operatorlearninghyperbolicpartial}, to the best of our knowledge, no theoretical results on general convergence for efficient Neural Operator architectures exist for hyperbolic problems (\cite{kovachki2021universalapproximationerrorbounds}~requires functions to lie in a Sobolev space, while~\cite{kovachki2021neural} guarantees existence of an approximator for a class more generic than those used in practice).

The flexibility of learning-based solvers for PDEs has also been increasingly used to recover the underlying dynamics of PDEs and potentially extrapolate it to other settings~\cite{raissi2017physics, li2021physics, lu2019deeponet}. 

In this work, we present a framework for learning scalar conservation laws by leveraging the structure of finite volume schemes and replacing the traditional numerical flux with a neural network. This flexible approach offers three key advantages: it eliminates the need to select an adapted numerical flux, a process that typically requires experimentation and domain expertise~\cite{NwaigweMungkasi2021}; it relaxes strong assumptions on the flux, such as monotonicity, thereby enabling improved performance; and it supports purely data-driven modeling. 

Recent works have explored similar directions, both in supervised and data-driven settings~\cite{morand2024deep,ChenEtAl2024}, as well as using Fourier Neural Operators (FNOs) to learn the flux from a spatial stencil~\cite{kim2025approximating}. 

Our approach builds on these efforts by generalizing the notion of the stencil to include temporal information, leading to significant improvements in model performance, as our goal is to provide a shelf-ready efficient solver.

Beyond this, we extend prior work by introducing an unsupervised training strategy tailored to the discontinuous nature of the solutions. We benchmark our model against established state-of-the-art numerical schemes, achieving higher accuracy for a given computational cost. Additionally, we provide theoretical insights, including convergence guarantees and upper bounds on model size. To our knowledge, this is the first work to offer such guarantees for neural architectures grounded in numerical methods. Finally, we evaluate the generalizability of the supervised method by learning to predict experimental highway data, comparing it to traditional schemes as a baseline.

\subsection{Contributions}
This work makes several key contributions to the field of neural network-based PDE solvers:

\begin{enumerate}
    \item We provide a flexible architecture to learn numerical fluxes for FVM schemes, as our method yields a conservative solution without enforcing strong constraints on the numerical flux.
    \item We demonstrate the effectiveness of our method by benchmarking it against state-of-the-art methods for hyperbolic PDEs such as WENO and DG, with improved performance for a given computation time.
    \item We introduce a learning framework based on the weak formulation of the PDE to learn a flux without generating ground truth data. We find experimentally that the structure of the numerical scheme does not require additional constraints to obtain the entropy solution.
    \item We provide the first theoretical guarantees, to our knowledge, for models based on numerical schemes. We prove that the method can achieve arbitrarily good performance with respect to given metrics, and provide upper bounds on the size of the required network, giving insight into the complexity of the learning task.
    \item We demonstrate the effectiveness of our framework on field experimental highway data using the Berkeley DeepDrive drone dataset~\cite{wu2022b3d}, effectively learning the system dynamics on limited data, and extrapolating them to accurately predict traffic evolution.
\end{enumerate}

The rest of the paper is organized as follows. In \Cref{sec:background}, we introduce necessary background on PDEs, conservation laws, and finite volume methods. In \Cref{sec:nfvm}, we introduce the proposed Neural FVM (NFVM) method, augmenting traditional FVMs using NNs. In \Cref{sec:theoretical_results}, we show theoretical results regarding the convergence and upper bounds on model complexity of the proposed approach. In \Cref{sec:results}, we conduct experiments comparing NFVM with classical schemes, both on synthetic data and experimental traffic data. Finally, \Cref{sec:conclusion} concludes this paper.

\section{Background}
\label{sec:background}

In this section, we first review the background on hyperbolic conservation laws and their solution spaces, followed by an overview of the finite volume methods used to solve them.

\subsection{PDEs and Conservation Laws}
\label{sec:pdes_conservation}

Let us consider a standard conservation law Cauchy problem of the form
\begin{equation}\label{E:main}
    \begin{cases}
       \partial_t u(t,x) + \partial_x\left[f(u(t,x))\right] = 0 &, (t,x) \in (0,T)\times\mathbb{R} \\
        u(0, x) = u_0(x) &, x \in \RR
    \end{cases}
\end{equation}
where \( u : (0,T) \times \mathbb{R} \to \mathbb{R} \) is the unknown density with \( T \in \R_{\ge 0}\). The initial condition is given by \( u_0\), and the flux function \( f \) is assumed to be concave or convex and sufficiently regular.

Let us recall a few definitions that will be used through out this paper.
\begin{definition}[Total Variation]\label{def:total_variation}
Let \( u \in L^1_{\loc}(\Omega) \), for an open set \( \Omega \subset \mathbb{R}^n \). The \emph{total variation} of \( u \) on \( \Omega \) is defined as
\begin{equation}
\abs{u}_{TV(\Omega)} \Def \sup \left\{ \int_\Omega u(x) \, \mathrm{div} \, \varphi(x) \, \dt x 
\ \middle| \ 
\varphi \in \mathcal C_c^1(\Omega; \mathbb{R}^n), \ \|\varphi\|_{L^\infty(\Omega)} \leq 1 
\right\},
\end{equation}
where $\mathrm{div}$ is the divergence operator.
Furthermore, a function \( u \in L^1(\Omega)\) is said to be of \emph{bounded variation}, denoted by \( u \in \mathrm{BV}(\Omega; \mathbb{R}) \) if $\norm{u}_{BV(\Omega)} < \infty$ where,
\begin{equation}\label{E:BV-norm}
    \norm{u}_{BV(\Omega)} \Def \norm{u}_{L^1(\Omega)} +\abs{u}_{TV(\Omega)}.
\end{equation}

\end{definition}


\begin{definition}[Weak solution]\label{def:weak_solution}
\text{A function}
\begin{equation}
u\in L^\infty((0,T); L^{\infty}(\mathbb{R})) \cap \mC([0, T]; L^1(\RR)) ,
\end{equation}
is a \emph{weak solution} of equation \eqref{E:main} if for all test functions $\varphi \in \Phi \overset{\mathrm{def}}{=} \mathcal C_c^1((-\infty,T) \times \mathbb{R}; \mathbb{R})$, the following integral identity holds~\cite{evans2022partial}:
\begin{equation}\label{E:weak_sol}
\int_0^T \!\! \int_\mathbb{R} \left( u \, \partial_t \varphi + f(u) \, \partial_x \varphi \right) + \int_\mathbb{R} u_0(\cdot) \varphi(0, \cdot) = 0.
\end{equation}

\end{definition}

While weak solutions for a given Cauchy problem may not be unique, uniqueness can be restored by selecting the physically admissible solution via entropy conditions. The notion of \emph{entropy solution}, introduced by Kružkov~\cite{kruvzkov1970first}, provides a framework ensuring both existence and uniqueness.

\subsection{Finite Volume Methods ($\text{FVM}$)} \label{S:FVM}
The solutions of the conservation law \eqref{E:main} are known to exhibit discontinuities which in turn leads to both theoretical and computational complexities \cite{leveque2002finite}. To this end, numerical schemes rely on different approximations in an attempt to improve the prediction accuracy. Classical Finite Difference (FD) methods employ approximation of the derivatives. However, such methods can fail near discontinuities where the differential equations are not well-defined \cite{leveque2002finite}. To overcome such issues, Finite Volume Methods (FVMs) instead use the following integral form:
\begin{equation}
    \label{E:intergral_form}
    \frac{\dt}{\dt t} \int_{a}^{b} u(t, x) \dt x = f(u(t, a)) - f(u(t, b)).
\end{equation}
In particular,  FVMs discretize space into a collection of control volumes (or grid cells), and the integral (or cell average) of \( u \) is approximated over each cell, rather than evaluating \( u \) pointwise at grid nodes. The approximate solution is updated locally in each cell, and the central challenge lies in defining a suitable \textit{numerical flux} at cell interfaces. This flux must accurately approximate the physical flux—which depends on the (generally unknown) solution at cell boundaries—using only the available cell average values.

To formalize this, let \( n \in \{0, \dots, \lfloor T/\Delta t \rfloor\} \) and define \( t^n = n\Delta t \), with the time step of the length \( \Delta t \). We consider a uniform spatial discretization with the cell length of \( \Delta x > 0 \), where cell centers are located at \( x_j = j \Delta x \), and cell boundaries at \( x_{j \pm 1/2} = (j \pm 1/2)\Delta x \). 

If \( u(\cdot, \cdot) \) is the solution to the conservation law, the average density in cell \( j \in \mathbb{Z} \) at timestep \( n \) is defined by
\begin{equation} \label{E:sol_cell_avg}
    u_j^{n} \Def \dfrac{1}{\Delta x} \int_{x_{j - 1/2}}^{x_{j + 1/2}} u(t^{n}, x) \, \mathrm{d}x.
\end{equation}

The conservative integral \eqref{E:intergral_form} can be written in the form of \begin{equation} \label{E:conservation_derivative}
    \frac{\dt}{\dt t} \int_{x_{j - 1/2}}^{x_{j + 1/2}} u(t,x) \dt x = f(u(t, x_{j - 1/2})) - f(u(t, x_{j + 1/2})).
\end{equation}
Integrating this equation from \(t^{n}\) to \(t^{n+1}\) and dividing by \(\Delta x\) yields
\begin{equation}\label{E:average_density_update}
    u_j^{n + 1} = u_j^{n} + \frac{\Delta t}{\Delta x} (\mathcal{F}_{j - 1/2}^n - \mathcal{F}_{j + 1/2}^n),
\end{equation} 
where $\mathcal{F}_{j - 1/2}^n$ is the average value of the flux at $x = x_{j-1/2}$ between time $t^{n}$ and time $t^{n + 1}$, defined as 
\begin{equation}\label{E:numerical_flux}
    \mathcal{F} _{j - 1/2}^n \Def \dfrac{1}{\Delta t} \int_{t^{n}}^{t^{n + 1}} f(u(t, x_{j - 1/2})) \dt t.
\end{equation}


    

Equation~\eqref{E:average_density_update} enables us to iteratively construct the numerical solution. However, $\mathcal{F} _{j - 1/2}^n $ cannot be directly computed as it requires knowledge of $u(\cdot, \cdot)$ beyond time $t^n$. Finite volume methods compute approximations of $\mathcal{F} _{j - 1/2}^n$ using only values at the current timestep:
\begin{equation} \label{E:numerical_flux_approx}
    F(u_{j - 1}^{n}, u_j^{n}) \approx \mathcal{F}_{j - 1/2}^n, 
\end{equation}
where $F$ is the numerical flux. A numerical approximation $\hat u_j^n$ of the solution \(u_j^n\) can then be obtained from \eqref{E:average_density_update} using  
\begin{equation} \label{E:FV-updates}
     \hat u_j^{n + 1} =  \hat u_j^{n} + \frac{\Delta t}{\Delta x} \left(F(\hat u_{j -1}^{n}, \hat u_j^{n}) - F(\hat u_{j}^{n}, \hat u_{j + 1}^{n})\right)  
\end{equation}
where the numerical flux $F$ can be calculated using various numerical schemes (see \Cref{app:numerical_schemes}). In the remainder of the manuscript, we denote $\hat u_j^n$ as simply $u_j^n$ when there is no ambiguity. 

\section{Neural Finite Volume Method (NFVM)}
\label{sec:nfvm}
In this section, we introduce our proposed method, NFVM, a NN-based numerical flux that aims to approximate the flux \eqref{E:numerical_flux}. We first define the NFVM architecture, which we generalize to an arbitrary stencil size $(a,b)$ and denote NFVM$_a^b$. We then introduce the two variants of the method: one based on supervised training (NFVM$_a^b$), when data is available and can be leveraged, another based on unsupervised training (UNFVM$_a^b$), when accurate solutions are unavailable or expensive to compute.
\subsection{NN-Flux}
\subsubsection{NN architecture}
Our approach consists of using a neural network to compute the numerical flux. 
We denote $\Theta$ the space of possible parameters, and define a standard fully-connected neural network \( F_{\theta} : [0,u_{\max}]^2 \to \mathbb{R} \), parametrized by weights \( \theta \in \Theta\), such that it approximates the numerical flux at the interface \( x_{j - 1/2} \):
\[ 
F_{\theta}(u_{j-1}^n, u_{j}^n) \approx \mathcal{F}^n_{j - 1/2}.
\]

\subsubsection{NFVM update}

An important property of many classical numerical fluxes is that they satisfy a discrete maximum principle: if the solution at time \( n \) satisfies \( u^n \in [0, u_{\max}]^{\mathbb{Z}} \), then the update preserves this bound, i.e., \( u^{n+1} \in [0, u_{\max}]^{\mathbb{Z}} \). However, this is not necessarily the case for our learned flux. To enforce stability at each iteration, we apply the correction 
\begin{equation} \label{E:NFV-updates}
    u_j^{n+1} = \max\left( \min\left( u_j^n + \frac{\Delta t}{\Delta x} \left( F(u_{j-1}^n, u_j^n) - F(u_j^n, u_{j+1}^n) \right), u_{\max} \right), 0 \right)
\end{equation}
after the update to ensure that the solution remains within the admissible range.
This in particular ensures that
\[
\forall n \in \mathbb{N},\ \forall j \in \mathbb{Z},\quad u_j^n \in [0, u_{\max}].
\]

\begin{remark} \label{non-conservativity}
The correction step is necessary in the theoretical analysis to ensure the solution remains within the admissible range.
To this end, the classical approach in enforcing the maximum principle while preserving the conservation laws for numerical schemes, is to consider proper slop limiters; see e.g. \cite{shu1997essentially, zhang2011maximum}. In this work, while the classic approaches could be applied, we consider the corrections in the form of \Cref{E:NFV-updates}. Nevertheless, the direct clipping violates the conservation of mass. However, in practice, omitting correction step does not significantly affect accuracy. Moreover, since the PDE solution stays within the range $[0, u_{\max}]$, any loss of conservation is confined to regions where the scheme fails, implying that the method remains conservative almost everywhere when accurate.
\end{remark}

\subsubsection{Stencil extension (NFVM$_a^b$)}

It is possible to leverage the flexibility of neural networks to compute the numerical flux using more than just the left and right values typically used by first order numerical schemes. We denote \( a \in 2\mathbb{N}^* \) the space stencil size, which we assume to be even for symmetry around the cell interface, and \( b \in \mathbb{N}^* \) the time stencil size. The \((a, b)\)-stencil is then defined as the following $a \times b$ rectangle of cells:

\[
\mathcal{U}^n_{j - 1/2}(a,b) = \left( u_i^m \right)^{m = n - b,\,\dots,\,n - 1}_{i = j  - \tfrac{a}{2} ,\,\dots,\,j - 1 + \tfrac{a}{2}}
\]

Then the neural network learns the mapping
\(
F_{\theta} : \mathbb{R}^{a\times b} \to \mathbb{R}
\)
such that
\[
F_{\theta}\left(\mathcal{U}^n_{j - 1/2}(a,b)\right) \approx \mathcal{F}_{j - 1/2}^n
\]

This significantly improves performance in practice, as it allows the network to encode non-local spatio-temporal information and alleviate error propagation.

\subsection{Supervised Loss}
\label{sec: sup loss}

\subsubsection{Notations}

Let us index the spatial and temporal discretization by $N\in \mathbb{N}^*$ : $\Delta x^N, \Delta t^N$.
The space of solutions of the finite volume method is the set of piecewise constant functions over space-time cells:
\[
S_N = \left\{ v : [0,T) \times \mathbb{R} \to \mathbb{R} \mid v \text{ is constant on each }  [t^n, t^{n+1}) \times [x_{j-\frac{1}{2}}, x_{j+\frac{1}{2}}) \right\}.
\]

Since we use the  update rule defined in \eqref{E:FV-updates}, it is convenient to define an operator \( H_N \) that maps any numerical flux \( F \) and initial condition to the corresponding discrete solution in \( S_N \).

Let \( u_\circ \in U_0 = (L^\infty \cap \mathrm{BV})(\mathbb{R}, [0, u_{\max}]) \) and define
\[
\mathcal F \Def  L^\infty([0, u_{\max}]^2; \RR) 
\]
Formally, we present the FV updates by: 
\begin{equation}\label{E:map_tildeh}
    \tilde h : \left\{ 
        \begin{array}{ccc}
            [0, u_{\max}]^{\mathbb{Z}} \times \mathcal{F} &\to& \RR^{\ZZ} \\
            (u, F) &\mapsto& \left( u_j 
- \dfrac{\Delta t^N}{\Delta x^N} \left[ F(u_j, u_{j+1}) - F(u_{j-1}, u_j) \right] \right)_{j \in \mathbb{Z}}.
        \end{array}
    \right.
\end{equation}
Next, define the clipping operator \(P\)
for bounded \( u \in \mathbb{R}^{\mathbb{Z}} \), by the Lipschitz map
\[
    P : \left\{ 
        \begin{array}{ccc}
            \mathbb{R}^{\mathbb{Z}} &\to& [0, u_{\max}]^{\mathbb{Z}} \\
            u &\mapsto& \left(\max\left( \min\left( u_j, u_{\max} \right), 0 \right)\right)_{j \in \mathbb{Z}}.
        \end{array}
    \right.
\]
Let us further define
\begin{equation}\label{E:h_function}
h : [0, u_{\max}]^{\mathbb{Z}} \times \mathcal{F} \to [0, u_{\max}]^{\mathbb{Z}}, \quad (u, F) \mapsto P\left( \tilde{h}(u, F) \right).
\end{equation}
This enables us to introduce the numerical solution operator \( H_N(F)\), mapping an initial condition \( u_0 \in U_0 \) to an approximate solution in \( S_N \):
\[
H_N(F) : U_0 \to S_N.
\]
defined by
\begin{equation}\label{E:H_N_function}
\forall n \in \left[0, \left\lfloor \nicefrac{T}{\Delta t^N} \right\rfloor \right] \quad (H_N(F)(u_0)(t^{n+1}, x_j))_{j \in \mathbb{Z}} = \left( h(H_N(F)(u_0)(t^n, x_j),F) \right)_{j \in \mathbb{Z}},
\end{equation}
with initial condition
\[
(H_N(F)(u_0)(0, x_j))_{j \in \mathbb{Z}} = \left( \frac{1}{\Delta x^N} \int_{x_j}^{x_{j+1}} u_0(x) \, \mathrm{d}x \right)_{j \in \mathbb{Z}}.
\]
In particular, for a given \( u_0 \), \( H_N(F) \) is constant on each cell.\\
We also define the entropy solution operator \( H^* : U_0 \to \mathbb{R}^{[0,T) \times \mathbb{R} } \) mapping each initial condition to the corresponding entropy solution introduced in  \cite{kruvzkov1970first}.  
\subsubsection{Supervised loss definition}\label{S:supervised_loss_defn}
In order to establish theoretical guarantees, we will evaluate the model on compact set $\mathcal{X} \subset \RR$. We define the supervised loss as
\begin{equation}\label{E:supervised_loss}
\mathcal{L}_{N, \sup}(F)
= \mathbb{E}_{u_0 \sim U_0}
\| H_N(F)(u_0) - H^*(u_0) \|_{L^1([0,T)\times \mathcal{X})} .
\end{equation}
where the expectation over \(u_0 \sim U_0\) refers to the average over finite subset of $U_0$. 
\subsection{Unsupervised loss definition}

We recall that \( \Phi = \mC_c^1(\mathcal U_T; \mathbb{R}) \), where $\mathcal U_T \Def (0, T) \times \RR$. We introduce for \( R > 0 \) the ball:
\begin{equation} \label{E:Phi_R}
\Phi_R = \left\{ \varphi \in \Phi \;\middle|\; \|\varphi\|_{\Phi} \le R \right\},
\end{equation}
with the norm \( \|\cdot\|_{\Phi} \) defined as
\begin{equation}
    \begin{split}
        \norm{\varphi}_{\Phi} \Def \norm{\varphi}_{W^{1,1}(\mathcal U_T)} + \norm{\varphi}_{W^{1, \infty}(\mathcal U_T)} + \norm{\varphi(0, \cdot)}_{L^1(\RR)}
    \end{split}
\end{equation}
where $\norm{\cdot}_{W^{k, p}}$ refers to the standard $(k, p)$-Sobolev norm. The exact formulation of \( \|\cdot\|_{\Phi} \) is tied to its use in the rest of the proof.

To quantify the violation of the weak formulation of the conservation law, we define the \emph{unsupervised loss}:
\begin{multline}\label{E:unsupervised_loss}
\mathcal{L}_{N,\mathrm{unsup}}^R(F) =\\
\ \underset{\substack{u_0 \sim U_0\\ \varphi \sim \Phi_R}}{\mathbb{E}}
\left[
\left(
\int_0^T \int_{\mathbb{R}} \left(H_N(F)(u_0)\, \partial_t \varphi
+ f(H_N(F)(u_0))\, \partial_x \varphi \right)
+ \int_{\mathbb{R}} u_0(\cdot)\, \varphi(0,\cdot) 
\right)^2
\right].
\end{multline}
Here, the expectation over \( u_0 \sim U_0 \) and \( \varphi \sim \Phi_R \) refer to discrete average sampling over finite subsets of $U_\circ$ and $\Phi_R$ which will be elaborated later in this section. 

\begin{remark}
    In the literature, several choices of test function families have been proposed for weak formulations of PINNs, including orthogonal polynomials \cite{kharazmi2019vpinn}.In our case, to preserve generality, we consider a distribution of test functions supported in a $R$-ball of the space \( \mC_c^1 \). The restriction to a ball can be justified by the linearity of the weak formulation~\eqref{E:weak_sol}.
\end{remark}

\section{Theoretical results}
\label{sec:theoretical_results}

The goal of this section is to prove that our NFVM$_2^1$ architecture can achieve arbitrarily small training error for a given $T > 0$. This is primarily achieved by controlling the error propagation within our numerical scheme, and we leverage this estimate to show that neural networks can match—and often outperform—traditional numerical methods, even with relatively small architectures. We also investigate the case where training is performed on a finite number of samples, and provide theoretical bounds on the training set size required to reach a desired level of accuracy. Our results cover both the supervised and unsupervised training settings.

\subsection{Setting}
\label{S:setting}
We make several technical assumptions in this section. For convenience, we assume $u_{\max} = 1$ (the results generalize via scaling) and that the $\|\cdot\|_{\mathrm{BV}}$ norm---see \Cref{E:BV-norm}---of initial conditions in $U_0 \subset L^\infty([0,1], \mathbb{R})$ is bounded by an arbitrary constant $I > 0$, i.e.,

\begin{equation}
\forall u_0 \in U_0, \quad \|u_0\|_{\mathrm{BV}}  \leq I.
\end{equation}

We choose mesh sizes such that $\Delta x^N \to 0$ as $N \to \infty$, $\Delta t^N = \alpha \Delta x^N$, and $\frac{1}{\alpha} \geq \sup_{u_{\min} \leq u \leq u_{\max}} |f'(u)|$, so that the discretization satisfies the CFL condition~\cite{leveque2002finite}.
In this work, we consider $f \in \mathcal{C}^1(\RR) \cap \mW^{k, \infty}([0, 1];\RR)$ for all $k \in \mathbb{N}$; it mainly simplifies the analysis and yields more tractable bounds on the required network size. For a broader treatment of how approximation rates depend on regularity, we refer to~\cite{Yarotsky16}.

Let \( F_{\mathrm{LxF}} \) denote the Lax-Friedrichs numerical flux~\cite{lax1955nonlinear}, defined by
\[
F_{\mathrm{LxF}}: (u^-, u^+) \mapsto \frac{1}{2}\left[f(u^-) + f(u^+)\right] - \frac{1}{2}(u^+ - u^-)\sup_{u \in [0, u_{\max}]} |f'(u)|.
\] 
This flux is consistent, monotone, under which the FVM is $L^1$-convergent to the entropy solution as $N \to \infty$; see ~\cite{godunov1959finite, Sanders1983}. Since \( f \) is Lipschitz, so is \( F_{\mathrm{LxF}} \) on \([0, u_{\max}]^2\), and we let \( M_{\mathrm{LxF}} > 0 \) denote its Lipschitz constant with respect to the $\|\cdot\|_\infty$ norm.

We define the multilayer perceptron
\begin{equation}
    F_\theta: \left\{
        \array{ccc}
            [0, u_{\max}]^2 &\to& \mathbb{R} \\
            z &\mapsto& W^{(L)} z^{(L-1)} + b^{(L)} 
        \endarray
    \right.
\end{equation}
with 
\begin{equation}
    \begin{cases}
        z^{(0)} = z \in [0, u_{\max}]^2\\
        z^{(l)} = \sigma(W^{(l)} z^{(l-1)} + b^{(l)}), &l = 1, \dots, L-1
    \end{cases}
\end{equation}
with \( \sigma(x) = \max\{0, x\} \) is applied element-wise and \( \theta = \{W^{(l)}, b^{(l)}\}_{l=1}^L \).

The output of \( F_\theta \) is clipped to the interval \( [-D, D] \), where
\[
D \geq \max_{(u^-, u^+) \in [0, u_{\max}]^2} |F_{\mathrm{LxF}}(u^-, u^+)|,
\]
so that the neural flux remains uniformly bounded, i.e., \( \|F_\theta\|_\infty \leq D \) for all \( \theta \in \Theta \).

As shown in Appendix~\ref{approx:clip}, the clipped network still approximates \( F_{\mathrm{LxF}} \) effectively. The clipping step guarantees uniform boundedness, which is required for the analytical results of the following sections. However, in practice, the outputs of \( F_\theta \) tend to remain within bounds naturally, so explicit clipping is optional.
\subsection{Error propagation lemma}
In this part, we establish a regularity result for the mapping between numerical fluxes in $\mathcal{F}$ and solutions produced by the numerical scheme. 
\begin{lemma}\label{lemme:h}
For any two flux functions \(F\), \(G\) $\in \mathcal{F}$ and for any $u \in [0,u_{\max}]^{\mathbb{Z}}$,
\[ \|\,h(u,F)-h(u,G)\|_{\infty}
\;\le\;
\frac{2\,\Delta t^N}{\Delta x^N}\,\|F-G\|_{\infty}.
\]
In other words, \(h(u,\cdot)\) is Lipschitz with constant \(2\Delta t^N/\Delta x^N\).
\end{lemma}

\begin{proof}
Let $u \in \mathbb{R}^{\mathbb{Z}}$ and, for each $j \in \mathbb{Z}$, denote $\tilde{h}_j(u,F) \Def [\tilde{h}(u,F)]_j$, the $j$-th component of the update. By subtracting the definition of the scheme \eqref{E:map_tildeh} for $F$ and $G$, we get
\[
\tilde{h}_j(u,F) - \tilde{h}_j(u,G) = -\frac{\Delta t^N}{\Delta x^N} \left[ (F - G)(u_j, u_{j+1}) - (F - G)(u_{j-1}, u_j) \right].
\]
Using the definition of the infinity norm,  $|F(u_j, u_{j+1}) - G(u_j, u_{j+1})| \leq \|F - G\|_{\infty}$ and $|F(u_{j-1}, u_j) - G(u_{j-1}, u_j)| \leq \|F - G\|_{\infty}$.
Therefore, 
$$|\tilde{h}_j(u,F) - \tilde{h}_j(u,G)| \leq 2\,\frac{\Delta t^N}{\Delta x^N} \, \|F - G\|_{\infty},$$
and taking the supremum over all $j \in \mathbb{Z}$ yields
\[
\| \tilde{h}(u, F) - \tilde{h}(u, G) \|_{\infty} \leq 2 \frac{\Delta t^N}{\Delta x^N} \| F - G \|_{\infty}.
\]
Considering the definition of operator $P$, we conclude that
\[
\|h(u, F) - h(u, G) \|_{\infty} \leq \| \tilde{h}(u, F) - \tilde{h}(u, G) \|_{\infty} \leq 2 \frac{\Delta t^N}{\Delta x^N} \| F - G \|_{\infty},
\]
which completes the proof.
\end{proof}
\normalsize
\begin{lemma}\label{lemma:H_Nlip}
Let $F \in \mathcal{F}$. 
Assume that F is Lipschitz continuous with constant $M_{F}$, i.e.
\begin{equation*}
\forall\, (x_1,x_2), (y_1, y_2) \in [0, u_{\max}]^2,\quad
|F(x_1, x_2)-F(y_1,y_2)| \leq M_F \| (x_1, x_2)-(y_1, y_2) \|_{\infty}.
\end{equation*}
Then, $\forall N  \in \mathbb{N^*}, \quad \forall G \in \mathcal{F}, \quad \forall u_0 \in U_0,$
\begin{equation*}
\|H_N(F)(u_0) - H_N(G)(u_0)\|_{\infty} 
\leq \frac{\|F - G\|_{\infty}}{M_F} 
\left[ \left(1 + \frac{2\Delta t^N}{\Delta x^N} M_F \right)^{\left\lfloor \nicefrac{T}{\Delta t^N} \right\rfloor} - 1 \right].
\end{equation*}
\end{lemma}
\begin{proof}
Let  \( u_0 \in U_0 \).  
Let \( F \),  \( G \) $\in \mathcal{F}$ such that $F $ is Lipschitz with constant $M_F$.
We denote \( u^n(F) = \bigl(u^n_j(F)\bigr)_{j \in \mathbb{Z}} \), the solution vector at timestep \(n\) obtained with flux \(F\), so that \( u^{n+1}(F) = h\bigl(u^n(F), F\bigr) \) and initial condition \( u_0 \).
In particular, for $n\in \left\{0,\cdots , \left\lfloor \nicefrac{T}{\Delta t^N} \right\rfloor-1 \right\}$, we have that
\begin{align*}
\|u^{n+1}(G) - u^{n+1}(F)\|_{\infty} 
&= \|h(u^{n}(G), G) - h(u^{n}(F), F)\|_{\infty} \\
&\hspace{-0.8cm}= \|h(u^{n}(F), F) - h(u^{n}(G), F) + h(u^{n}(G), F) - h(u^{n}(G), G)\|_{\infty}.
\end{align*}
Using the triangle inequality, we get that
\begin{align} \label{eq:step1}
\|u^{n+1}(F) - u^{n+1}(G)\|_{\infty} 
&\leq \|h(u^{n}(G), F) - h(u^{n}(G), G)\|_{\infty} \notag \\
&\quad + \|h(u^{n}(F), F) - h(u^{n}(G), F)\|_{\infty} .
\end{align}

On the other hand, by Lemma~\ref{lemme:h}
\begin{equation} \label{eq:step2}
\|h(u^{n}(G), F) - h(u^{n}(G), G)\|_{\infty} 
\leq \frac{2\,\Delta t^N}{\Delta x^N} \|F - G\|_{\infty}.
\end{equation}
Moreover, for \( j \in \mathbb{Z} \), and arbitrary  \( a, b \in [0,u_{\max}]^{\mathbb{Z} } \),
\begin{align*}
[\tilde{h}(a, F) - \tilde{h}(b, F)]_j
&= (a_j - b_j) + \frac{\Delta t^N}{\Delta x^N} \big[ F(a_j, a_{j+1}) - F(b_j, b_{j+1}) \\
&\hphantom{= (a_j - b_j) + \frac{\Delta t^N}{\Delta x^N} [} - F(a_{j-1}, a_j) + F(b_{j-1}, b_j) \big].
\end{align*}

Therefore, by \eqref{E:h_function} and considering:
\begin{equation*}
\begin{split}
|[h(a, F) - h(b, F)]_j| 
&\leq |[\tilde{h}(a, F) - \tilde{h}(b, F)]_j| \\
&\leq |a_j - b_j| + \frac{\Delta t^N}{\Delta x^N} \big(
    |F(a_j, a_{j+1}) - F(b_j, b_{j+1})| \\
&\hphantom{\leq |a_j - b_j| + \frac{\Delta t^N}{\Delta x^N} \big(} + |F(a_{j-1}, a_j) - F(b_{j-1}, b_j)|
\big).
\end{split}
\end{equation*}
Additionally, using that $F$ is Lipschitz continuous and summing over  \( j \) we get
\begin{equation} \label{eq:step3}
\|h(a, F) - h(b, F)\|_{\infty} 
\leq \|a - b\|_{\infty} 
+ 2\frac{\Delta t^N}{\Delta x^N} M_F \|a - b\|_{\infty}.
\end{equation}
By applying \eqref{eq:step3} to $a=u^n(F), b=u^n(G)$, and applying \eqref{eq:step2} in \eqref{eq:step1}, we get

\begin{equation} \label{eq:recursive}
\|u^{n+1}(F) - u^{n+1}(G)\|_{\infty}
\le
\underbrace{\frac{2\Delta t^N}{\Delta x^N}\|F-G\|_{\infty}}_{\text{by~\eqref{eq:step2}}}
+
\underbrace{\left(1 + \frac{2\Delta t^N}{\Delta x^N} M_F\right)
\|u^{n}(F) - u^{n}(G)\|_{\infty}}_{\text{by~\eqref{eq:step3}}}.
\end{equation}
We recognize an arithmetic-geometric sequence in  \eqref{eq:recursive}, we rewrite
\begin{align*}
&\|u^{n+1}(F) - u^{n+1}(G)\|_{\infty} 
+ \frac{\|F - G\|_{\infty}}{M_F} \\
&\quad \quad\leq \left(1 + \frac{2\Delta t^N}{\Delta x^N} M_F \right) 
\left( \|u^n(F) - u^n(G)\|_{\infty} + \frac{\|F - G\|_{\infty}}{M_F} \right)
\end{align*}
By iterating and since the terms are positive, we have:  
\begin{equation} 
\|u^{n}(F) - u^{n}(G)\|_{\infty} + \frac{\|F - G\|_{\infty}}{M_F} 
\leq  \left(1 + \frac{2\Delta t^N}{\Delta x^N} M_F \right)^n \left(\|u^{0} - u^{0}\|_{\infty} + \frac{\|F - G\|_{\infty}}{M_F}\right).
\end{equation}
Yet $u^{0}(F) = (u^0_j)_{ j \in \mathbb{Z}} = u^{0}(G)$, thus,
\begin{equation} \label{eq:closed-form}
\|u^{n}(F) - u^{n}(G)\|_{\infty} 
\leq \frac{\|F - G\|_{\infty}}{M_F} 
\left[ \left(1 + \frac{2\Delta t^N}{\Delta x^N} M_F \right)^n - 1 \right].
\end{equation}
Let \( x \in \mathbb{R} \) and \( t \in [0,T) \). 
The solution \( H_N(F)(u_0) \) is piecewise constant,  
thus, for each \( (t, x) \in  [0,T) \times \mathbb{R} \), there exist unique indices \( j \) and \( n \) such that
\[
(t, x) \in [t^n, t^{n+1}) \times [x_{j - 1/2}, x_{j+1/2}),
\quad \text{then} \quad
H_N(F)(u_0)(t, x) = u^n(F)_j.
\]
Using \eqref{eq:closed-form}, we have that 
\begin{equation*}
\begin{split}
|H_N(F)(u_0)(t, x) - H_N(G)(u_0)(t, x)|
&= |u^n(F)_j - u^n(G)_j|  \\ 
&\leq \frac{\|F - G\|_{\infty}}{M_F} 
\left[ \left(1 + \frac{2\Delta t^N}{\Delta x^N} M_F \right)^n - 1 \right]  \\
&\leq \frac{\|F - G\|_{\infty}}{M_F} 
\left[ \left(1 + \frac{2\Delta t^N}{\Delta x^N} M_F \right)^{\left\lfloor \frac{T}{\Delta t^N}\right\rfloor} - 1 \right],
\label{eq:HN-bound}
\end{split}
\end{equation*}
and hence the claim follows. 
\end{proof}
\begin{remark}[\textbf{Stability With Respect to Numerical Flux}]
    Lemma \ref{lemma:H_Nlip} provides the stability bound of the solution with respect to the numerical fluxes. More precisely, for a fixed $N \in \NN$ and over a finite time horizon $[0, T]$, $T >0$, for sufficiently close $F, G \in \mathcal F$, the associated solutions $H_N(F)$ and $H_N(G)$ remain close. 
\end{remark}
\begin{definition}[\textbf{$\eps$-{Proximity}}]
Let $N \in \NN$ and $\eps >0$ be fixed scalars. In addition, let $F \in \mathcal F$ with associate solution $H_N(F)$. We say the solution is weakly $\eps$-consistent with respect to $F$ over a finite time interval $[0, T]$, if
\begin{equation}
    \mathcal L_{N, \text{unsup}}^R(F) \le \eps, 
\end{equation}
for the unsupervised loss \eqref{E:unsupervised_loss} over this time interval, and strongly $\eps$-consistent with respect to $F$ over the time interval $[0, T]$ if
\begin{equation}
    \mathcal L_{N, \text{sup}}(F) \le \eps, 
\end{equation}
for the supervised loss \eqref{E:supervised_loss} over this time interval. 
\end{definition}
Our goal is to show that our proposed schemes create $\eps$-consistent solutions. 
\subsection{Unsupervised Learning}
\label{subsec: unsupervised}
This section is devoted to establishing that minimizing a loss with respect to the weak formulation yields numerical solutions with arbitrarily small residuals (in the same weak sense).
The first result shows that our proposed method performs at least as well as known schemes (which is considered to be the Lax-Friedrichs numerical flux in this work). \\
We denote the parameter that minimizes the unsupervised loss function defined in \eqref{E:unsupervised_loss} by
\begin{equation}\label{E:minimizer}
\theta^{*,R}_{N,\text{unsup}} = \operatorname*{arg\,min}_{\theta \in \Theta} \mathcal{L}_{N,\text{unsup}}^R(F_\theta).
\end{equation}
\begin{theorem}\label{T:unsup}
For any \(\epsilon>0\), \(R> 0 \), \(N \in \mathbb{N}^*\), and for a neural network of depth and width at least 
\[
\mathcal O\left(\log\left(\frac{1}{\epsilon} \right)+ \left\lfloor \frac{T}{\Delta t^N} \right\rfloor+\log(R)\right),
\]
where $R$ is defined in \eqref{E:Phi_R}, we have that 
\begin{equation} \label{E:eps_close_optimal_unsepervised}
\mathcal{L}^R_{N, \mathrm{unsup}}\big(F_{\theta^{*,R}_{N,\mathrm{unsup}}}\big) \leq \epsilon + \mathcal{L}^R_{N, \mathrm{unsup}}(F_{\mathrm{LxF}}). 
\end{equation} 
\end{theorem}

\begin{proof}
Let \( \epsilon > 0 \), \( R > 0 \), and \( N \in \mathbb{N} \). Under the assumption that $f \in \mathcal{C}^1 \cap W^{k,\infty}([0,1]; [0, f_{\max}])$ (see Section~4.1), the Lax-Friedrichs numerical flux function \( F_{\mathrm{LxF}} \) is Lipschitz continuous and sufficiently smooth to apply the approximation bounds from~\cite{Yarotsky16}.\\
By the \textit{universal approximation theorem}~\cite{hornik1991approximation}, for any \( \varepsilon_0 > 0 \), there exists a parameter configuration \( \theta_{\mathrm{LxF}} \in \Theta \) such that $\|F_{\mathrm{LxF}} - F_{\theta_{\mathrm{LxF}}}\|_\infty \leq \varepsilon_0$, where $F_{\theta_{\mathrm{LxF}}}$ denotes a neural network approximation of $F_{\mathrm{LxF}}$.
In addition, the width and depth of the neural network required to achieve this approximation can be chosen as
$\mathcal O\left( \log\left({\varepsilon_0^{-1}} \right) \right)$ for a fully connected ReLU network \cite{Yarotsky16}. Therefore, by \eqref{E:minimizer}, we have $\mathcal{L}_{N,\mathrm{unsup}}^R(F_{\theta^{*,R}_{N,\mathrm{unsup}}}) \leq \mathcal{L}_{N,\mathrm{unsup}}^R(F_{\theta_{\mathrm{LxF}}})$.

We now introduce the following notations,
\begin{itemize}
  \item \( A_N \Def \left({ \left(1 + {2\Delta t^N}{\Delta x^{-N})} M_{\text{LxF}} \right)^{\left\lfloor{T}{\Delta t^{-N}} \right\rfloor} - 1 }\right){M_{\text{LxF}}^{-1}} \)  from Lemma~\ref{lemma:H_Nlip};
  \item \( C \Def \max \set{u_{\max},\ \max_{u \in [0, u_{\max}]} |f(u)|} \);
  \item \( L \Def \max \set{1, W} \), where \( W \) is a Lipschitz constant of \( f \).
\end{itemize}
Using the definition \eqref{E:H_N_function} of the function $H_N$, we have that
\begin{equation}
\|H_N(F_{\mathrm{LxF}})(u_0)\|_{\infty} \leq C,
\quad
\|f(H_N(F_{\mathrm{LxF}})(u_0))\|_{\infty} \leq C.
\end{equation}
Let us consider $\varphi \in \Phi_R, u_0 \in U_0$.
Furthermore, considering \eqref{E:unsupervised_loss}, we define
\begin{equation} \label{E:l_N_function}
l_N(F,u_0,\varphi) =
\left(
\int_0^T \int_{\mathbb{R}} H_N(F)\, \partial_t \varphi
+
\int_0^T \int_{\mathbb{R}} f(H_N(F))\, \partial_x \varphi
+\int_{\mathbb{R}}u_0(x) \varphi(x,0) \mathrm{d}x
\right)^2.
\end{equation}
For simplicity, since the initial condition is fixed, we denote $H_N(F)(u_0) = H_N(F)$. Using \eqref{E:l_N_function}, we can write

\begin{equation}
\label{eq:product_form}
\begin{split}
&\left| l_N(F_{\theta_{\mathrm{LxF}}}, u_0, \varphi) 
- l_N(F_{\mathrm{LxF}}, u_0, \varphi) \right|  \\
&=
\Bigg|
\left(
\int_0^T \!\! \int_{\mathbb{R}} 
H_N(F_{\theta_{\mathrm{LxF}}})\, \partial_t \varphi +
\int_0^T \!\! \int_{\mathbb{R}} 
f(H_N(F_{\theta_{\mathrm{LxF}}}))\, \partial_x \varphi 
+ \int_{\mathbb{R}} u_0 \varphi(0,\cdot)
\right)^{2}  \\
&\qquad -
\left(
\int_0^T \!\! \int_{\mathbb{R}} 
H_N(F_{\mathrm{LxF}})\, \partial_t \varphi +
\int_0^T \!\! \int_{\mathbb{R}} 
f(H_N(F_{\mathrm{LxF}}))\, \partial_x \varphi 
+ \int_{\mathbb{R}} u_0 \varphi(0,\cdot)
\right)^{2}
\Bigg|  \\
&= 
\Bigg|
\left(
\int_0^T \!\! \int_{\mathbb{R}} 
\left[ H_N(F_{\theta_{\mathrm{LxF}}}) 
- H_N(F_{\mathrm{LxF}}) \right] \partial_t \varphi  +
\left[ f(H_N(F_{\theta_{\mathrm{LxF}}})) 
- f(H_N(F_{\mathrm{LxF}})) \right] \partial_x \varphi 
\right)  \\
&\qquad \times
\left(
\int_0^T \!\! \int_{\mathbb{R}} 
\left[ H_N(F_{\theta_{\mathrm{LxF}}}) 
+ H_N(F_{\mathrm{LxF}}) \right] \partial_t \varphi 
\right. \\
&\left. \qquad \qquad +
\int_0^T \!\! \int_{\mathbb{R}} 
\left[ f(H_N(F_{\theta_{\mathrm{LxF}}})) 
+ f(H_N(F_{\mathrm{LxF}})) \right] \partial_x \varphi 
+ 2 \int_{\mathbb{R}} u_0 \varphi(0,\cdot)
\right)
\Bigg| \\
& = \mathcal I_1 \times \mathcal I_2.
\end{split}
\end{equation}
Next, we can show the following bounds on $\mathcal I_1$ and $\mathcal I_2$,
\begin{equation}
\begin{aligned}
\mathcal I_1  
&\leq 
\int_0^T\!\!\int_{\mathbb{R}} 
\left| H_N(F_{\theta_{\mathrm{LxF}}}) - H_N(F_{\mathrm{LxF}}) \right|\, |\partial_t \varphi| 
+ 
\left| f(H_N(F_{\theta_{\mathrm{LxF}}})) - f(H_N(F_{\mathrm{LxF}})) \right|\, |\partial_x \varphi| \\
&\leq 
\| H_N(F_{\theta_{\mathrm{LxF}}}) - H_N(F_{\mathrm{LxF}}) \|_{\infty} 
\left( 
\int_0^T\!\!\int_{\mathbb{R}} |\partial_t \varphi| 
+ 
W \int_0^T\!\!\int_{\mathbb{R}} |\partial_x \varphi| 
\right) \\
&\leq 
L A_N 
\left( 
\int_0^T\!\!\int_{\mathbb{R}} |\partial_t \varphi| 
+ 
\int_0^T\!\!\int_{\mathbb{R}} |\partial_x \varphi| 
\right)
\|F_{\theta_{\mathrm{LxF}}} - F_{\mathrm{LxF}}\|_{\infty} \\
&\leq 
L A_N R\, \|F_{\theta_{\mathrm{LxF}}} - F_{\mathrm{LxF}}\|_{\infty},
\end{aligned}
\label{eq:I_1}
\end{equation}
\normalsize
where in the last inequality we are using the fact that $\varphi \in \Phi_R$. 
For the second term we have that:
\begin{multline*}
\mathcal I_2 \leq 
\int_0^T \!\! \int_{\mathbb{R}} 
\left| H_N(F_{\theta_{\mathrm{LxF}}}) 
+ H_N(F_{\mathrm{LxF}}) \right| 
|\partial_t \varphi|\\
+
\int_0^T \!\! \int_{\mathbb{R}} 
\left| f(H_N(F_{\theta_{\mathrm{LxF}}})) 
+ f(H_N(F_{\mathrm{LxF}})) \right| 
|\partial_x \varphi| + 
2 u_{\max} R
\end{multline*}
First, let's note that
\begin{equation} \label{ineq1}
\begin{aligned}
\int_0^T \!\! \int_{\mathbb{R}} 
\left| H_N(F_{\theta_{\mathrm{LxF}}}) 
+ H_N(F_{\mathrm{LxF}}) \right| |\partial_t \varphi|
\leq\;
&\int_0^T \!\! \int_{\mathbb{R}} 
\left| H_N(F_{\theta_{\mathrm{LxF}}}) 
- H_N(F_{\mathrm{LxF}}) \right| |\partial_t \varphi| \\
&\qquad \qquad + 
2 \int_0^T \!\! \int_{\mathbb{R}} 
\left| H_N(F_{\mathrm{LxF}}) \right| |\partial_t \varphi| \\
\leq\;
&\left( A_N L \| F_{\theta_{\mathrm{LxF}}} 
- F_{\mathrm{LxF}} \|_{\infty} + 2C \right)
\int_0^T \!\! \int_{\mathbb{R}} |\partial_t \varphi|.
\end{aligned}
\end{equation}
Similarly, we can rewrite and bound
\begin{align*}
&\int_0^T \!\! \int_{\mathbb{R}} 
\left| f(H_N(F_{\theta_{\mathrm{LxF}}})) 
+ f(H_N(F_{\mathrm{LxF}})) \right| 
|\partial_x \varphi| 
\\
&\hfill \leq\;
\int_0^T \!\! \int_{\mathbb{R}} 
\left| f(H_N(F_{\theta_{\mathrm{LxF}}})) 
- f(H_N(F_{\mathrm{LxF}})) \right| 
|\partial_x \varphi|
+ 2 \int_0^T \!\! \int_{\mathbb{R}} 
\left| f(H_N(F_{\mathrm{LxF}})) \right| 
|\partial_x \varphi| 
\\
& \hfill \leq\; \left( A_N L \| F_{\theta_{\mathrm{LxF}}} 
- F_{\mathrm{LxF}} \|_{\infty} + 2C \right)
\int_0^T \!\! \int_{\mathbb{R}} 
|\partial_x \varphi|.
\end{align*}
Therefore, putting all together, we have that
\begin{equation} \label{eq:I_2}
\begin{aligned}
\mathcal I_2 \leq\; &\left( 2C + A_N L \| F_{\theta_{\mathrm{LxF}}} - F_{\mathrm{LxF}} \|_{\infty} \right) \int_0^T \!\! \int_{\mathbb{R}} \left( |\partial_t \varphi| + |\partial_x \varphi| \right) + 2 u_{\max} R \\
\leq\; &\left( 2C + A_N L \| F_{\theta_{\mathrm{LxF}}} - F_{\mathrm{LxF}} \|_{\infty} \right) R + 2 u_{\max} R \\
=\; &\left( 4C + A_N L \| F_{\theta_{\mathrm{LxF}}} - F_{\mathrm{LxF}} \|_{\infty} \right) R, 
\end{aligned}
\end{equation}
where we are using the fact that $C \geq u_{\max}$. Then, using the decomposition given in \eqref{eq:product_form}, together with the bounds from \eqref{eq:I_1} and \eqref{eq:I_2}, we obtain
\begin{multline}
\left| l_N(F_{\theta_{\mathrm{LxF}}},u_0 , \varphi) 
- l_N(F_{\mathrm{LxF}}, u_0, \varphi) \right| \\
 \leq 
L A_N \left( 
    4C 
    + A_N L \| F_{\theta_{\mathrm{LxF}}} - F_{\mathrm{LxF}} \|_{\infty} 
\right) 
\left( R  \right)^2 
\cdot \| F_{\theta_{\mathrm{LxF}}} - F_{\mathrm{LxF}} \|_{\infty}.
\label{sansexpxt}
\end{multline}
We deduce by taking the expectation of \eqref{sansexpxt} over $\varphi \sim \Phi_R$ and $u_0 \sim U_0$
\begin{align*}
\mathcal{L}_{N,\mathrm{unsup}}^R(F_{\theta^{*,R}_{N,\mathrm{unsup}}}) 
&\leq \mathcal{L}_{N,\mathrm{unsup}}^R(F_{\theta_{\mathrm{LxF}}}) \\
&\leq \mathcal{L}_{N,\mathrm{unsup}}^R(F_{\mathrm{LxF}}) 
+ 4 L A_N C R^2 \cdot \| F_{\theta_{\mathrm{LxF}}} - F_{\mathrm{LxF}} \|_{\infty} \\
&\quad + (A_N L)^2 R^2 \cdot \| F_{\theta_{\mathrm{LxF}}} - F_{\mathrm{LxF}} \|_{\infty}^2.
\end{align*}

Moreover, by taking $N$ large enough, as $A_N \xrightarrow{N \to \infty} +\infty$, we can make the following assumption (this is not restrictive but gives a simpler upper bound)
\[
 4C+ A_N L \| F_{\theta_{\mathrm{LxF}}} - F_{\mathrm{LxF}} \|_{\infty} 
\leq 2 A_N L \| F_{\theta_{\mathrm{LxF}}} - F_{\mathrm{LxF}} \|_{\infty}.
\]
Therefore, we have
\begin{equation} \label{eq:unsup_bound_LxF}
\begin{aligned}
\mathcal{L}_{N,\mathrm{unsup}}^R(F_{\theta^{*,R}_{N,\mathrm{unsup}}}) 
    &\leq \mathcal{L}_{N,\mathrm{unsup}}^R(F_{\theta_{\mathrm{LxF}}}) \\
    &\leq \mathcal{L}_{N,\mathrm{unsup}}^R(F_{\mathrm{LxF}})
    + 2  \left(L  A_N  R  
     \| F_{\theta_{\mathrm{LxF}}} - F_{\mathrm{LxF}} \|_{\infty}\right)^2.
\end{aligned}
\end{equation}
This implies that it suffices to choose $\varepsilon_0 \leq \epsilon^{1/2}/(\sqrt{2} L A_N R)$ to ensure \eqref{E:eps_close_optimal_unsepervised} holds.

Consequently, the required depth and width of the neural network is given by 
\[
\mathcal O \left( \log(\varepsilon_0^{-1}) \right) = \mathcal O \left( \log\left(\frac{1}{\epsilon} \right)+ \left\lfloor \frac{T}{\Delta t^N} \right\rfloor +\log(R)\right).
\]
More importantly, this shows that the growth order of the network with respect to the $\epsilon$ is logarithmic. 
The $\left\lfloor \frac{T}{\Delta t^N} \right\rfloor$ term comes from the fact that errors propagate as the numerical scheme is iterated. Finally, the dependency in $R$ corresponds to the magnitude of the test functions, as $\mathcal{L}^R_{N, \mathrm{unsup}}$ increases with $R$.
\end{proof}

The next result follows from the convergence of the Lax–Friedrichs scheme and shows that the associated loss vanishes as the resolution increases.

\begin{lemma}{\label{F_L:unsup}}{}
For all \( \epsilon > 0 \) and \(R>0\), there exists \( N \in \mathbb{N} \) such that
\[
\mathcal{L}^R_{N, \mathrm{unsup}}(F_{\mathrm{LxF}}) \leq \epsilon.
\]
\end{lemma}
The proof can be found in \Cref{proof: Lemma unsupervised} in the Appendix.

\begin{corollary}\label{cor:unsup-univ}
For all \( \epsilon > 0 \) and \( R > 0 \), there exists \( N \in \mathbb{N} \) such that for a neural network of depth and width at most
\[
O \left( \log\left(\frac{1}{\epsilon} \right)+ \left\lfloor \frac{T}{\Delta t^N} \right\rfloor +\log(R)\right)
\]
we have \( \mathcal{L}_{N,\mathrm{unsup}}^R(F_{\theta^{*,R}_{N,\mathrm{unsup}}}) \leq \epsilon \).

\end{corollary}

\begin{proof}
Let us fix $\epsilon > 0$ and $R > 0$. We choose $N$ sufficiently large such that by Lemma~\ref{F_L:unsup} we have $\mathcal{L}_{N,\mathrm{unsup}}^R(F_{\mathrm{LxF}}) \leq \epsilon/2$ and hence by Theorem \ref{T:unsup} the claim follows.
\end{proof}

So far, we have shown that the solution of the neural network can get arbitrarily close to a weak solution under proper depth and width. In addition, it is crucial to understand the complexity of the samplings to achieve the desired accuracy. 
To do so, we start by defining the value of the loss for a given flux function and initial condition. More precisely, for any \( u_0 \in U_0 \) and a numerical flux $F$, we define
\begin{multline}\label{E:m_emp}
m_{N,\mathrm{unsup}}^R(F, u_0) 
\Def \\
\mathbb{E}_{\varphi \sim \Phi_R}
\Bigg[
\Bigg(
\int_0^T \!\! \int_{\mathbb{R}} 
H_N(F)(u_0)\, \partial_t \varphi 
+ f(H_N(F)(u_0))\, \partial_x \varphi
 + \int_{\mathbb{R}} u_0 \varphi(0,\cdot) 
\Bigg)^2
\Bigg].
\end{multline}
Therefore, the empirical loss and the associated parameter can be defined by
\begin{equation} \label{E:empirical}
\begin{cases}
\widehat{\mathcal{L}}_{N, \mathrm{unsup}}^{R, n_s}(F) \Def \frac{1}{n_s} \sum_{i=1}^{n_s} m_{N,\mathrm{unsup}}^R(F, u_0^i),\\
\widehat{\theta}^{R,n_s}_{N,\mathrm{unsup}} = \operatorname*{arg\,min}_{\theta \in \Theta} \widehat{\mathcal{L}}_{N, \mathrm{unsup}}^{R, n_s}(F_\theta)
\end{cases}
\end{equation} 
where $u_0^i \sim U_0$ are i.i.d. samples.

Employing \eqref{E:empirical}, the following results show that with high probability, arbitrary precisions can be achieved by training over empirical samples.
\begin{theorem}\label{conceration:unsup}
    Let us fix \(N\in \mathbb{N}\), \( \delta > 0 \), \( \epsilon > 0 \), and \( R > 0 \). Then, there exists an \( n_s \in \mathbb{N} \) such that, with probability $1 - \delta$,
\begin{equation*} 
\mathcal{L}_{N,\mathrm{unsup}}^R(F_{\theta^{*,R}_{N,\mathrm{unsup}}}) \le \mathcal{L}_{N,\mathrm{unsup}}^R(F_{\widehat{\theta}^{R,n_s}_{N,\mathrm{unsup}}})   + \epsilon. 
\end{equation*}  
\end{theorem}
\begin{proof}
We first note that
\begin{align}
\mathcal{L}_{N,\mathrm{unsup}}^R(F_{\theta^{*,R}_{N,\mathrm{unsup}}}) 
- \mathcal{L}_{N,\mathrm{unsup}}^R\big(F_{\widehat{\theta}^{R,n_s}_{N,\mathrm{unsup} }}\big) 
= \mathcal{L}_{N,\mathrm{unsup}}^R(F_{\theta^{*,R}_{N,\mathrm{unsup}}}) 
- \widehat{\mathcal{L}}_{N, \mathrm{unsup}}^{R, n_s}\big(F_{\widehat{\theta}^{R,n_s}_{N,\mathrm{unsup} }}\big) \nonumber \\ 
+ \widehat{\mathcal{L}}_{N, \mathrm{unsup}}^{R, n_s}\big(F_{\widehat{\theta}^{R,n_s}_{N,\mathrm{unsup} }}\big) 
- \mathcal{L}_{N,\mathrm{unsup}}^R\big(F_{\widehat{\theta}^{R,n_s}_{N,\mathrm{unsup} }}\big).
\label{decompo}
\end{align}
We should note that, for any $n_s$, 
\begin{equation} \label{E:expectation} 
\mathbb{E}_{u_0 \sim U_0}\left[ m_{N,\mathrm{unsup}}^R(F, u_0) \right] = \mathcal{L}_{N,\mathrm{unsup}}^R(F). 
\end{equation}

Furthermore, we can show that \( m_{N,\mathrm{unsup}}^R(F_{\theta}, u_0) \) is uniformly bounded for all \( u_0 \in U_0\) and \( \theta \in \Theta\). Indeed, for all \( u_0 \in U_0 \) and $N$ large enough, using the same calculations used to derive equation \eqref{eq:unsup_bound_LxF},
\begin{equation*}
\begin{aligned}
\left| m_{N,\mathrm{unsup}}^R(F_{\theta}, u_0) \right| 
- \left| m_{N,\mathrm{unsup}}^R(F_{\mathrm{LxF}}, u_0) \right| 
&\leq 
\left| m_{N,\mathrm{unsup}}^R(F_{\theta}, u_0) 
- m_{N,\mathrm{unsup}}^R(F_{\mathrm{LxF}}, u_0) \right| \\
&\leq 
\left( 2 L A_N R 
\| F_{\theta} - F_{\mathrm{LxF}} \|_{\infty} \right)^2
\end{aligned}
\end{equation*}
Using the definition of $D$, we have that
\begin{equation*}
\begin{aligned}
 |m_{N,\mathrm{unsup}}^R(F_{\theta}, u_0)|
&\leq  2 (L A_N R  (\| F_{\theta}\|_{\infty}+\|F_{\mathrm{LxF}} \|_{\infty}))^2 +|m_{N,\mathrm{unsup}}^R(F_{\mathrm{LxF}}, u_0)|\\
& \leq 2 (2 D L A_N R )^2 +|m_{N,\mathrm{unsup}}^R(F_{\mathrm{LxF}}, u_0)|.
\end{aligned}
\end{equation*}
Moreover, using the notation of Theorem \ref{T:unsup} and \eqref{E:m_emp}, we have \ \begin{equation*}|m_{N,\mathrm{unsup}}^R(F_{\mathrm{LxF}}, u_0)| \leq (C R)^2. 
\end{equation*}
This implies that \( m_{N,\mathrm{unsup}}^R(F_{\theta}, u_0) \) is uniformly bounded with respect to $\theta \in \Theta$ and $u_0 \in U_0$. Considering \eqref{E:expectation}, the Hoeffding's Lemma~\cite{vershynin2018high} for the random variables 
\begin{equation*} 
\begin{cases}
\big[m_{N,\mathrm{unsup}}^R(F_\theta, \cdot)\big]_i: U_0 \to \RR, \quad  i \in \set{1, \cdots, n_s} \\
u_0 \mapsto m_{N,\mathrm{unsup}}^R(F_\theta, u_0), 
\end{cases}
\end{equation*}
implies that for
\begin{equation*} n_s \geq \frac{B_{\mathrm{unsup}}^2}{2 (\epsilon/2)^2} \log\left( \frac{2}{\delta} \right), \qquad B_{\mathrm{unsup}} \Def 8 ( D L A_N R )^2 + (C R)^2
\end{equation*} 
with probability $1 - \delta$ and letting  $\theta=\widehat{\theta}^{R,n_s}_{N,\mathrm{unsup}}$, we have
\begin{equation} \label{Hoeffding}
\left| \widehat{\mathcal{L}}_{N, \mathrm{unsup}}^{R, n_s}(F_{\widehat{\theta}^{R,n_s}_{N,\mathrm{unsup}}}) - \mathcal{L}_{N,\mathrm{unsup}}^R(F_{\widehat{\theta}^{R,n_s}_{N,\mathrm{unsup}}}) \right| \leq \epsilon / 2.
\end{equation}
This implies that $\mathcal{L}_{N,\mathrm{unsup}}^R(F_{\widehat{\theta}^{R,n_s}_{N,\mathrm{unsup}}}) \leq \widehat{\mathcal{L}}_{N, \mathrm{unsup}}^{R, n_s}(F_{\widehat{\theta}^{R,n_s}_{N,\mathrm{unsup}}}) + \epsilon / 2$.
In addition, by definition of $\theta^*$, we have that
\begin{equation} \label{def_min}
\mathcal{L}_{N,\mathrm{unsup}}^R(F_{\theta^{*,R}_{N,\mathrm{unsup}}})\le \mathcal{L}_{N,\mathrm{unsup}}^R(F_{\widehat{\theta}^{R}_{N,\mathrm{unsup}}}) .
\end{equation}
Putting \eqref{decompo}, \eqref{Hoeffding}, and \eqref{def_min} together, the claim follows.
\end{proof}

\begin{corollary}{\label{cor:unsup}}
For all \( \epsilon > 0 \), \( R > 0 \), and \( \delta > 0 \),  
there exist \( N \in \mathbb{N} \), \( n_s \in \mathbb{N} \) such that,  
for a network of depth and width at most
\[
\mathcal O\left(
\log(R) + \log\left(\frac{1}{\epsilon} \right) +
\left\lfloor \frac{T}{\Delta t^N} \right\rfloor
\right),
\]
with probability \( 1 - \delta \), we have $\mathcal{L}_{N,\mathrm{unsup}}^R(F_{\widehat{\theta}^{R,n_s}_{N,\mathrm{unsup}}}) \leq \epsilon$.

\end{corollary}

\begin{proof}
    The proof follows the same arguments as Corollary~\ref{cor:unsup-univ}; see Corollary~\ref{proof:*unusp} in the Appendix for details
\end{proof}
\subsection{Supervised Learning}
In this section, we extend our analysis to the supervised setting, where the model is trained on known entropy solutions by minimizing the discrepancy between the neural network approximation and the ground truth. The results in this section are analogous to those of Section~\ref{subsec: unsupervised}. For consistency, we use the notations defined in Section~\ref{sec: sup loss}.

The optimal parameter of the supervised loss is defined by
$$
\theta^*_{N, \sup}
\Def \operatorname*{arg\,min}_{\theta \in \Theta}
\mathcal{L}_{N, \sup}(F_\theta).
$$
\begin{theorem}{\label{ineq:sup}}{}
For any fixed \( \epsilon > 0 \), \( N \in \NN \), and a neural network of depth and width at most
\[
\mathcal O \left(\log\left(\frac{1}{\epsilon} \right)+ \left\lfloor \frac{T}{\Delta t^N} \right\rfloor \right),
\]
we have that $\mathcal{L}_{N, \sup}(F_{\theta^*_{N, \sup}}) \leq \epsilon + \mathcal{L}_{N, \sup}(F_{\mathrm{LxF}})$.
\end{theorem}
\begin{proof}
From \textit{universal approximation theorem}, for any $\varepsilon_0 > 0$, there exists $\theta_{\mathrm{LxF}} \in \Theta$ such that $\|F_{\theta_{\mathrm{LxF}}} - F_{\mathrm{LxF}}\|_{\infty} \leq \varepsilon_0$.
For \( u_0 \in U_0 \), by definition of the minimizer and defining $\mathcal X_T \Def [0, T] \times \mathcal X$, for a compact set $\mathcal X \subset \RR$ (see Section \ref{S:supervised_loss_defn} for notations), we have 
\begin{equation} \label{sup loss without Exp}
\begin{aligned}
\|H_N(F_{\theta^*_{N, \sup}})(u_0) - H^*(u_0)\|_{L^1(\mathcal X_T)} 
&\leq \|H_N(F_{\theta_{\mathrm{LxF}}})(u_0) - H^*(u_0)\|_{L^1(\mathcal X_T)}  \\
&\leq  \|H_N(F_{\theta_{\mathrm{LxF}}})(u_0) \!- \!H_N(F_{\mathrm{LxF}})(u_0)\|_{L^1(\mathcal X_T)} \\
&\quad\quad + 
  \|H_N(F_{\small{\mathrm{LxF}}})(u_0) - H^*(u_0)\|_{L^1(\mathcal X_T)}.
\end{aligned}
\end{equation}
Moreover, we can write
\begin{equation*}
\begin{split}
\| H_N(F_{\theta_{\mathrm{LxF}}})(u_0) 
- H_N(F_{\mathrm{LxF}})(u_0) \|_{L^1(\mathcal X_T)} &=\int_0^T \!\! \int_{\mathcal{X}} 
\left| H_N(F_{\theta_{\mathrm{LxF}}})(u_0) 
- H_N(F_{\mathrm{LxF}})(u_0) \right| \\
&\leq T |\mathcal{X}| 
\left\| H_N(F_{\theta_{\mathrm{LxF}}})(u_0) 
\!-\! H_N(F_{\small{\mathrm{LxF}}})(u_0) \right\|_{\infty} \\
&\leq T |\mathcal{X}| A_N  \| F_{\theta_{\mathrm{LxF}}} - F_{\mathrm{LxF}} \|_{\infty},
\end{split}
\end{equation*}
where the last inequality is by Lemma \ref{lemma:H_Nlip}. Here, $\abs{\mathcal X}$ denotes the Lebesgue measure of set $\mathcal X$. Taking the expectation over \( u_0 \sim U_0 \) in \eqref{sup loss without Exp}, we obtain
\[
\mathcal{L}_{N, \sup}\big(F_{\theta^*_{N, \sup}}\big) 
\leq \mathcal{L}_{N, \sup}(F_{\mathrm{LxF}}) + T |\mathcal{X}| A_N \| F_{\theta_{\mathrm{LxF}}} - F_{\mathrm{LxF}} \|_{\infty}.
\]
Thus, it suffices to choose $\varepsilon_0 \leq \epsilon / (A_N T |\mathcal{X}|)$.
Accordingly, the required depth and width of the neural network is 
\[
\mathcal O \left( \log(\varepsilon_0^{-1}) \right) 
= \mathcal O\left( \log\left(\frac{1}{\epsilon} \right) + \left\lfloor \frac{T}{\Delta t^N} \right\rfloor \right).
\]
\end{proof}

The following result shows that the solution computed from the scheme using the neural network flux can be arbitrarily close to the entropy solution. 

\begin{theorem}
For any $\epsilon > 0$, there exists $N \in \mathbb{N}$ such that for a neural network of depth and width at most $\mathcal O\left( \log\left(1/\epsilon \right) + \left\lfloor T / \Delta t^N \right\rfloor \right)$,
\begin{equation*} \mathcal{L}_{N, \sup}(F_{\theta^*_{N, \mathrm{sup}}}) \leq \epsilon. 
\end{equation*}
\end{theorem}
\begin{proof}
Using notations of Lemma \ref{F_L:unsup} and computations from Appendix \ref{bound:K}, there exists a constant $K_{f, \alpha}$, \cite{godounov1959difference}, such that for any $u_0 \in U_0$
\[\|H_N(F_{\text{LxF}})(u_0) - H^*(u_0)\|_{L^1([0,T)\times \mathbb{R})} \leq 2 K_{f,\alpha} (\Delta x^N)^{1/2} T^{3/2} I.\] 
Furthermore, by monotonicity of Lebesgue measure
$$\|H_N(F_{\text{LxF}})(u_0) - H^*(u_0)\|_{L^1([0,T)\times \mathcal{X})}\leq \|H_N(F_{\text{LxF}})(u_0) - H^*(u_0)\|_{L^1([0,T)\times \mathbb{R})}.$$
By taking the expectation over \( u_0 \sim U_0 \), we obtain 
\[
\mathcal{L}_{N, \sup}(F_{\text{LxF}})
\leq 2 K_{f,\alpha} (\Delta x^N)^{1/2} T^{3/2} I.
\]
Hence, by choosing $N$ such that $\Delta x^N \leq {\epsilon^2}/{(4 K_{f,\alpha} T^{3/2} I)^2} $, we get $\mathcal{L}_{N, \sup}(F_{\text{LxF}}) \leq \epsilon / 2$ and hence, considering Theorem \ref{ineq:sup}, the claim follows.
\end{proof}

Similar to the calculations of Theorem \ref{concentration:sup}, we show the $\eps$-approximability of the by considering a finite number of initial conditions. In particular, we define
\[
m_{N, \sup}(F_\theta, u_0) \Def \|H_N(F_{\theta})(u_0)-H^{*}(u_0)\|_{L^1([0,T)\times \mathcal{X})},
\]
and for i.i.d. \( u_0^i \sim U_0 \), the average loss function is then denoted by
\[
\widehat{\mathcal{L}}_{N,\sup}^{n_s}(F_\theta) \Def \frac{1}{n_s} \sum_{i=1}^{n_s} m_{N, \sup}(F_\theta, u_0^i).
\]
We further define the minimizer as
\[
\widehat{\theta}^{n_s}_{N, \sup} \Def \operatorname*{arg\,min}_{\theta \in \Theta} \widehat{\mathcal{L}}_{N,\sup}^{n_s}(F_\theta).
\]
The following results show that it suffices to train on a finite number of initial conditions to achieve arbitrary precision with high probability. 

\begin{theorem}\label{concentration:sup}
    For any \(N\in \mathbb{N}\), \( \delta > 0 \), and \( \epsilon > 0 \), there exists an \( n_s \in \mathbb{N} \) such that with probability \( 1 - \delta \), 
\begin{equation*} 
\mathcal{L}_N^{\sup}(F_{\widehat{\theta}^{n_s}_{N, \sup}}) \leq \mathcal{L}_N^{\sup}(F_{\theta^*_{N, \sup}}) + \epsilon .
\end{equation*}
\end{theorem}
\begin{proof}
See Theorem \ref{proof: concentration:sup} in the Appendix. 
\end{proof}

\begin{corollary}{}{}
For any \( \epsilon > 0 \) and \( \delta > 0 \),  
there exist \( N \in \mathbb{N} \) and \( n_s \in \mathbb{N} \) such that,  
for a network of depth and width at least
\[
\mathcal O\left(
 \log\left(\frac{1}{\epsilon} \right) +
\left\lfloor \frac{T}{\Delta t^N} \right\rfloor
\right),
\]
and with probability \( 1 - \delta \), we have
  \begin{equation*} 
  \mathcal{L}_{N, \sup}(F_{\widehat{\theta}^{n_s}_{N, \sup}}) \leq \epsilon.
  \end{equation*}
\end{corollary}
\begin{proof}
    The proof follows from that of Corollary \ref{cor:unsup}.
\end{proof}

\section{Numerical Experiments} \label{sec:results}

We evaluate the performance of our method, both supervised (NFVM$_a^b$) and unsupervised (UNFVM$_a^b$), as introduced in \Cref{sec:nfvm}, for various stencil sizes $(a,b)$, against classical numerical schemes. As benchmarks to our approach, we consider a variety of first-order FVMs (\emph{Godunov}~\cite{godunov1959finite}, \emph{Lax-Friedrichs}~\cite{lax1954initial}, \emph{Engquist-Osher}~\cite{engquistOsher1981}), higher-order FVMs (\emph{Essentially Non-Oscillatory (ENO)}~\cite{shu1999high} and \emph{Weighted Essentially Non-Oscillatory (WENO)}~\cite{shu1999high}), and \emph{Discontinuous Galerkin (DG)}~\cite{first_paper_DG}, a finite-element method. These are classical numerical methods that are widely used in the literature to compute approximate solutions to hyperbolic conservation laws.

We consider two different conservation equations that can both model traffic flow: the \emph{Lighthill-Whitham-Richards} (LWR)~\cite{lighthill1955kinematic, R56} model
\begin{equation}
    \partial_t \rho + \partial_x f(\rho) = 0, \label{eq:lwr}
\end{equation}
where \(\rho\), a function of space and time, denotes the traffic density; and the \emph{inviscid Burgers' equation}~\cite{cameron2011notes} (which we will refer to as Burgers')
\begin{equation}
    \partial_t u + u\partial_x u = 0, \label{eq:burgers}
\end{equation}
where 
$u$ is also a function of time and space, denoting the traffic density. These two equations are instances of the general conservation law~\eqref{E:main} and notations are defined in \Cref{sec:pdes_conservation}.

In the case of LWR, $f$ is the flux function, which is a concave function of density. In this work, six different models have been considered for the flow function $f$: Greenshields'~\cite{greenshields1935study}, Triangular~\cite{geroliminis2008existence} (both symmetrical and skewed, later referred to as ``Triangular Sym'' and ``Triangular Skw''), Trapezoidal~\cite{geroliminis2011properties}, Greenberg~\cite{greenberg1959analysis} and Underwood~\cite{Underwood1961}. These models behave \textit{very} differently, as shown in~\Cref{fig:lwr_heatmaps}, and can be considered as six different equations. Analytical formulations of these models are given in~\Cref{app:lwr_flows}. Burgers' equation can actually be formulated as a special case of LWR, but with the strictly convex flux \[
    f: u \mapsto \frac{1}{2} u^2.
\]

\subsection{Training details}

The following applies to all six LWR models as well as to Burgers'. Training and evaluation data are independent for each of those seven models.

\paragraph{Training set} The training set is composed of on the order of 2000 randomly-generated Riemann initial problems, each one parametrized by densities $(\rho_1, \rho_2)$ and discontinuity location $x_0$, for which the exact solution is known (see \Cref{app:lwr_riemann}, or~\cite{leveque2002finite} for more details). The quantity of training data is chosen to be large enough to ensure that the initial conditions $(\rho_1, \rho_2)$ of the Riemann problems are sufficiently dense in the parameter space $[0, \rho_{\max}]^2$. In addition, for certain LWR equations, additional training samples are generated near the critical density, where the flow exhibits sharp changes (see \Cref{app:lwr_flows}). The data are discretized using $\Delta t \in [10^{-4}, 10^{-3}]$ (depending on equations) and $\Delta x = 10^{-3}$, with up to 100 spatial cells and 250 timesteps.  To promote generalization, a small Gaussian perturbation is added to the discontinuity location $x_0$, ensuring variation within the discretized grid. The overall objective is to train the model to solve these fundamental Riemann problem building blocks, and have it generalize to be able to solve more complex initial conditions.

\paragraph{Evaluation set} The evaluation set is composed of 1000 randomly generated complex piecewise constant initial conditions with 10 pieces each (see~\Cref{fig:lwr_heatmaps}) on a finer grid. Since each pair of adjacent pieces in the initial condition consists of a Riemann problem, one can build up the exact solution of any piecewise constant initial problem with a finite number of pieces (or infinite and of measure zero); this is known as the Lax-Hopf algorithm~\cite{SimoniClaudel2017}. The discretization for the evaluation set is $\Delta x = 5\mathrm{e}^{-3}$ and $\Delta t = 5\mathrm{e}^{-4}$, with 200 spatial cells and 1000 timesteps.

\FloatBarrier 

\paragraph{Boundary conditions} In this work, we provide the model with the true left and right boundary conditions, computed from the exact solution of Riemann problems during training, or using the Lax-Hopf scheme during evaluation. During training, the domain is chosen to be large enough so that waves do not reach the boundaries, meaning boundary conditions do not come into play and constant ghost cell \cite[Chapter~4]{leveque2002finite} values suffice for forward prediction. At evaluation time, however, we simulate a larger domain with piecewise constant initial conditions and then crop the interior to capture waves that may enter the system from the boundaries. This ensures that the boundary data supplied to the model accurately reflects the true dynamics of the problem, improving its ability to generalize and handle realistic scenarios. Indeed, if a wave enters the domain from the boundaries after the initial condition, the model would not be able to predict it correctly without appropriate boundary information. Another possible approach could have been to use periodic boundary conditions; however, we chose not to pursue this option since it does not align well with the physical nature of the problem.

\paragraph{Autoregressive prediction} The model is trained to predict $n_T$ timesteps autoregressively. Since this process amplifies small errors, we increase $n_T$ from 10 to 250 over the course of training. Initially, the model is only trained to predict $n_T=10$ timesteps. We find that this accounts for most of the training, as the model already often outperforms the Godunov scheme on the evaluation set after this first phase. However we keep fine-tuning the model with longer horizons (50, then 100, then 250) to further improve long-term prediction accuracy. During evaluation, the prediction is also autoregressive: the model knows only the (local) initial condition and the two boundary conditions, and predicts everything in-between, using its previous output as its next input. 

\paragraph{Model architecture} The model is applied locally on each group of cells to estimate the corresponding numerical flux, as detailed in \Cref{sec:nfvm}. The model is implemented as a one-dimensional CNN with six convolution layers. For NFVM$_a^b$, the first layer uses a kernel size of $a$ and outputs $b$ channels, while the five remaining convolutional layers use kernel size 1 with 15 channels each, except for the final layer, which outputs a single channel. This design is equivalent to sliding a fully-connected NN along the cells, but leverages the efficiency of CNN vectorization for significantly faster inference. NFVM$_2^1$ takes a single input timestep, whereas NFVM$_4^5$ uses five input timesteps, each provided as a separate channel. The total number of trainable parameters is approximately 1{,}200 for NFVM$_2^1$ and 1{,}800 for NFVM$_4^5$. Models used either ReLU or ELU activation functions, empirically chosen depending on the equations.

\paragraph{Training hyperparameters} The learning rate progressively decreases from $10^{-4}$ to $10^{-7}$ throughout the different stages of training introduced above ($10^{-4}$ while training on $n_T=10$, $10^{-5}$ while training on $n_T=50$, $10^{-6}$ while training on $n_T=100$, and $10^{-7}$ while training on $n_T=250$). Optimization is done using Adam with a batch size ranging between 256 and 2000 depending on equations. Even with proper weight initialization, model weights occasionally diverge from the start and models are unable to learn; in that case training is stopped early and restarted. 

Each training is performed on a single RTX5000 GPU and takes 5--10 minutes to reach 90\% of the final performance, and another 30--40 minutes of fine-tuning to reach 99\%, since the predictions on longer horizons take longer to compute.

\subsection{LWR Experiments}

\input{tex/lwr_heatmaps.tex}

For each of the six LWR models considered, we trained two first-order models: a supervised variant (NFVM$_2^1$) and an unsupervised variant (UNFVM$_2^1$). This section compares their performance with that of standard classical numerical schemes. \Cref{fig:lwr_heatmaps} presents the predictions of the NFVM$_2^1$ models on the evaluation set for each of the six LWR models, alongside the exact solutions obtained using the Lax-Hopf scheme. Note that a separate neural network was trained for each flow. Although the learned models exhibit some diffusion over time, which is a typical behavior of first-order methods, they appear to accurately capture the overall dynamics of the system.

\input{tex/lwr_metrics.tex}

\input{tex/lwr_l1_error_bar_plots.tex}
\vspace{-0.8cm}
\input{tex/lwr_l2_error_bar_plots.tex}

More precisely, \Cref{tab:lwr_metrics} compares the accuracy of all twelve NFVM$_2^1$ and UNFVM$_2^1$ models with that of the baseline numerical schemes across each of the six LWR models, highlighting their ability to generalize from training on coarse Riemann data to evaluation on fine-mesh complex initial conditions. The same results are also visualized as error bar plots in \Cref{fig:lwr_l1_error_bar_plots} and \Cref{fig:lwr_l2_error_bar_plots}, for the L1 and L2 errors, respectively. Overall, both NFVM$_2^1$ and UNFVM$_2^1$ consistently outperform the other first-order schemes considered (Godunov, Lax-Friedrichs, and Engquist-Osher) across all LWR models. Remarkably, the neural networks even surpass the higher-order ENO and WENO schemes in $L_1$ and $L_2$ errors for certain equations, including the widely used Greenshields and Triangular models, while achieving lower standard deviations. As expected, DG yields the best results in terms of accuracy and stability, being the only finite-element method evaluated.

\input{tex/lwr_winrates.tex}

To demonstrate that our models perform robustly beyond average accuracy, we also report the proportion of test cases in which the models outperform the baselines, as shown in \Cref{fig:lwr_winrates}. The learned models surpass all other first-order methods in every case for most fluxes, and in more than 90\% of cases for the remaining fluxes. Although they never outperform DG, they do outperform ENO and WENO on certain cases while being outperformed on others, which is a promising result given that ENO and WENO are higher-order methods compared to NFVM$_2^1$ and UNFVM$_2^1$.

In \Cref{fig:lwr_convergence_plot}, we further analyze the convergence properties of the learned models by evaluating their accuracy on different mesh sizes. Both models consistently achieve lower errors while maintaining a convergence rate similar to Godunov's method, which is known to converge to the entropy-satisfying solution as the discretization is refined. This suggests that the neural network-based schemes converge to the correct entropy solution as well. This result is particularly noteworthy for UNFVM$_2^1$, which was optimized only with respect to the weak formulation, a setting that typically admits many non-entropic solutions.

\input{tex/lwr_convergence_plot.tex}

One of the key properties required of a numerical flux $F_{\theta}$ in most convergence proofs is monotonicity: it must be non-decreasing in its first argument and non-increasing in its second argument~\cite{Bertoluzza2009}. \Cref{fig:lwr_numerical_flux} illustrates that the Godunov numerical flux strictly satisfies this condition in an explicit, constructed manner while also verifying the consistency property, namely $F_{\theta}(u,u) = f(u)$ for all $u$. Our models produce fluxes with similar shapes and preserve the consistency property, but do not satisfy monotonicity. Although this implies that standard convergence proofs do not apply, we believe this relaxation of monotonicity helps explain the improved performance of our models.

\input{tex/lwr_numerical_flux.tex}

Furthermore, we also consider a higher-order model, NFVM$_4^5$, which takes five timesteps as input while introducing minimal implementation or computational overhead. This demonstrates the flexibility of the NFVM$_a^b$ framework: whereas traditional higher-order numerical schemes typically require significant modifications to handle multiple inputs and often cannot extend beyond a single timestep, our method can easily incorporate additional spatial and temporal stencils by simply adjusting the neural network input size.

As reported in \Cref{tab:lwr_2d_metrics}, NFVM$_4^5$ achieves performance that surpasses higher-order schemes like WENO by an order of magnitude in both mean and standard deviation, and it approaches the performance of DG. Importantly, NFVM$_4^5$ remains nearly as simple as NFVM$_2^1$ in terms of implementation and model size, while delivering significantly improved results. Indeed, thanks to its vectorized architecture, NFVM$_4^5$ achieves inference speeds comparable to NFVM$_2^1$. In contrast, the DG algorithm is far more complex and, despite considerable effort to vectorize our implementation, remains orders of magnitude slower than NFVM.

\input{tex/lwr_2d_metrics.tex}

Additional results and movies are available at \href{https://www.nathanlichtle.com/research/nfv}{nathanlichtle.com/research/nfv}. This page presents the behavior of each numerical scheme on each flow model for various initial conditions, both as heatmaps and as videos showing the evolution of the solution over time. The final movies in particular highlight how closely NFVM$_4^5$ predictions track the true solution over time, compared to Godunov's predictions. 

\FloatBarrier

\subsection{Burgers' experiments}

\input{tex/burgers_heatmaps.tex}

As with the LWR experiments, we trained two first-order models on Burgers' equation: a supervised variant (NFVM$_2^1$) and an unsupervised variant (UNFVM$_2^1$). \Cref{fig:burgers_heatmaps} illustrates the exact solutions for Burgers' equation across different initial conditions. Overall, the performance of these models relative to baseline numerical schemes closely aligns with the patterns observed with the LWR models. As reported in \Cref{tab:burgers_metrics}, NFVM$_2^1$ achieves the lowest errors among first-order methods. Remarkably, despite being first-order methods, both NFVM$_2^1$ and UNFVM$_2^1$ also outperform higher-order schemes such as ENO and WENO in terms of $L_2$ and relative errors.

Furthermore, the win rates displayed in \Cref{fig:burgers_winrates} confirm that NFVM$_2^1$ and UNFVM$_2^1$ outperform Godunov's method, achieving lower $L_2$ errors than Godunov on 100\% and 93\% of the evaluation test cases, respectively. The consistent superiority across different conservation laws further highlights the robustness and versatility of the proposed approach.

\input{tex/burgers_metrics.tex}

\input{tex/burgers_winrates.tex}

We also consider a higher-order model, NFVM$_4^5$, which takes five timesteps as input while introducing only a few hundred additional trainable parameters compared to NFVM$_2^1$. \Cref{tab:burgers_2d_metrics} shows that NFVM\(_4^5\) again outperforms all other FVMs by an order of magnitude and get results close to DG. As was the case on LWR, this level of accuracy is higher precision than we thought possible from a finite volume scheme. 

\input{tex/burgers_2d_metrics.tex}

Furthermore, we compare the last timestep produced by the different schemes against the PDE solution in \Cref{fig:burgers_2d_final_density}. Since the methods are autoregressive, this last step reflects the maximum accumulated error, making it easier to see the discrepancy than on heatmap plots. The prediction of the learned models is significantly more accurate than Godunov's, particularly in regions with shocks. Moreover, unlike the trained models, Godunov's method exhibits a sonic glitch whenever the flux's derivative is zero, which is a known issue~\cite{van1989sonic}.  Those solutions accurately capture shocks whereas classical FVMs tend to smear them.

\input{tex/burgers_2d_final_density.tex}

The results presented for both the LWR and Burgers' equations show only a fraction of the potential of the proposed method. In addition to adding inputs in space and time, other problem-specific variables can also be incorporated at essentially no additional implementation cost, further improving precision.

\FloatBarrier

\subsection{Modeling experimental highway data}

A key advantage of NFVM is its ability to learn directly from data. While the conservation of mass remains a hard constraint, the flux function is flexible and can be inferred from empirical observations. Field experimental traffic data is often noisy and deviates from idealized conditions. In this context, learning the flux function offers a promising alternative to hand-designed models, potentially capturing complex or counter-intuitive behaviors that arise in practice.

\subsubsection{Dataset}

\input{tex/drone_highway.tex}

\input{tex/drone_tsd.tex}

The highway data used in this study comes from the Berkeley DeepDrive Drone Dataset~\cite{wu2022b3d}, which provides high-resolution aerial footage of highway traffic. The dataset is collected using drones flying over highways, capturing detailed vehicle trajectories and traffic patterns. As shown in~\Cref{fig:drone_highway}, the drone footage allows for precise tracking of individual vehicles over time, using computer vision algorithms.

\input{tex/drone_dataset_heatmap.tex}

We consider segment \texttt{highway\_4} from the dataset, as it contains the majority of traffic waves. The segment spans 400 meters over four lanes, on which vehicle trajectory data is recorded for 15 minutes. We compute the average density over cells of 10 meters long (spanning the whole width of the four lanes), over 1-second intervals. The resulting time-space diagram is shown in~\Cref{fig:drone_dataset_heatmap}, where four congestion waves propagating along the highway can be observed.

\input{tex/fundamental_diagram.tex}

The corresponding fundamental diagram is shown in~\Cref{fig:fundamental_diagram}. For reduced noise, we used a larger sliding window of 100 meters and compute density by counting the number of vehicles in that window. We then estimate flow by counting the number of vehicles passing through the middle of the window over an interval of 10 seconds.

\subsubsection{Calibrated flux functions}
\label{sec:calibrated_flux_functions}

We first evaluate the performance of classical numerical schemes on the drone dataset, which does not necessarily satisfy the LWR PDE or conservation properties. We consider five variants of LWR with different flow functions (Greenshields', Triangular, Trapezoidal, Greenberg, and Underwood, detailed in \Cref{app:lwr_flows}) and attempt to fit them to the experimental data. To this end, we consider the following two different approaches.

\paragraph{Calibration on scheme prediction} We calibrate the parameters of the flow functions by minimizing the mean-squared error (MSE) between predicted and true densities. Calibration is performed on the last 70 seconds of data, the same wave used to train the NFVM model, shown in~\Cref{fig:drone_dataset_heatmap}. Parameters are optimized with \texttt{scipy}'s differential-evolution solver~\cite{storn1997differential} within broad bounds around physically plausible values (e.g., the maximum speed $v_{\max}$ can range between 60 and 150km/h, while the wave propagation speed $-w$ can range from 10 to 70km/h). While that process might output a non-physical flow function, this relaxation seems to give the fairest comparison possible. For each candidate parameter set, we roll out the finite-volume scheme autoregressively across the entire spatio-temporal training domain; the resulting MSE is the objective minimized to get the optimal flow parameters $\theta^*$:
\begin{equation}
  \theta^* = \arg\min_{\theta} \| u_\text{fv} - u_\text{gt} \|_{L_2}, 
  \label{eq:mse_pred}
\end{equation}
where $u_\text{gt}$ is the ground truth density (on the training set) and $u_\text{fv}$ is the density predicted by the finite-volume scheme (using initial and boundary conditions from $u_\text{gt}$, as done previously on the synthetic data).

\paragraph{Calibration on fundamental diagram} We calibrate the parameters of the flow functions by minimizing the MSE between the flow function and the fundamental diagram computed on the experimental data (shown in \Cref{fig:fundamental_diagram}). This second approach is perhaps the most natural way to calibrate the flow function. In this case, we do not impose any bounds on the parameters, aside from them having the correct sign, as fitting the true fundamental diagram should yield physically realistic values. Parameters are again optimized with \texttt{scipy}'s differential-evolution solver. Optimal flow parameters $\theta^*$ is defined as:
\begin{equation}
  \theta^* = \arg\min_{\theta} \; \mathbb{E}_{(k, q) \in \text{FD}} \left[ (q - \max(f_\theta(k), 0))^2 \right],
  \label{eq:mse_fd}
\end{equation}
 where $\text{FD}$ represents all the (density, flow) points from the fundamental diagram, and the maximum is taken to ensure non-negative flow values.

For both calibration methods, and for all five fitted flows over five different FV schemes, \Cref{fig:drone_fitted_all} shows the predictions and fitted flows along with the MSE with respect to the fundamental diagram; \Cref{tab:drone_fitted_mse_pred} reports the MSE with respect to the autoregressive prediction; and \Cref{tab:drone_fitted_params} lists the parameters of the fitted flows. 

\textbf{Sanity check:} Flows calibrated on the fundamental diagram indeed achieve a lower MSE on the fundamental diagram (\Cref{fig:drone_fitted_all}, bottom) than flows calibrated on the autoregressive prediction (\Cref{fig:drone_fitted_all}, top), as shown in the figure. Reversely, flows calibrated on the autoregressive prediction achieve a lower MSE on the autoregressive prediction (\Cref{tab:drone_fitted_mse_pred}, top) than flows calibrated on the fundamental diagram (\Cref{tab:drone_fitted_mse_pred}, bottom), as evidenced in the table. One exception to the latter check is Greenshields' flow, for which the difference is small but can be explained by the fact that the flows are calibrated on the training set consisting of a single wave, while evaluations are done on the entire test set consisting of all four waves.

\input{tex/drone_fitted_fd.tex}  

\subsubsection{Results and comparison}

\input{tex/drone_dtw.tex}

\input{tex/drone_heatmap_best.tex}

Prediction results are displayed in~\Cref{fig:drone_heatmap_best}. Numerical schemes tend to suffer from glitching, as can be seen in \Cref{fig:drone_fitted_all}, especially at the $x=0$ boundary. This is an issue that does not happen with the learned models. 
A visual inspection in~\Cref{fig:drone_heatmap_best} shows the incredible performance of NFVM\(_4^5\); it gives the only solution that takes into account that the middle part of waves one, two and four are denser than their boundaries. Moreover, NFVM\(_4^5\) is the only model that correctly captures the fine-grained directional patterns of traffic density: the thin lines between the main waves are oriented along the true traffic flow, unlike in the other solutions where they are misaligned or absent. Its solution is also the only one that does not transform the flux coming at the $x=400m$ boundary at the beginning into a wave. This result is impressive considering the very limited amount of training data (only 70 timesteps) and the presence of noise. Besides, recall that all schemes predict the whole dataset autoregressively at once, only being given as input the initial condition ($t=0$) and the boundary conditions ($x=0$ and $x=400m$). 

However, in terms of traditional error metrics such as L1 and L2, due to averaging effects, the best numerical schemes outperform the learned model, despite producing visibly worse results. To address this discrepancy, we considered Dynamic Time Warping (DTW) \cite{muller2007dtw}, a standard algorithm for comparing time series. For each position $x$ along the highway, we compute the time series of traffic density for the ground truth, the FVM fit with the best $L_2$ score, and NFVM$_4^5$. As shown in~\Cref{fig:drone_dtw}, the NFVM consistently achieves lower DTW distance than the FVM. We hypothesize that this metric better captures temporal alignment of the wave patterns: even small misalignments in the fine-grained wave structure can lead to poor L1 and L2 scores, while DTW remains sensitive to the overall temporal dynamics. In contrast, the FVM produces smoother but less realistic predictions, while it performes better under L1 and L2 norms, this is largely due to the absence of such fine-scale variations. However, obtaining more robust metrics to better assess the predictive performance of different models in a way that matches visual expectations remains challenging.

\section{Conclusion} \label{sec:conclusion}

We illustrated that NN-based numerical solvers can provide an efficient and theoretically justified framework to solve PDEs, which can be adapted to data-driven approaches, even with limited data.

For future work, several key directions remain open. First, extending these methods to higher-dimensional and time-dependent PDEs would test their scalability and flexibility. Second, enforcing entropy conditions, in the same line as \cite{de2024wpinns}, could improve the theoretical guarantees. Moreover, integrating learning-based solvers with hybrid approaches, like combining data-driven techniques with classical numerical methods, may improve accuracy and robustness. Third, adapting the framework to incorporate advanced computational techniques, such as adaptive meshes and time stepping, could further enhance performance.


\bibliographystyle{siamplain}
\newpage
\bibliography{references}

\newpage

\appendix

\section{Proofs of Section 3} \label{app:proofs}

\begin{lemma}\label{approx:clip}
Let \(F\in\mathcal{F}\) be bounded and fix \(\epsilon>0\).
Let \(\theta\in\Theta\). If an (unclipped) neural network \(\tilde F_{\theta}\) satisfies
\begin{itemize}
\item[(i)] \(\|F-\tilde F_{\theta}\|_{\infty}\le \epsilon\).
\item[(ii)] \(D \geq \max_{\substack{(u^-, u^+) \in [0, u_{\max}]^2}} \left| F(u^-, u^+) \right|\).
\item[(iii)] \(\forall u^-,u^+\in[0,u_{\max}]\quad
    F_{\theta}(u^-,u^+)\;=\;
    \operatorname{clip}\!\left(
    \tilde F_{\theta}(u^-,u^+),\,
    -D,\,D \right)\).
\end{itemize}
\noindent Then
\[
\|F-F_{\theta}\|_{\infty}\leq \epsilon .
\]
\end{lemma}

\begin{proof}
Let \((u^-,u^+)\in[0,u_{\max}]^{2}\).  
Using the first assumption, we have
\[
\left|\tilde F_{\theta}(u^-,u^+)-F(u^-,u^+)\right|\le\epsilon .
\]

\paragraph{Case 1: No clipping}
If
\(
\tilde F_{\theta}(u^-,u^+)\in[-D,D],
\)
then
\(
F_{\theta}(u^-,u^+)=\tilde F_{\theta}(u^-,u^+)
\). Thus,
\[
\left|F_{\theta}(u^-,u^+)-F(u^-,u^+)\right|
= \left|\tilde F_{\theta}(u^-,u^+)-F(u^-,u^+)\right|
\le\epsilon .
\]

\paragraph{Case 2: Upper clipping.}
If
\(
\tilde F_{\theta}(u^-,u^+)>D\) then \(
F_{\theta}(u^-,u^+)=D.
\) which gives,
\begin{align*}
\left|F_{\theta}(u^-,u^+)-F(u^-,u^+)\right|
&= D - F(u^-,u^+) \\
&\le \tilde F_{\theta}(u^-,u^+) - F(u^-,u^+) \\
&= \left|\tilde F_{\theta}(u^-,u^+)-F(u^-,u^+)\right| \\
&\le \epsilon.
\end{align*}

\paragraph{Case 3: Lower clipping.}
If
\(
\tilde F_{\theta}(u^-,u^+)< -D
\) then \(
F_{\theta}(u^-,u^+)= -D
\), giving,
\begin{align*}
\left|F_{\theta}(u^-,u^+)-F(u^-,u^+)\right|
&= F(u^-,u^+) + D \\
&\le F(u^-,u^+) - \tilde F_{\theta}(u^-,u^+) \\
&= \left|\tilde F_{\theta}(u^-,u^+)-F(u^-,u^+)\right| \\
&\le \epsilon.
\end{align*}

In all three cases, the pointwise error is at most \(\epsilon\).  
Taking the supremum over \((u^-,u^+)\in[0,u_{\max}]^{2}\) finally yields
\[
\left\|F - F_{\theta}\right\|_{\infty} \le \epsilon,
\]
\end{proof}

\begin{lemma}{\label{proof: Lemma unsupervised}}{}

For all \( \epsilon, R \in \R_+^\star \), there exists \( N \in \mathbb{N} \) such that
\[
    \mathcal{L}^R_{N, \mathrm{unsup}}(F_{\mathrm{LxF}}) \leq \epsilon.
\]
\end{lemma}{}{}
\begin{proof} Let \(W\) a Lipschitz constant of \(f\) and \(u_{\max}\) the maximum value of \(u\). Let \(L = \max\{1, W\}\). We have 
\[
C = \max \set{u_{\max}, \max_{u \in [0, u_{\max}]} f(u)}.
\]

For all \( u \in (L^\infty\cap \mathrm{BV})\left(\mathbb{R}\times[0,T); [0, u_{\max}]\right) \), we define:
\[
W^R(u) \coloneqq \mathbb{E}_{\varphi \sim \Phi_R} \left[
\left(
\int_0^T \int_{\mathbb{R}} u\, \partial_t \varphi
+ \int_0^T \int_{\mathbb{R}} f(u)\, \partial_x \varphi
+ \int_{\mathbb{R}} u(0,\cdot) \varphi(0,\cdot)
\right)^2
\right].
\]
Let \( u, v \in (L^\infty \cap \mathrm{BV})\left(\mathbb{R}\times[0,T); [0, u_{\max}]\right) \). 
We bound the difference:
\begin{align*}
|W^R(u) - W^R(v)| 
&\leq \mathbb{E}_{\varphi \sim \Phi_R}
\Bigg|
\Bigg(
\int_0^T \!\!\int_{\mathbb{R}} \left( u\, \partial_t \varphi
+ f(u)\, \partial_x \varphi\right)
+ \int_{\mathbb{R}} u(0,\cdot) \varphi(0,\cdot)
\Bigg)^2 \\
&\qquad\qquad\qquad -
\Bigg(
\int_0^T \!\!\int_{\mathbb{R}}\left( v\, \partial_t \varphi
+ f(v)\, \partial_x \varphi \right)
+ \int_{\mathbb{R}} v(0,\cdot) \varphi(0,\cdot)\, 
\Bigg)^2
\Bigg|
 \\
&\leq \mathbb{E}_{\varphi \sim \Phi_R} \Bigg[
\Bigg(
\int_0^T \!\!\int_{\mathbb{R}} 
\left|(u + v)\partial_t \varphi + (f(u) + f(v))\partial_x \varphi\right| \\
&\qquad\qquad\qquad + \int_{\mathbb{R}} \left|u(0,\cdot)+v(0,\cdot)\right|\times\left|\varphi(0,\cdot)\right|
\Bigg) \\
&\qquad\qquad\qquad\qquad\times
\Bigg(
\int_0^T \!\!\int_{\mathbb{R}} 
|(u - v)\partial_t \varphi + (f(u) - f(v))\partial_x \varphi|
\Bigg)
\Bigg].
\end{align*}

Let \( \varphi \sim \Phi_R \). We bound the first term:
\[
\int_0^T \int_{\mathbb{R}} |(u + v)\partial_t \varphi + (f(u) + f(v))\partial_x \varphi|
\leq 2 C R \  \text{and}\ 
\int_{\mathbb{R}} |u(0,\cdot)+v(0,\cdot)|\, |\varphi(0,\cdot)| \leq 2 C R.
\]

\noindent For the second term:
\begin{align*}
\int_0^T \int_{\mathbb{R}} |(u - v)\partial_t \varphi + (f(u) - f(v))\partial_x \varphi|
&\leq \int_0^T \int_{\mathbb{R}} \left( |u - v||\partial_t \varphi| + |f(u) - f(v)||\partial_x \varphi| \right) \\
&\leq \int_0^T \int_{\mathbb{R}} |u - v| \left(|\partial_t \varphi| + W|\partial_x \varphi|\right) \\
&\leq R L \int_0^T \int_{\mathbb{R}} |u - v| \\
&= R L\|u - v\|_{L^1(\mathbb{R} \times [0,T))}.
\end{align*}
\noindent Hence,
\begin{equation} \label{eq:WR-diff}
|W^R(u) - W^R(v)| \leq 4CL R^2 \|u - v\|_{L^1(\mathbb{R} \times [0,T))}.
\end{equation}

Let \( u_0 \in U_0 \). We apply \eqref{eq:WR-diff} using \( u = H_N(F_{\mathrm{LxF}})(u_0) \), \( v = H^*(u_0) \), noting their values lie in \( [0, u_{\max}] \). Then we have:
\[
|W^R(H_N(F_{\mathrm{LxF}})(u_0)) - W^R(H^*(u_0))| \leq 4CL R^2 \|H_N(F_{\mathrm{LxF}})(u_0) - H^*(u_0)\|_{L^1(\mathbb{R} \times [0,T))}.
\]

At this point, it is important to observe that if we denote by \( \tilde{H}_N(F) \) the solution obtained using the update \( \tilde{h} \), we have \( \tilde{H}_N(F) = H_N(F) \). Indeed, as shown in \cite{barth2004finite}, since $F_{\mathrm{LxF}}$ is both monotonous and consistent, we have \( u_0 \in U_0 \) and all \( (t,x) \in [0,T) \times \mathbb{R} \),
\[
\forall u_0 \in U_0 \forall t \in [0,T) \forall x \in \mathbb{R}\quad 0 \leq \tilde{H}_N(F_{\mathrm{LxF}})(u_0)(t,x) \leq u_{\max}.
\]
Thus, the clipping step has no effect on the update and we apply a standard finite volume update.

As shown in \cite{godunov1959finite} and \cite{Sanders1983}, there exists a constant \( K_{f,\alpha} \) such that for any consistent and monotone flux, for example \( F_{\mathrm{LxF}} \), for any \(t \in [0,T) \):
\begin{align*}
\|H_N(F_{\mathrm{LxF}})(u_0)(t,\cdot) - H^*(u_0)(t,\cdot)\|_{L^1(\mathbb{R})}
&= \|\tilde{H}_N(F_{\mathrm{LxF}})(u_0)(t,\cdot) - H^*(u_0)(t,\cdot)\|_{L^1(\mathbb{R})} \\
&\leq K_{f,\alpha} \left(\Delta x^N + (\Delta x^N T)^{1/2}\right) \|u_0\|_{\mathrm{BV}}.
\end{align*}
We recall that since \( H^*(u_0) \) is the entropy solution, \( W^R(H^*(u_0)) = 0 \).

For large \( N \), we may assume (the proof still holds without this, but the expression of the upper bound is more complex):
\[
\left\| H_N(F_{\mathrm{LxF}})(u_0)(t,\cdot) - H^*(u_0)(t,\cdot) \right\|_{L^1(\mathbb{R})}
\leq 2K_{f,\alpha} \sqrt{\Delta x^N T} \, \|u_0\|_{\mathrm{BV}}.
\]
Let \( I = \sup_{u_0 \in U_0} \|u_0\|_{\mathrm{BV}} \). Integrating over time yields:
\begin{align} \label{bound:K}
\|H_N(F_{\mathrm{LxF}})(u_0) - H^*(u_0)\|_{L^1(\mathbb{R} \times [0,T))}
&\leq 2 K_{f,\alpha} (\Delta x^N)^{1/2} T^{3/2} \|u_0\|_{\mathrm{BV}} \notag \\
&\leq 2 K_{f,\alpha} (\Delta x^N)^{1/2} T^{3/2} I, 
\end{align}
Taking the expectation over \( u_0 \sim U_0 \) in \eqref{eq:WR-diff} and using \eqref{bound:K}, we obtain:
\[
\mathcal{L}_{N,\mathrm{unsup}}^R(F_{\mathrm{LxF}}) \leq 8CL R^2 K_{f,\alpha} (\Delta x^N)^{1/2} T^{3/2} I.
\]
There exist \( N \) such that:
\[
\Delta x^N \leq \left( \frac{\epsilon}{8CL R^2 K_{f,\alpha} T^{3/2} I} \right)^2,
\]
Finally, we have:
\[
\mathcal{L}_{N,\mathrm{unsup}}^R(F_{\mathrm{LxF}}) \leq \epsilon.
\]
\end{proof}
\begin{corollary}{\label{proof:*unusp}}{}
For any fixed \( \epsilon, R \in \R_+^\star \), there exists \( N \in \mathbb{N} \) such that for a neural network of depth and width at least
\[
\mathcal O\left(
\log(R) + \log\left(\frac{1}{\epsilon} \right) +
\left\lfloor \frac{T}{\Delta t^N} \right\rfloor
\right),
\]
we have
\[
\mathcal{L}_{N,\mathrm{unsup}}^R(F_{\theta^{*,R}_{N,\mathrm{unsup}}}) \leq \epsilon.
\]
\end{corollary}{}{}
\begin{proof}
Considering \eqref{F_L:unsup}, we may choose \( N \) such that
\[
\mathcal{L}_{N,\mathrm{unsup}}^R(F_{\text{LxF}}) \leq \frac{\epsilon}{2}.
\]
Employing Theorem \ref{T:unsup}, for a network of depth and width at least
\[
\mathcal O\left(
\log(R) + \log\left(\frac{1}{\epsilon} \right) +
\left\lfloor \frac{T}{\Delta t^N} \right\rfloor
\right),
\]
we have
\[
\mathcal{L}_{N,\mathrm{unsup}}^R(F_{\theta^{*,R}_{N,\mathrm{unsup}}}) \leq \frac{\epsilon}{2} + \mathcal{L}_{N,\mathrm{unsup}}^R(F_{\text{LxF}}).
\]
Hence,
\[
\mathcal{L}_{N,\mathrm{unsup}}^R(F_{\theta^{*,R}_{N,\mathrm{unsup}}}) \leq \epsilon.
\]
\end{proof}

\begin{corollary}{\label{proof:unsup }}
For all \( \epsilon, R \in \R_+^\star \), and \( \delta \in \R_+^\star \),  
there exist \( N \in \mathbb{N} \), \( n_s \in \mathbb{N} \) such that,  
for a network of depth and width at least
\[
O\left(
\log(R) + \log\left(\frac{1}{\epsilon} \right) +
\left\lfloor \frac{T}{\Delta t^N} \right\rfloor
\right),
\]
with probability \( 1 - \delta \), we have
\[
\mathcal{L}_{N,\mathrm{unsup}}^R\left(F_{\widehat{\theta}^{R,n_s}_{N,\mathrm{unsup} }}\right) \leq \epsilon.
\]
    
\end{corollary}

\begin{proof}
From Corollary \ref{cor:unsup-univ}, for a network of depth and width at least
\[
\mathcal O\left(
    \log(R) + \log\left(\frac{1}{\epsilon} \right) +
    \left\lfloor \frac{T}{\Delta t^N} \right\rfloor 
\right),
\]
there exists \( N \in \mathbb{N} \) such that
\[
\mathcal{L}_{N,\mathrm{unsup}}^R\left(F_{\theta^{*,R}_{N,\mathrm{unsup}}}\right) \leq \frac{\epsilon}{2}.
\]

Moreover, using Theorem \ref{conceration:unsup}, there exists \( n_s \in \mathbb{N} \) such that, with probability  \( 1 - \delta \),
\[
\mathcal{L}_{N,\mathrm{unsup}}^R\left(F_{\widehat{\theta}^{R,n_s}_{N,\mathrm{unsup} }} \right)
\leq \mathcal{L}_{N,\mathrm{unsup}}^R\left(F_{\theta^{*,R}_{N,\mathrm{unsup}}} \right) + \frac{\epsilon}{2}.
\]

Hence,
\[
\mathcal{L}_{N,\mathrm{unsup}}^R\left(F_{\widehat{\theta}^{R,n_s}_{N,\mathrm{unsup} }} \right) \leq \epsilon.
\]
\end{proof}

\begin{theorem}{\label{proof: concentration:sup}}{}
    For all \(N\in \mathbb{N}\), \( \delta > 0 \), \( \epsilon > 0 \), there exists \( n_s \in \mathbb{N} \) such that with probability  \( 1 - \delta \),
\[
\mathcal{L}_N^{\sup }(F_{\widehat{\theta}^{n_s}_{N, \sup}}) \leq \mathcal{L}_N^{\sup }(F_{\theta^*_{N, \sup}}) + \epsilon.
\]
\end{theorem}
\begin{proof}

We start by decomposing the quantity of interest:
\begin{equation}
\begin{aligned}\label{decompoS}
&\mathcal{L}_{N, \sup}(F_{\theta^*_{N, \sup}}) - \mathcal{L}_{N, \sup}(F_{\widehat{\theta}^{n_s}_{N,\sup}}) \\
&\quad= \mathcal{L}_{N, \sup}(F_{\theta^*_{N, \sup}}) - \widehat{\mathcal{L}}_{N,\sup}^{n_s}(F_{\widehat{\theta}^{n_s}_{N,\sup}}) + \widehat{\mathcal{L}}_{N,\sup}^{n_s}(F_{\widehat{\theta}^{n_s}_{N,\sup}}) - \mathcal{L}_{N, \sup}(F_{\widehat{\theta}^{n_s}_{N,\sup}}).
\end{aligned}
\end{equation}
Note that
\[
\mathbb{E}_{u_0 \sim U_0}\left[ \widehat{\mathcal{L}}_{N,\sup}^{n_s}(F) \right] = \mathcal{L}_{N, \sup}(F).
\]

Furthermore, we can show that \( m_{N, \sup}(F_{\theta}, u_0) \) is uniformly bounded for all \( u_0 \) and \(\theta\). Indeed, for all \( u_0 \in U_0 \) ,

We recall that : 
\[
m_{N, \sup}(F_\theta, u_0)= \|H_N(F_{\theta})(u_0)-H^{*}(u_0)\|_{L^1([0,T)\times \mathcal{X})}
\]

\[
\begin{aligned}
|m_{N, \sup}(F_{\theta}, u_0)| - |m_{N, \sup}(F_{\mathrm{LxF}}, u_0)|
&\leq |m_{N, \sup}(F_{\theta}, u_0) - m_{N, \sup}(F_{\mathrm{LxF}}, u_0)| \\
&\leq \|H_N(F_{\theta}) - H_N(F_{\mathrm{LxF}})\|_{L^1([0,T)\times \mathcal{X})} \\
&\leq T |\mathcal{X}|\|H_N(F_{\theta}) - H_N(F_{\mathrm{LxF}})\|_{\infty} \\
&\leq A_N T |\mathcal{X}| \cdot \|F_{\theta} - F_{\mathrm{LxF}}\|_{\infty}.
\end{aligned}
\]

Since
\begin{align*}
|m_{N, \sup}(F_{\mathrm{LxF}}, u_0)|
&= \|H_N(F_{\mathrm{LxF}})(u_0) - H^{*}(u_0)\|_{L^1([0,T)\times \mathcal{X})} \\
&\leq \|H_N(F_{\mathrm{LxF}})(u_0)\|_{L^1([0,T)\times \mathcal{X})} + \|H^{*}(u_0)\|_{L^1([0,T)\times \mathcal{X})} \\
&\leq 2|\mathcal{X}| T u_{\max},
\end{align*}
we will have that
\begin{align*}
|m_{N, \sup}(F_{\theta}, u_0)| 
&\leq 2|\mathcal{X}| T u_{\max} 
+ A_N T |\mathcal{X}| \left( \|F_{\theta}\|_{\infty} + \|F_{\mathrm{LxF}}\|_{\infty} \right) \\
&\leq 2|\mathcal{X}| T u_{\max} + 2 A_N T |\mathcal{X}| D.
\end{align*}

Hence, we conclude that \( m_n^{\sup}(F_{\theta}, u_0) \) is uniformly bounded. We denote this upper bound by \( B_{\sup}=2|\mathcal{X}|Tu_{\max} + 2A_N T |\mathcal{X}| D \).
Applying the Hoeffding's Lemma to the i.i.d. sequence
\begin{equation*} 
\begin{cases}
\big[m_{N,\sup}(F_\theta, \cdot)\big]_i: U_0 \to \RR, \quad  i \in \set{1, \cdots, n_s} \\
u_0 \mapsto m_{N,\sup}(F_\theta, u_0), 
\end{cases}
\end{equation*}
we may conclude that for
\[
n_s \geq \frac{B_{\sup}^2}{2(\epsilon/2)^2} \log\left( \frac{2}{\delta} \right),
\]
with probability  \( 1 - \delta \),
\begin{equation} \label{HoeffdingS}
\left| \widehat{\mathcal{L}}_{N,\sup}^{n_s}(F_{\widehat{\theta}^{n_s}_{N,\sup}}) - \mathcal{L}_{N, \sup}(F_{\widehat{\theta}^{n_s}_{N,\sup}}) \right| \leq \epsilon / 2.
\end{equation}
    
This implies
\[
\mathcal{L}_{N, \sup}(F_{\widehat{\theta}^{n_s}_{N,\sup }}) \leq \widehat{\mathcal{L}}_{N,\sup}^{n_s}(F_{\widehat{\theta}^{n_s}_{N,\sup }}) + \epsilon / 2,
\]
and additionally, by definition of the minimizer,
\begin{equation} \label{def_minS}
\mathcal{L}_{N, \sup}(F_{\theta^*_{N, \sup}}) \leq \widehat{\mathcal{L}}_{N,\sup}^{n_s}(F_{\widehat{\theta}^{n_s}_{N,\sup }}) + \epsilon / 2.
\end{equation}

By combining \eqref{decompoS}, \eqref{HoeffdingS}, and \eqref{def_minS}, we conclude:
\[
\mathcal{L}_{N, \sup}(F_{\theta^*_{N, \sup}}) - \mathcal{L}_{N, \sup}(F_{\widehat{\theta}^{n_s}_{N,\sup }}) \leq \epsilon.
\]
\end{proof}
\section{Flow Models: LWR PDE Model}\label{app:lwr_flows}

\input{tex/flows.tex}

This work considers six distinct traffic flow models from the literature for the LWR equation~\eqref{eq:lwr}, resulting in six LWR models with significantly different behaviors. Depending on the application and the specific dataset or road under consideration, the flux function can take various forms~\cite{ardekani2011macroscopic}. We summarize below the different flux functions considered in this study, which include the most commonly used models. All of them map the vehicle density $\rho \in [0, \rho_{\max}]$ to the flow $f \in [0, f_{\max}]$, where $\rho_{\max}$ and $f_{\max}$ denote the maximum possible density (road capacity) and flow, respectively.

In many classical conservation law models, the flow function is assumed to be concave, satisfying $f(0) = 0$ (no vehicles means no flow) and $f(\rho_{\max}) = 0$ (full road implies no movement), or $\lim_{\rho \to \rho_{\max}} f(\rho) = 0$ if $\rho_{\max}$ is infinite. The maximum flow is typically attained at a critical density $\rho_c \in (0, \rho_{\max})$, with densities below $\rho_c$ corresponding to free-flow and densities above $\rho_c$ to congestion. However, a few of the flux functions considered here depart from these classical assumptions, for example by lacking concavity or a unique critical density.

The models considered in this work can be seen in~\Cref{fig:flows}, and their closed-form formulations are given as follows.

\paragraph{Greenshields' flow~\cite{greenshields1935study}}
The Greenshields flux is given by
\[
f(\rho) = v_{\max} \, \rho \left( 1 - \frac{\rho}{\rho_{\max}} \right),
\]
with default parameters $v_{\max} = 1.0$ and $\rho_{\max} = 1.0$. This classical parabolic flow is concave and achieves its maximum at critical density $\rho_c = \rho_{\max}/2$.

\paragraph{Triangular flow~\cite{geroliminis2008existence}}
The triangular flux takes the form
\[
f(\rho) = 
\begin{cases}
v_{\max} \rho & \text{if } 0 \leq \rho \leq \rho_c, \\
w \, (\rho - \rho_{\max}) & \text{if } \rho_c < \rho \leq \rho_{\max},
\end{cases}
\]
with critical density
\[
\rho_c = \frac{\rho_{\max} w}{w - v_{\max}}.
\]
We consider two variants of this flux: a symmetric form (“Triangular Sym”) with $v_{\max} = 1.0$, $w = -1.0$, $\rho_{\max} = 1.0$, and a skewed form (“Triangular Skw”) with a higher maximum speed $v_{\max} = 2.0$, and other parameters unchanged.

\paragraph{Trapezoidal flow~\cite{geroliminis2011properties}}
The trapezoidal flux can be expressed as
\[
f(\rho) = 
\begin{cases}
v_{\max} \rho & 0 \leq \rho \leq \rho_{c_1}, \\
v_{\max} \rho_{c_1} + s (\rho - \rho_{c_1}) & \rho_{c_1} < \rho \leq \rho_{c_2}, \\
w (\rho - \rho_{\max}) & \rho_{c_2} < \rho \leq \rho_{\max},
\end{cases}
\]
where the slope of the transition region is
\[
s = \frac{w (\rho_{c_2} - \rho_{\max}) - v_{\max} \rho_{c_1}}{\rho_{c_2} - \rho_{c_1}}.
\]
By default, we set $v_{\max} = 1.0$, $w = -1.5$, $\rho_{c_1} = 0.2$, $\rho_{c_2} = 0.8$, and $\rho_{\max} = 1.0$. When $s=0$, the transition is flat with maximum flow $q_{\max} = v_{\max} \rho_{c_1}$; however, we also allow for a more general form with nonzero slope in the transition region, which is therefore not strictly a trapezoid.

\paragraph{Greenberg flow~\cite{greenberg1959analysis}}
The Greenberg flux is given by
\[
f(\rho) = c_0 \, \rho \, \log \left( \frac{\rho_{\max}}{\rho} \right),
\]
with default parameters $c_0 = 2.0$ and $\rho_{\max} = 1.0$. This logarithmic function is concave on $(0, \rho_{\max}]$.

\paragraph{Underwood flow~\cite{Underwood1961}}
The Underwood flux takes the form
\[
f(\rho) = c_1 \, \rho \, \exp \left( -c_2 \rho + 1 \right),
\]
with default parameters $c_1 = 0.25$, $c_2 = 1.0$, and $\rho_{\max} = 1.0$. This model represents an exponentially decaying flow with respect to density beyond a certain threshold and is not necessarily concave depending on the parameters.

\section{Riemann Problems and Solutions: LWR PDE Model} \label{app:lwr_riemann}

A Riemann problem is a special type of initial value problem for a hyperbolic conservation law (here for the LWR equation~\eqref{eq:lwr}) where the initial condition is piecewise constant with a single discontinuity. Formally, the initial data is given by
\[\rho_0(x) = \begin{cases}
  \rho_1, & x < x_0 \\
  \rho_2, & x \geq x_0
\end{cases}\]
where $\rho_1$ and $\rho_2$ are constant states to the left and right of the discontinuity located at $x_0$. The Riemann problem is fundamental in the theory of conservation laws because its solution describes the evolution of a single jump discontinuity, and provides the building block for understanding more general initial conditions by superposition or approximate methods. 

For LWR, the entropy solution depends on whether a shock or a rarefaction develops:
\begin{itemize}
  \item \textbf{Shock case:} If $\rho_1 < \rho_2$, a shock wave forms, propagating at the Rankine-Hugoniot speed
        \[
           s = \frac{f(\rho_2)-f(\rho_1)}{\rho_2-\rho_1}.
        \]
        The solution is
        \[
          \rho(t,x)=
          \begin{cases}
            \rho_1, & (x-x_0)/t < s,\\
            \rho_2, & (x-x_0)/t \ge s.
          \end{cases}
        \]

  \item \textbf{Rarefaction case:} If $\rho_1 > \rho_2$, a rarefaction (expansion) wave forms, with left and right characteristic speeds
        \[
          \lambda_1 = f'(\rho_1), 
          \qquad
          \lambda_2 = f'(\rho_2).
        \]
        The solution is
        \[
          \rho(t,x)=
          \begin{cases}
            \rho_1, & (x-x_0)/t < \lambda_1,\\
            (f')^{-1}((x-x_0)/t), & (x-x_0)/t \in [\lambda_1,\lambda_2),\\
            \rho_2, & (x-x_0)/t \ge \lambda_2,
          \end{cases}
        \]
        where $(f')^{-1}$ is the inverse of the derivative of the flux function $f$ from the LWR model~\eqref{eq:lwr}.
\end{itemize}


\section{Numerical schemes}\label{app:numerical_schemes}

\subsection{Finite Volume Methods}

\paragraph{Godunov}
The Godunov Scheme is a first-order method which is widely accepted for numerical solutions of the hyperbolic conservation laws. The Godunov scheme considers the following numerical flux scheme to update the solution in each timestep. 
\begin{equation}\label{E:N_Godunov}
F_\text{G}\left(\rho_1, \rho_2\right)= 
\begin{cases}
    f\left(\rho_2\right) & \text { if } \rho_{\text{c}}<\rho_2<\rho_1 \\
    f\left(\rho_{\text{c}}\right) & \text { if } \rho_2<\rho_{\text{c}}<\rho_1 \\ 
    f\left(\rho_1\right) & \text { if } \rho_2<\rho_1<\rho_{\text{c}} \\
    \min \left(f\left(\rho_1\right), f\left(\rho_2\right)\right) & \text { if } \rho_1 \leq \rho_2.
\end{cases}
\end{equation}

\paragraph{Enquist Osher}
\begin{equation}
F_{\text{EO}}\left(\rho_1, \rho_2\right)= 
\begin{cases}
    f\left(\rho_1\right) & \text { if } \rho_2,\rho_1 \leq \rho_{\text{c}} \\
    f\left(\rho_2\right) & \text { if } \rho_2,\rho_1 > \rho_{\text{c}} \\
    f\left(\rho_{\text{c}}\right) & \text { if } \rho_2 \leq \rho_{\text{c}} <\rho_1 \\
    f(\rho_1)+f(\rho_2) - f(\rho_{\text{c}}) & \text { if } \rho_1 \leq \rho_{\text{c}} <\rho_2 .
\end{cases}
\end{equation}

\paragraph{Lax-Friedrichs}
\begin{equation}\label{E:N_flux_lax}
F_{\text{LxF}} : (\rho_{1},\rho_{2}) \mapsto \tfrac{1}{2}\big( f(\rho_{1}) + f(\rho_{2}) - \tfrac{\Delta x}{\Delta t}(\rho_{2} - \rho_{1})\big)
\end{equation}

\paragraph{Essentially non-oscillatory (ENO)} In the ENO scheme, the goal is to construct local polynomial interpolations using adaptive stencils. And the reconstruction of $\rho_{1},\rho_{2}$ relies on choosing a single ``best" stencil in a way that avoids discontinuities. Here we use a five-cell case as an example:

\paragraph{Weighted essentially non-oscillatory (WENO)} WENO is similar in spirit to ENO but combines multiple stencils in a weighted fashion instead of choosing just one. 

\subsection{Discontinuous Galerkin Method: Semi-Discretized Scheme}
We consider the LWR model in the conservative form:
\[
\partial_t \rho + \partial_x [f(\rho)] = 0,
\]
where $(x, t) \mapsto \rho(x, t)$ is the traffic density and $\rho \mapsto f(\rho)$ is the flux function. To approximate this equation using the Discontinuous Galerkin (DG) method, we partition the domain into a mesh of elements $(I_j)_j = ([x_{j-\frac{1}{2}}, x_{j+\frac{1}{2}}])_j$ and seek a piecewise polynomial approximation $\rho_h$ in each element.

Let $V_h^k$ be the space of discontinuous piecewise polynomials of degree at most $k$ on each element. The semi-discrete DG formulation reads: find $\rho_h \in V_h^k$ such that for each element $I_j$ and for all test functions $v \in V_h^k$,
\[
\int_{I_j} \partial_t \rho_h \, v \, \dt x 
- \int_{I_j} f(\rho_h) \, \partial_x v \, \dt x 
+ \mathcal F_{j+\frac{1}{2}} v^{-}_{j+\frac{1}{2}} 
- \mathcal F_{j-\frac{1}{2}} v^{+}_{j-\frac{1}{2}} = 0,
\]
where $\mathcal F_{j\pm \frac{1}{2}} = F(\rho_h^-, \rho_h^+)$ is a numerical flux function consistent with $f$, and $\rho_h^{\pm}$ denote the values of $\rho_h$ at the interfaces from the left and right, respectively.

This semi-discrete system of ODEs in time is advanced using a strong stability-preserving (SSP) Runge–Kutta method. For example, using a third-order SSP RK scheme:
\begin{align*}
\rho_h^{(1)} &= \rho_h^n + \Delta t \, L(\rho_h^n), \\
\rho_h^{(2)} &= \frac{3}{4} \rho_h^n + \frac{1}{4} \left( \rho_h^{(1)} + \Delta t \, L(\rho_h^{(1)}) \right), \\
\rho_h^{n+1} &= \frac{1}{3} \rho_h^n + \frac{2}{3} \left( \rho_h^{(2)} + \Delta t \, L(\rho_h^{(2)}) \right),
\end{align*}
where $L(\rho_h)$ denotes the spatial DG discretization operator defined by the weak form above.

\section{Lax-Hopf Exact Solution}\label{app:lax_hopf}

To generalize the explicit Riemann solution (which will be employed for training the NN) to any problem with piecewise constant initial conditions, we will use the Lax-Hopf algorithmic approach \cite{mazare2011analytical, claudel2010lax_a, claudel2010lax_b}; such explicit solutions will be used for the evaluation of the performance of the NN-solver. While from the theoretical point of view we are interested in the Cauchy problem \eqref{E:main}, from the practical point of view we need to define the problem by incorporating (numerical) boundary conditions. More precisely, let $T >0$ and $a, b\in \RR$ and in practice we are interested in the solution $\rho(t,x)$ of \eqref{E:main} over the domain $[0,T) \times [a, b]$. To avoid including the boundary values $x = a$ and $x = b$ in the calculation of solutions, we consider an extended domain $[a - \underline a, b + \bar b]$ for some $\underline a$ and $\bar b >0$ (the details of such extension will be elaborated in what follows). Denoting by $x_\circ = a - \underline a$ and $x_n = b + \bar b$, we consider the initial and boundary value (IBV) problem 
\begin{equation} \label{E:LWR}
    \begin{cases}
        \rho_t + f(\rho)_x = 0  \\
        \rho(0, x) = \phi_\circ(x)\\
        \rho(t, x_\ell) = \phi_L(t) \\
        \rho(t, x_u) = \phi_U(t) 
    \end{cases}
\end{equation}
over $[0, t_m] \times [x_\circ, x_n]$. Here, $\set{x_\circ, \cdots, x_n}$ and $\set{0, \cdots, t_m}$ are defined based on the piecewise initial and boundary functions $\phi_\circ$, $\phi_L$ and $\phi_U$. In particular, 
we consider the partition $\set{x_\circ, \cdots, x_n}$ for space and $\set{t_\circ, \cdots, t_m}$ for time. The initial densities $\set{\rho_\circ^{(i)}: i = 0, \cdots, n -1}$. The upstream flow $\set{\phi_{L}^{(j)}: j = 0, \cdots, m-1}$ and the downstream flow $\set{\phi_{U}^{(j)}: j = 0, \cdots, m - 1}$ are known. Using this approach the initial value can be presented as a piecewise constant function 
\begin{equation*}
    \rho_\circ(x) = \rho_\circ^{(i)}, \quad \forall x \in [x_i, x_{i + 1}],  
\end{equation*}
and similarly for the upstream and downstream boundaries. 
Note that the grid in both space and time is not necessarily uniform. In other words, $x_k - x_{k -1}$ and $t_r - t_{r - 1}$ are not necessarily uniform over $k$ and $r$. 


It is customary to address this problem in the context of Moskowitz framework (see e.g. \cite{claudel2010lax_a}). Let $M(t,x)$ define the Moskowitz function. Then we can retrieve the following information 
\begin{equation}
    \rho(t,x ) = \frac{-  \partial M(t, x)}{\partial x}, \quad f(\rho(t,x)) \Def \frac{\partial M(t, x)}{\partial t}
\end{equation}
and hence we have the Hamilton-Jacobi PDE for which $f$ is Hamiltonian:
\begin{equation}
    M_t(t,x) - f(- M_x(t,x )) = 0.
\end{equation}
We can then consider the following Dirichlet problem
\begin{equation} \label{E:Moskowitz}
    \begin{cases}
        M_t - f(-M_x) = 0  &, (t, x) \in [0, t_m] \times [x_\circ, x_n] \\
        M(0, x) = \varphi_\circ(x)\\
        M(t, x_\ell) = \varphi_L(t) \\
        M(t, x_u) = \varphi_U(t) 
    \end{cases}
\end{equation}
and the boundary conditions should be understood in the weak sense. 

The initial condition on the Moskowitz PDE can be obtained by integrating the initial condition of the LWR PDE with the assumption that $M(0, x_\circ) = 0$. In particular, for $x \in [x_k , x_{k + 1})$ we have
\begin{equation}\label{E:initial_M}
    M(0, x) = - \int_{x_\circ}^{x} \rho(0, \xi) d \xi = - \sum_{r =0}^{k -1}(x_{r+ 1} - x_r) \rho^{(r)}_\circ - (x - x_{k}) \rho_\circ^{(k)}
\end{equation}
Similarly, we can present the boundary functions $M(t, x_\circ)$ and $M(t, x_n)$ for any $t \in [t_j, t_{j + 1}]$. 
\subsection{Extension of the domain and the boundary data.}\label{S:Bounded_Domain}
We will use the Lax-Hopf method to calculate the solution $\rho(t, x)$ for any $(t, x) \in [0,T) \times [a, b]$. The main concept in such calculations is that depending on the desired point $(t,x)$ we recognize the initial and boundary conditions that are influenced by the maximum and minimum wave speed. 
To this end, if $\underline a = \underline a(a, v_{\max}, w)$ and $\bar b = \bar b(b, v_{\max}, w)$, which are defined above to extend the interval of interest, are sufficiently large then for any $(t, x) \in [0,T) \times [a, b]$ we have that (see \cite[Algorithm 1]{mazare2011analytical} for the details)
\begin{equation}\label{E:bdry_1}
   \min \set{i - 1: \frac{x - x_\circ}{v_{\max}} \ge  t - t_i} = 0
\end{equation}
and
\begin{equation}\label{E:bdry_2}
    \min \set{j - 1: \frac{x_n - x}{w} \le t_j - t} = - \infty
\end{equation}
with the convention that $\min \varnothing = - \infty$. Using \eqref{E:bdry_1} and \eqref{E:bdry_2}, by selecting $\underline a$ and $\bar b$ properly, we can only focus on the initial condition as the boundary values at $x = x_\circ$ and $x = x_n$ will not affect the solution $\rho(t,x)$ in the domain $[0,T) \times [a, b]$. Collecting all together, we will employ the following algorithmic scheme to find the exact solution of the initial value Cauchy problem \eqref{E:main}.  
\subsection{Lax-Hopf exact calculation of the solution: algorithmic approach}
In this section, we will discuss the main algorithm that will be used to calculate the explicit solution of the conservation law with any piecewise continuous flux function. 

\begin{itemize}
    \item[] \textbf{Input.} $x \in [x_\circ, x_n]$, and $t \in [0,T)$ 
    \item[] \textbf{Output.} Solution $\rho(t, x)$
     \item Let $J_\ell \Def  \max \set{0, \min(i -1) \mid x_i  \ge x - v_{\max} t} $; min space index for the initial data 
     \item Let $J_u \Def  \min \set{n -1, \max(i) \mid x_i  \le x + w t} $;  max space index for the initial data 
     \item $M =  \infty$; Initialization of Moskowitz function
     \item \textbf{For} $i = J_\ell$ to $J_u$ \textbf{do} 
     \begin{itemize}
         \item [] Compute $M_{c_\circ^{(i)}}(t, x)$ from \eqref{E:Mosk}
         \item[] \textbf{if} $M_{c_\circ^{(i)}}(t, x) < M$ \textbf{then}
         \begin{itemize}
            \item[$\circ$] Let $M =M_{c_\circ^{(i)}}(t, x)$ 
            \item[$\circ$] Compute $\rho = \rho_{c_\circ^{(i)}}(t,x)$ from \eqref{E:density}
        \end{itemize}
         \item[] \textbf{end if} 
     \end{itemize}
     \item \textbf{end For}
\end{itemize}
Let 
\begin{equation}
    x \in [x_i, x_{i + 1}], \quad c_\circ^{(i)} \Def M_\circ(x)
\end{equation}
where, $M_\circ(x)$ is defined as of \eqref{E:initial_M}.
Let us recall the convex transformation 
\begin{equation}
    R(u) \Def \sup_{\rho \in [0, \rhomax]}(f(\rho) - u \rho)
\end{equation}

The goal is to calculate the solution associated with an affine, locally defined initial condition, indexed by $i$:
\begin{equation}
   x \in [x_i, x_{i + 1}], \quad M_\circ^{(i)} = - \rho_\circ^{(i)} x + b_i,
\end{equation}
where, $b_i = \rho_\circ^{(i)} x_i - \sum_{l = 0}^{i - 1}(x_{l + 1} - x_l) \rho_\circ^{(l)}$ which allows for the continuity of the initial conditions in $(0, x_i)$. 
The analytical solution of the problem with the initial condition $M_\circ^{(i)}$ will be (See \cite{claudel2010lax_b})
\begin{equation}\label{E:Mosk}
    M_{c_\circ^{(i)}}(t, x) = \begin{cases}
        t f(\rho_\circ^{(i)}) - \rho_\circ^{(i)} x + b_i &, x_i + t f'(\rho_\circ^{(i)}) \le x \le x_{i + 1} + t f'(\rho_\circ^{(i)}) \\
        t R\left(\frac{x - x_i}{t}\right) - \rho_\circ^{(i)} x_i + b_i &, x_i - t w \le x \le x_i + t f'(\rho_\circ^{(i)}) \\
        t R \left(\frac{x - x_{i+ 1}}{t}\right) - \rho_\circ^{(i)} x_{i + 1} + b_i &, x_{i + 1} + t f'(\rho_\circ^{(i)}) \le x \le x_{i + 1} + t v_{\max}. 
    \end{cases}
\end{equation}
where by $f'$ we refer to the derivative as well as $R'$. It should be noted that if the derivative does not exist at some point, then the left or right derivative will be used. In addition, let us recall that $f$ has a left derivative $v_{\max}$ at $\rho = 0$ and $w$ at $\rho = \rhomax$. 
Then, the density function will be updated by 
\begin{equation}\label{E:density}
    \rho_{c_\circ^{(i)}}(t, x) \Def 
        - \frac{\partial M_{c_\circ^{(i)}}}{\partial x} = \begin{cases}
        \rho_{\circ}^{(i)} &, x_i + t f'(\rho_\circ^{(i)}) \le x \le x_{i + 1} + t f'(\rho_\circ^i) \\
        - R' \left(\frac{x - x_i}{t}\right) &, x_i - t w \le x \le x_i + t f'(\rho_\circ^{(i)}) \\
        - R' \left(\frac{x - x_{i + 1}}{t}\right) &, x _{i + 1} + t f'(\rho_\circ^{(i)}) \le x \le x_{i + 1} + t v_{\max}.
    \end{cases}
\end{equation}

For the triangular flux, this explicit form of solution can be simplified.

\end{document}

%% file: tex/lwr_heatmaps.tex
\begin{figure}
    \centering
        \begingroup
\setlength{\tabcolsep}{1pt}
    \begin{tabularx}{\textwidth}{m{.5cm}CCCCCC}
    & {\footnotesize Greenshields'} & {\footnotesize Triang. Sym} & {\footnotesize Triang. Skw} & {\footnotesize Trapezoidal} &  {\footnotesize Greenberg} & {\footnotesize Underwood} \\
    
    \(\vcenter{\hbox{\rotatebox{90}{\strut LxH}}}\) & \(\vcenter{\hbox{\includegraphics[width=.15\textwidth]{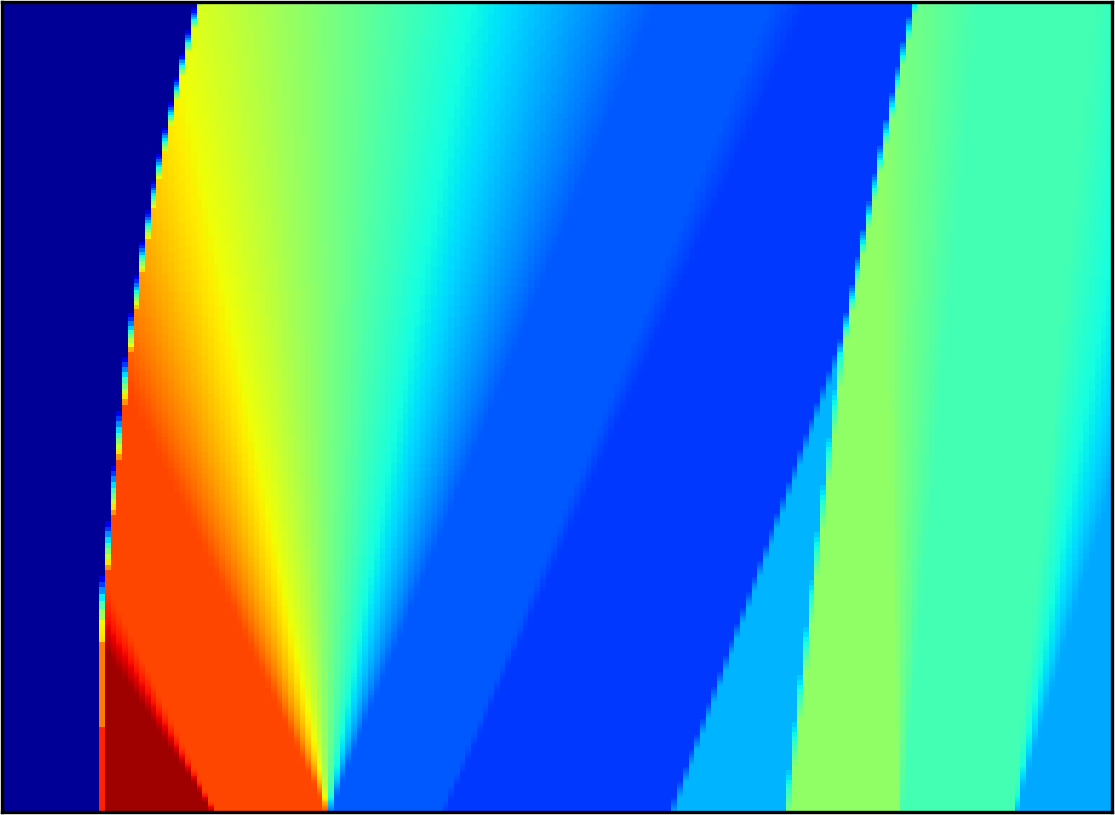}}}\) &
    \(\vcenter{\hbox{\includegraphics[width=.15\textwidth]{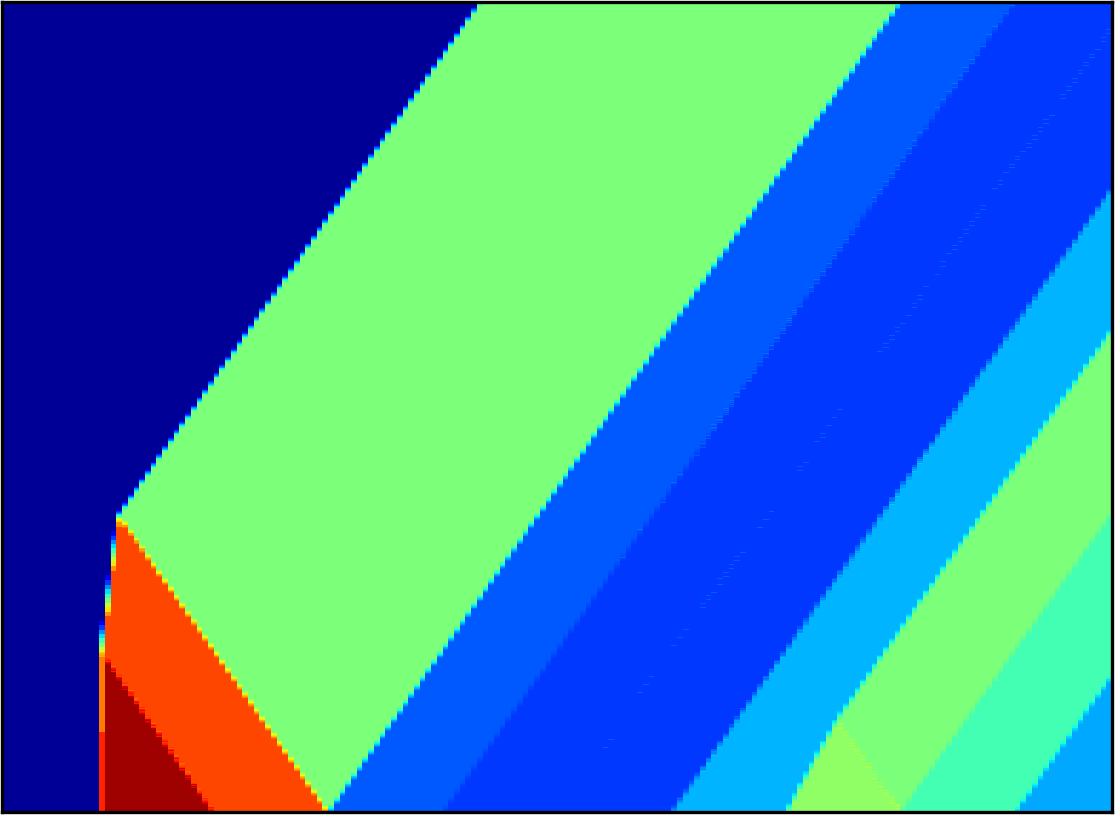}}}\) &
    \(\vcenter{\hbox{\includegraphics[width=.15\textwidth]{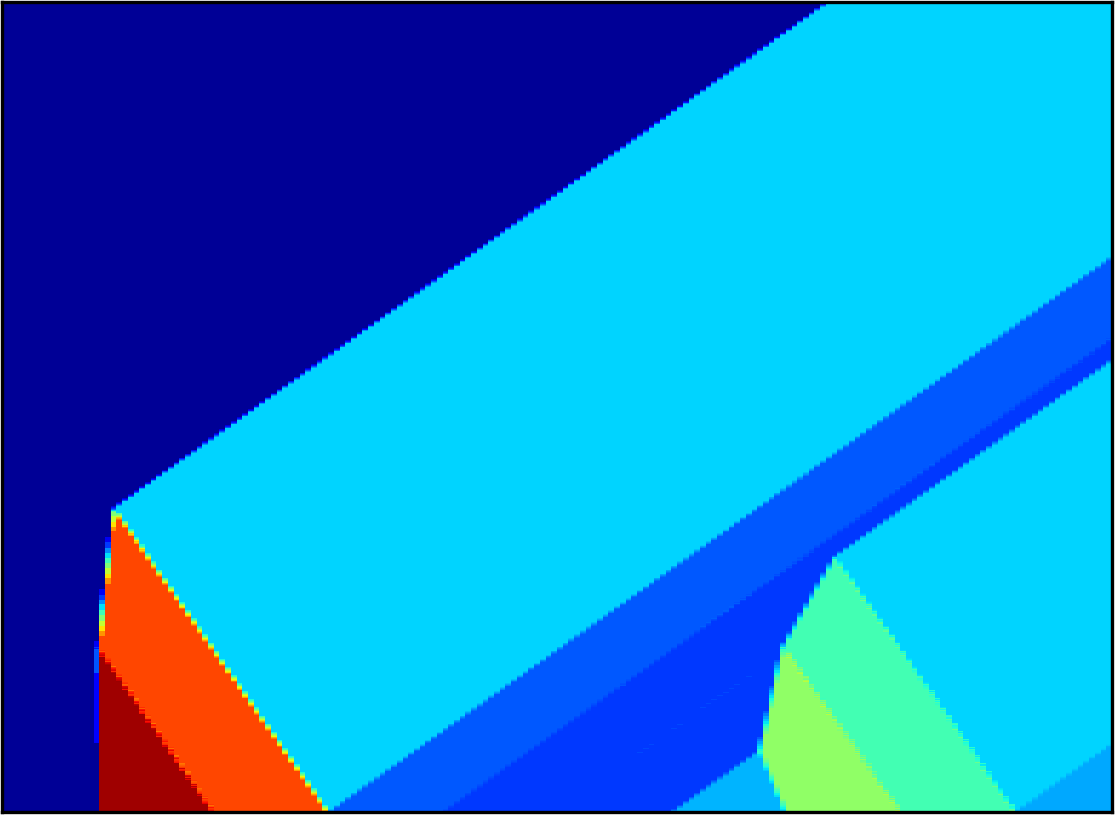}}}\) &
    \(\vcenter{\hbox{\includegraphics[width=.15\textwidth]{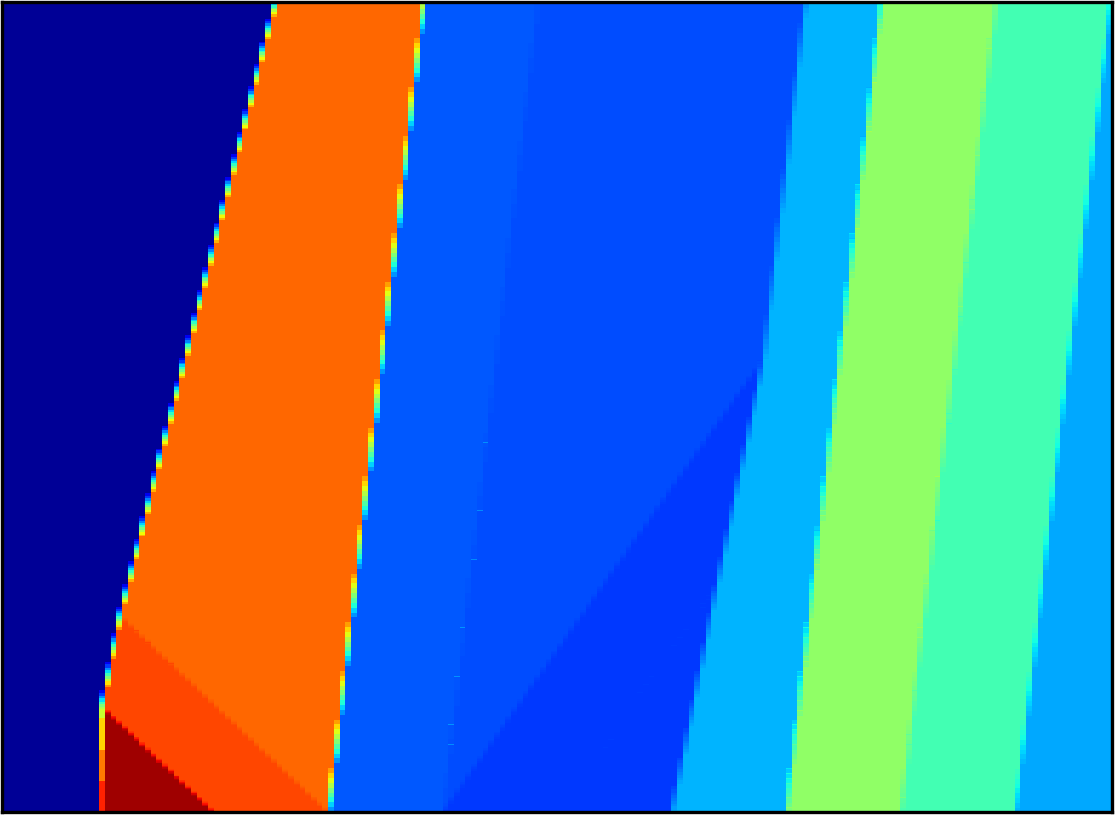}}}\) &
    \(\vcenter{\hbox{\includegraphics[width=.15\textwidth]{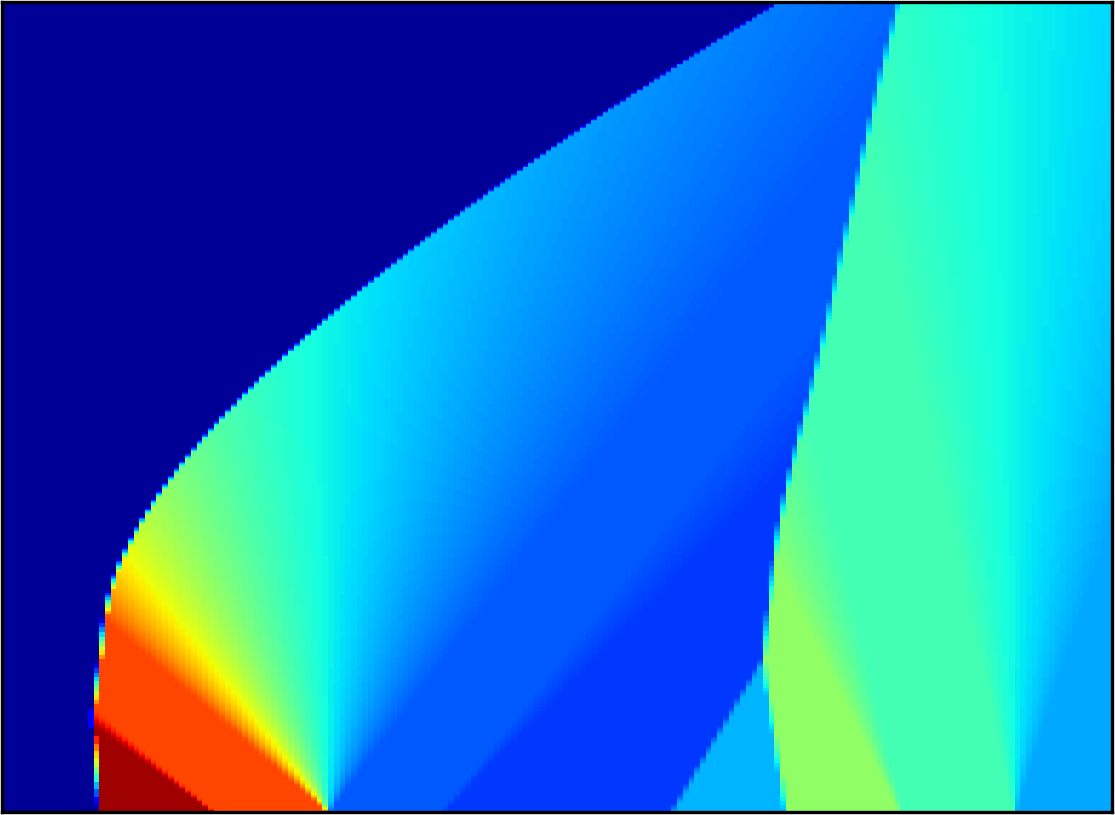}}}\) &
    \(\vcenter{\hbox{\includegraphics[width=.15\textwidth]{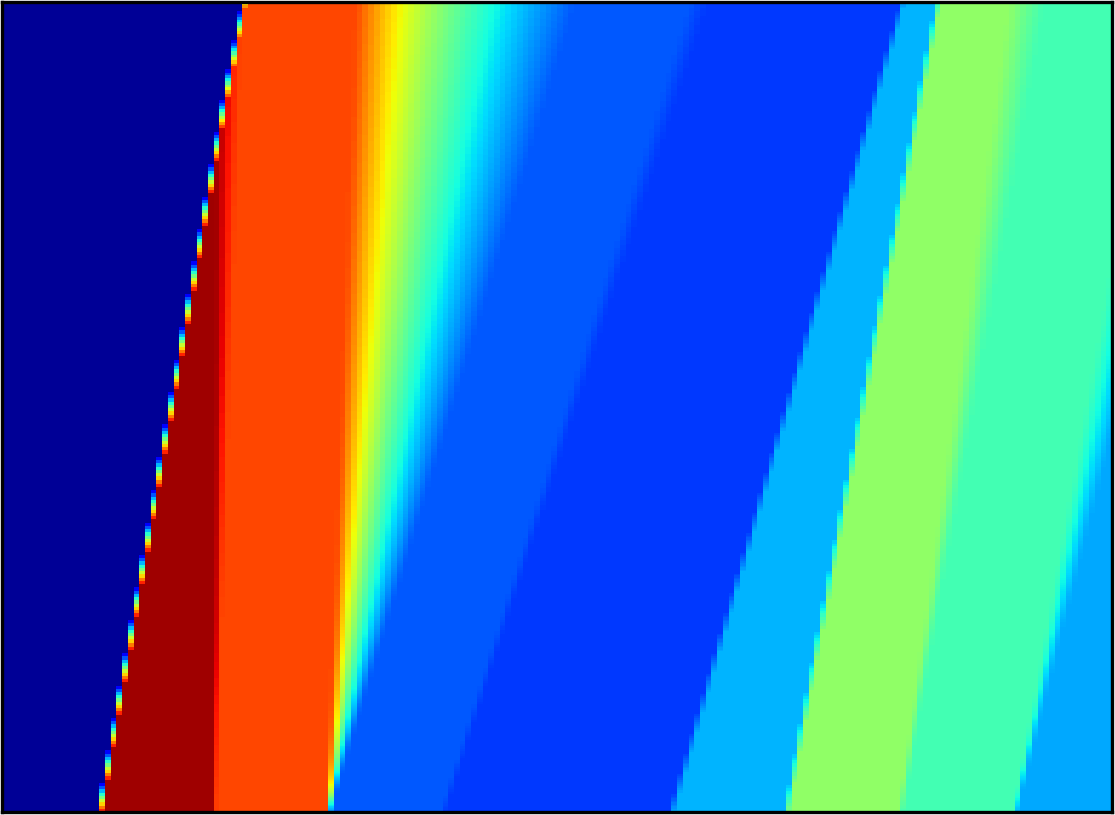}}}\) \\
    
    \(\vcenter{\hbox{\rotatebox{90}{\strut{NFVM$_2^1$}}}}\) & \(\vcenter{\vspace{2pt}\hbox{\includegraphics[width=.15\textwidth]{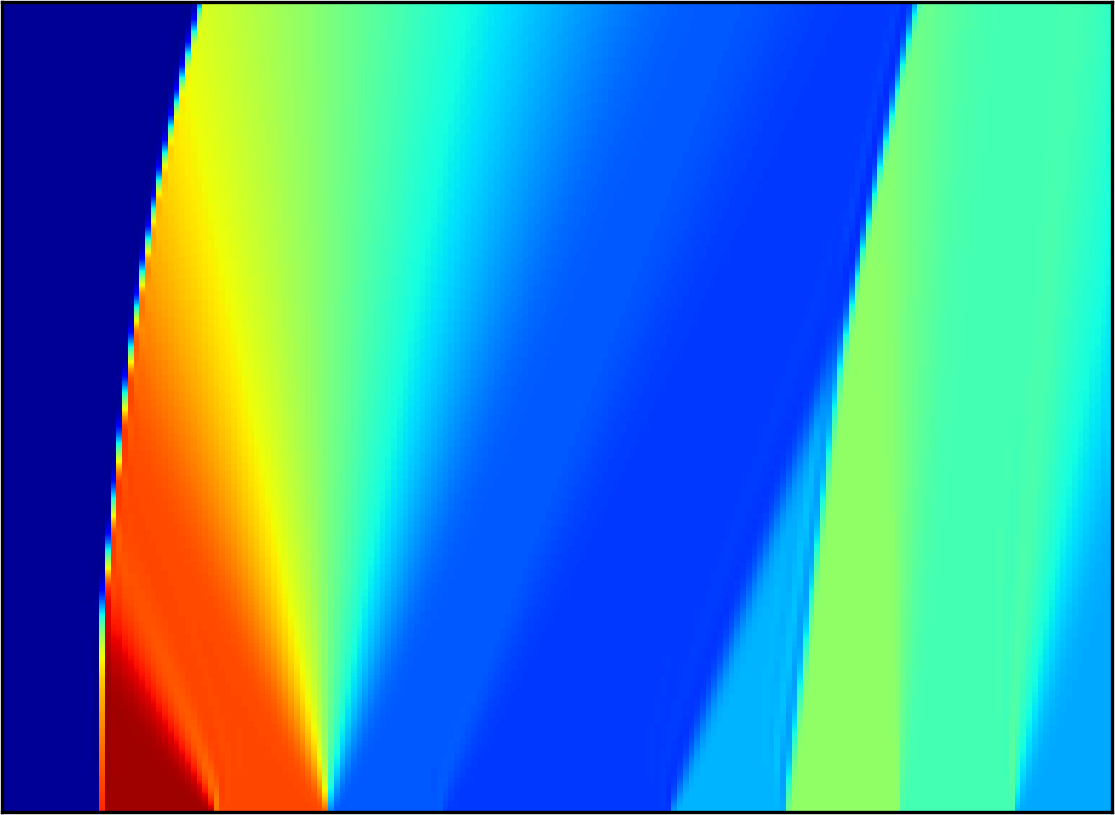}}\vspace{2pt}}\) &
    \(\vcenter{\vspace{2pt}\hbox{\includegraphics[width=.15\textwidth]{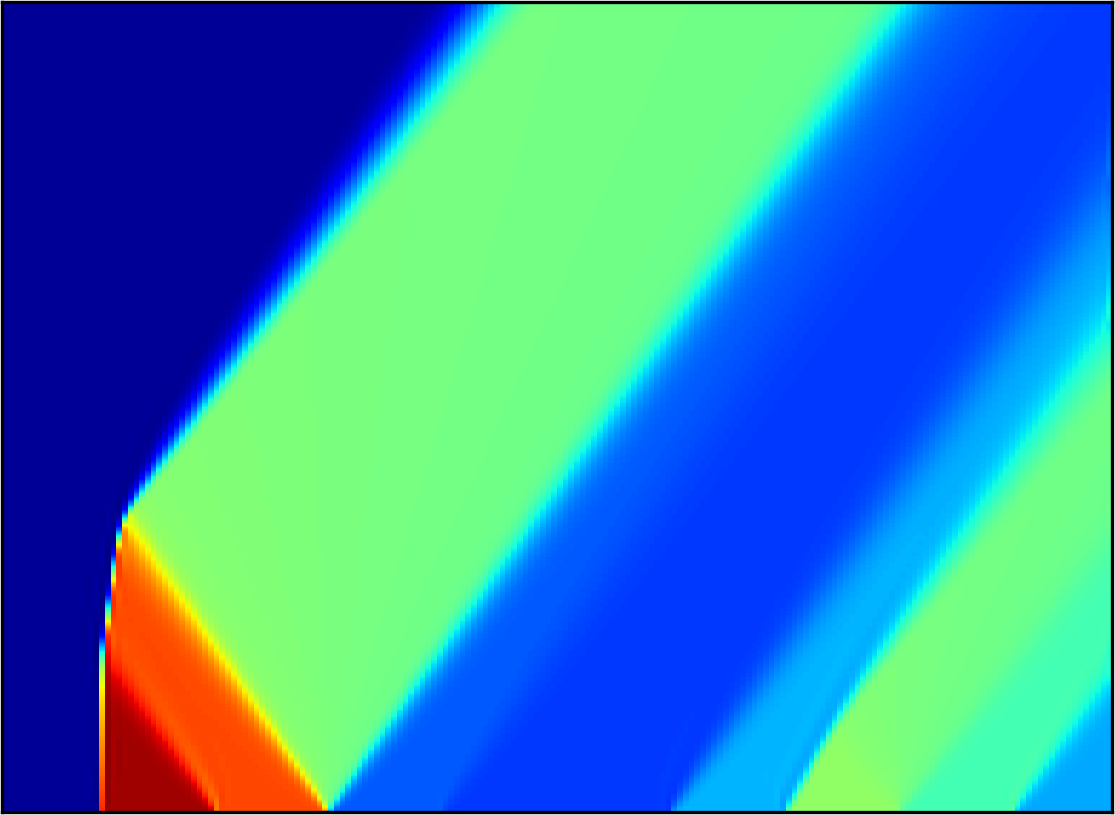}}\vspace{2pt}}\) &
    \(\vcenter{\vspace{2pt}\hbox{\includegraphics[width=.15\textwidth]{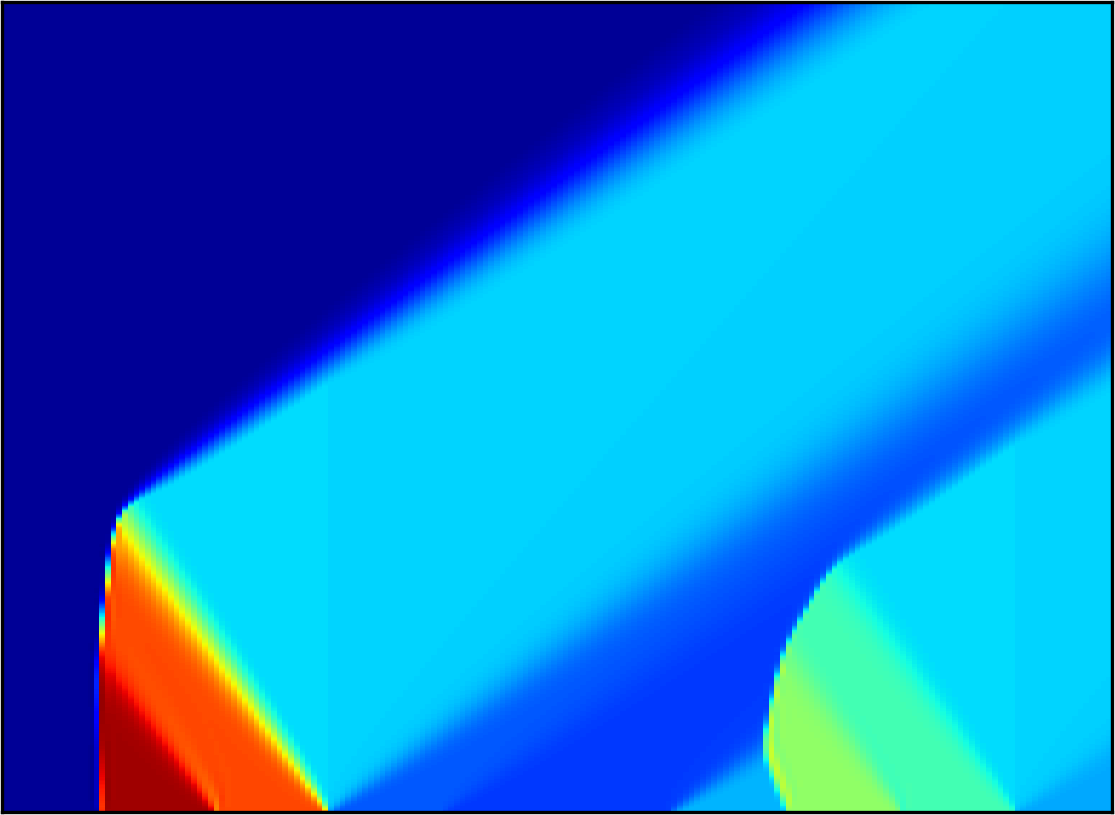}}\vspace{2pt}}\) &
    \(\vcenter{\vspace{2pt}\hbox{\includegraphics[width=.15\textwidth]{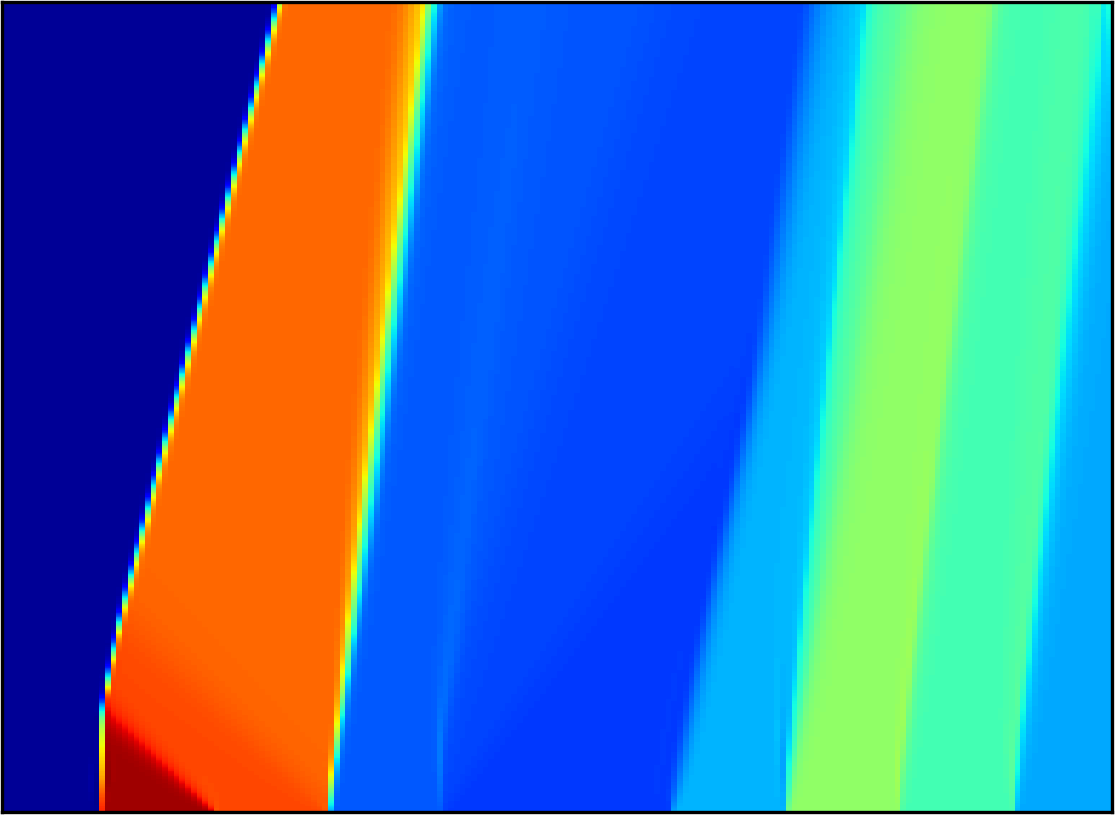}}\vspace{2pt}}\) &
    \(\vcenter{\vspace{2pt}\hbox{\includegraphics[width=.15\textwidth]{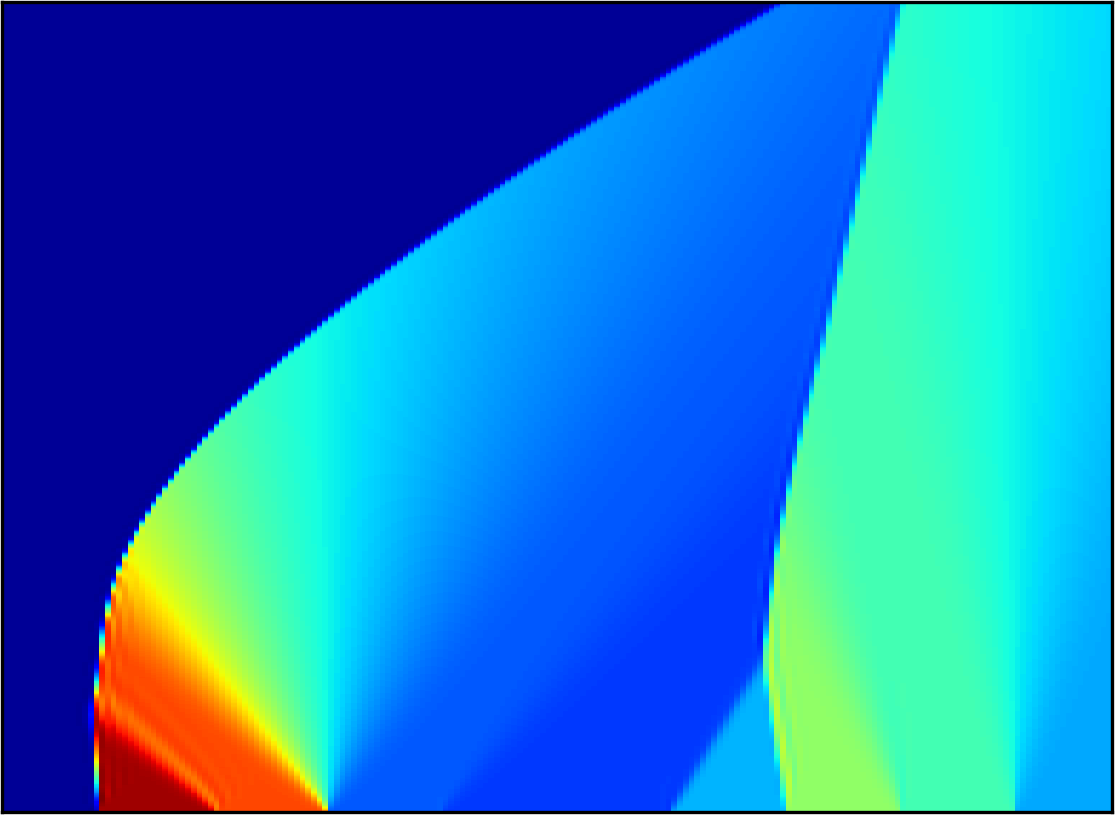}}\vspace{2pt}}\) &
    \(\vcenter{\vspace{2pt}\hbox{\includegraphics[width=.15\textwidth]{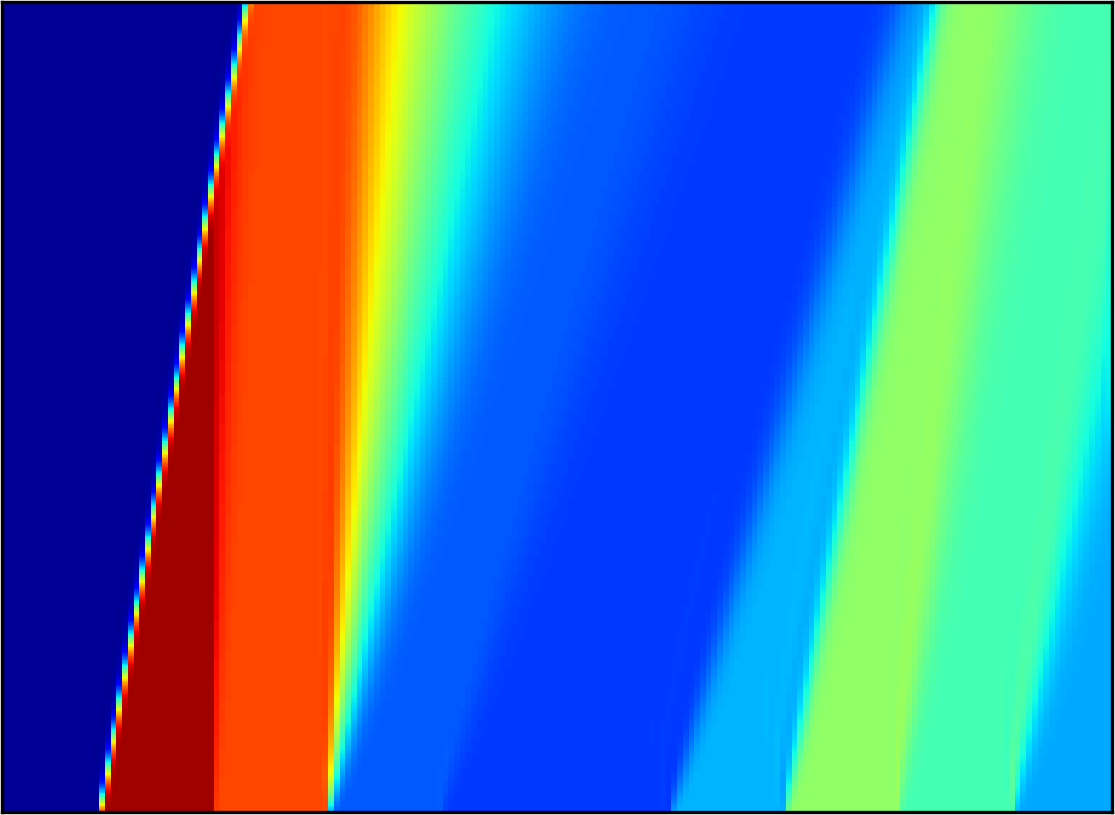}}\vspace{2pt}}\) \\

    \Xhline{1pt}
    
   \(\vcenter{\hbox{\rotatebox{90}{\strut{LxH}}}}\) & \(\vcenter{\vspace{2pt}\hbox{\includegraphics[width=.15\textwidth]{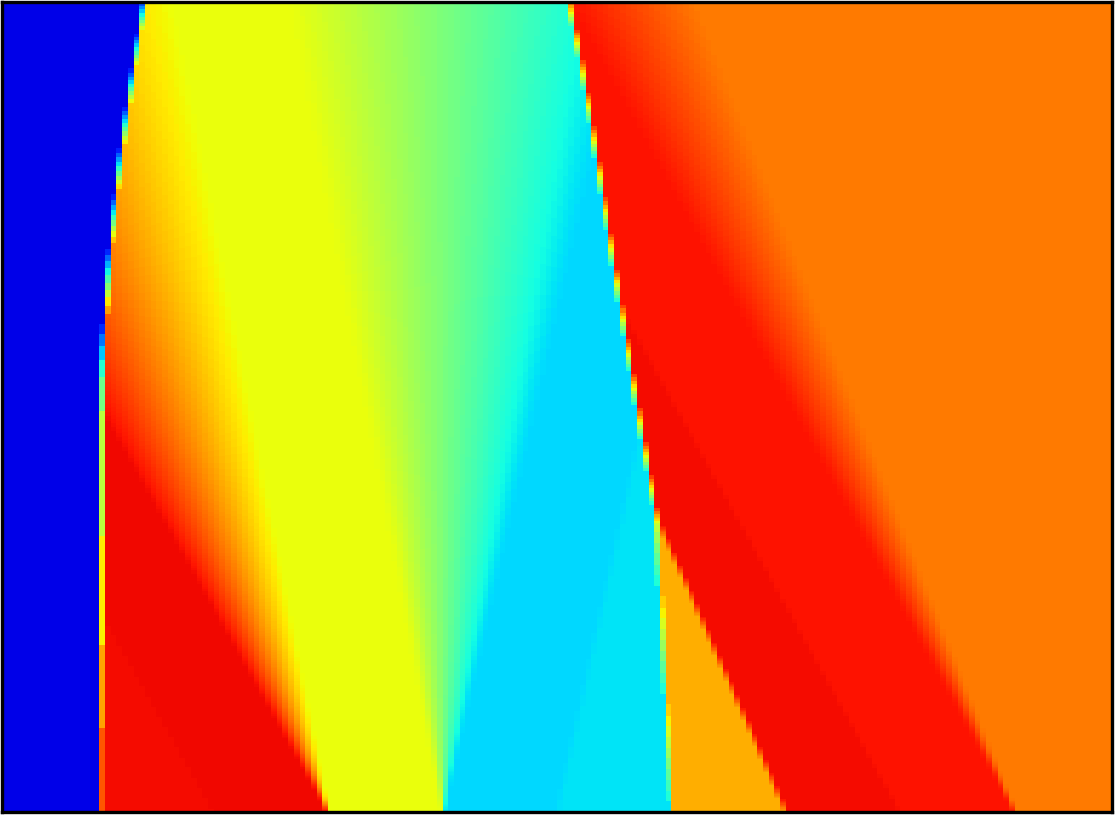}}}\) &
   \(\vcenter{\vspace{2pt}\hbox{\includegraphics[width=.15\textwidth]{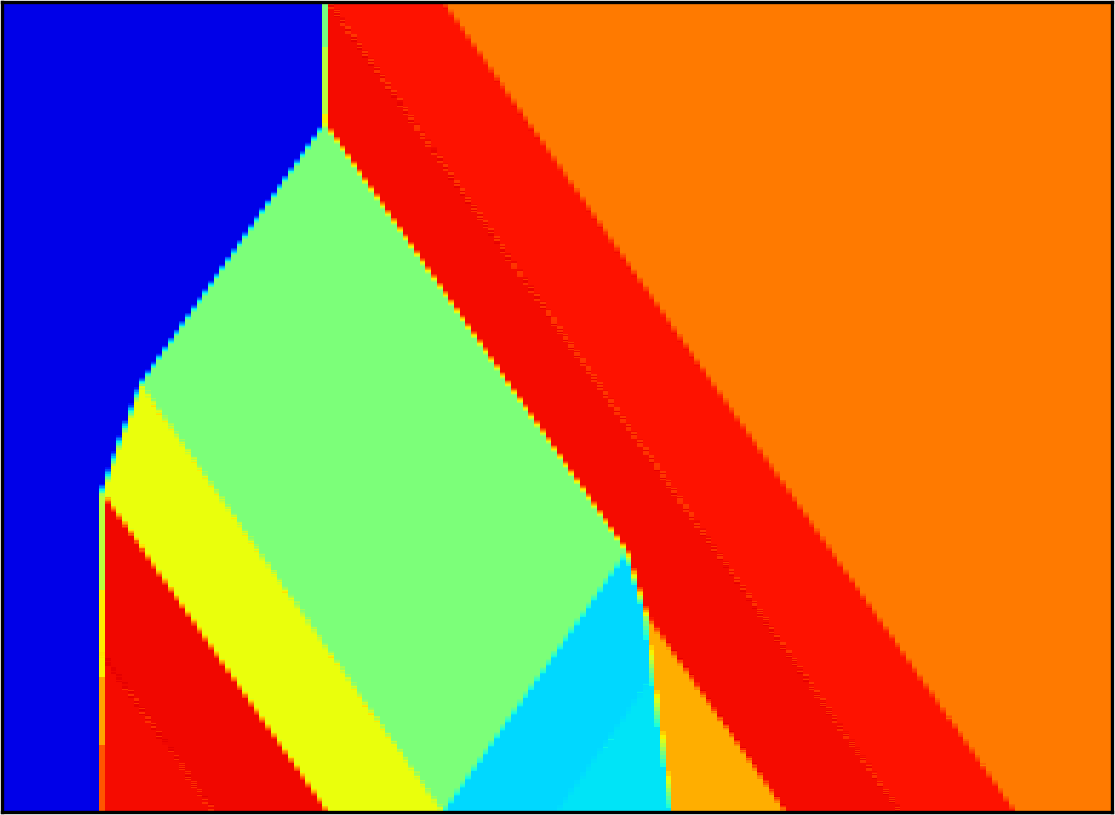}}}\) &
   \(\vcenter{\vspace{2pt}\hbox{\includegraphics[width=.15\textwidth]{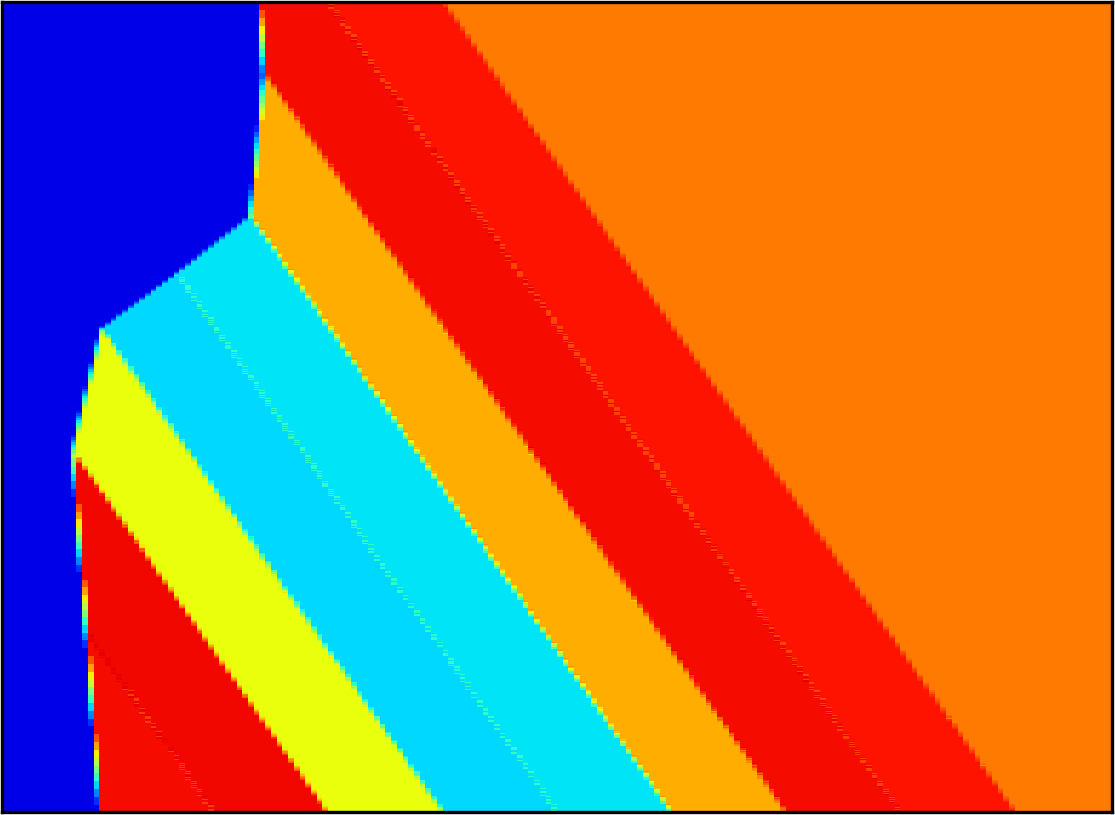}}}\) &
   \(\vcenter{\vspace{2pt}\hbox{\includegraphics[width=.15\textwidth]{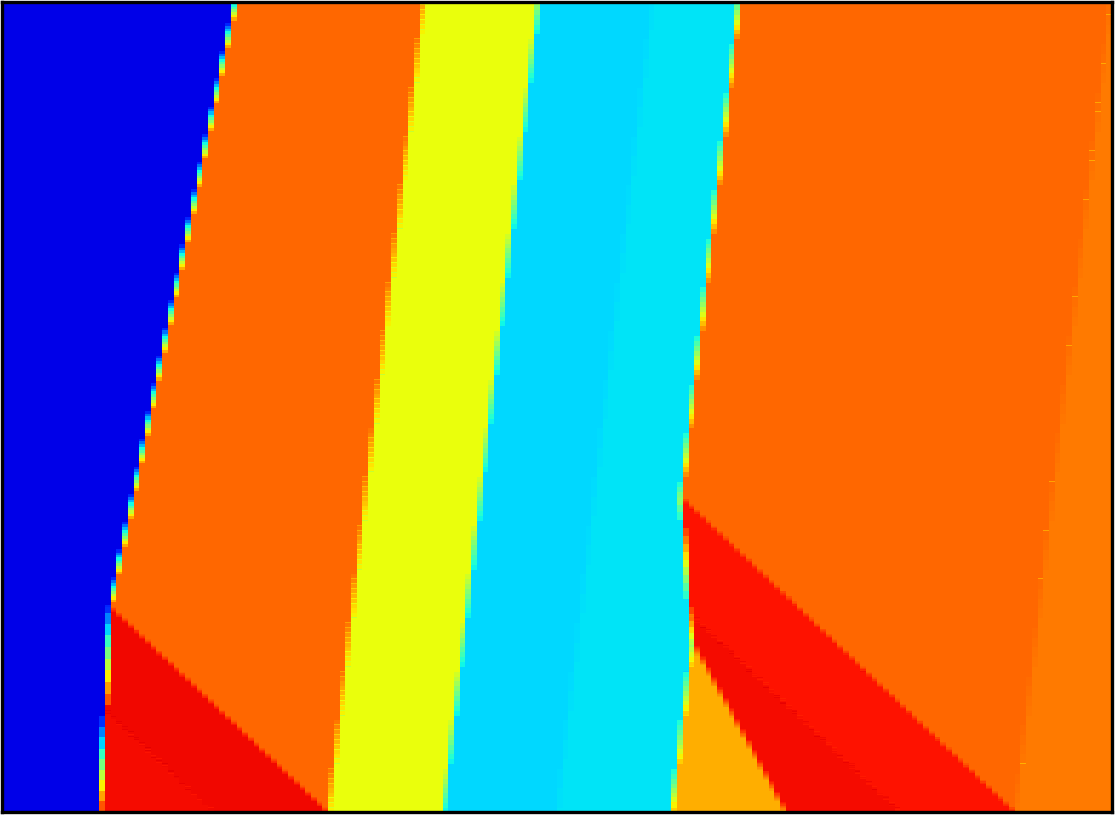}}}\) &
   \(\vcenter{\vspace{2pt}\hbox{\includegraphics[width=.15\textwidth]{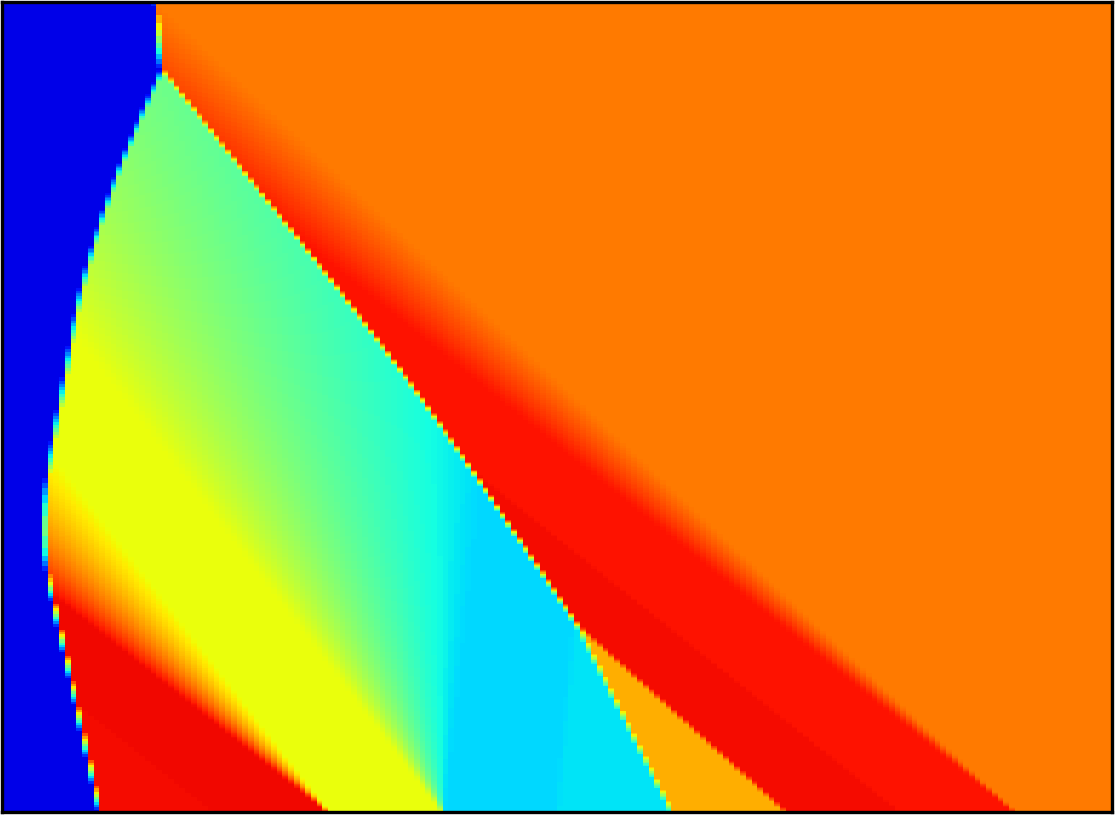}}}\) &
   \(\vcenter{\vspace{2pt}\hbox{\includegraphics[width=.15\textwidth]{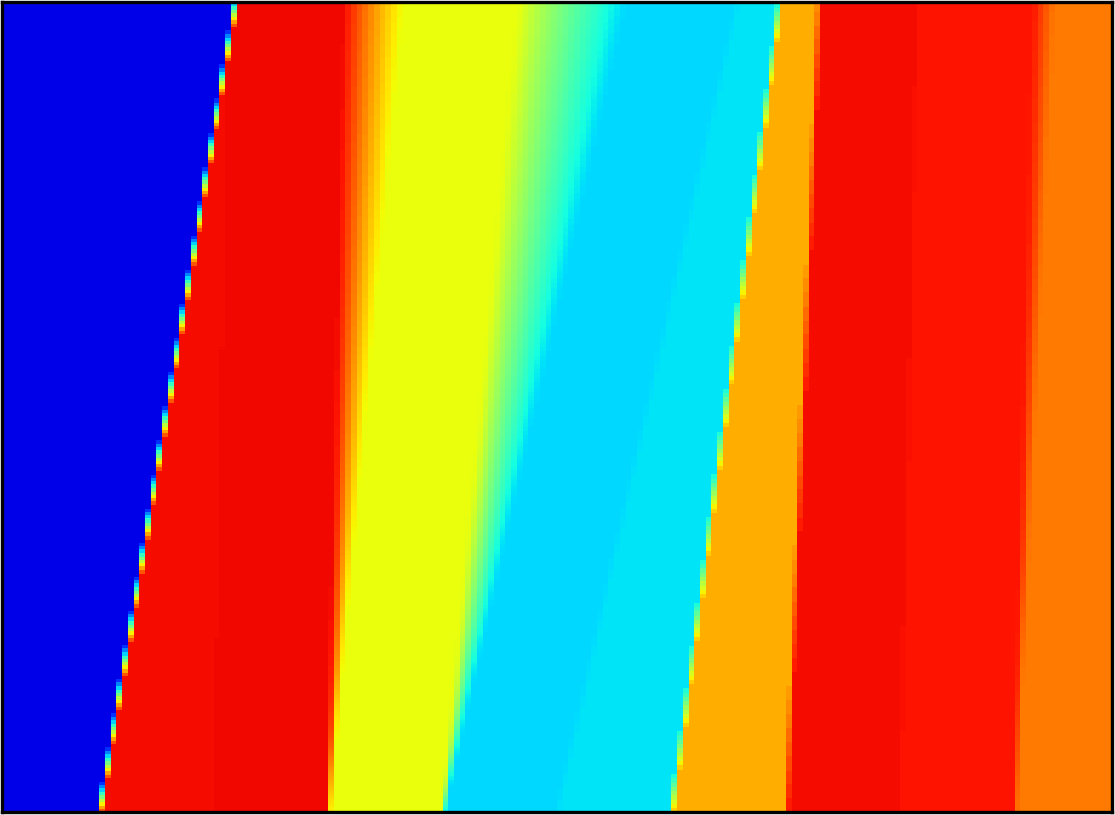}}}\) \\
   
   \(\vcenter{\hbox{\rotatebox{90}{\strut NFVM$_2^1$}}}\) & \(\vcenter{\vspace{2pt}\hbox{\includegraphics[width=.15\textwidth]{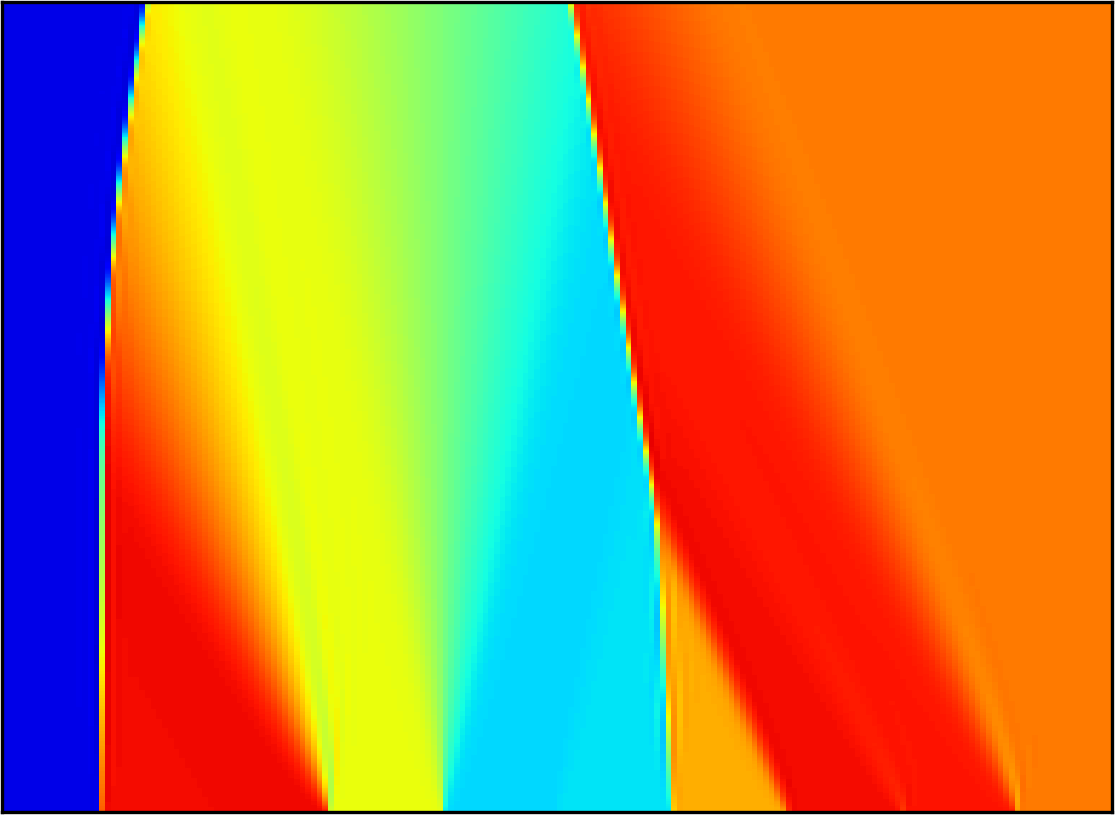}}}\) &
   \(\vcenter{\vspace{2pt}\hbox{\includegraphics[width=.15\textwidth]{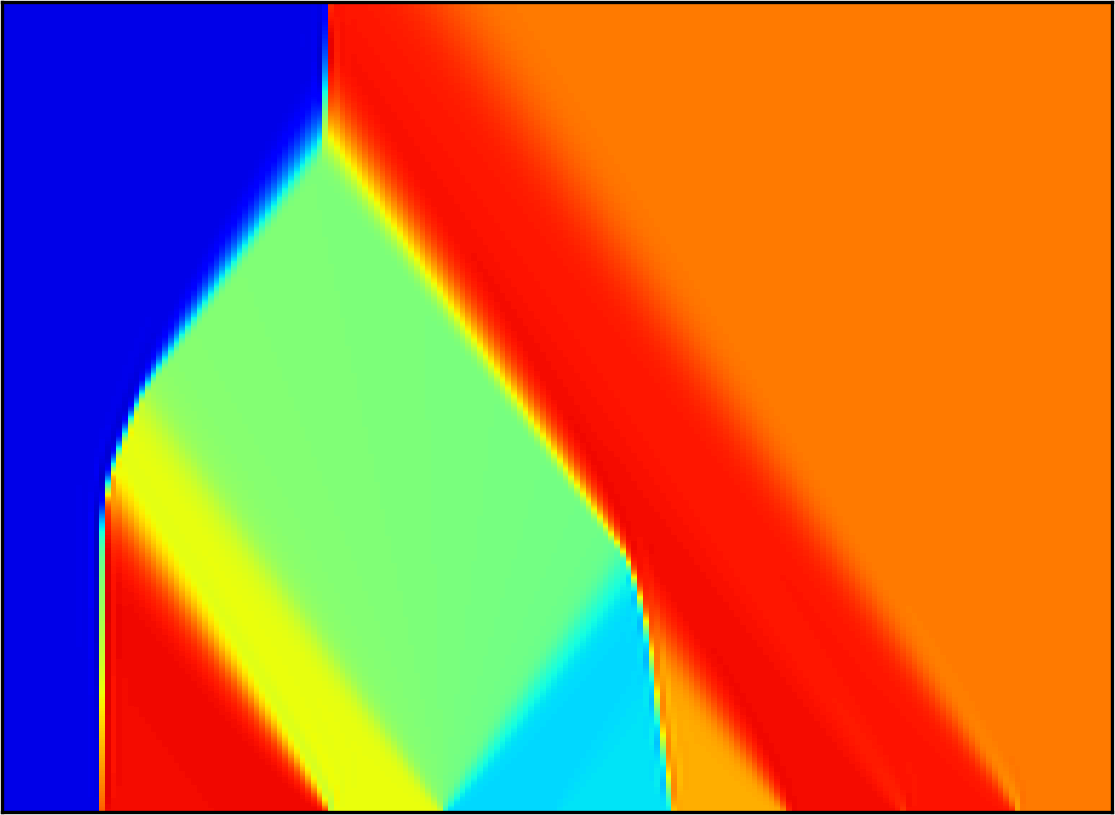}}}\) &
   \(\vcenter{\vspace{2pt}\hbox{\includegraphics[width=.15\textwidth]{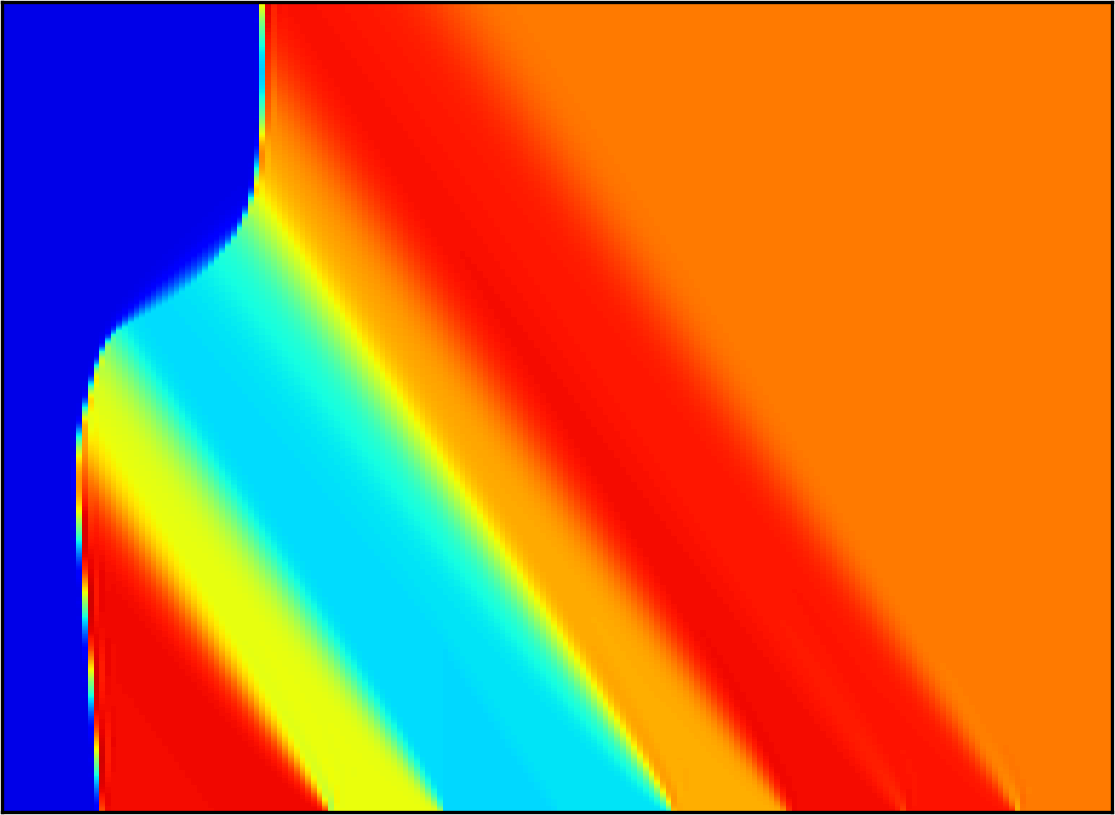}}}\) &
   \(\vcenter{\vspace{2pt}\hbox{\includegraphics[width=.15\textwidth]{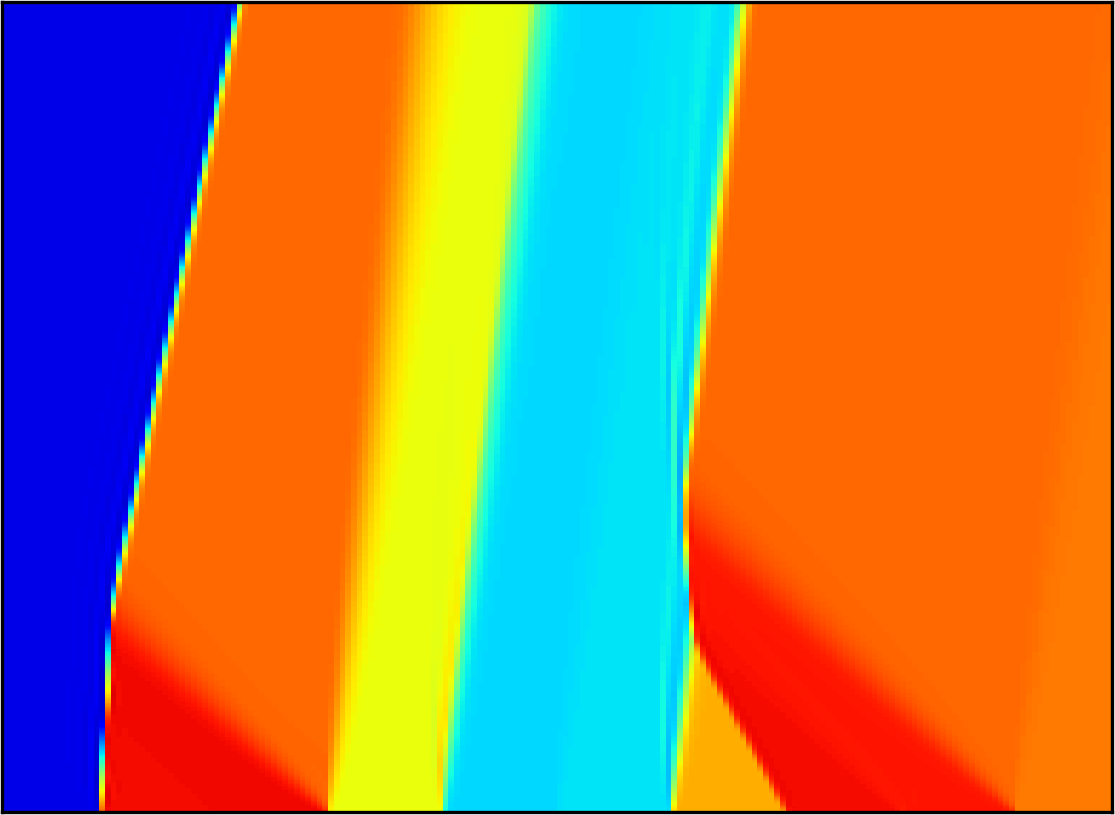}}}\) &
   \(\vcenter{\vspace{2pt}\hbox{\includegraphics[width=.15\textwidth]{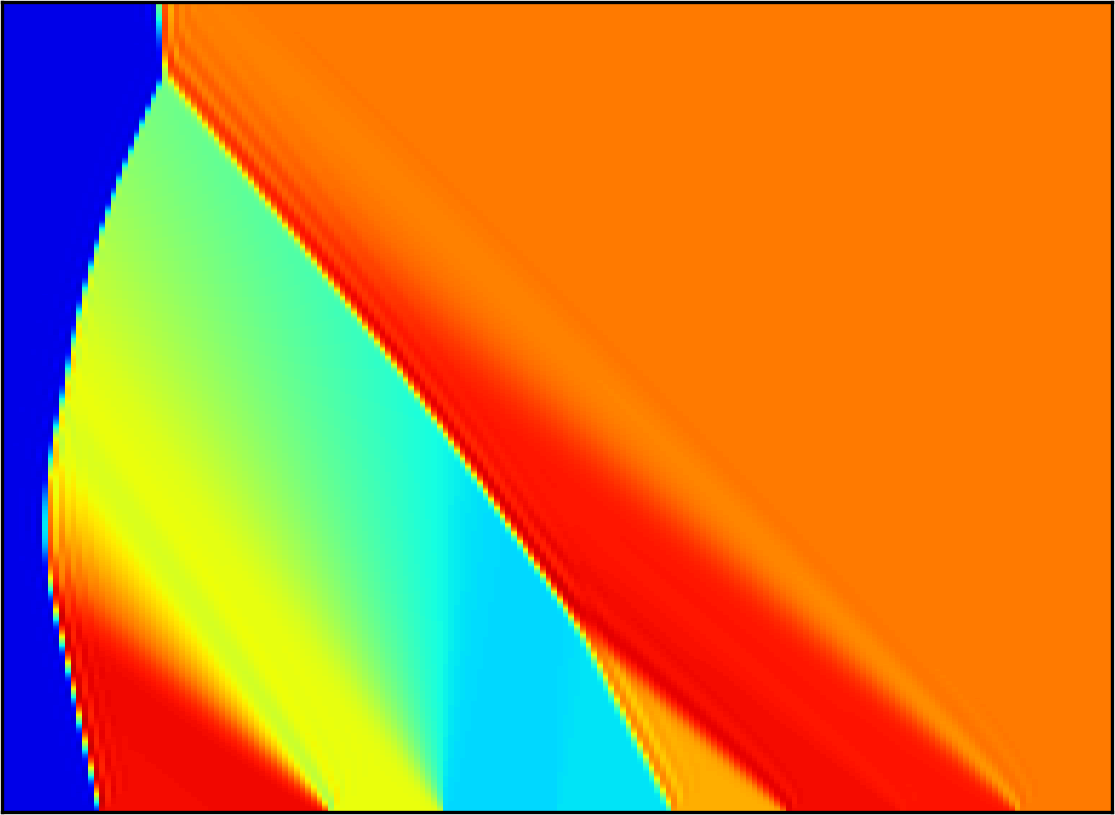}}}\) &
   \(\vcenter{\vspace{2pt}\hbox{\includegraphics[width=.15\textwidth]{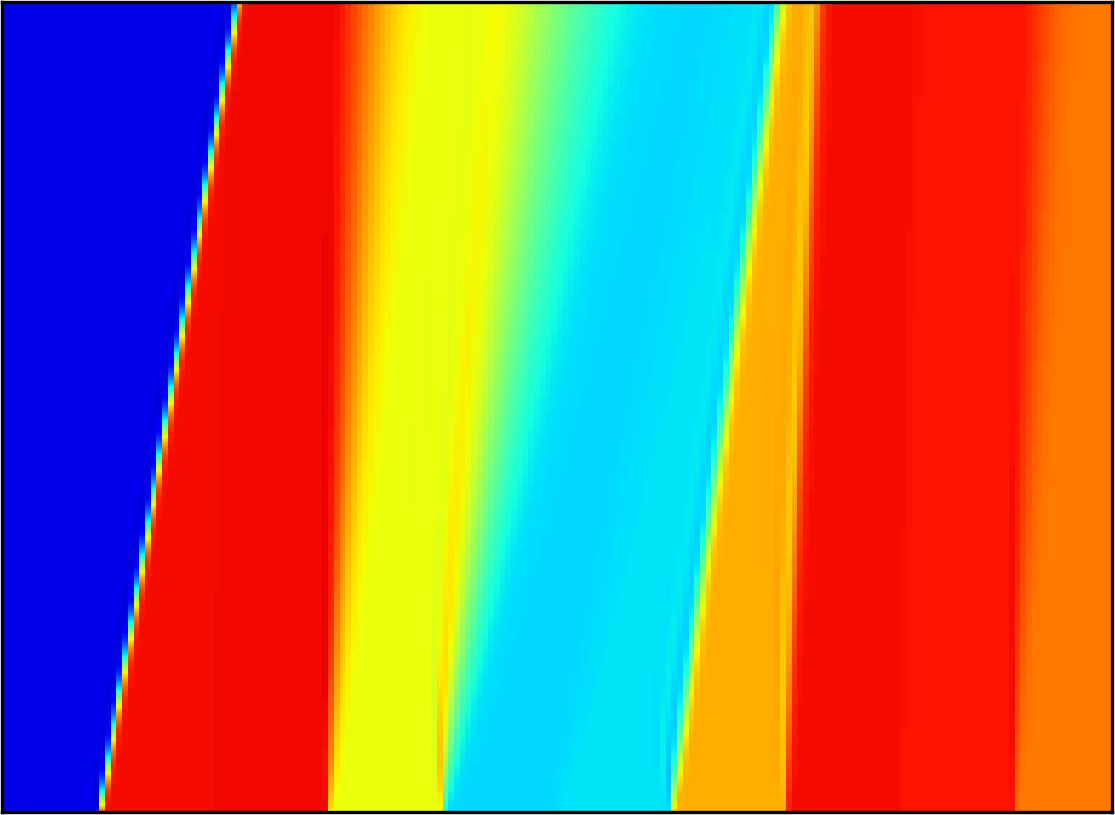}}}\) \\
    
    \end{tabularx}
    \endgroup
    \caption{\small Lax-Hopf ground truth vs. NFVM$_2^1$ prediction on two different initial conditions. Each flow (column) uses a different trained model. See \ref{fig:burgers_heatmaps} for axes and legendt}
    \label{fig:lwr_heatmaps}
\end{figure}

%% file: tex/lwr_metrics.tex
\afterpage{%
\clearpage
\begin{landscape}
\begin{table}
    \centering
    \begin{tabularx}{\linewidth}{cc||cCCCC||CC||C}
    \hline
         \Xhline{1.pt} 
         & & \multicolumn{5}{c||}{\textbf{1\textsuperscript{st} order FVM}} & \multicolumn{2}{c||}{\textbf{Higher order FVM}} & \multicolumn{1}{c}{\textbf{FEM}}\\
         \Xhline{1.pt} 
         & & NFVM$_2^1$ & UNFVM$_2^1$ & GD & LF & EO & ENO & WENO & DG\\
         \Xhline{1.pt} 
        \input{figs/lwr_metrics}
    \end{tabularx}
    \caption{Performance comparison between neural network models and classical numerical schemes. 
    Results are computed over the evaluation set of 1000 initial conditions (see~\Cref{sec:results}). For each method, we report mean and standard deviation for three metrics: L1 error (\(\mean(|u-\hat u|)\)), L2 error (\(\mean((u-\hat u)^2)\)), and relative error (\(\mean(|u-\hat u|/|\max\{\varepsilon, u\}|)\)).}
    \label{tab:lwr_metrics}
\end{table}
\end{landscape}\clearpage}

%% file: figs/lwr_metrics.tex
\multirow{3}{*}{Greenshields'} & L1 & \(\bm{5.3\mathrm{e}^{\shortminus 3}}\)\std{\(8\mathrm{e}^{\shortminus 4}\)} & \(6.3\mathrm{e}^{\shortminus 3}\)\std{\(8\mathrm{e}^{\shortminus 4}\)} & \(1.1\mathrm{e}^{\shortminus 2}\)\std{\(2\mathrm{e}^{\shortminus 3}\)} & \(8.3\mathrm{e}^{\shortminus 2}\)\std{\(1\mathrm{e}^{\shortminus 2}\)} & \(1.1\mathrm{e}^{\shortminus 2}\)\std{\(2\mathrm{e}^{\shortminus 3}\)} & \(9.8\mathrm{e}^{\shortminus 3}\)\std{\(4\mathrm{e}^{\shortminus 3}\)} & \(9.7\mathrm{e}^{\shortminus 3}\)\std{\(4\mathrm{e}^{\shortminus 3}\)} & \(2.1\mathrm{e}^{\shortminus 3}\)\std{\(3\mathrm{e}^{\shortminus 4}\)} \\
 & L2 & \(\bm{1.3\mathrm{e}^{\shortminus 4}}\)\std{\(4\mathrm{e}^{\shortminus 5}\)} & \(2.0\mathrm{e}^{\shortminus 4}\)\std{\(6\mathrm{e}^{\shortminus 5}\)} & \(4.5\mathrm{e}^{\shortminus 4}\)\std{\(2\mathrm{e}^{\shortminus 4}\)} & \(1.3\mathrm{e}^{\shortminus 2}\)\std{\(4\mathrm{e}^{\shortminus 3}\)} & \(4.5\mathrm{e}^{\shortminus 4}\)\std{\(2\mathrm{e}^{\shortminus 4}\)} & \(6.4\mathrm{e}^{\shortminus 4}\)\std{\(4\mathrm{e}^{\shortminus 4}\)} & \(6.4\mathrm{e}^{\shortminus 4}\)\std{\(4\mathrm{e}^{\shortminus 4}\)} & \(3.1\mathrm{e}^{\shortminus 5}\)\std{\(1\mathrm{e}^{\shortminus 5}\)} \\
 & Rel. & \(\bm{2.1\mathrm{e}^{\shortminus 2}}\)\std{\(1\mathrm{e}^{\shortminus 2}\)} & \(2.5\mathrm{e}^{\shortminus 2}\)\std{\(2\mathrm{e}^{\shortminus 2}\)} & \(4.1\mathrm{e}^{\shortminus 2}\)\std{\(3\mathrm{e}^{\shortminus 2}\)} & \(3.6\mathrm{e}^{\shortminus 1}\)\std{\(2\mathrm{e}^{\shortminus 1}\)} & \(4.2\mathrm{e}^{\shortminus 2}\)\std{\(3\mathrm{e}^{\shortminus 2}\)} & \(1.9\mathrm{e}^{\shortminus 2}\)\std{\(6\mathrm{e}^{\shortminus 3}\)} & \(1.9\mathrm{e}^{\shortminus 2}\)\std{\(5\mathrm{e}^{\shortminus 3}\)} & \(7.4\mathrm{e}^{\shortminus 3}\)\std{\(4\mathrm{e}^{\shortminus 3}\)} \\
	\Xhline{1.pt}
\multirow{3}{*}{Triangular Sym} & L1 & \(\bm{1.7\mathrm{e}^{\shortminus 2}}\)\std{\(5\mathrm{e}^{\shortminus 3}\)} & \(1.9\mathrm{e}^{\shortminus 2}\)\std{\(6\mathrm{e}^{\shortminus 3}\)} & \(2.2\mathrm{e}^{\shortminus 2}\)\std{\(7\mathrm{e}^{\shortminus 3}\)} & \(6.1\mathrm{e}^{\shortminus 2}\)\std{\(2\mathrm{e}^{\shortminus 2}\)} & \(2.2\mathrm{e}^{\shortminus 2}\)\std{\(7\mathrm{e}^{\shortminus 3}\)} & \(1.8\mathrm{e}^{\shortminus 2}\)\std{\(1\mathrm{e}^{\shortminus 2}\)} & \(1.7\mathrm{e}^{\shortminus 2}\)\std{\(1\mathrm{e}^{\shortminus 2}\)} & \(4.5\mathrm{e}^{\shortminus 3}\)\std{\(1\mathrm{e}^{\shortminus 3}\)} \\
 & L2 & \(\bm{1.4\mathrm{e}^{\shortminus 3}}\)\std{\(6\mathrm{e}^{\shortminus 4}\)} & \(1.9\mathrm{e}^{\shortminus 3}\)\std{\(9\mathrm{e}^{\shortminus 4}\)} & \(2.3\mathrm{e}^{\shortminus 3}\)\std{\(1\mathrm{e}^{\shortminus 3}\)} & \(9.6\mathrm{e}^{\shortminus 3}\)\std{\(4\mathrm{e}^{\shortminus 3}\)} & \(2.3\mathrm{e}^{\shortminus 3}\)\std{\(1\mathrm{e}^{\shortminus 3}\)} & \(2.0\mathrm{e}^{\shortminus 3}\)\std{\(2\mathrm{e}^{\shortminus 3}\)} & \(1.9\mathrm{e}^{\shortminus 3}\)\std{\(2\mathrm{e}^{\shortminus 3}\)} & \(2.6\mathrm{e}^{\shortminus 4}\)\std{\(1\mathrm{e}^{\shortminus 4}\)} \\
 & Rel. & \(\bm{9.7\mathrm{e}^{\shortminus 2}}\)\std{\(1\mathrm{e}^{\shortminus 1}\)} & \(1.1\mathrm{e}^{\shortminus 1}\)\std{\(1\mathrm{e}^{\shortminus 1}\)} & \(1.4\mathrm{e}^{\shortminus 1}\)\std{\(2\mathrm{e}^{\shortminus 1}\)} & \(4.1\mathrm{e}^{\shortminus 1}\)\std{\(5\mathrm{e}^{\shortminus 1}\)} & \(1.4\mathrm{e}^{\shortminus 1}\)\std{\(2\mathrm{e}^{\shortminus 1}\)} & \(5.2\mathrm{e}^{\shortminus 2}\)\std{\(3\mathrm{e}^{\shortminus 2}\)} & \(4.5\mathrm{e}^{\shortminus 2}\)\std{\(3\mathrm{e}^{\shortminus 2}\)} & \(2.5\mathrm{e}^{\shortminus 2}\)\std{\(3\mathrm{e}^{\shortminus 2}\)} \\
	\Xhline{1.pt}
\multirow{3}{*}{Triangular Skw} & L1 & \(\bm{2.2\mathrm{e}^{\shortminus 2}}\)\std{\(9\mathrm{e}^{\shortminus 3}\)} & \(2.6\mathrm{e}^{\shortminus 2}\)\std{\(1\mathrm{e}^{\shortminus 2}\)} & \(3.0\mathrm{e}^{\shortminus 2}\)\std{\(1\mathrm{e}^{\shortminus 2}\)} & \(7.6\mathrm{e}^{\shortminus 2}\)\std{\(3\mathrm{e}^{\shortminus 2}\)} & \(3.0\mathrm{e}^{\shortminus 2}\)\std{\(1\mathrm{e}^{\shortminus 2}\)} & \(4.0\mathrm{e}^{\shortminus 2}\)\std{\(2\mathrm{e}^{\shortminus 2}\)} & \(4.0\mathrm{e}^{\shortminus 2}\)\std{\(2\mathrm{e}^{\shortminus 2}\)} & \(5.6\mathrm{e}^{\shortminus 3}\)\std{\(2\mathrm{e}^{\shortminus 3}\)} \\
 & L2 & \(\bm{2.4\mathrm{e}^{\shortminus 3}}\)\std{\(1\mathrm{e}^{\shortminus 3}\)} & \(3.1\mathrm{e}^{\shortminus 3}\)\std{\(2\mathrm{e}^{\shortminus 3}\)} & \(3.8\mathrm{e}^{\shortminus 3}\)\std{\(2\mathrm{e}^{\shortminus 3}\)} & \(1.4\mathrm{e}^{\shortminus 2}\)\std{\(8\mathrm{e}^{\shortminus 3}\)} & \(3.8\mathrm{e}^{\shortminus 3}\)\std{\(2\mathrm{e}^{\shortminus 3}\)} & \(5.8\mathrm{e}^{\shortminus 3}\)\std{\(4\mathrm{e}^{\shortminus 3}\)} & \(5.8\mathrm{e}^{\shortminus 3}\)\std{\(4\mathrm{e}^{\shortminus 3}\)} & \(4.1\mathrm{e}^{\shortminus 4}\)\std{\(2\mathrm{e}^{\shortminus 4}\)} \\
 & Rel. & \(\bm{6.2\mathrm{e}^{\shortminus 2}}\)\std{\(4\mathrm{e}^{\shortminus 2}\)} & \(6.8\mathrm{e}^{\shortminus 2}\)\std{\(4\mathrm{e}^{\shortminus 2}\)} & \(7.7\mathrm{e}^{\shortminus 2}\)\std{\(4\mathrm{e}^{\shortminus 2}\)} & \(2.2\mathrm{e}^{\shortminus 1}\)\std{\(1\mathrm{e}^{\shortminus 1}\)} & \(7.8\mathrm{e}^{\shortminus 2}\)\std{\(4\mathrm{e}^{\shortminus 2}\)} & \(8.5\mathrm{e}^{\shortminus 2}\)\std{\(3\mathrm{e}^{\shortminus 2}\)} & \(8.5\mathrm{e}^{\shortminus 2}\)\std{\(4\mathrm{e}^{\shortminus 2}\)} & \(1.5\mathrm{e}^{\shortminus 2}\)\std{\(8\mathrm{e}^{\shortminus 3}\)} \\
	\Xhline{1.pt}
\multirow{3}{*}{Trapezoidal} & L1 & \(\bm{1.4\mathrm{e}^{\shortminus 2}}\)\std{\(2\mathrm{e}^{\shortminus 3}\)} & \(1.6\mathrm{e}^{\shortminus 2}\)\std{\(4\mathrm{e}^{\shortminus 3}\)} & \(2.0\mathrm{e}^{\shortminus 2}\)\std{\(5\mathrm{e}^{\shortminus 3}\)} & \(1.2\mathrm{e}^{\shortminus 1}\)\std{\(3\mathrm{e}^{\shortminus 2}\)} & \(2.0\mathrm{e}^{\shortminus 2}\)\std{\(5\mathrm{e}^{\shortminus 3}\)} & \(8.5\mathrm{e}^{\shortminus 3}\)\std{\(2\mathrm{e}^{\shortminus 3}\)} & \(7.6\mathrm{e}^{\shortminus 3}\)\std{\(2\mathrm{e}^{\shortminus 3}\)} & \(4.9\mathrm{e}^{\shortminus 3}\)\std{\(1\mathrm{e}^{\shortminus 3}\)} \\
 & L2 & \(\bm{1.1\mathrm{e}^{\shortminus 3}}\)\std{\(4\mathrm{e}^{\shortminus 4}\)} & \(1.6\mathrm{e}^{\shortminus 3}\)\std{\(7\mathrm{e}^{\shortminus 4}\)} & \(2.1\mathrm{e}^{\shortminus 3}\)\std{\(8\mathrm{e}^{\shortminus 4}\)} & \(2.5\mathrm{e}^{\shortminus 2}\)\std{\(1\mathrm{e}^{\shortminus 2}\)} & \(2.1\mathrm{e}^{\shortminus 3}\)\std{\(8\mathrm{e}^{\shortminus 4}\)} & \(6.2\mathrm{e}^{\shortminus 4}\)\std{\(2\mathrm{e}^{\shortminus 4}\)} & \(5.3\mathrm{e}^{\shortminus 4}\)\std{\(2\mathrm{e}^{\shortminus 4}\)} & \(2.9\mathrm{e}^{\shortminus 4}\)\std{\(1\mathrm{e}^{\shortminus 4}\)} \\
 & Rel. & \(\bm{4.6\mathrm{e}^{\shortminus 2}}\)\std{\(2\mathrm{e}^{\shortminus 2}\)} & \(5.6\mathrm{e}^{\shortminus 2}\)\std{\(3\mathrm{e}^{\shortminus 2}\)} & \(7.3\mathrm{e}^{\shortminus 2}\)\std{\(4\mathrm{e}^{\shortminus 2}\)} & \(4.9\mathrm{e}^{\shortminus 1}\)\std{\(3\mathrm{e}^{\shortminus 1}\)} & \(7.3\mathrm{e}^{\shortminus 2}\)\std{\(4\mathrm{e}^{\shortminus 2}\)} & \(2.6\mathrm{e}^{\shortminus 2}\)\std{\(1\mathrm{e}^{\shortminus 2}\)} & \(2.2\mathrm{e}^{\shortminus 2}\)\std{\(8\mathrm{e}^{\shortminus 3}\)} & \(1.6\mathrm{e}^{\shortminus 2}\)\std{\(8\mathrm{e}^{\shortminus 3}\)} \\
	\Xhline{1.pt}
\multirow{3}{*}{Greenberg} & L1 & \(\bm{4.0\mathrm{e}^{\shortminus 3}}\)\std{\(1\mathrm{e}^{\shortminus 3}\)} & \(7.3\mathrm{e}^{\shortminus 3}\)\std{\(1\mathrm{e}^{\shortminus 3}\)} & \(9.3\mathrm{e}^{\shortminus 3}\)\std{\(2\mathrm{e}^{\shortminus 3}\)} & \(4.5\mathrm{e}^{\shortminus 2}\)\std{\(9\mathrm{e}^{\shortminus 3}\)} & \(9.4\mathrm{e}^{\shortminus 3}\)\std{\(2\mathrm{e}^{\shortminus 3}\)} & \(1.4\mathrm{e}^{\shortminus 2}\)\std{\(5\mathrm{e}^{\shortminus 3}\)} & \(1.4\mathrm{e}^{\shortminus 2}\)\std{\(5\mathrm{e}^{\shortminus 3}\)} & \(2.3\mathrm{e}^{\shortminus 3}\)\std{\(4\mathrm{e}^{\shortminus 3}\)} \\
 & L2 & \(\bm{1.4\mathrm{e}^{\shortminus 4}}\)\std{\(9\mathrm{e}^{\shortminus 5}\)} & \(3.8\mathrm{e}^{\shortminus 4}\)\std{\(2\mathrm{e}^{\shortminus 4}\)} & \(4.9\mathrm{e}^{\shortminus 4}\)\std{\(2\mathrm{e}^{\shortminus 4}\)} & \(5.3\mathrm{e}^{\shortminus 3}\)\std{\(2\mathrm{e}^{\shortminus 3}\)} & \(4.9\mathrm{e}^{\shortminus 4}\)\std{\(2\mathrm{e}^{\shortminus 4}\)} & \(1.1\mathrm{e}^{\shortminus 3}\)\std{\(6\mathrm{e}^{\shortminus 4}\)} & \(1.2\mathrm{e}^{\shortminus 3}\)\std{\(9\mathrm{e}^{\shortminus 4}\)} & \(3.4\mathrm{e}^{\shortminus 4}\)\std{\(2\mathrm{e}^{\shortminus 3}\)} \\
 & Rel. & \(\bm{8.4\mathrm{e}^{\shortminus 3}}\)\std{\(1\mathrm{e}^{\shortminus 3}\)} & \(1.8\mathrm{e}^{\shortminus 2}\)\std{\(4\mathrm{e}^{\shortminus 3}\)} & \(2.0\mathrm{e}^{\shortminus 2}\)\std{\(4\mathrm{e}^{\shortminus 3}\)} & \(1.2\mathrm{e}^{\shortminus 1}\)\std{\(3\mathrm{e}^{\shortminus 2}\)} & \(2.1\mathrm{e}^{\shortminus 2}\)\std{\(4\mathrm{e}^{\shortminus 3}\)} & \(3.0\mathrm{e}^{\shortminus 2}\)\std{\(7\mathrm{e}^{\shortminus 2}\)} & \(4.9\mathrm{e}^{\shortminus 2}\)\std{\(2\mathrm{e}^{\shortminus 1}\)} & \(5.9\mathrm{e}^{\shortminus 2}\)\std{\(4\mathrm{e}^{\shortminus 1}\)} \\
	\Xhline{1.pt}
\multirow{3}{*}{Underwood} & L1 & \(\bm{9.3\mathrm{e}^{\shortminus 3}}\)\std{\(2\mathrm{e}^{\shortminus 3}\)} & \(1.4\mathrm{e}^{\shortminus 2}\)\std{\(3\mathrm{e}^{\shortminus 3}\)} & \(1.6\mathrm{e}^{\shortminus 2}\)\std{\(3\mathrm{e}^{\shortminus 3}\)} & \(1.3\mathrm{e}^{\shortminus 1}\)\std{\(3\mathrm{e}^{\shortminus 2}\)} & \(1.6\mathrm{e}^{\shortminus 2}\)\std{\(3\mathrm{e}^{\shortminus 3}\)} & \(4.0\mathrm{e}^{\shortminus 3}\)\std{\(6\mathrm{e}^{\shortminus 4}\)} & \(3.6\mathrm{e}^{\shortminus 3}\)\std{\(5\mathrm{e}^{\shortminus 4}\)} & \(3.2\mathrm{e}^{\shortminus 3}\)\std{\(6\mathrm{e}^{\shortminus 4}\)} \\
 & L2 & \(\bm{3.8\mathrm{e}^{\shortminus 4}}\)\std{\(1\mathrm{e}^{\shortminus 4}\)} & \(6.9\mathrm{e}^{\shortminus 4}\)\std{\(2\mathrm{e}^{\shortminus 4}\)} & \(9.2\mathrm{e}^{\shortminus 4}\)\std{\(3\mathrm{e}^{\shortminus 4}\)} & \(2.7\mathrm{e}^{\shortminus 2}\)\std{\(1\mathrm{e}^{\shortminus 2}\)} & \(9.2\mathrm{e}^{\shortminus 4}\)\std{\(3\mathrm{e}^{\shortminus 4}\)} & \(1.1\mathrm{e}^{\shortminus 4}\)\std{\(3\mathrm{e}^{\shortminus 5}\)} & \(9.8\mathrm{e}^{\shortminus 5}\)\std{\(2\mathrm{e}^{\shortminus 5}\)} & \(5.9\mathrm{e}^{\shortminus 5}\)\std{\(2\mathrm{e}^{\shortminus 5}\)} \\
 & Rel. & \(\bm{5.0\mathrm{e}^{\shortminus 2}}\)\std{\(4\mathrm{e}^{\shortminus 2}\)} & \(6.7\mathrm{e}^{\shortminus 2}\)\std{\(5\mathrm{e}^{\shortminus 2}\)} & \(8.4\mathrm{e}^{\shortminus 2}\)\std{\(6\mathrm{e}^{\shortminus 2}\)} & \(6.6\mathrm{e}^{\shortminus 1}\)\std{\(5\mathrm{e}^{\shortminus 1}\)} & \(8.4\mathrm{e}^{\shortminus 2}\)\std{\(6\mathrm{e}^{\shortminus 2}\)} & \(1.8\mathrm{e}^{\shortminus 2}\)\std{\(1\mathrm{e}^{\shortminus 2}\)} & \(1.6\mathrm{e}^{\shortminus 2}\)\std{\(9\mathrm{e}^{\shortminus 3}\)} & \(1.2\mathrm{e}^{\shortminus 2}\)\std{\(7\mathrm{e}^{\shortminus 3}\)} \\
	\Xhline{1.pt}

%% file: tex/lwr_l1_error_bar_plots.tex
\begin{figure}[H]
    \centering
    \begin{subfigure}[t]{.35\textwidth}
        \begin{tikzpicture}
            \begin{axis}[
                ymin=0,
                tick label style={font=\scriptsize},
                xticklabel style={rotate=-45, anchor=west, xshift=-5pt, yshift=-5pt},
                ylabel={\(L_1\) error},
                xtick={1, 2, 3, 4, 5, 6, 7},
                height=4cm,
                xticklabels={NFVM$_2^1$, UNFVM$_2^1$, GD, EO, ENO, WENO, DG},
                width=1\textwidth,
            ]
            \addplot[darkmagenta, 
                    only marks, 
                    mark=x,
                    mark options={ultra thick, scale=1.5},
                    error bars/.cd,
                    y dir=both, 
                    y explicit,
                    error bar style={very thick, color=pink}
            ] coordinates { 
            (1, 5.3e-3) +- (8e-4, 8e-4)
            (2, 6.3e-3) +- (8e-4, 8e-4)
            (3, 1.1e-2) +- (2e-3, 2e-3)
            (4, 1.1e-2) +- (2e-3, 2e-3)
            (5, 9.8e-3) +- (4e-3, 4e-3)
            (6, 9.7e-3) +- (4e-3, 4e-3)
            (7, 2.1e-3) +- (3e-4, 3e-4)
            };

            \draw[thick, dashed] (axis cs:4.5,\pgfkeysvalueof{/pgfplots/ymin}) -- (axis cs:4.5,\pgfkeysvalueof{/pgfplots/ymax});
            \draw[thick, dashed] (axis cs:6.5,\pgfkeysvalueof{/pgfplots/ymin}) -- (axis cs:6.5,\pgfkeysvalueof{/pgfplots/ymax});
            \end{axis}
        \end{tikzpicture}
        \vspace{-0.2cm}
        \caption{Greenshields'.}
    \end{subfigure}
    \begin{subfigure}[t]{.31\textwidth}
    \centering
        \begin{tikzpicture}
            \begin{axis}[
                ymin=0,
                tick label style={font=\scriptsize},
                xticklabel style={rotate=-45, anchor=west, xshift=-5pt, yshift=-5pt},
                ylabel style={yshift=-20pt},
                xtick={1, 2, 3, 4, 5, 6, 7},
                height=4cm,
                xticklabels={NFVM$_2^1$, UNFVM$_2^1$, GD, EO, ENO, WENO, DG},
                width=1.2\textwidth,
            ]
            \addplot[darkmagenta, 
                    only marks, 
                    mark=x,
                    mark options={ultra thick, scale=1.5},
                    error bars/.cd,
                    y dir=both, 
                    y explicit,
                    error bar style={very thick, color=pink}
            ] coordinates { 
                (1, 1.7e-2) +- (5e-3, 5e-3)
                (2, 1.9e-2) +- (6e-3, 6e-3)
                (3, 2.2e-2) +- (7e-3, 7e-3)
                (4, 2.2e-2) +- (7e-3, 7e-3)
                (5, 1.8e-2) +- (1e-2, 1e-2)
                (6, 1.7e-2) +- (1e-2, 1e-2)
                (7, 4.5e-3) +- (1e-3, 1e-3)
            };

            \draw[thick, dashed] (axis cs:4.5,\pgfkeysvalueof{/pgfplots/ymin}) -- (axis cs:4.5,\pgfkeysvalueof{/pgfplots/ymax});
            \draw[thick, dashed] (axis cs:6.5,\pgfkeysvalueof{/pgfplots/ymin}) -- (axis cs:6.5,\pgfkeysvalueof{/pgfplots/ymax});
            \end{axis}
        \end{tikzpicture}
        \vspace{-0.2cm}
        \caption{Triangular Sym.}
    \end{subfigure}
    \begin{subfigure}[t]{.31\textwidth}
    \centering
        \begin{tikzpicture}
            \begin{axis}[
                ymin=0,
                tick label style={font=\scriptsize},
                xticklabel style={rotate=-45, anchor=west, xshift=-5pt, yshift=-5pt},
                ylabel style={yshift=-20pt},
                xtick={1, 2, 3, 4, 5, 6, 7},
                height=4cm,
                xticklabels={NFVM$_2^1$, UNFVM$_2^1$, GD, EO, ENO, WENO, DG},
                width=1.2\textwidth,
            ]
            \addplot[darkmagenta, 
                    only marks, 
                    mark=x,
                    mark options={ultra thick, scale=1.5},
                    error bars/.cd,
                    y dir=both, 
                    y explicit,
                    error bar style={very thick, color=pink}
            ] coordinates { 
                (1, 2.2e-2) +- (9e-3, 9e-3)
                (2, 2.6e-2) +- (1e-2, 1e-2)
                (3, 3.0e-2) +- (1e-2, 1e-2)
                (4, 3.0e-2) +- (1e-2, 1e-2)
                (5, 4.0e-2) +- (2e-2, 2e-2)
                (6, 4.0e-2) +- (2e-2, 2e-2)
                (7, 5.6e-3) +- (2e-3, 2e-3)
            };

            \draw[thick, dashed] (axis cs:4.5,\pgfkeysvalueof{/pgfplots/ymin}) -- (axis cs:4.5,\pgfkeysvalueof{/pgfplots/ymax});
            \draw[thick, dashed] (axis cs:6.5,\pgfkeysvalueof{/pgfplots/ymin}) -- (axis cs:6.5,\pgfkeysvalueof{/pgfplots/ymax});
            \end{axis}
        \end{tikzpicture}
        \vspace{-0.2cm}
        \caption{Triangular Skw.}
    \end{subfigure}
            \begin{subfigure}[t]{.35\textwidth}
            \centering
        \begin{tikzpicture}
            \begin{axis}[
                ymin=0,
                tick label style={font=\scriptsize},
                xticklabel style={rotate=-45, anchor=west, xshift=-5pt, yshift=-5pt},
                ylabel={\(L_1\) error},
                xtick={1, 2, 3, 4, 5, 6, 7},
                height=4cm,
                xticklabels={NFVM$_2^1$, UNFVM$_2^1$, GD, EO, ENO, WENO, DG},
                width=1\textwidth,
            ]
            \addplot[darkmagenta, 
                    only marks, 
                    mark=x,
                    mark options={ultra thick, scale=1.5},
                    error bars/.cd,
                    y dir=both, 
                    y explicit,
                    error bar style={very thick, color=pink}
            ] coordinates { 
                (1, 1.4e-2) +- (2e-3, 2e-3)
                (2, 1.6e-2) +- (4e-3, 4e-3)
                (3, 2.0e-2) +- (5e-3, 5e-3)
                (4, 2.0e-2) +- (5e-3, 5e-3)
                (5, 8.5e-3) +- (2e-3, 2e-3)
                (6, 7.6e-3) +- (2e-3, 2e-3)
                (7, 4.9e-3) +- (1e-3, 1e-3)
            };

            \draw[thick, dashed] (axis cs:4.5,\pgfkeysvalueof{/pgfplots/ymin}) -- (axis cs:4.5,\pgfkeysvalueof{/pgfplots/ymax});
            \draw[thick, dashed] (axis cs:6.5,\pgfkeysvalueof{/pgfplots/ymin}) -- (axis cs:6.5,\pgfkeysvalueof{/pgfplots/ymax});
            \end{axis}
        \end{tikzpicture}
        \vspace{-0.2cm}
        \caption{Trapezoidal.}
    \end{subfigure}
    \begin{subfigure}[t]{.31\textwidth}
    \centering
        \begin{tikzpicture}
            \begin{axis}[
                ymin=0,
                tick label style={font=\scriptsize},
                xticklabel style={rotate=-45, anchor=west, xshift=-5pt, yshift=-5pt},
                ylabel style={yshift=-20pt},
                xtick={1, 2, 3, 4, 5, 6, 7},
                height=4cm,
                xticklabels={NFVM$_2^1$, UNFVM$_2^1$, GD, EO, ENO, WENO, DG},
                width=1.2\textwidth,
            ]
            \addplot[darkmagenta, 
                    only marks, 
                    mark=x,
                    mark options={ultra thick, scale=1.5},
                    error bars/.cd,
                    y dir=both, 
                    y explicit,
                    error bar style={very thick, color=pink}
            ] coordinates { 
                (1, 4.0e-3) +- (1e-3, 1e-3)
                (2, 7.3e-3) +- (1e-3, 1e-3)
                (3, 9.3e-3) +- (2e-3, 2e-3)
                (4, 9.4e-3) +- (2e-3, 2e-3)
                (5, 1.4e-2) +- (5e-3, 5e-3)
                (6, 1.4e-2) +- (5e-3, 5e-3)
                (7, 2.3e-3) +- (4e-3, 4e-3)
            };

            \draw[thick, dashed] (axis cs:4.5,\pgfkeysvalueof{/pgfplots/ymin}) -- (axis cs:4.5,\pgfkeysvalueof{/pgfplots/ymax});
            \draw[thick, dashed] (axis cs:6.5,\pgfkeysvalueof{/pgfplots/ymin}) -- (axis cs:6.5,\pgfkeysvalueof{/pgfplots/ymax});
            \end{axis}
        \end{tikzpicture}
        \vspace{-0.62cm}
        \caption{Greenberg.}
    \end{subfigure}
    \begin{subfigure}[t]{.31\textwidth}
    \centering
        \begin{tikzpicture}
            \begin{axis}[
                ymin=0,
                tick label style={font=\scriptsize},
                xticklabel style={rotate=-45, anchor=west, xshift=-5pt, yshift=-5pt},
                ylabel style={yshift=-20pt},
                xtick={1, 2, 3, 4, 5, 6, 7},
                height=4cm,
                xticklabels={NFVM$_2^1$, UNFVM$_2^1$, GD, EO, ENO, WENO, DG},
                width=1.2\textwidth,
            ]
            \addplot[darkmagenta, 
                    only marks, 
                    mark=x,
                    mark options={ultra thick, scale=1.5},
                    error bars/.cd,
                    y dir=both, 
                    y explicit,
                    error bar style={very thick, color=pink}
            ] coordinates { 
            (1, 9.3e-3) +- (2e-3, 2e-3)
            (2, 1.4e-2) +- (3e-3, 3e-3)
            (3, 1.6e-2) +- (3e-3, 3e-3)
            (4, 1.6e-2) +- (3e-3, 3e-3)
            (5, 4.0e-3) +- (6e-4, 6e-4)
            (6, 3.6e-3) +- (5e-4, 5e-4)
            (7, 3.2e-3) +- (6e-4, 6e-4)
            };

            \draw[thick, dashed] (axis cs:4.5,\pgfkeysvalueof{/pgfplots/ymin}) -- (axis cs:4.5,\pgfkeysvalueof{/pgfplots/ymax});
            \draw[thick, dashed] (axis cs:6.5,\pgfkeysvalueof{/pgfplots/ymin}) -- (axis cs:6.5,\pgfkeysvalueof{/pgfplots/ymax});
            \end{axis}
        \end{tikzpicture}
        \vspace{-0.62cm}
        \caption{Underwood.}
    \end{subfigure}
    \vspace{-0.5cm}
    \caption{Performance of learned models against baselines in $L_1$ error on evaluation set, with standard deviations reported as error bars. Exact values are reported in \Cref{tab:lwr_metrics}.}
	\label{fig:lwr_l1_error_bar_plots}
\end{figure}
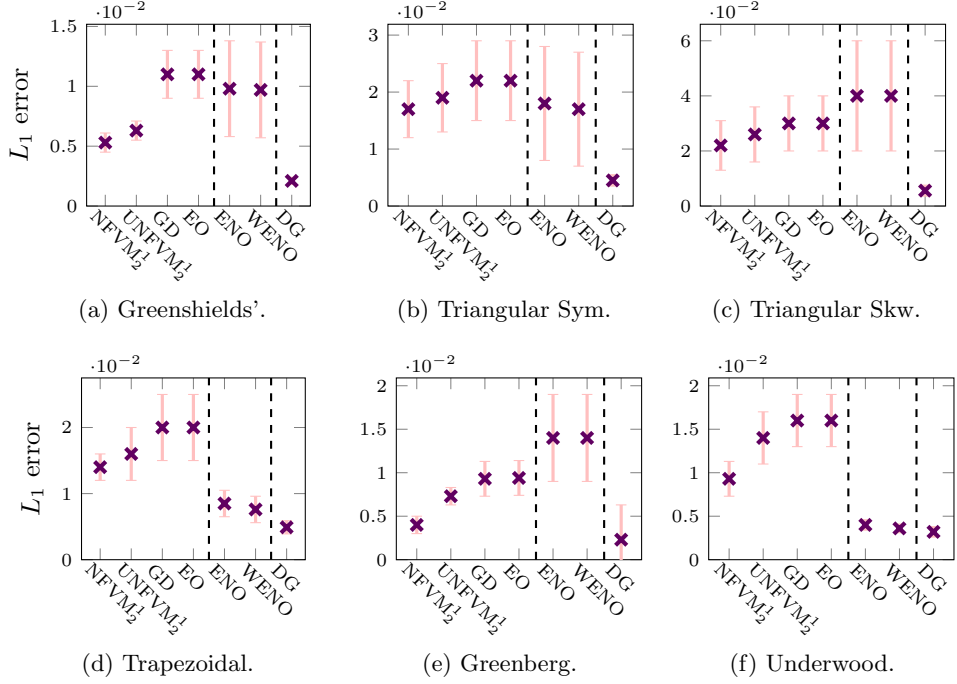

%% file: tex/lwr_l2_error_bar_plots.tex
\begin{figure}[H]
    \centering
    \begin{subfigure}[t]{.35\textwidth}
        \begin{tikzpicture}
            \begin{axis}[
                ymin=0,
                tick label style={font=\scriptsize},
                xticklabel style={rotate=-45, anchor=west, xshift=-5pt, yshift=-5pt},
                ylabel={\(L_2\) error},
                xtick={1, 2, 3, 4, 5, 6, 7},
                height=4cm,
                xticklabels={NFVM$_2^1$, UNFVM$_2^1$, GD, EO, ENO, WENO, DG},
                width=1\textwidth,
            ]
            \addplot[darkmagenta, 
                    only marks, 
                    mark=x,
                    mark options={ultra thick, scale=1.5},
                    error bars/.cd,
                    y dir=both, 
                    y explicit,
                    error bar style={very thick, color=pink}
            ] coordinates { 
                (1, 1.3e-4) +- (4e-5, 4e-5)
                (2, 2.0e-4) +- (6e-5, 6e-5)
                (3, 4.5e-4) +- (2e-4, 2e-4)
                (4, 4.5e-4) +- (2e-4, 2e-4)
                (5, 6.4e-4) +- (4e-4, 4e-4)
                (6, 6.4e-4) +- (4e-4, 4e-4)
                (7, 3.1e-5) +- (1e-5, 1e-5)
            };

            \draw[thick, dashed] (axis cs:4.5,\pgfkeysvalueof{/pgfplots/ymin}) -- (axis cs:4.5,\pgfkeysvalueof{/pgfplots/ymax});
            \draw[thick, dashed] (axis cs:6.5,\pgfkeysvalueof{/pgfplots/ymin}) -- (axis cs:6.5,\pgfkeysvalueof{/pgfplots/ymax});
            \end{axis}
        \end{tikzpicture}
        \vspace{-0.2cm}
        \caption{Greenshields'.}
    \end{subfigure}
    \begin{subfigure}[t]{.31\textwidth}
    \centering
        \begin{tikzpicture}
            \begin{axis}[
                ymin=0,
                tick label style={font=\scriptsize},
                xticklabel style={rotate=-45, anchor=west, xshift=-5pt, yshift=-5pt},
                ylabel style={yshift=-20pt},
                xtick={1, 2, 3, 4, 5, 6, 7},
                height=4cm,
                xticklabels={NFVM$_2^1$, UNFVM$_2^1$, GD, EO, ENO, WENO, DG},
                width=1.2\textwidth,
            ]
            \addplot[darkmagenta, 
                    only marks, 
                    mark=x,
                    mark options={ultra thick, scale=1.5},
                    error bars/.cd,
                    y dir=both, 
                    y explicit,
                    error bar style={very thick, color=pink}
            ] coordinates { 
            (1, 1.4e-3) +- (6e-4, 6e-4)
            (2, 1.9e-3) +- (9e-4, 9e-4)
            (3, 2.3e-3) +- (1e-3, 1e-3)
            (4, 2.3e-3) +- (1e-3, 1e-3)
            (5, 2.0e-3) +- (2e-3, 2e-3)
            (6, 1.9e-3) +- (2e-3, 2e-3)
            (7, 2.6e-4) +- (1e-4, 1e-4)
            };

            \draw[thick, dashed] (axis cs:4.5,\pgfkeysvalueof{/pgfplots/ymin}) -- (axis cs:4.5,\pgfkeysvalueof{/pgfplots/ymax});
            \draw[thick, dashed] (axis cs:6.5,\pgfkeysvalueof{/pgfplots/ymin}) -- (axis cs:6.5,\pgfkeysvalueof{/pgfplots/ymax});
            \end{axis}
        \end{tikzpicture}
        \vspace{-0.2cm}
        \caption{Triangular Sym.}
    \end{subfigure}
    \begin{subfigure}[t]{.31\textwidth}
    \centering
        \begin{tikzpicture}
            \begin{axis}[
                ymin=0,
                tick label style={font=\scriptsize},
                xticklabel style={rotate=-45, anchor=west, xshift=-5pt, yshift=-5pt},
                ylabel style={yshift=-20pt},
                xtick={1, 2, 3, 4, 5, 6, 7},
                height=4cm,
                xticklabels={NFVM$_2^1$, UNFVM$_2^1$, GD, EO, ENO, WENO, DG},
                width=1.2\textwidth,
            ]
            \addplot[darkmagenta, 
                    only marks, 
                    mark=x,
                    mark options={ultra thick, scale=1.5},
                    error bars/.cd,
                    y dir=both, 
                    y explicit,
                    error bar style={very thick, color=pink}
            ] coordinates { 
            (1, 2.4e-3) +- (1e-3, 1e-3)
            (2, 3.1e-3) +- (2e-3, 2e-3)
            (3, 3.8e-3) +- (2e-3, 2e-3)
            (4, 3.8e-3) +- (2e-3, 2e-3)
            (5, 5.8e-3) +- (4e-3, 4e-3)
            (6, 5.8e-3) +- (4e-3, 4e-3)
            (7, 4.1e-4) +- (2e-4, 2e-4)
            };

            \draw[thick, dashed] (axis cs:4.5,\pgfkeysvalueof{/pgfplots/ymin}) -- (axis cs:4.5,\pgfkeysvalueof{/pgfplots/ymax});
            \draw[thick, dashed] (axis cs:6.5,\pgfkeysvalueof{/pgfplots/ymin}) -- (axis cs:6.5,\pgfkeysvalueof{/pgfplots/ymax});
            \end{axis}
        \end{tikzpicture}
        \vspace{-0.62cm}
        \caption{Triangular Skw.}
    \end{subfigure}
            \begin{subfigure}[t]{.35\textwidth}
            \centering
        \begin{tikzpicture}
            \begin{axis}[
                ymin=0,
                tick label style={font=\scriptsize},
                xticklabel style={rotate=-45, anchor=west, xshift=-5pt, yshift=-5pt},
                ylabel={\(L_2\) error},
                xtick={1, 2, 3, 4, 5, 6, 7},
                height=4cm,
                xticklabels={NFVM$_2^1$, UNFVM$_2^1$, GD, EO, ENO, WENO, DG},
                width=1\textwidth,
            ]
            \addplot[darkmagenta, 
                    only marks, 
                    mark=x,
                    mark options={ultra thick, scale=1.5},
                    error bars/.cd,
                    y dir=both, 
                    y explicit,
                    error bar style={very thick, color=pink}
            ] coordinates { 
            (1, 1.1e-3) +- (4e-4, 4e-4)
            (2, 1.6e-3) +- (7e-4, 7e-4)
            (3, 2.1e-3) +- (8e-4, 8e-4)
            (4, 2.1e-3) +- (8e-4, 8e-4)
            (5, 6.2e-4) +- (2e-4, 2e-4)
            (6, 5.3e-4) +- (2e-4, 2e-4)
            (7, 2.9e-4) +- (1e-4, 1e-4)
            };

            \draw[thick, dashed] (axis cs:4.5,\pgfkeysvalueof{/pgfplots/ymin}) -- (axis cs:4.5,\pgfkeysvalueof{/pgfplots/ymax});
            \draw[thick, dashed] (axis cs:6.5,\pgfkeysvalueof{/pgfplots/ymin}) -- (axis cs:6.5,\pgfkeysvalueof{/pgfplots/ymax});
            \end{axis}
        \end{tikzpicture}
        \vspace{-0.2cm}
        \caption{Trapezoidal.}
    \end{subfigure}
    \begin{subfigure}[t]{.31\textwidth}
    \centering
        \begin{tikzpicture}
            \begin{axis}[
                ymin=0,
                tick label style={font=\scriptsize},
                xticklabel style={rotate=-45, anchor=west, xshift=-5pt, yshift=-5pt},
                ylabel style={yshift=-20pt},
                xtick={1, 2, 3, 4, 5, 6, 7},
                height=4cm,
                xticklabels={NFVM$_2^1$, UNFVM$_2^1$, GD, EO, ENO, WENO, DG},
                width=1.2\textwidth,
            ]
            \addplot[darkmagenta, 
                    only marks, 
                    mark=x,
                    mark options={ultra thick, scale=1.5},
                    error bars/.cd,
                    y dir=both, 
                    y explicit,
                    error bar style={very thick, color=pink}
            ] coordinates { 
            (1, 1.4e-4) +- (9e-5, 9e-5)
            (2, 3.8e-4) +- (2e-4, 2e-4)
            (3, 4.9e-4) +- (2e-4, 2e-4)
            (4, 4.9e-4) +- (2e-4, 2e-4)
            (5, 1.1e-3) +- (6e-4, 6e-4)
            (6, 1.2e-3) +- (9e-4, 9e-4)
            (7, 3.4e-4) +- (2e-3, 2e-3)
            };

            \draw[thick, dashed] (axis cs:4.5,\pgfkeysvalueof{/pgfplots/ymin}) -- (axis cs:4.5,\pgfkeysvalueof{/pgfplots/ymax});
            \draw[thick, dashed] (axis cs:6.5,\pgfkeysvalueof{/pgfplots/ymin}) -- (axis cs:6.5,\pgfkeysvalueof{/pgfplots/ymax});
            \end{axis}
        \end{tikzpicture}
        \vspace{-0.2cm}
        \caption{Greenberg.}
    \end{subfigure}
    \begin{subfigure}[t]{.31\textwidth}
    \centering
        \begin{tikzpicture}
            \begin{axis}[
                ymin=0,
                tick label style={font=\scriptsize},
                xticklabel style={rotate=-45, anchor=west, xshift=-5pt, yshift=-5pt},
                ylabel style={yshift=-20pt},
                xtick={1, 2, 3, 4, 5, 6, 7},
                height=4cm,
                xticklabels={NFVM$_2^1$, UNFVM$_2^1$, GD, EO, ENO, WENO, DG},
                width=1.2\textwidth,
            ]
            \addplot[darkmagenta, 
                    only marks, 
                    mark=x,
                    mark options={ultra thick, scale=1.5},
                    error bars/.cd,
                    y dir=both, 
                    y explicit,
                    error bar style={very thick, color=pink}
            ] coordinates { 
            (1, 3.8e-4) +- (1e-4, 1e-4)
            (2, 6.9e-4) +- (2e-4, 2e-4)
            (3, 9.2e-4) +- (3e-4, 3e-4)
            (4, 9.2e-4) +- (3e-4, 3e-4)
            (5, 1.1e-4) +- (3e-5, 3e-5)
            (6, 9.8e-5) +- (2e-5, 2e-5)
            (7, 5.9e-5) +- (2e-5, 2e-5)
            };

            \draw[thick, dashed] (axis cs:4.5,\pgfkeysvalueof{/pgfplots/ymin}) -- (axis cs:4.5,\pgfkeysvalueof{/pgfplots/ymax});
            \draw[thick, dashed] (axis cs:6.5,\pgfkeysvalueof{/pgfplots/ymin}) -- (axis cs:6.5,\pgfkeysvalueof{/pgfplots/ymax});
            \end{axis}
        \end{tikzpicture}
        \vspace{-0.62cm}
        \caption{Underwood.}
    \end{subfigure}
    \vspace{-0.5cm}
    \caption{Performance of learned models against baselines in $L_2$ error on evaluation set, with standard deviations reported as error bars. Exact values are reported in \Cref{tab:lwr_metrics}.}
    \label{fig:lwr_l2_error_bar_plots}
\end{figure}

%% file: tex/lwr_winrates.tex
\newcommand{\Fig}[5]{%
    \input{#2}
    \begin{subfigure}[t]{#5 \textwidth}
        \centering
        \tiny
        \begin{tikzpicture}[scale=.45]
              \foreach \y [count=\n] in \mydata {
                \foreach \x [count=\m] in \y {
                  \ifnum\n=\m
                    \node[fill=black, draw=none, minimum size=4.5mm, text=black] at (\m,-\n) {};
                  \else
                    \node[fill=teal!\x!purple!80!white, minimum size=4.5mm, text=white, inner sep=0] at (\m,-\n) {\x};
                  \fi
            }
                }

          \draw[thick] (0.5,-2.5) -- (8.5,-2.5); 
          \draw[thick, dashed] (5.5,-.5) -- (5.5,-8.5); 
          \draw[thick, dashed] (7.5,-.5) -- (7.5,-8.5); 

        \ifthenelse{\equal{#4}{true}}{%
            \foreach \a [count=\i] in {NFVM$_2^1$,UNFVM$_2^1$,GD,LF,EO,ENO,WENO,DG}  {
                \node[inner sep=0, anchor=south] at (\i + .2, -.4) {\rotatebox{45}{\a}};
              }%
        }{}
        \ifthenelse{\equal{#3}{true}}{%
              \foreach \a [count=\i] in {NFVM$_2^1$,UNFVM$_2^1$,GD,LF,EO,ENO,WENO,DG} {
                \node[] at (-.5,-\i+.2) {\rotatebox{-45}{\a}};
              }%
        }{}
        \end{tikzpicture}
        \subcaption{#1}
    \end{subfigure}%
}

\begin{figure}
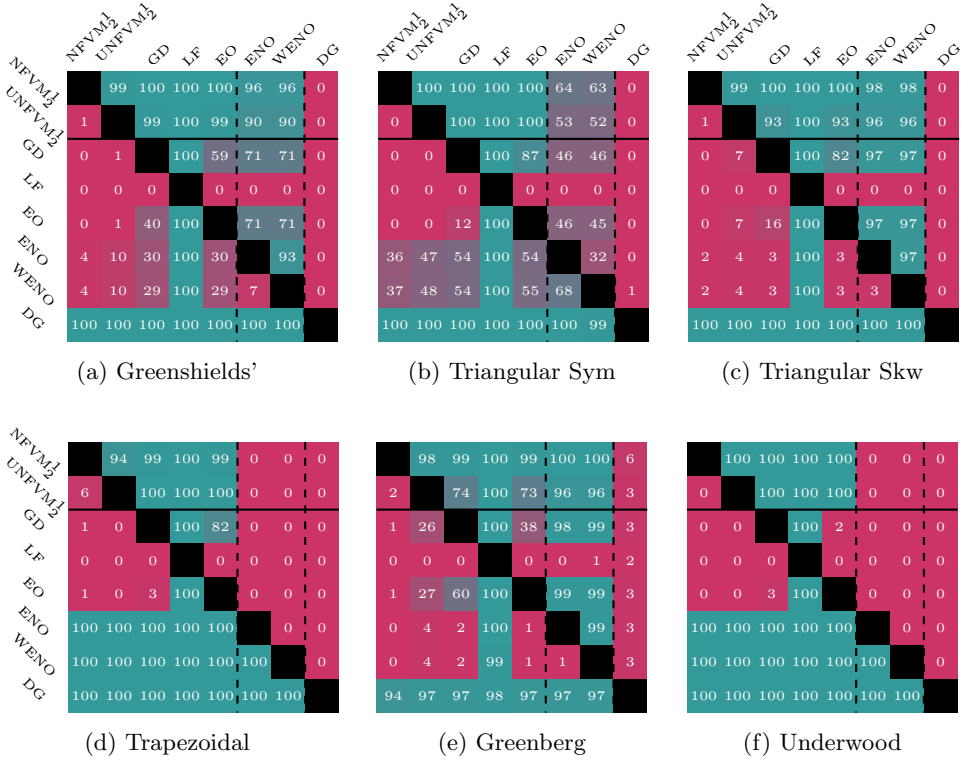

    \centering
    \Fig{Greenshields'}{figs/lwr_winrates/winrates_greenshield.tex}{true}{true}{.36}
    \Fig{Triangular Sym}{figs/lwr_winrates/winrates_triangular.tex}{false}{true}{.29}
    \Fig{Triangular Skw}{figs/lwr_winrates/winrates_triangular_skewed.tex}{false}{true}{.29}
    \Fig{Trapezoidal}{figs/lwr_winrates/winrates_trapezoidal}{true}{false}{.36}
    \Fig{Greenberg}{figs/lwr_winrates/winrates_greenberg}{false}{false}{.29}
    \Fig{Underwood}{figs/lwr_winrates/winrates_underwood}{false}{false}{.29}
    \caption{Winrates of schemes against each other, computed over the evaluation set of 1000 initial conditions (see~\Cref{sec:results}) for each of the 6 flux functions. A winrate of $N$ for row-scheme $X$ against column-scheme $Y$ means that scheme $X$ has a lower $L_2$ error than scheme $Y$ on $N$\% of the initial conditions for that flux.}
    \label{fig:lwr_winrates}
\end{figure}

%% file: tex/lwr_convergence_plot.tex
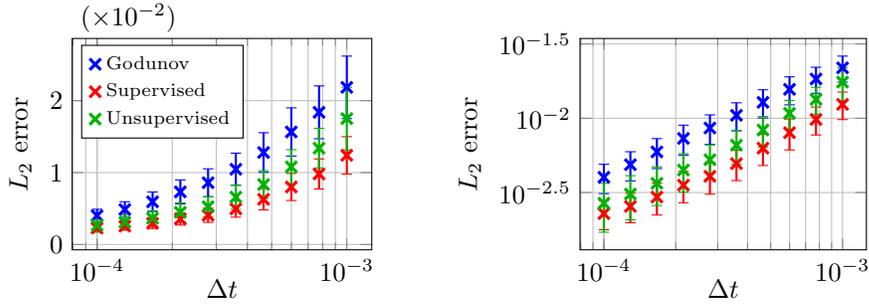
\begin{figure}
    \centering
    \begin{subfigure}[t]{.45\textwidth}
        \begin{tikzpicture}
            \begin{axis}[
                xlabel={\(\Delta t\)},
                ylabel={\(L_2\) error},
                grid=both,   
                xmode=log,
                x label style={at={(axis description cs:.5,-.1)}},
                legend style={at={(.27,0.84)}, anchor=center},
                legend cell align=left,
                width=.95\textwidth,
                height=.75\textwidth,
                enlargelimits=0.1,
                y tick label style={/pgf/number format/fixed},
                ytick scale label code/.code={($\times 10^{#1}$)},
                legend style={font=\scriptsize},
                legend pos = north west,
            ]
            \addplot[blue, 
                     only marks, 
                     mark=x,
                     mark options={very thick, scale=1.5},
                     error bars/.cd,
                     y dir=both, 
                     y explicit,
                     error bar style={thick, color=blue}
            ] coordinates { 
(0.001, 0.02184995597091381) +- (0.0043518786623826635, 0.0043518786623826635)
(0.0007742636826811271, 0.01838306431386583) +- (0.0036746066675479644, 0.0036746066675479644)
(0.000599484250318941, 0.015653029710674806) +- (0.003349639050582098, 0.003349639050582098)
(0.0004641588833612779, 0.012793922524939815) +- (0.0027536709411917237, 0.0027536709411917237)
(0.00035938136638046273, 0.010466471035864865) +- (0.002237362854026388, 0.002237362854026388)
(0.0002782559402207124, 0.00859585272119339) +- (0.0019355612437412593, 0.0019355612437412593)
(0.0002154434690031884, 0.00731524856816584) +- (0.0016680487247411687, 0.0016680487247411687)
(0.0001668100537200059, 0.005937887644113963) +- (0.0013556247170399092, 0.0013556247170399092)
(0.00012915496650148844, 0.0048693872096709015) +- (0.0010829769972929286, 0.0010829769972929286)
(9.999999999999998e-05, 0.004004110215715265) +- (0.0008950858175905402, 0.0008950858175905402)
            };
            \addlegendentry{Godunov};

            \addplot[red, 
                    only marks, 
                    mark=x,
                    mark options={very thick, scale=1.5},
                    error bars/.cd,
                    y dir=both, 
                    y explicit,
                    error bar style={thick, color=red}
            ] coordinates { 
(0.001, 0.012398871031558099) +- (0.002606770378204459, 0.002606770378204459)
(0.0007742636826811271, 0.009808386748395601) +- (0.0020948545447002634, 0.0020948545447002634)
(0.000599484250318941, 0.008000624403867375) +- (0.0018837330042691983, 0.0018837330042691983)
(0.0004641588833612779, 0.006292124619445425) +- (0.0014867072580761918, 0.0014867072580761918)
(0.00035938136638046273, 0.004953528186474813) +- (0.0011505203252874792, 0.0011505203252874792)
(0.0002782559402207124, 0.0040670904715640055) +- (0.0009853362593965048, 0.0009853362593965048)
(0.0002154434690031884, 0.0035447599931533852) +- (0.0008459611690404069, 0.0008459611690404069)
(0.0001668100537200059, 0.002947867099663713) +- (0.0007179783456420399, 0.0007179783456420399)
(0.00012915496650148844, 0.0025506941116867107) +- (0.0005716555331476384, 0.0005716555331476384)
(9.999999999999998e-05, 0.0022822297822621355) +- (0.0005064092785838558, 0.0005064092785838558)
            };
            \addlegendentry{Supervised};

            \addplot[green!70!black, 
                    only marks, 
                    mark=x,
                    mark options={very thick, scale=1.5},
                    error bars/.cd,
                    y dir=both, 
                    y explicit,
                    error bar style={thick, color=green!70!black}
            ] coordinates { 
(0.001, 0.01751107666113228) +- (0.004354443993719672, 0.004354443993719672)
(0.0007742636826811271, 0.013397033639360269) +- (0.0027598736739124496, 0.0027598736739124496)
(0.000599484250318941, 0.010798895899196886) +- (0.0023861777480975507, 0.0023861777480975507)
(0.0004641588833612779, 0.008338250248229498) +- (0.0019417621218139338, 0.0019417621218139338)
(0.00035938136638046273, 0.0065911226747773655) +- (0.0016938655926175478, 0.0016938655926175478)
(0.0002782559402207124, 0.005256896162823991) +- (0.001359921357957374, 0.001359921357957374)
(0.0002154434690031884, 0.004480195132407675) +- (0.0013006906911855578, 0.0013006906911855578)
(0.0001668100537200059, 0.003647601879135498) +- (0.0010632675179153542, 0.0010632675179153542)
(0.00012915496650148844, 0.0030884805502485115) +- (0.0010175122636898468, 0.0010175122636898468)
(9.999999999999998e-05, 0.0026790538326641154) +- (0.0009695418685616329, 0.0009695418685616329)
            };
            \addlegendentry{Unsupervised};

            \end{axis}
        \end{tikzpicture}
    \end{subfigure}
    \begin{subfigure}[t]{.45\textwidth}
        \begin{tikzpicture}
            \begin{axis}[
                xlabel={\(\Delta t\)},
                ylabel={\(L_2\) error},
                grid=both,   
                xmode=log,
                ymode=log,
                x label style={at={(axis description cs:.5,-.1)}},
                width=.92\textwidth,
                height=.75\textwidth,
                enlargelimits=0.1,
            ]
        
            \addplot[blue, 
            only marks, 
            mark=x,
            mark options={very thick, scale=1.5},
            error bars/.cd,
            y dir=both, 
            y explicit,
            error bar style={thick, color=blue}
            ] coordinates { 
            (0.001, 0.02184995597091381) +- (0.0043518786623826635, 0.0043518786623826635)
            (0.0007742636826811271, 0.01838306431386583) +- (0.0036746066675479644, 0.0036746066675479644)
            (0.000599484250318941, 0.015653029710674806) +- (0.003349639050582098, 0.003349639050582098)
            (0.0004641588833612779, 0.012793922524939815) +- (0.0027536709411917237, 0.0027536709411917237)
            (0.00035938136638046273, 0.010466471035864865) +- (0.002237362854026388, 0.002237362854026388)
            (0.0002782559402207124, 0.00859585272119339) +- (0.0019355612437412593, 0.0019355612437412593)
            (0.0002154434690031884, 0.00731524856816584) +- (0.0016680487247411687, 0.0016680487247411687)
            (0.0001668100537200059, 0.005937887644113963) +- (0.0013556247170399092, 0.0013556247170399092)
            (0.00012915496650148844, 0.0048693872096709015) +- (0.0010829769972929286, 0.0010829769972929286)
            (9.999999999999998e-05, 0.004004110215715265) +- (0.0008950858175905402, 0.0008950858175905402)
            };

            \addplot[red, 
            only marks, 
            mark=x,
            mark options={very thick, scale=1.5},
            error bars/.cd,
            y dir=both, 
            y explicit,
            error bar style={thick, color=red}
            ] coordinates { 
            (0.001, 0.012398871031558099) +- (0.002606770378204459, 0.002606770378204459)
            (0.0007742636826811271, 0.009808386748395601) +- (0.0020948545447002634, 0.0020948545447002634)
            (0.000599484250318941, 0.008000624403867375) +- (0.0018837330042691983, 0.0018837330042691983)
            (0.0004641588833612779, 0.006292124619445425) +- (0.0014867072580761918, 0.0014867072580761918)
            (0.00035938136638046273, 0.004953528186474813) +- (0.0011505203252874792, 0.0011505203252874792)
            (0.0002782559402207124, 0.0040670904715640055) +- (0.0009853362593965048, 0.0009853362593965048)
            (0.0002154434690031884, 0.0035447599931533852) +- (0.0008459611690404069, 0.0008459611690404069)
            (0.0001668100537200059, 0.002947867099663713) +- (0.0007179783456420399, 0.0007179783456420399)
            (0.00012915496650148844, 0.0025506941116867107) +- (0.0005716555331476384, 0.0005716555331476384)
            (9.999999999999998e-05, 0.0022822297822621355) +- (0.0005064092785838558, 0.0005064092785838558)
            };

            \addplot[green!70!black, 
            only marks, 
            mark=x,
            mark options={very thick, scale=1.5},
            error bars/.cd,
            y dir=both, 
            y explicit,
            error bar style={thick, color=green!70!black}
            ] coordinates { 
            (0.001, 0.01751107666113228) +- (0.004354443993719672, 0.004354443993719672)
            (0.0007742636826811271, 0.013397033639360269) +- (0.0027598736739124496, 0.0027598736739124496)
            (0.000599484250318941, 0.010798895899196886) +- (0.0023861777480975507, 0.0023861777480975507)
            (0.0004641588833612779, 0.008338250248229498) +- (0.0019417621218139338, 0.0019417621218139338)
            (0.00035938136638046273, 0.0065911226747773655) +- (0.0016938655926175478, 0.0016938655926175478)
            (0.0002782559402207124, 0.005256896162823991) +- (0.001359921357957374, 0.001359921357957374)
            (0.0002154434690031884, 0.004480195132407675) +- (0.0013006906911855578, 0.0013006906911855578)
            (0.0001668100537200059, 0.003647601879135498) +- (0.0010632675179153542, 0.0010632675179153542)
            (0.00012915496650148844, 0.0030884805502485115) +- (0.0010175122636898468, 0.0010175122636898468)
            (9.999999999999998e-05, 0.0026790538326641154) +- (0.0009695418685616329, 0.0009695418685616329)
            };

            \end{axis}
        \end{tikzpicture}
    \end{subfigure}
    \caption{\textbf{Convergence plots on Greenshields' flux.} The average \(L_2\) error is computed against the exact solution on the evaluation set, with standard deviation shown as error bars. The ratio \(\Delta t / \Delta x = 0.1\) remains constant as the mesh is refined.}
    \label{fig:lwr_convergence_plot}
\end{figure}

%% file: tex/lwr_numerical_flux.tex
\begin{figure}
    \centering
    
    \begin{subfigure}{\textwidth}
      \centering
      \begin{tikzpicture}
        \begin{axis}[
          hide axis,
          legend columns=3,
          legend style={
            draw=none,
            /tikz/every even column/.append style={column sep=0.5cm}
          },
          width=0.8\textwidth,
        ]
          \addplot[line width=1pt,mark=none,color={rgb,1:red,0.173;green,0.627;blue,0.173}] coordinates {(0,0) (0,0)};
          \addlegendentry{\footnotesize Godunov};
          \addplot[line width=1pt,mark=none,color={rgb,1:red,0.122;green,0.467;blue,0.706}] coordinates {(0,0) (0,0)};
          \addlegendentry{\footnotesize Supervised (NFVM$_2^1$)};
          \addplot[line width=1pt,mark=none,color={rgb,1:red,1.000;green,0.498;blue,0.055}] coordinates {(0,0) (0,0)};
          \addlegendentry{\footnotesize Unsupervised (UNFVM$_2^1$)};
        \end{axis}
      \end{tikzpicture}
    \end{subfigure}
    \vspace{-3.3cm}
  
    \begin{subfigure}{0.33\textwidth}
      \centering
      \begin{tikzpicture}
        \begin{axis}[
          width=\textwidth,
          height=0.75\textwidth,
          xlabel={$\rho$},
          ylabel={$\mathcal{F}(0.25, \rho)$},
          ymin=0, ymax=0.3,
          xmin=0, xmax=1,
          grid=major,
          label style={font=\footnotesize},
          tick label style={font=\footnotesize},
        ]
          \addplot[line width=1pt,mark=none,color={rgb,1:red,0.173;green,0.627;blue,0.173}] table[x expr=\coordindex/200, y index=0] {figs/lwr_numerical_flux/greenshields/flux_cross_greenshield_godunov_0.25_x.csv};
          \addplot[line width=1pt,mark=none,color={rgb,1:red,0.122;green,0.467;blue,0.706}] table[x expr=\coordindex/200, y index=0] {figs/lwr_numerical_flux/greenshields/flux_cross_greenshield_supervised_0.25_x.csv};
          \addplot[line width=1pt,mark=none,color={rgb,1:red,1.000;green,0.498;blue,0.055}] table[x expr=\coordindex/200, y index=0] {figs/lwr_numerical_flux/greenshields/flux_cross_greenshield_unsupervised_0.25_x.csv};
        \end{axis}
      \end{tikzpicture}
    \end{subfigure}%
    \begin{subfigure}{0.33\textwidth}
      \centering
      \begin{tikzpicture}
        \begin{axis}[
          width=\textwidth,
          height=0.75\textwidth,
          xlabel={$\rho$},
          ylabel={$\mathcal{F}(0.50, \rho)$},
          ymin=0, ymax=0.3,
          xmin=0, xmax=1,
          grid=major,
          label style={font=\footnotesize},
          tick label style={font=\footnotesize},
        ]
          \addplot[line width=1pt,mark=none,color={rgb,1:red,0.173;green,0.627;blue,0.173}] table[x expr=\coordindex/200, y index=0] {figs/lwr_numerical_flux/greenshields/flux_cross_greenshield_godunov_0.50_x.csv};
          \addplot[line width=1pt,mark=none,color={rgb,1:red,0.122;green,0.467;blue,0.706}] table[x expr=\coordindex/200, y index=0] {figs/lwr_numerical_flux/greenshields/flux_cross_greenshield_supervised_0.50_x.csv};
          \addplot[line width=1pt,mark=none,color={rgb,1:red,1.000;green,0.498;blue,0.055}] table[x expr=\coordindex/200, y index=0] {figs/lwr_numerical_flux/greenshields/flux_cross_greenshield_unsupervised_0.50_x.csv};
        \end{axis}
      \end{tikzpicture}
    \end{subfigure}%
    \begin{subfigure}{0.33\textwidth}
      \centering
      \begin{tikzpicture}
        \begin{axis}[
          width=\textwidth,
          height=0.75\textwidth,
          xlabel={$\rho$},
          ylabel={$\mathcal{F}(0.75, \rho)$},
          ymin=0, ymax=0.3,
          xmin=0, xmax=1,
          grid=major,
          label style={font=\footnotesize},
          tick label style={font=\footnotesize},
        ]
          \addplot[line width=1pt,mark=none,color={rgb,1:red,0.173;green,0.627;blue,0.173}] table[x expr=\coordindex/200, y index=0] {figs/lwr_numerical_flux/greenshields/flux_cross_greenshield_godunov_0.75_x.csv};
          \addplot[line width=1pt,mark=none,color={rgb,1:red,0.122;green,0.467;blue,0.706}] table[x expr=\coordindex/200, y index=0] {figs/lwr_numerical_flux/greenshields/flux_cross_greenshield_supervised_0.75_x.csv};
          \addplot[line width=1pt,mark=none,color={rgb,1:red,1.000;green,0.498;blue,0.055}] table[x expr=\coordindex/200, y index=0] {figs/lwr_numerical_flux/greenshields/flux_cross_greenshield_unsupervised_0.75_x.csv};
        \end{axis}
      \end{tikzpicture}
    \end{subfigure}
    
    \begin{subfigure}{0.33\textwidth}
      \centering
      \begin{tikzpicture}
        \begin{axis}[
          width=\textwidth,
          height=0.75\textwidth,
          xlabel={$\rho$},
          ylabel={$\mathcal{F}(\rho, 0.25)$},
          ymin=0, ymax=0.3,
          xmin=0, xmax=1,
          grid=major,
          label style={font=\footnotesize},
          tick label style={font=\footnotesize},
        ]
          \addplot[line width=1pt,mark=none,color={rgb,1:red,0.173;green,0.627;blue,0.173}] table[x expr=\coordindex/200, y index=0] {figs/lwr_numerical_flux/greenshields/flux_cross_greenshield_godunov_x_0.25.csv};
          \addplot[line width=1pt,mark=none,color={rgb,1:red,0.122;green,0.467;blue,0.706}] table[x expr=\coordindex/200, y index=0] {figs/lwr_numerical_flux/greenshields/flux_cross_greenshield_supervised_x_0.25.csv};
          \addplot[line width=1pt,mark=none,color={rgb,1:red,1.000;green,0.498;blue,0.055}] table[x expr=\coordindex/200, y index=0] {figs/lwr_numerical_flux/greenshields/flux_cross_greenshield_unsupervised_x_0.25.csv};
        \end{axis}
      \end{tikzpicture}
    \end{subfigure}%
    \begin{subfigure}{0.33\textwidth}
      \centering
      \begin{tikzpicture}
        \begin{axis}[
          width=\textwidth,
          height=0.75\textwidth,
          xlabel={$\rho$},
          ylabel={$\mathcal{F}(\rho, 0.50)$},
          ymin=0, ymax=0.3,
          xmin=0, xmax=1,
          grid=major,
          label style={font=\footnotesize},
          tick label style={font=\footnotesize},
        ]
          \addplot[line width=1pt,mark=none,color={rgb,1:red,0.173;green,0.627;blue,0.173}] table[x expr=\coordindex/200, y index=0] {figs/lwr_numerical_flux/greenshields/flux_cross_greenshield_godunov_x_0.50.csv};
          \addplot[line width=1pt,mark=none,color={rgb,1:red,0.122;green,0.467;blue,0.706}] table[x expr=\coordindex/200, y index=0] {figs/lwr_numerical_flux/greenshields/flux_cross_greenshield_supervised_x_0.50.csv};
          \addplot[line width=1pt,mark=none,color={rgb,1:red,1.000;green,0.498;blue,0.055}] table[x expr=\coordindex/200, y index=0] {figs/lwr_numerical_flux/greenshields/flux_cross_greenshield_unsupervised_x_0.50.csv};
        \end{axis}
      \end{tikzpicture}
    \end{subfigure}%
    \begin{subfigure}{0.33\textwidth}
      \centering
      \begin{tikzpicture}
        \begin{axis}[
          width=\textwidth,
          height=0.75\textwidth,
          xlabel={$\rho$},
          ylabel={$\mathcal{F}(\rho, 0.75)$},
          ymin=0, ymax=0.3,
          xmin=0, xmax=1,
          grid=major,
          label style={font=\footnotesize},
          tick label style={font=\footnotesize},
        ]
          \addplot[line width=1pt,mark=none,color={rgb,1:red,0.173;green,0.627;blue,0.173}] table[x expr=\coordindex/200, y index=0] {figs/lwr_numerical_flux/greenshields/flux_cross_greenshield_godunov_x_0.75.csv};
          \addplot[line width=1pt,mark=none,color={rgb,1:red,0.122;green,0.467;blue,0.706}] table[x expr=\coordindex/200, y index=0] {figs/lwr_numerical_flux/greenshields/flux_cross_greenshield_supervised_x_0.75.csv};
          \addplot[line width=1pt,mark=none,color={rgb,1:red,1.000;green,0.498;blue,0.055}] table[x expr=\coordindex/200, y index=0] {figs/lwr_numerical_flux/greenshields/flux_cross_greenshield_unsupervised_x_0.75.csv};
        \end{axis}
      \end{tikzpicture}
    \end{subfigure}
    
    \caption{Comparison of the learned supervised and unsupervised numerical fluxes for NFVM$_2^1$ and UNFVM$_2^1$, respectively, against Godunov's for different left and right states, using the Greenshields model. Top row: Fixed left state $\rho_1$. Bottom row: Fixed right state $\rho_2$.}
    \label{fig:lwr_numerical_flux}
  \end{figure}
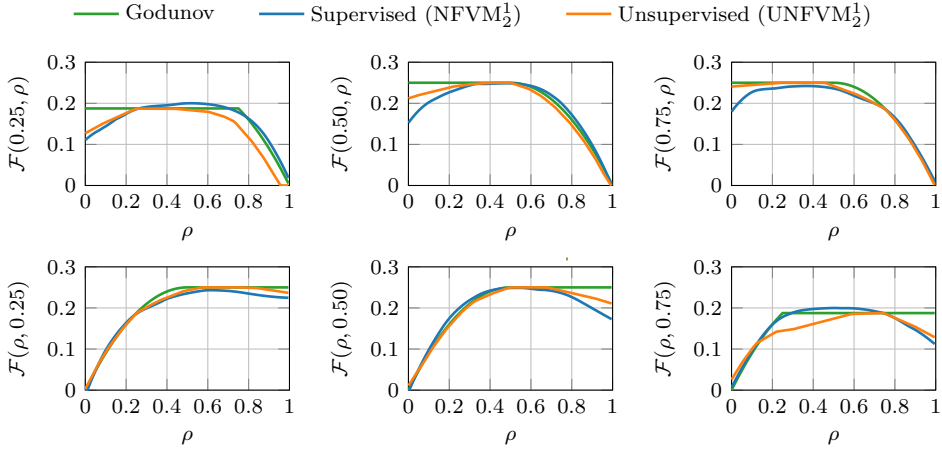

%% file: tex/lwr_2d_metrics.tex
\begin{table}
    \definecolor{deepcarmine}{rgb}{0.66, 0.13, 0.24}
    \setlength{\tabcolsep}{1pt}
    \centering
        \begin{tabularx}{\textwidth}{c|CCCCC}
        \Xhline{1.pt} 
        & Godunov & WENO & \textbf{NFVM}\(\bm{^1_2}\) & \textbf{NFVM}\(\bm{^5_4}\) & DG\\
        \Xhline{1.pt} 
        Greenshields' (L2) & \(4.1\color{deepcarmine}\mathrm{e}^{-4}\)\std{\(1\mathrm{e}^{-4}\)} & \(6.9\color{deepcarmine}\mathrm{e}^{-4}\)\std{\(4\mathrm{e}^{-4}\)}& \(1.2\color{deepcarmine}\mathrm{e}^{-4}\)\std{\(4\mathrm{e}^{-5}\)} & \(4.6\color{deepcarmine}\mathrm{e}^{-5}\)\std{\(3\mathrm{e}^{-5}\)} & \(\bm{4.2}\color{deepcarmine}\bm{\mathrm{e}^{-5}}\)\std{\(2\mathrm{e}^{-5}\)}\\
        Triangular Sym (L2) & \(2.2\color{deepcarmine}\mathrm{e}^{-3}\)\std{\(1\mathrm{e}^{-3}\)} & \(2.0\color{deepcarmine}\mathrm{e}^{-3}\)\std{\(2\mathrm{e}^{-3}\)} & \(1.3\color{deepcarmine}\mathrm{e}^{-3}\)\std{\(6\mathrm{e}^{-4}\)} & \(2.9\color{deepcarmine}\mathrm{e}^{-4}\)\std{\(2\mathrm{e}^{-4}\)} & \(\bm{2.7}\color{deepcarmine}\bm{\mathrm{e}^{-4}}\)\std{\(1\mathrm{e}^{-4}\)}\\
        \Xhline{1.pt} 
    \end{tabularx}
    \caption{Evaluation of NFVM\(_4^5\) on LWR using piecewise constant initial conditions.}
	\label{tab:lwr_2d_metrics}
\end{table}

%% file: tex/burgers_heatmaps.tex
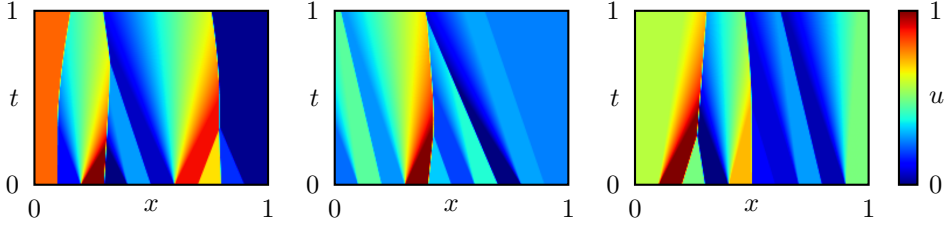
\begin{figure}
    \centering
    \begin{subfigure}[t]{.295\textwidth}
        \begin{tikzpicture}[scale=1., clip=false]
            \node[anchor=south west, inner sep=0] (img) at (0,0) {\includegraphics[height=2.3cm]{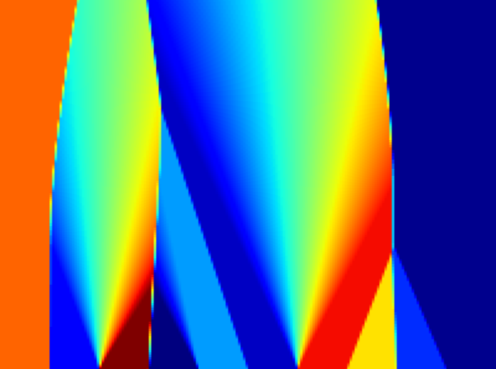}};
            \draw[black, thick] (img.south west) rectangle (img.north east);
            
            \path (img.south west) -- (img.south east) node[midway, below, yshift=-0.08cm] {\(x\)};
            \path (img.south west) -- (img.north west) node[midway, left, xshift=-0.08cm] {\(t\)};
            
            \node[below, yshift=-0.08cm] at (img.south west) {\(0\)};
            \node[below, yshift=-0.08cm] at (img.south east) {\(1\)};
            \node[left, xshift=-0.08cm] at (img.south west) {\(0\)};
            \node[left, xshift=-0.08cm] at (img.north west) {\(1\)};
        \end{tikzpicture}
    \end{subfigure}
    \hfill
    \begin{subfigure}[t]{.295\textwidth}
        \begin{tikzpicture}[scale=1., clip=false]
            \node[anchor=south west, inner sep=0] (img) at (0,0) {\includegraphics[height=2.3cm]{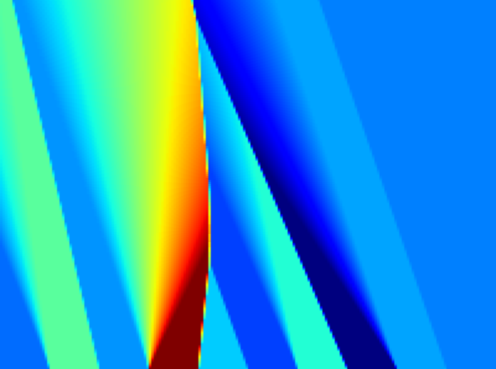}};
            \draw[black, thick] (img.south west) rectangle (img.north east);
            
            \path (img.south west) -- (img.south east) node[midway, below, yshift=-0.08cm] {\(x\)};
            \path (img.south west) -- (img.north west) node[midway, left, xshift=-0.08cm] {\(t\)};
            
            \node[below, yshift=-0.08cm] at (img.south west) {\(0\)};
            \node[below, yshift=-0.08cm] at (img.south east) {\(1\)};
            \node[left, xshift=-0.08cm] at (img.south west) {\(0\)};
            \node[left, xshift=-0.08cm] at (img.north west) {\(1\)};
        \end{tikzpicture}
    \end{subfigure}
    \hfill
    \begin{subfigure}[t]{.295\textwidth}
        \begin{tikzpicture}[scale=1., clip=false]
            \node[anchor=south west, inner sep=0] (img) at (0,0) {\includegraphics[height=2.3cm]{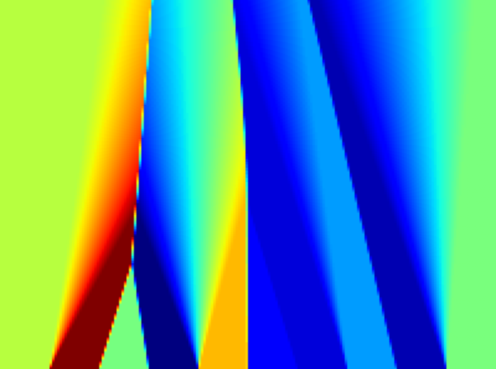}};
            \draw[black, thick] (img.south west) rectangle (img.north east);
            
            \path (img.south west) -- (img.south east) node[midway, below, yshift=-0.08cm] {\(x\)};
            \path (img.south west) -- (img.north west) node[midway, left, xshift=-0.08cm] {\(t\)};
            
            \node[below, yshift=-0.08cm] at (img.south west) {\(0\)};
            \node[below, yshift=-0.08cm] at (img.south east) {\(1\)};
            \node[left, xshift=-0.08cm] at (img.south west) {\(0\)};
            \node[left, xshift=-0.08cm] at (img.north west) {\(1\)};
        \end{tikzpicture}
    \end{subfigure}
    \hfill
    \begin{subfigure}[t]{.08\textwidth}
        \begin{tikzpicture}[scale=1., clip=false]
            \node[anchor=south west, inner sep=0] (img) at (0,0) {\includegraphics[height=2.3cm]{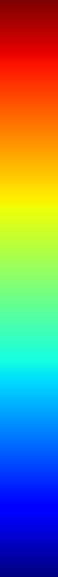}};
            \draw[black, thick] (img.south west) rectangle (img.north east);
            \path (img.south east) -- (img.north east) node
            [midway, right, xshift=+0.03cm] {\(u\)};
            \node[right, xshift=+0.03cm] at (img.south east) {\(0\)};
            \node[right, xshift=+0.03cm] at (img.north east) {\(1\)};
            \node[below, yshift=-0.08cm] at (img.south east) {\phantom{0}};
        \end{tikzpicture}
    \end{subfigure}

\caption{Solutions to Burgers' equation on piecewise constant initial conditions.}
	\label{fig:burgers_heatmaps}
\end{figure}

%% file: tex/burgers_metrics.tex
\begingroup
\setlength{\tabcolsep}{1pt}
\begin{table}
    \centering
    \begin{tabularx}{\linewidth}{c||cccc||CC||c}
        \hline
         \Xhline{1.pt} 
         & \multicolumn{4}{c||}{\textbf{1\textsuperscript{st} order FVM}} & \multicolumn{2}{c||}{\textbf{Higher order FVM}} & FEM\\
         \Xhline{1.pt} 
          & \textbf{NFVM}\(\bm{^1_2}\)\ & \textbf{UNFVM}\(\bm{^1_2}\)\ & GD & EO & ENO & WENO & DG\\
         \Xhline{1.pt}
         {L1} & \(\bm{1.2\mathrm{e}^{\shortminus2}}\)\std{\(2\mathrm{e}^{\shortminus3}\)} & \(1.6\mathrm{e}^{\shortminus2}\)\std{\(3\mathrm{e}^{\shortminus3}\)} & \(2.1\mathrm{e}^{\shortminus2}\)\std{\(4\mathrm{e}^{\shortminus3}\)} & \(1.9\mathrm{e}^{\shortminus2}\)\std{\(7\mathrm{e}^{\shortminus3}\)} & \(1.9\mathrm{e}^{\shortminus2}\)\std{\(7\mathrm{e}^{\shortminus3}\)} & \(2.0\mathrm{e}^{\shortminus4}\)\std{\(7\mathrm{e}^{\shortminus3}\)} & \(1.9\mathrm{e}^{\shortminus3}\)\std{\(4\mathrm{e}^{\shortminus4}\)}\\
         {L2} & \(\bm{8.5\mathrm{e}^{\shortminus4}}\)\std{\(3\mathrm{e}^{\shortminus4}\)} & \(1.3\mathrm{e}^{\shortminus3}\)\std{\(6\mathrm{e}^{\shortminus4}\)} & \(1.9\mathrm{e}^{\shortminus3}\)\std{\(7\mathrm{e}^{\shortminus4}\)} & \(2.6\mathrm{e}^{\shortminus3}\)\std{\(1\mathrm{e}^{\shortminus3}\)} & \(2.7\mathrm{e}^{\shortminus3}\)\std{\(1\mathrm{e}^{\shortminus3}\)} & \(2.8\mathrm{e}^{\shortminus3}\)\std{\(1\mathrm{e}^{\shortminus3}\)} & \(1.0\mathrm{e}^{\shortminus4}\)\std{\(4\mathrm{e}^{\shortminus5}\)}\\
         {Rel.} & \(\bm{5.9\mathrm{e}^{\shortminus2}}\)\std{\(1\mathrm{e}^{\shortminus2}\)} & \(8.5\mathrm{e}^{\shortminus2}\)\std{\(2\mathrm{e}^{\shortminus2}\)} & \(1.6\mathrm{e}^{\shortminus1}\)\std{\(4\mathrm{e}^{\shortminus2}\)} & \(1.5\mathrm{e}^{\shortminus1}\)\std{\(5\mathrm{e}^{\shortminus2}\)} & \(1.5\mathrm{e}^{\shortminus1}\)\std{\(5\mathrm{e}^{\shortminus2}\)} & \(1.5\mathrm{e}^{\shortminus1}\)\std{\(5\mathrm{e}^{\shortminus2}\)} & \(1.1\mathrm{e}^{\shortminus2}\)\std{\(5\mathrm{e}^{\shortminus3}\)}\\
         \Xhline{1.pt} 
    \end{tabularx}
    \caption{Results on Burgers' equation using piecewise constant initial conditions.}
    \label{tab:burgers_metrics}
\end{table}
\endgroup

%% file: tex/burgers_winrates.tex
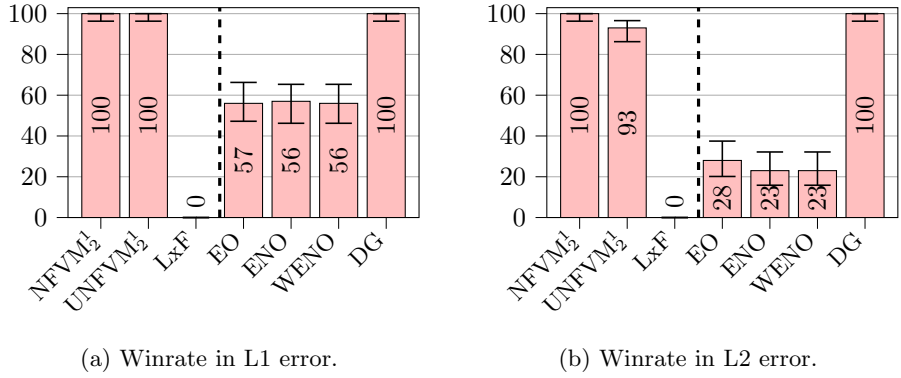
\begin{figure}
  \centering
  \begin{subfigure}[t]{.48\textwidth}
      \centering
      \begin{tikzpicture}
          \definecolor{darkgrey176}{RGB}{176,176,176}
          \definecolor{skyblue}{RGB}{135,206,235}
          \begin{axis}[
              tick align=outside,
              tick pos=left,
              x grid style={darkgrey176},
              xmin=-0.69, xmax=6.69,
              xtick style={color=black},
              xtick={0,1,2,3,4,5,6},
              xticklabels={NFVM$_2^1$,UNFVM$_2^1$,LxF, EO, ENO,WENO, DG},
              y grid style={darkgrey176},
              ymajorgrids,
              ymin=0, ymax=102.927592607118,
              ytick style={color=black},
              width=\textwidth,
              height=.7\textwidth,
              xticklabel style={rotate=45, anchor=east, font=\small},
          ]
              \draw[draw=black,fill=pink] (axis cs:-0.4,0) rectangle (axis cs:0.4,100);
              \draw[draw=black,fill=pink] (axis cs:0.6,0) rectangle (axis cs:1.4,100);
              \draw[draw=black,fill=pink] (axis cs:1.6,0) rectangle (axis cs:2.4,0);
              \draw[draw=black,fill=pink] (axis cs:2.6,0) rectangle (axis cs:3.4,56);
              \draw[draw=black,fill=pink] (axis cs:3.6,0) rectangle (axis cs:4.4,57);
              \draw[draw=black,fill=pink] (axis cs:4.6,0) rectangle (axis cs:5.4,56);
              \draw[draw=black,fill=pink] (axis cs:5.6,0) rectangle (axis cs:6.4,100);
              \path [draw=black, semithick]
              (axis cs:0,96.30051925)
              --(axis cs:0,100);
              
              \path [draw=black, semithick]
              (axis cs:1,96.30051925)
              --(axis cs:1,100);
              
              \path [draw=black, semithick]
              (axis cs:2,0)
              --(axis cs:2,0);
              
              \path [draw=black, semithick]
              (axis cs:3,47.2152117)
              --(axis cs:3,66.266861);
              
              \path [draw=black, semithick]
              (axis cs:4,46.22792771)
              --(axis cs:4,65.3281346);
                              
              \path [draw=black, semithick]
              (axis cs:5,46.22792771)
              --(axis cs:5,65.3281346);
              
              \path [draw=black, semithick]
              (axis cs:6,96.30051925)
              --(axis cs:6,100);
              
              \addplot [semithick, black, mark=-, mark size=5, mark options={solid}, only marks]
              table {%
                  0 96.30051925
                  1 96.30051925
                  2 0
                  3 47.2152117
                  4 46.22792771
                  5 46.22792771
                  6 96.30051925
              };
              \addplot [semithick, black, mark=-, mark size=5, mark options={solid}, only marks]
              table {%
                  0 100
                  1 100
                  2 0
                  3 66.266861
                  4 65.3281346
                  5 65.3281346
                  6 100
              };

              \draw (axis cs:0,50) node[
                anchor=center,
                text=black,
                rotate=90.0
              ]{100};
              \draw (axis cs:1,50) node[
                anchor=center,
                text=black,
                rotate=90.0
              ]{100};
              \draw (axis cs:2,8) node[
                anchor=center,
                text=black,
                rotate=90.0
              ]{0};
              \draw (axis cs:3,28) node[
                anchor=center,
                text=black,
                rotate=90.0
              ]{57};
              \draw (axis cs:4,28) node[
                anchor=center,
                text=black,
                rotate=90.0
              ]{56};
              \draw (axis cs:5,28) node[
                anchor=center,
                text=black,
                rotate=90.0
              ]{56};
              \draw (axis cs:6,50) node[
                anchor=center,
                text=black,
                rotate=90.0
              ]{100};

              \path [draw=black, very thick, dashed]
              (axis cs:2.5,0)
              --(axis cs:2.5,102.927592607118);
          \end{axis}
      \end{tikzpicture}
      \subcaption{Winrate in L1 error.}
  \end{subfigure}
  \begin{subfigure}[t]{.48\textwidth}
      \centering
      \begin{tikzpicture}
          \definecolor{darkgrey176}{RGB}{176,176,176}
          \definecolor{skyblue}{RGB}{135,206,235}
          \begin{axis}[
              tick align=outside,
              tick pos=left,
              x grid style={darkgrey176},
              xmin=-0.69, xmax=6.69,
              xtick style={color=black},
              xtick={0,1,2,3,4,5,6},
              xticklabels={NFVM$_2^1$,UNFVM$_2^1$,LxF, EO, ENO,WENO, DG},
              y grid style={darkgrey176},
              ymajorgrids,
              ymin=0, ymax=102.927592607118,
              ytick style={color=black},
              width=\textwidth,
              height=.7\textwidth,
              xticklabel style={rotate=45, anchor=east, font=\small},
          ]
              \draw[draw=black,fill=pink] (axis cs:-0.4,0) rectangle (axis cs:0.4,100);
              \draw[draw=black,fill=pink] (axis cs:0.6,0) rectangle (axis cs:1.4,93);
              \draw[draw=black,fill=pink] (axis cs:1.6,0) rectangle (axis cs:2.4,0);
              \draw[draw=black,fill=pink] (axis cs:2.6,0) rectangle (axis cs:3.4,28);
              \draw[draw=black,fill=pink] (axis cs:3.6,0) rectangle (axis cs:4.4,23);
              \draw[draw=black,fill=pink] (axis cs:4.6,0) rectangle (axis cs:5.4,23);
              \draw[draw=black,fill=pink] (axis cs:5.6,0) rectangle (axis cs:6.4,100);
              \path [draw=black, semithick]
              (axis cs:0,96.30051925)
              --(axis cs:0,100);
              
              \path [draw=black, semithick]
              (axis cs:1,86.25032899)
              --(axis cs:1,96.56811756);
              
              \path [draw=black, semithick]
              (axis cs:2,0)
              --(axis cs:2,0);
              
              \path [draw=black, semithick]
              (axis cs:3,20.13955917)
              --(axis cs:3,37.48821236);

              \path [draw=black, semithick]
              (axis cs:4,15.84315475)
              --(axis cs:4,32.15456486);
              
              \path [draw=black, semithick]
              (axis cs:5,15.84315475)
              --(axis cs:5,32.15456486);
              
              \path [draw=black, semithick]
              (axis cs:6,96.30051925)
              --(axis cs:6,100);
              
              \addplot [semithick, black, mark=-, mark size=5, mark options={solid}, only marks]
              table {%
                  0 96.30051925
                  1 86.25032899
                  2 0
                  3 20.13955917
                  4 15.84315475
                  5 15.84315475
                  6 96.30051925
              };
              \addplot [semithick, black, mark=-, mark size=5, mark options={solid}, only marks]
              table {%
                  0 100
                  1 96.56811756
                  2 0
                  3 37.48821236
                  4 32.15456486
                  5 32.15456486
                  6 100
              };

              \draw (axis cs:0,50) node[
                anchor=center,
                text=black,
                rotate=90.0
              ]{100};
              \draw (axis cs:1,46.5) node[
                anchor=center,
                text=black,
                rotate=90.0
              ]{93};
              \draw (axis cs:2,8) node[
                anchor=center,
                text=black,
                rotate=90.0
              ]{0};
              \draw (axis cs:3,9.5) node[
                anchor=center,
                text=black,
                rotate=90.0
              ]{28};
              \draw (axis cs:4,9.5) node[
                anchor=center,
                text=black,
                rotate=90.0
              ]{23};
              \draw (axis cs:5,9.5) node[
                anchor=center,
                text=black,
                rotate=90.0
              ]{23};
              \draw (axis cs:6,50) node[
                anchor=center,
                text=black,
                rotate=90.0
              ]{100};

          \path [draw=black, very thick, dashed]
              (axis cs:2.5,0)
              --(axis cs:2.5,102.927592607118);
          \end{axis}
      \end{tikzpicture}
      \subcaption{Winrate in L2 error.}
  \end{subfigure}
  \caption{Proportion of test cases where every considered method beats the Godunov scheme. Confidence interval are Wilson 95\% confidence interval~\cite{Wilson01061927}.}
  \label{fig:burgers_winrates}
\end{figure}

%% file: tex/burgers_2d_metrics.tex
\begin{table}
    \definecolor{deepcarmine}{rgb}{0.66, 0.13, 0.24}
    \setlength{\tabcolsep}{1pt}
    \centering
        \begin{tabularx}{\textwidth}{c|CCCCC}
        \Xhline{1.pt} 
        & Godunov & WENO & \textbf{NFVM}\(\bm{^1_2}\) & \textbf{NFVM}\(\bm{^5_4}\) & DG\\
        \Xhline{1.pt} 
        Burgers' (L2) & \(1.8\color{deepcarmine}\mathrm{e}^{-3}\)\std{\(6\mathrm{e}^{-4}\)} & \(2.6\color{deepcarmine}\mathrm{e}^{-3}\)\std{\(1\mathrm{e}^{-3}\)} & \(8.3\color{deepcarmine}\mathrm{e}^{-4}\)\std{\(3\mathrm{e}^{-4}\)} & \(2.2\color{deepcarmine}\mathrm{e}^{-4}\)\std{\(1\mathrm{e}^{-4}\)} & \(\bm{1.0}\color{deepcarmine}\bm{\mathrm{e}^{-4}}\)\std{\(4\mathrm{e}^{-5}\)}\\
        \Xhline{1.pt} 
    \end{tabularx}
    \caption{Evaluation of NFVM\(_4^5\) on Burgers' using piecewise constant initial conditions.}
	\label{tab:burgers_2d_metrics}
\end{table}

%% file: tex/burgers_2d_final_density.tex
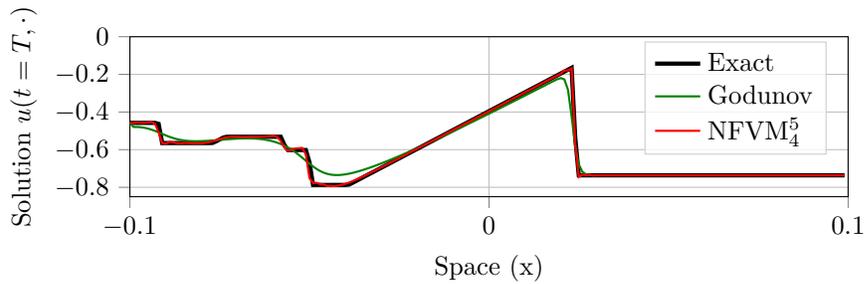
\begin{figure}
    \centering
    \begin{subfigure}[t]{.95\textwidth}
        \begingroup
        \centering
        \begin{tikzpicture}[baseline]
            \definecolor{darkgray176}{RGB}{176,176,176}
            \definecolor{darkorange25512714}{RGB}{255,0,0}
            \definecolor{forestgreen4416044}{RGB}{44,160,44}
            \definecolor{lightgray204}{RGB}{204,204,204}
            \definecolor{steelblue31119180}{RGB}{100,0,255}
            
            \begin{axis}[
                legend pos=north east,
                legend cell align={left},
                legend style={fill opacity=0.8, draw opacity=1, text opacity=1, draw=lightgray204},
                tick align=outside,
                tick pos=left,
                width=0.9\textwidth,
                height=.3\textwidth,
                xlabel={Space (x)},
                xmajorgrids,
                ymin=-.85, ymax=0,
                xmin=0-.1, xmax=.1, 
                ylabel={Solution \(\displaystyle u(t=T, \cdot)\)},
                ymajorgrids,
                legend columns=1,
                xtick={-.1, 0, .1},
            ]
            \addplot [ultra thick, black]
            table[y index=0,  x expr=\coordindex*0.001 - .1] {figs/2d_final_density/burgers_lax_hopf.tex};
            \addlegendentry{Exact};
            
            \addplot [thick, green!50!black]  
            table[y index=0, x expr=\coordindex*0.001 - .1] {figs/2d_final_density/burgers_godunov.tex};
            \addlegendentry{Godunov}

            \addplot [thick, red] 
            table[y index=0, x expr=\coordindex*0.001 - .1] {figs/2d_final_density/burgers_ours.tex};
            \addlegendentry{NFVM\(^5_4\)}

            \end{axis}
        \end{tikzpicture}
        \endgroup
    \end{subfigure}
    \caption{Comparison of the final state of $\text{NFVM}_4^5$ against Godunov on Burgers' equation, after 200 timesteps of autoregressive prediction.}
    \label{fig:burgers_2d_final_density}
\end{figure}

%% file: tex/drone_highway.tex
\begin{figure}[H]
	\centering
	\includegraphics[width=0.99\linewidth]{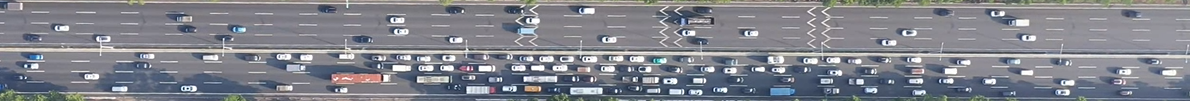}
	\caption{Snapshot from the drone-recorded highway video data (corresponding to a column in~\Cref{fig:drone_dataset_heatmap}).}
	\label{fig:drone_highway}
\end{figure}

%% file: tex/drone_tsd.tex
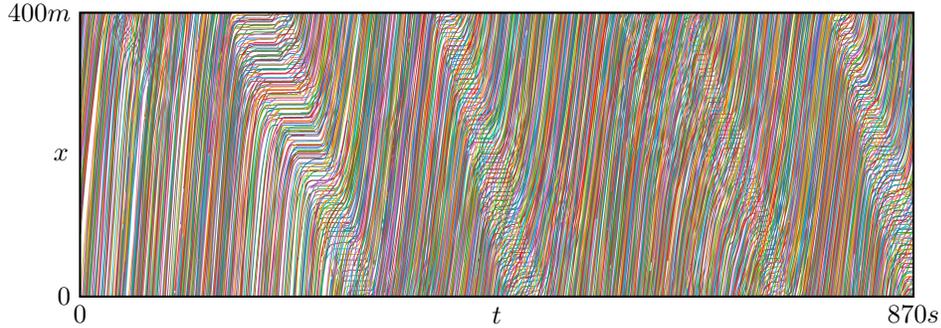
\begin{figure}[H]
	\centering
	\begin{tikzpicture}[scale=0.9, clip=false]
    \node[anchor=south west, inner sep=0] (img) at (0,0)
      {\includegraphics[width=0.85\textwidth]{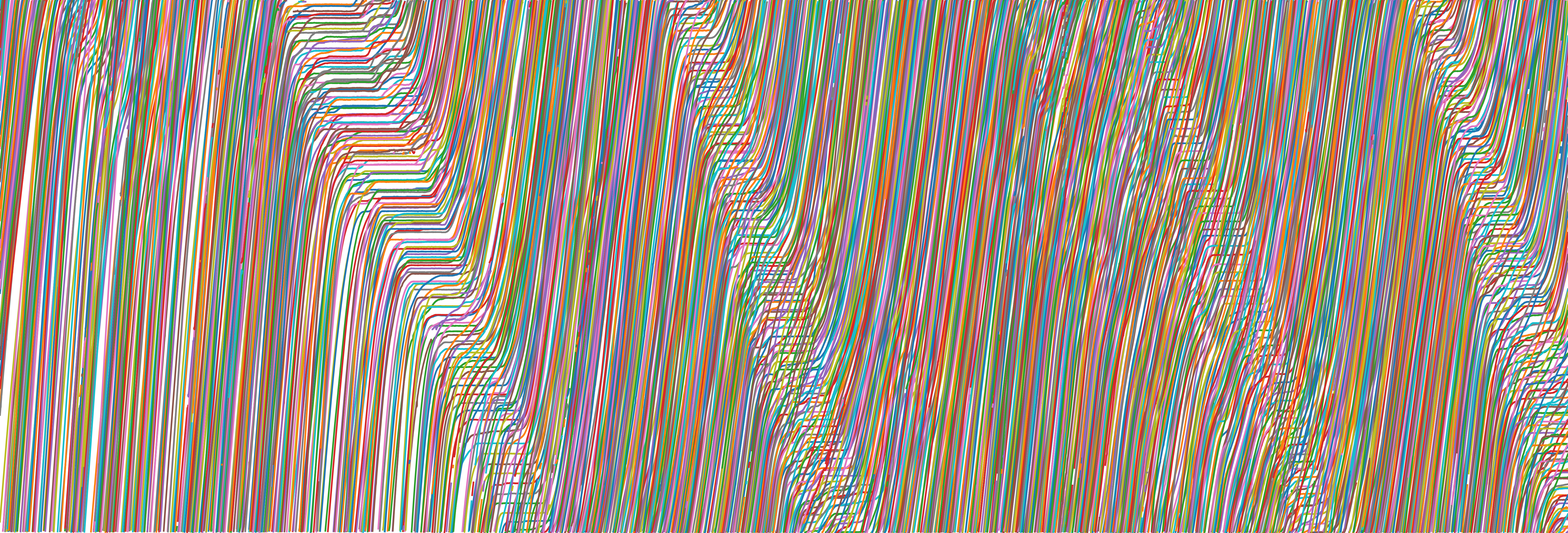}};

		\draw[black, thick] (img.south west) rectangle (img.north east);

		\path (img.south west) -- (img.south east) node[midway, anchor=center, yshift=-0.25cm] {\(t\)};
		\path (img.north west) -- (img.south west) node[midway, anchor=center, xshift=-0.25cm] {\(x\)};
		\node[below] at (img.south east) {\(870s\)};
		\node[left] at (img.north west) {\(400m\)};
		
		\node[below] at (0, 0) {0};

		\node[left] at (0, 0) {0};
	\end{tikzpicture}
	\caption{Space-time diagram of the drone dataset. Each line corresponds to one vehicle trajectory, where the $x$-axis spans along the highway. Note that we consider 4 lanes, explaining overlapping $x$-trajectories.}
	\label{fig:drone_tsd}
\end{figure}

%% file: tex/drone_dataset_heatmap.tex
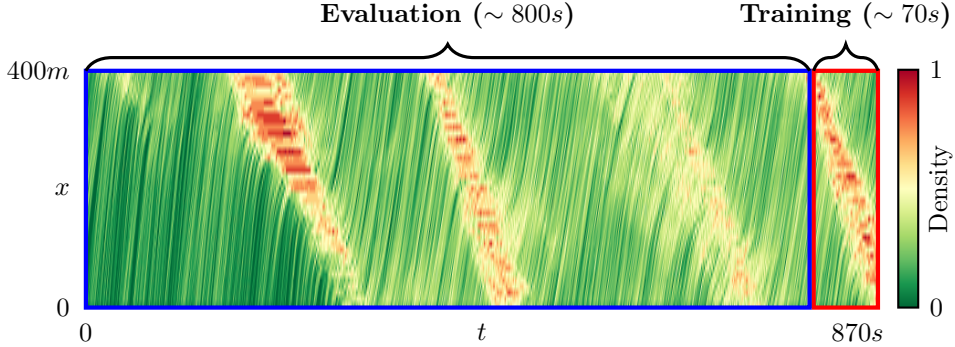
\begin{figure}
	\centering
	\begin{tikzpicture}[scale=0.9, clip=false]
		\node[anchor=south west, inner sep=0] (img) at (0,0) {\includegraphics[width=.81\textwidth]{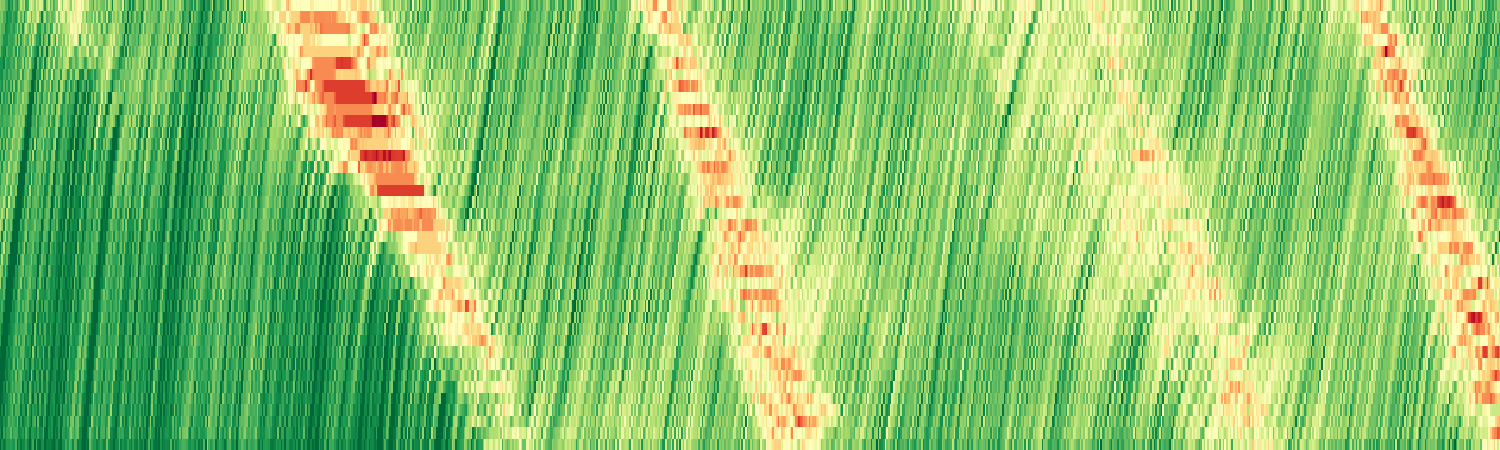}};
		
		\fill[red!50!black, opacity=0.0, draw=red, ultra thick] (10.75, 0) rectangle (11.7, 3.5);
		\draw[draw=red, ultra thick] (10.75, 0) rectangle (11.7, 3.5)
			node[red!30!black, midway, rotate=90] {}; 
		\fill[blue!50!black, opacity=0.0, draw=blue, ultra thick] (0, 0) rectangle (10.69, 3.5);
		\draw[draw=blue, ultra thick] (0, 0) rectangle (10.69, 3.5)
			node[blue!30!black, midway] {}; 
		\draw [
			very thick,
			decorate,
			decoration={brace,amplitude=10pt}
		] (0, 3.5) -- (10.69, 3.5)
			node [midway, above=12pt] {\textbf{Evaluation (\(\sim 800s\))}};
		\draw [
			very thick,
			decorate,
			decoration={brace,amplitude=10pt}
		] (10.75, 3.5) -- (11.7, 3.5)
			node [midway, above=12pt, xshift=-0.0cm] {\textbf{Training (\(\sim 70s\))}};

		\node[left] at (-0.1, 0) {0};
		\node[left] at (-0.1, 1.75) {\(x\)};
		\node[left] at (-0.1, 3.5) {\(400m\)};

		\node[below] at (0, -0.1) {0};
		\node[below] at (5.85, -0.1) {\(t\)};
		\node[below] at (11.4, -0.1) {\(870s\)};

		\node[anchor=south west, inner sep=0] (colorbar) at (12.0,0) {\includegraphics[height=3.15cm]{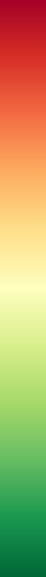}};
		\draw[black, thick] (colorbar.south west) rectangle (colorbar.north east);
		\path (colorbar.south east) -- (colorbar.north east) node[midway, rotate=90, anchor=center, xshift=0.0cm, yshift=-0.3cm] {Density};
		\node[right, xshift=+0.03cm] at (colorbar.south east) {0};
		\node[right, xshift=+0.03cm] at (colorbar.north east) {1};
		
	\end{tikzpicture}
	\caption{Time-space diagram of vehicle trajectories extracted from the drone video, color-coded by density.}
    \label{fig:drone_dataset_heatmap}
\end{figure}

%% file: tex/fundamental_diagram.tex
\begin{figure}[t]
    \centering
    \begin{tikzpicture}
        \begin{axis}[
            xlabel={Density (vehicules/miles/lane)},
            ylabel={Flow (vehicules/hour/lane)},
            grid=both,
            only marks, 
            width=11cm,
            height=6cm,
            scaled y ticks=false,
          ]
            \addplot[
              mark=*,
              mark size=1pt,
              blue,
              opacity=0.3
            ] table[
              col sep=comma,
              x expr=\thisrow{density}/4,
              y expr=\thisrow{flow}/4
            ] {figs/fundamental_diagram/data.csv};
          \end{axis}   
    \end{tikzpicture}
    \caption{Fundamental diagram of the traffic flow for the highway dataset.}
    \label{fig:fundamental_diagram}
\end{figure}
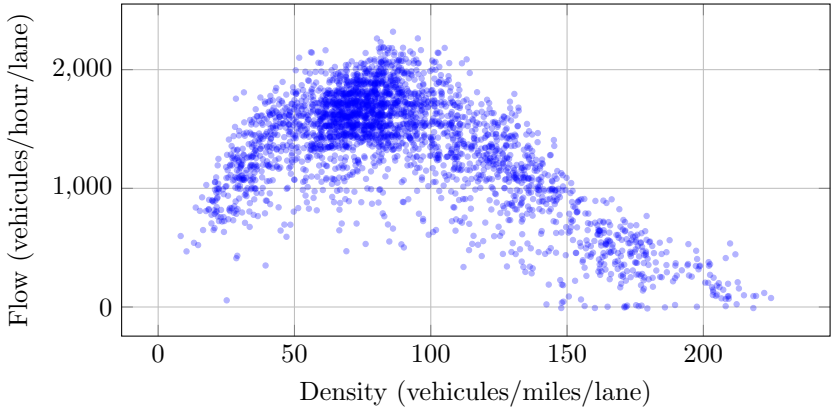

%% file: tex/drone_fitted_fd.tex
\begin{figure}
    \centering
    
    \vspace{-0.4cm}
    \caption*{\textbf{Calibration on scheme prediction}}
    \vspace{0.1cm}
    \small
    \raisebox{0.61cm}{\rotatebox[origin=c]{90}{GD}}
    \begin{subfigure}[t]{.193\textwidth}
      \centering
      \includegraphics[width=\textwidth]{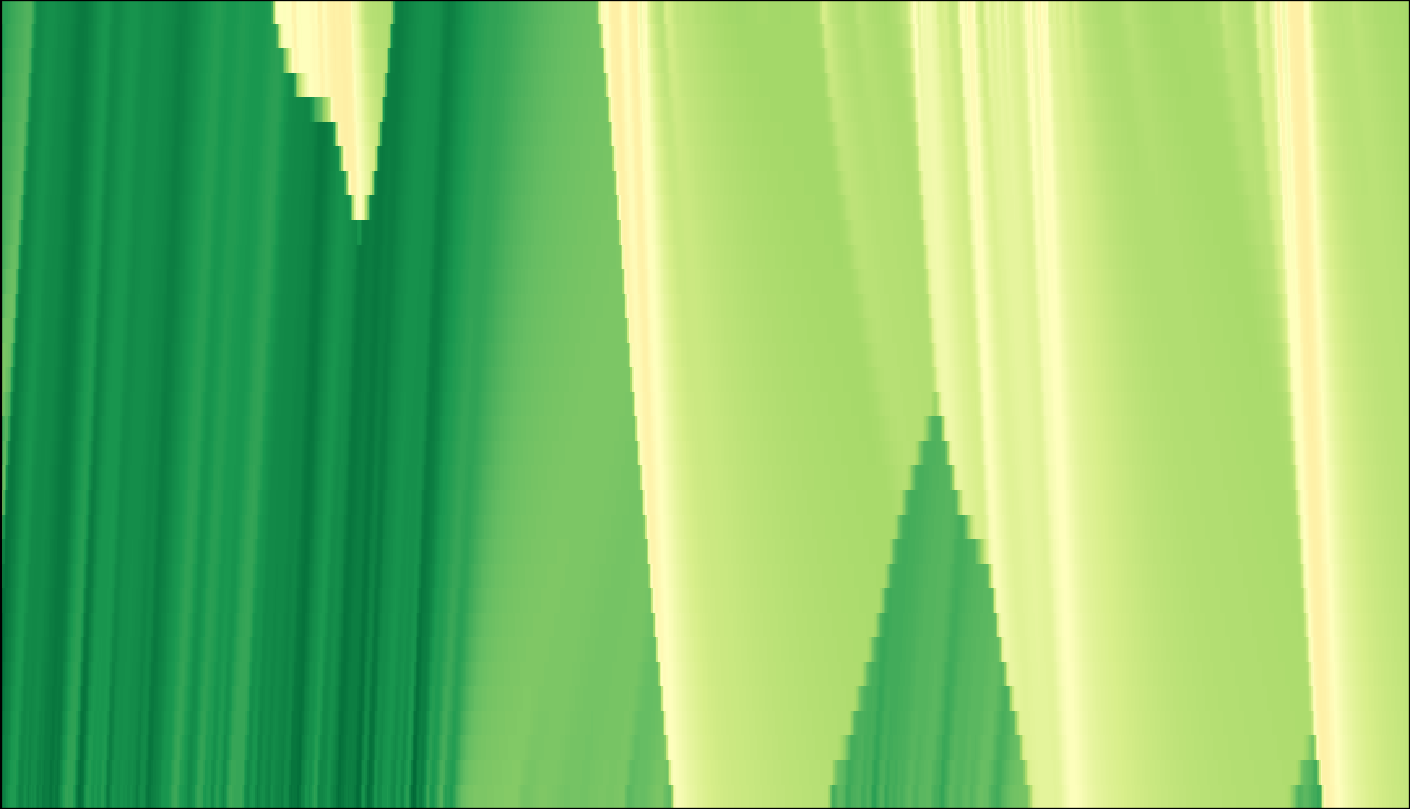}
    \end{subfigure}\hfill
    \begin{subfigure}[t]{.193\textwidth}
      \centering
      \includegraphics[width=\textwidth]{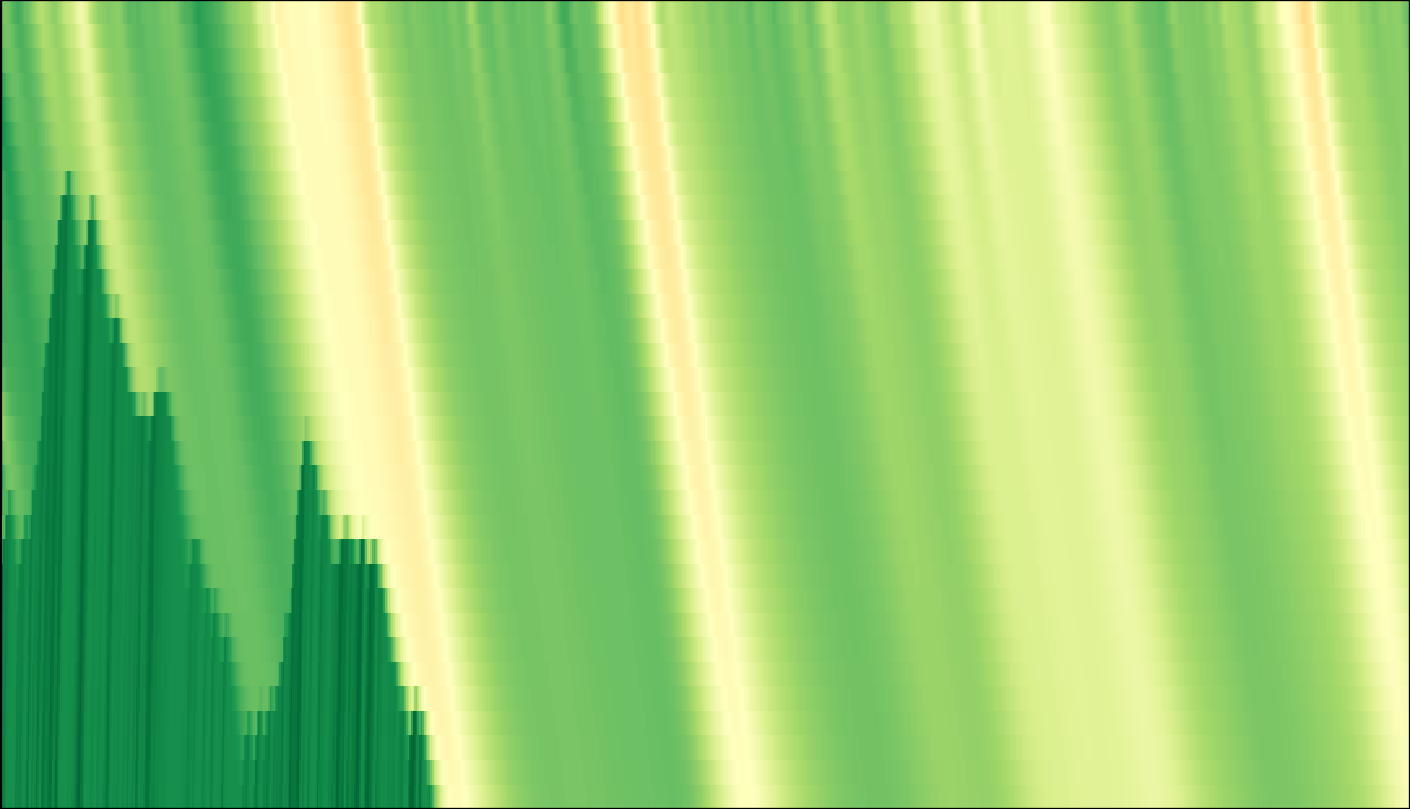}
    \end{subfigure}\hfill
    \begin{subfigure}[t]{.193\textwidth}
      \centering
      \includegraphics[width=\textwidth]{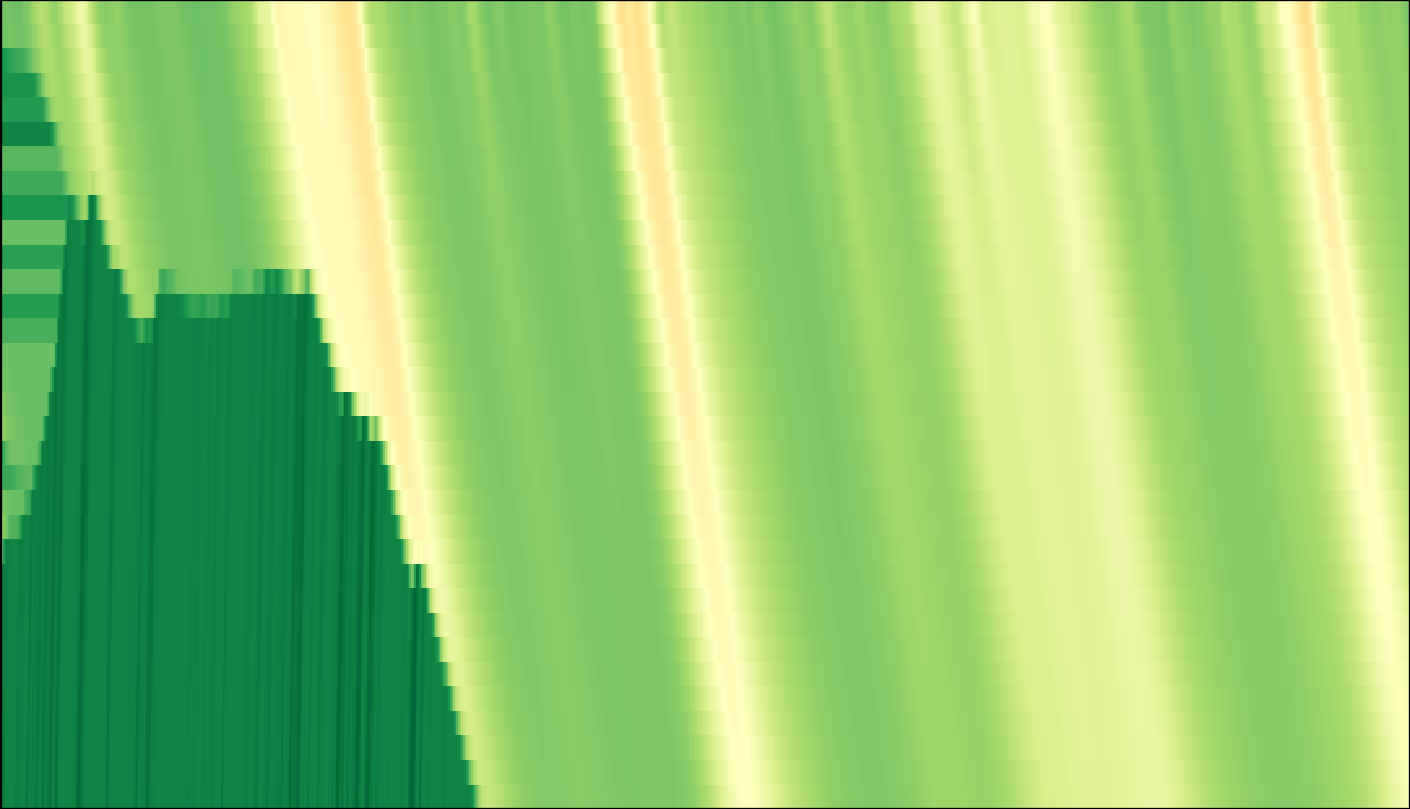}
    \end{subfigure}\hfill
    \begin{subfigure}[t]{.193\textwidth}
      \centering
      \includegraphics[width=\textwidth]{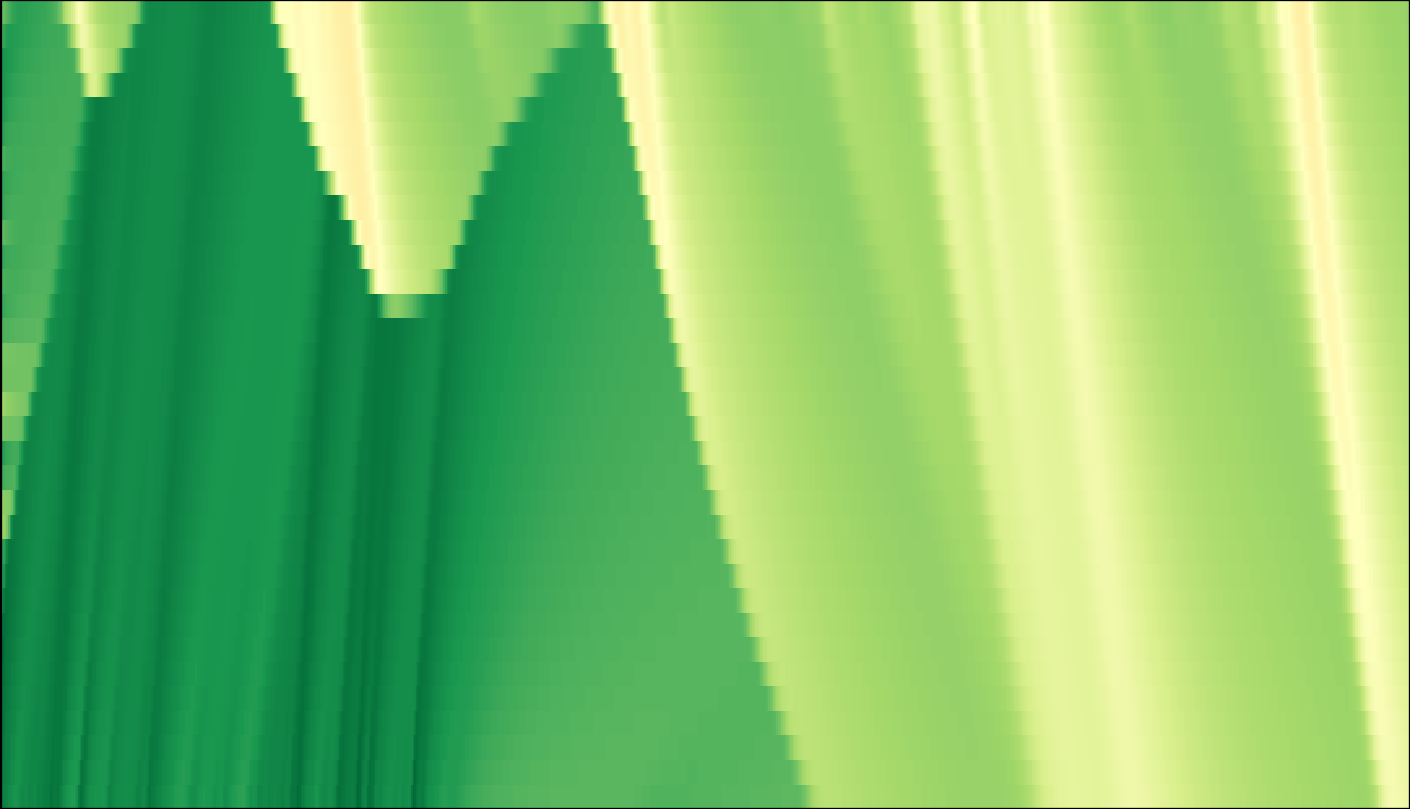}
    \end{subfigure}\hfill
    \begin{subfigure}[t]{.193\textwidth}
      \centering
      \includegraphics[width=\textwidth]{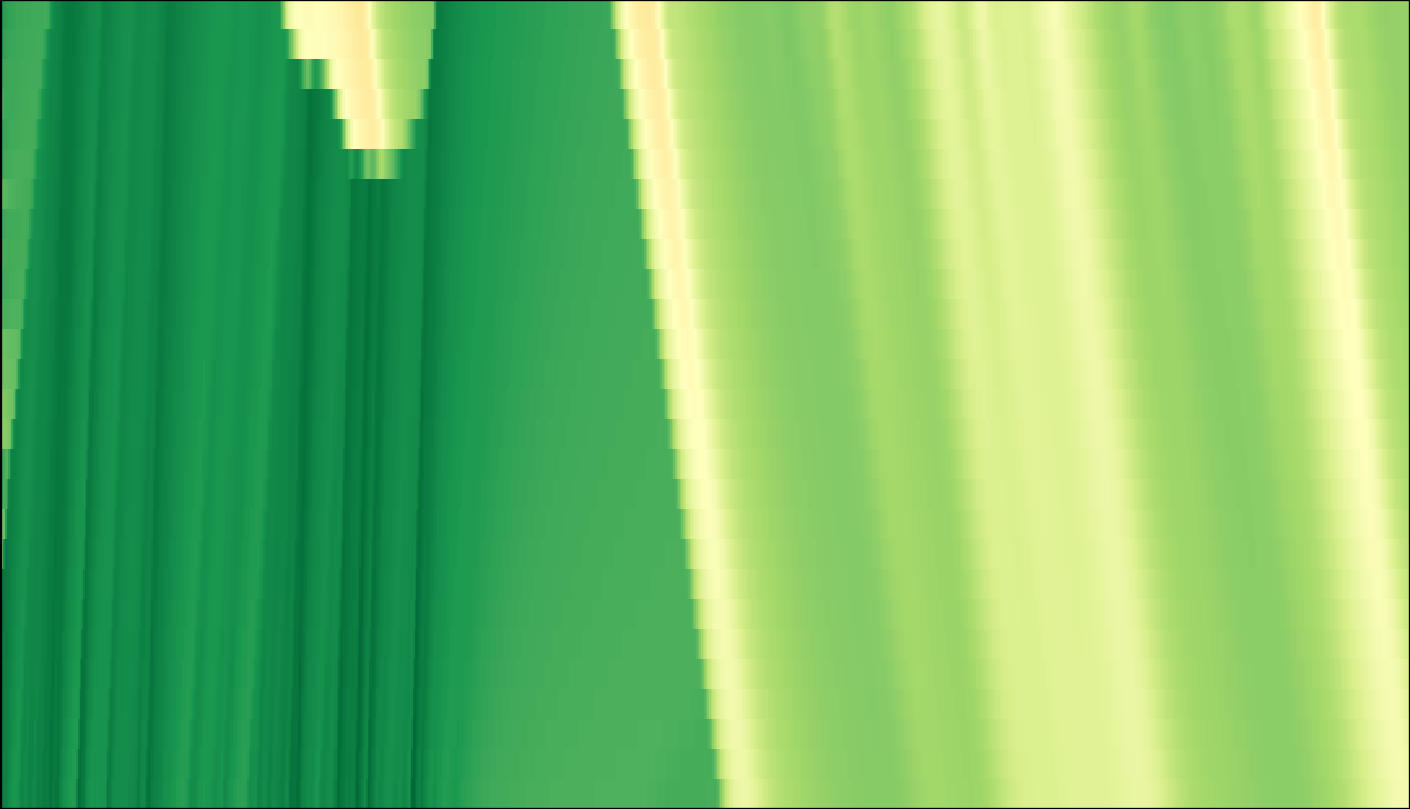}
    \end{subfigure}
  
    \raisebox{0.61cm}{\rotatebox[origin=c]{90}{LxF}}
    \begin{subfigure}[t]{.193\textwidth}
      \centering
      \includegraphics[width=\textwidth]{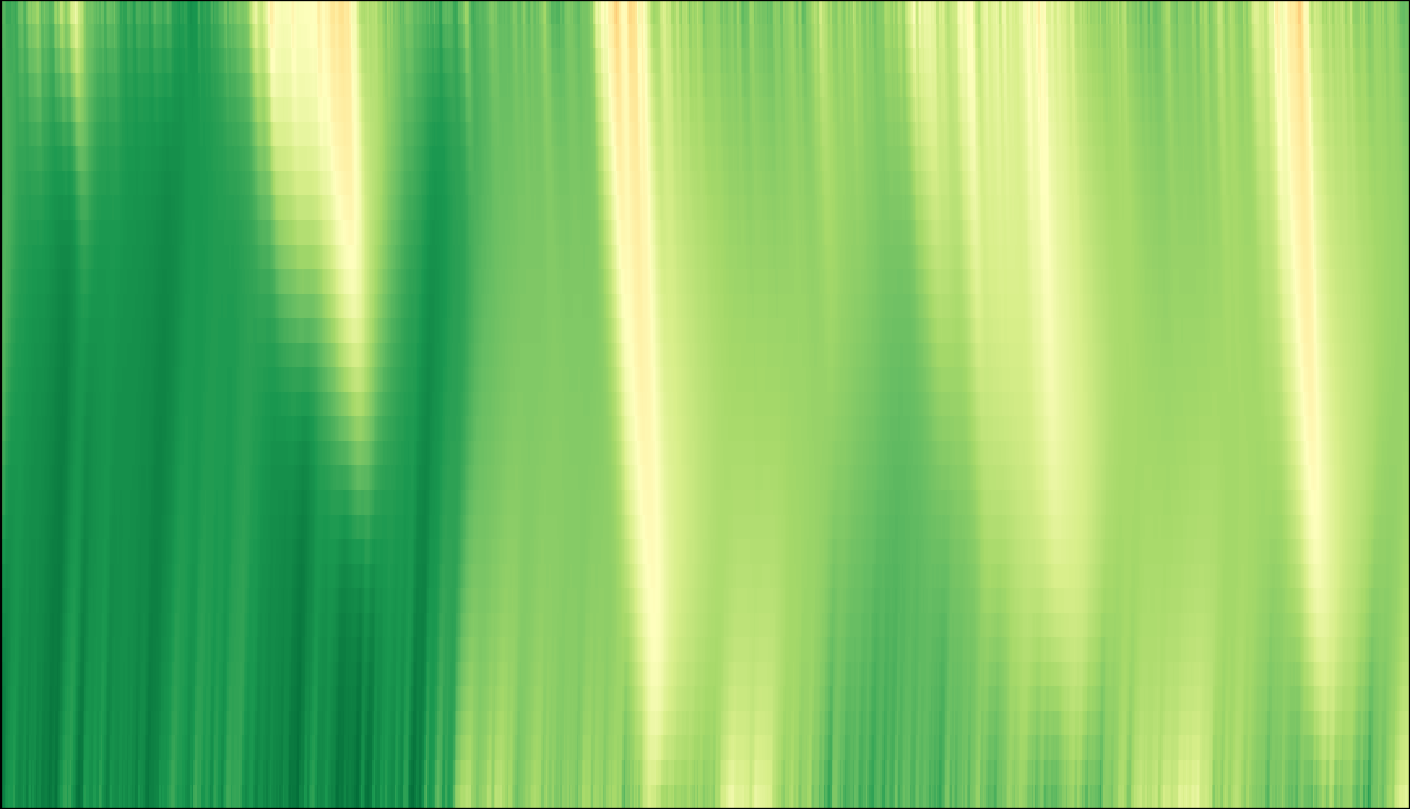}
    \end{subfigure}\hfill
    \begin{subfigure}[t]{.193\textwidth}
      \centering
      \includegraphics[width=\textwidth]{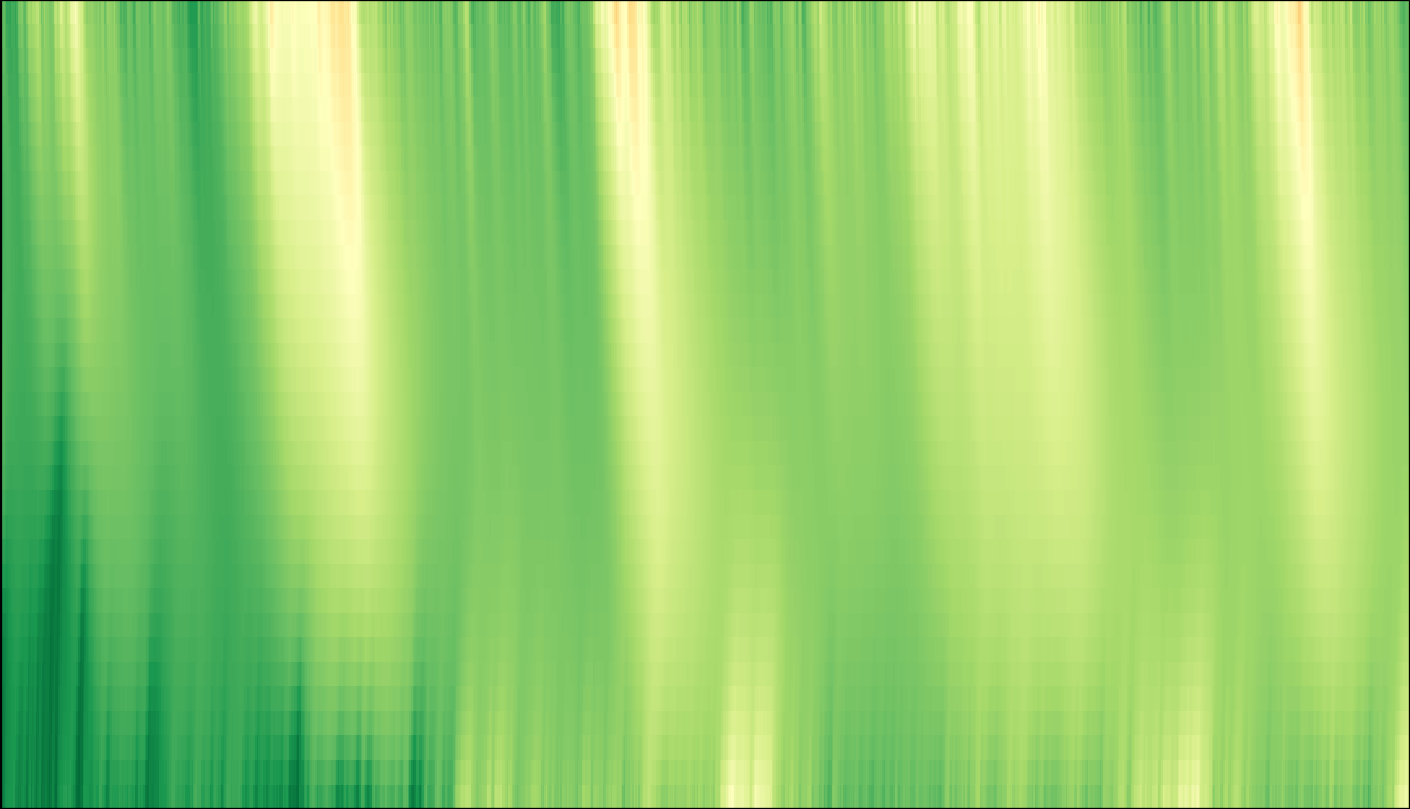}
    \end{subfigure}\hfill
    \begin{subfigure}[t]{.193\textwidth}
      \centering
      \includegraphics[width=\textwidth]{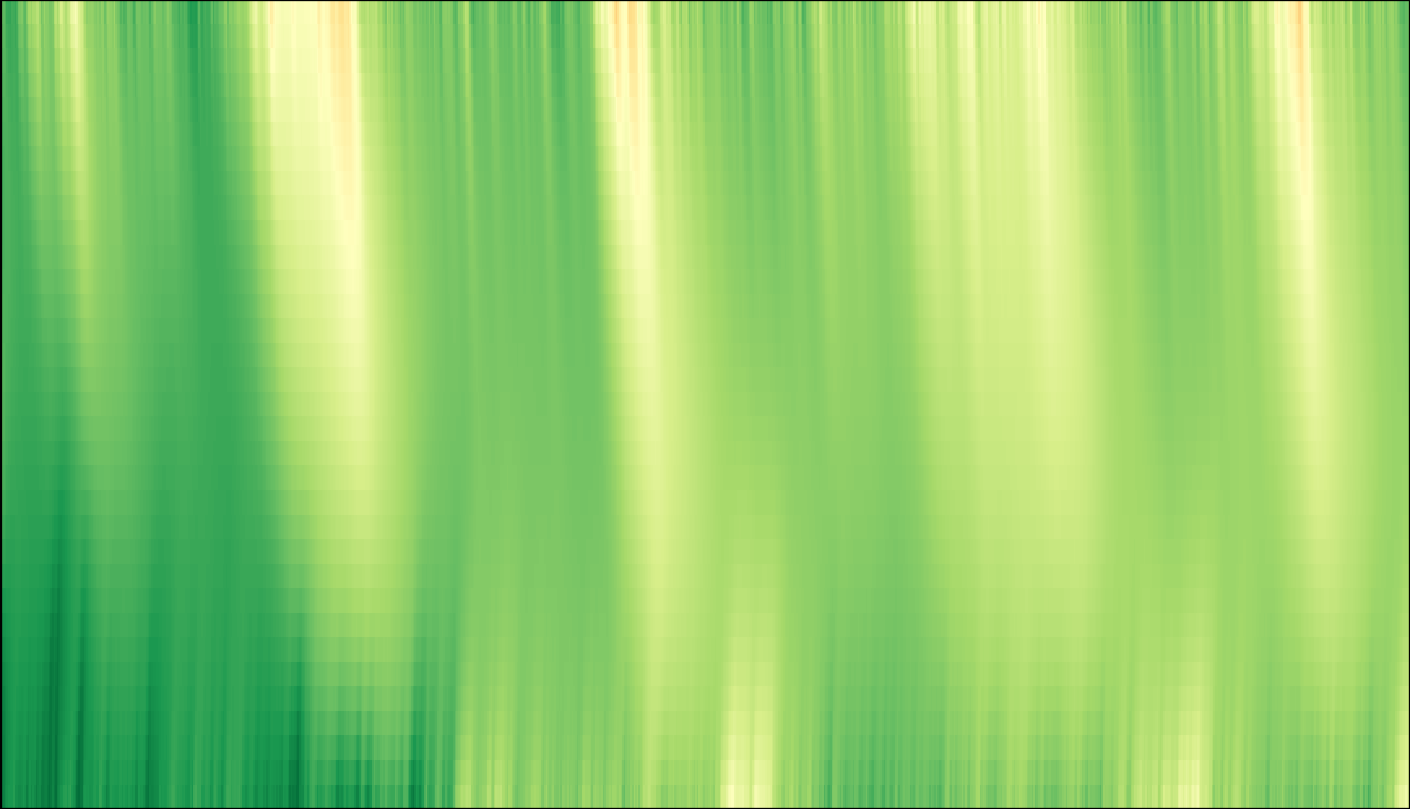}
    \end{subfigure}\hfill
    \begin{subfigure}[t]{.193\textwidth}
      \centering
      \includegraphics[width=\textwidth]{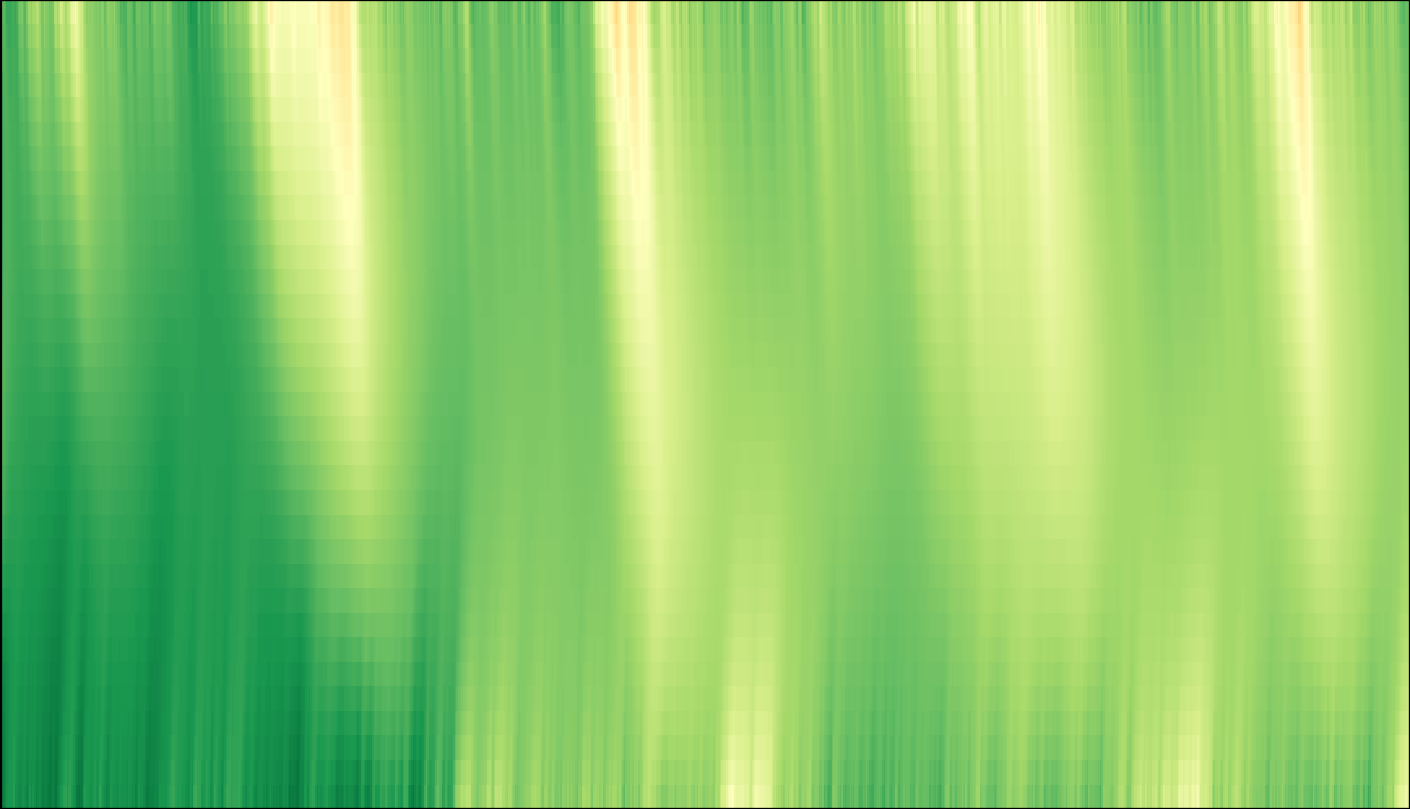}
    \end{subfigure}\hfill
    \begin{subfigure}[t]{.193\textwidth}
      \centering
      \includegraphics[width=\textwidth]{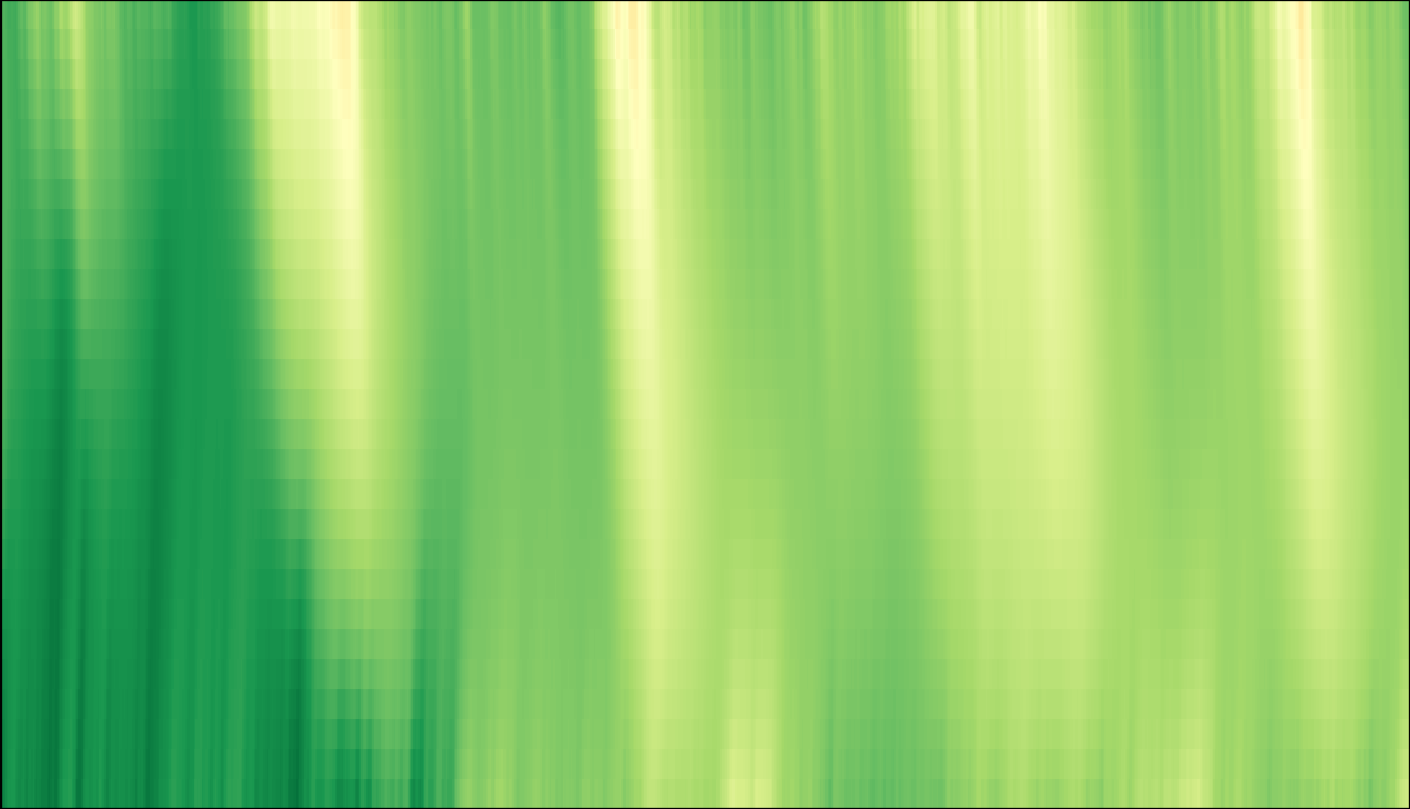}
    \end{subfigure}
  
    \raisebox{0.61cm}{\rotatebox[origin=c]{90}{EO}}
    \begin{subfigure}[t]{.193\textwidth}
      \centering
      \includegraphics[width=\textwidth]{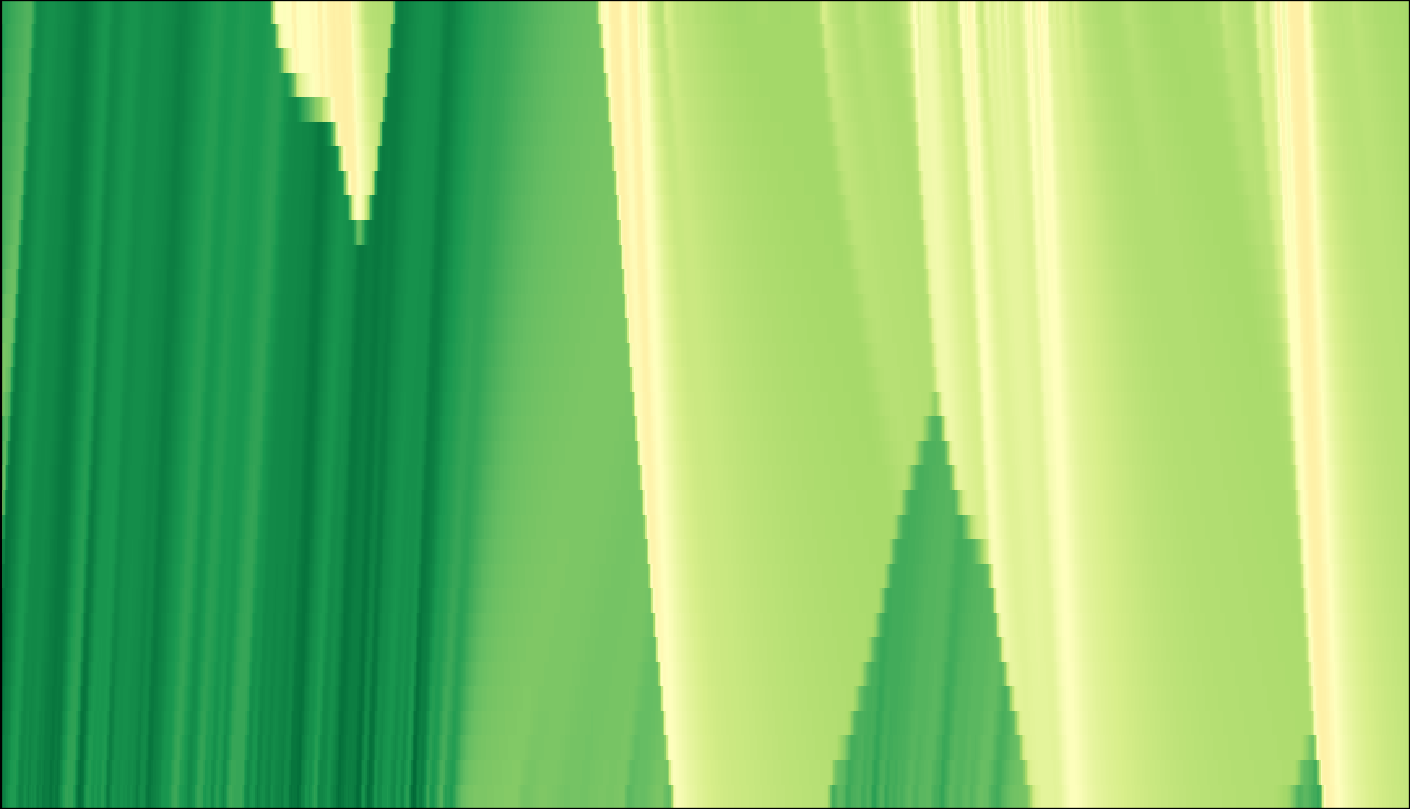}
    \end{subfigure}\hfill
    \begin{subfigure}[t]{.193\textwidth}
      \centering
      \includegraphics[width=\textwidth]{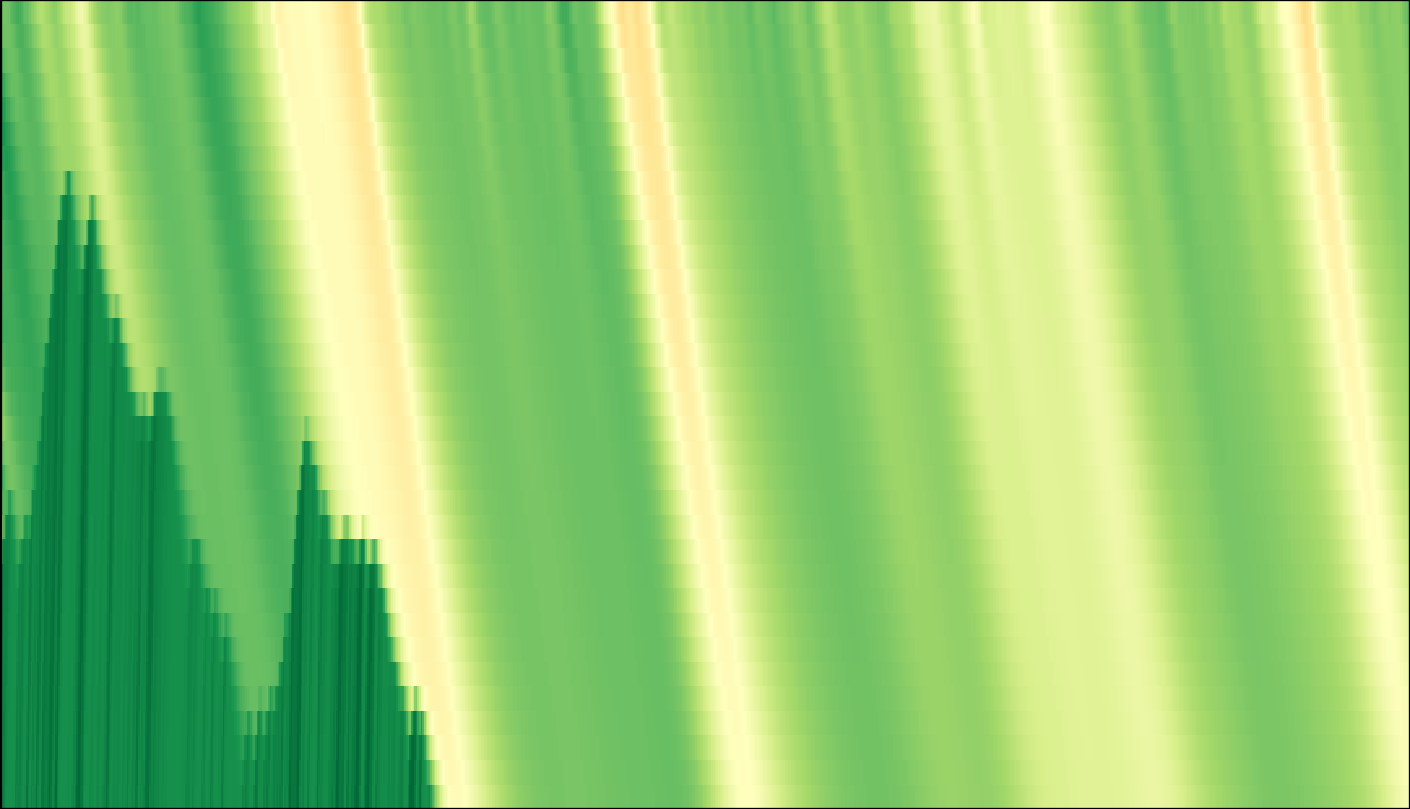}
    \end{subfigure}\hfill
    \begin{subfigure}[t]{.193\textwidth}
      \centering
      \includegraphics[width=\textwidth]{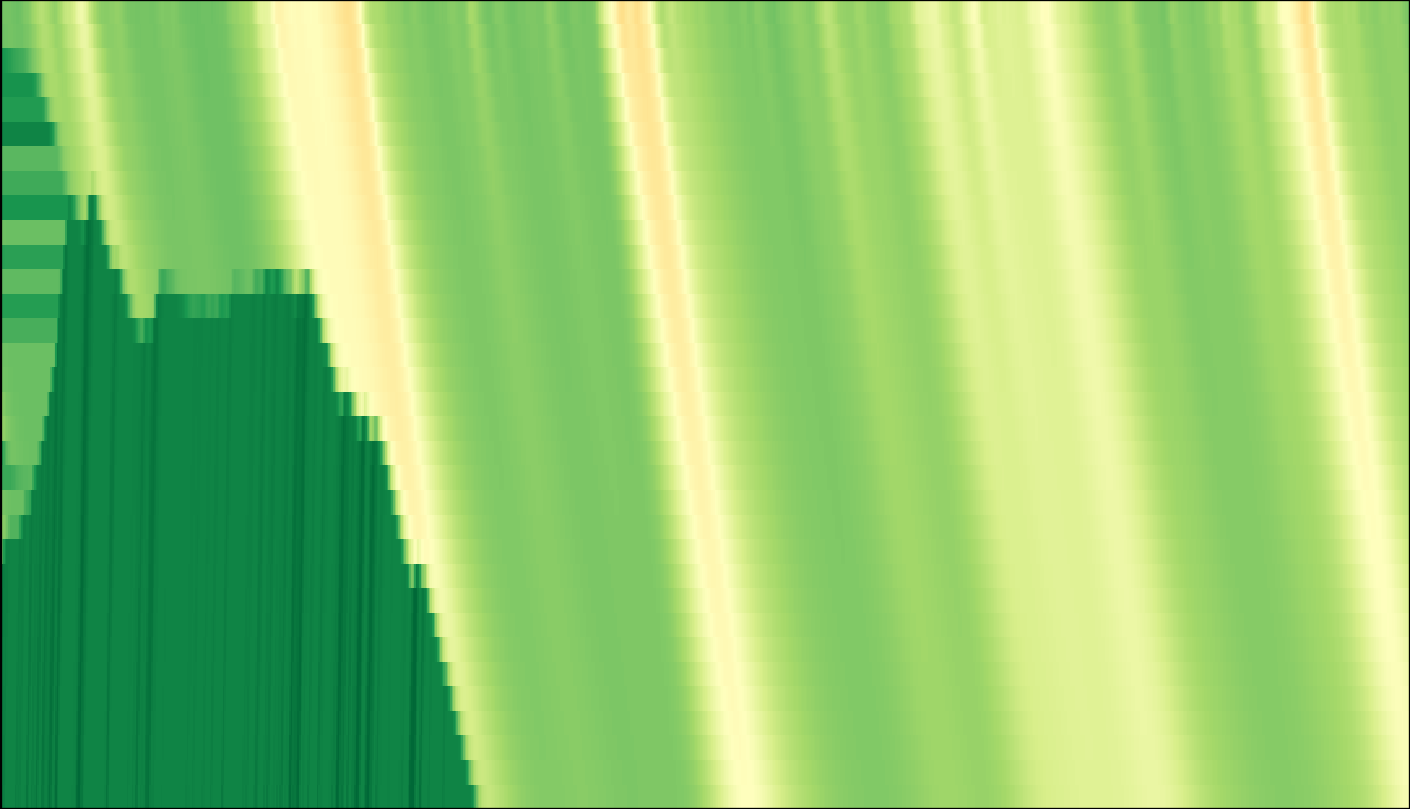}
    \end{subfigure}\hfill
    \begin{subfigure}[t]{.193\textwidth}
      \centering
      \includegraphics[width=\textwidth]{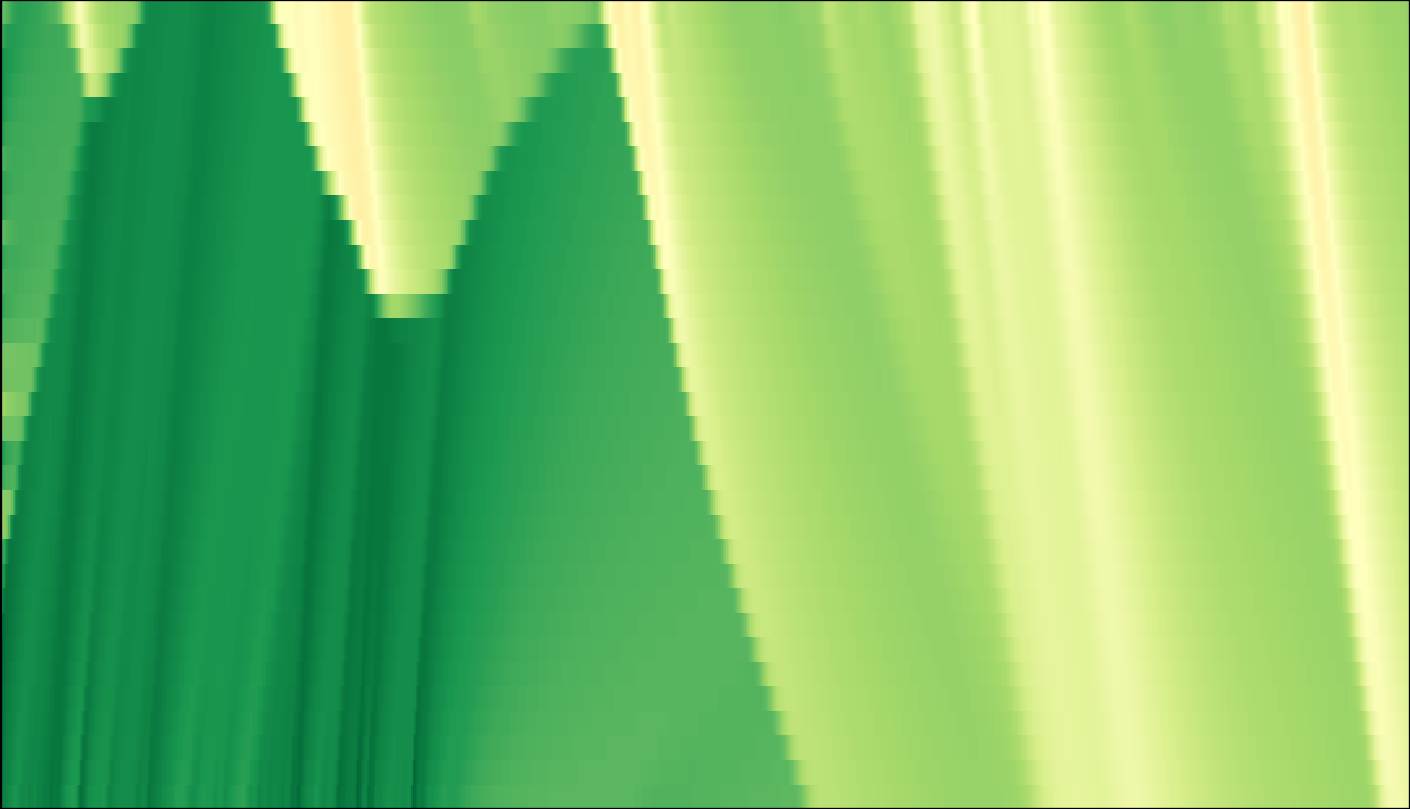}
    \end{subfigure}\hfill
    \begin{subfigure}[t]{.193\textwidth}
      \centering
      \includegraphics[width=\textwidth]{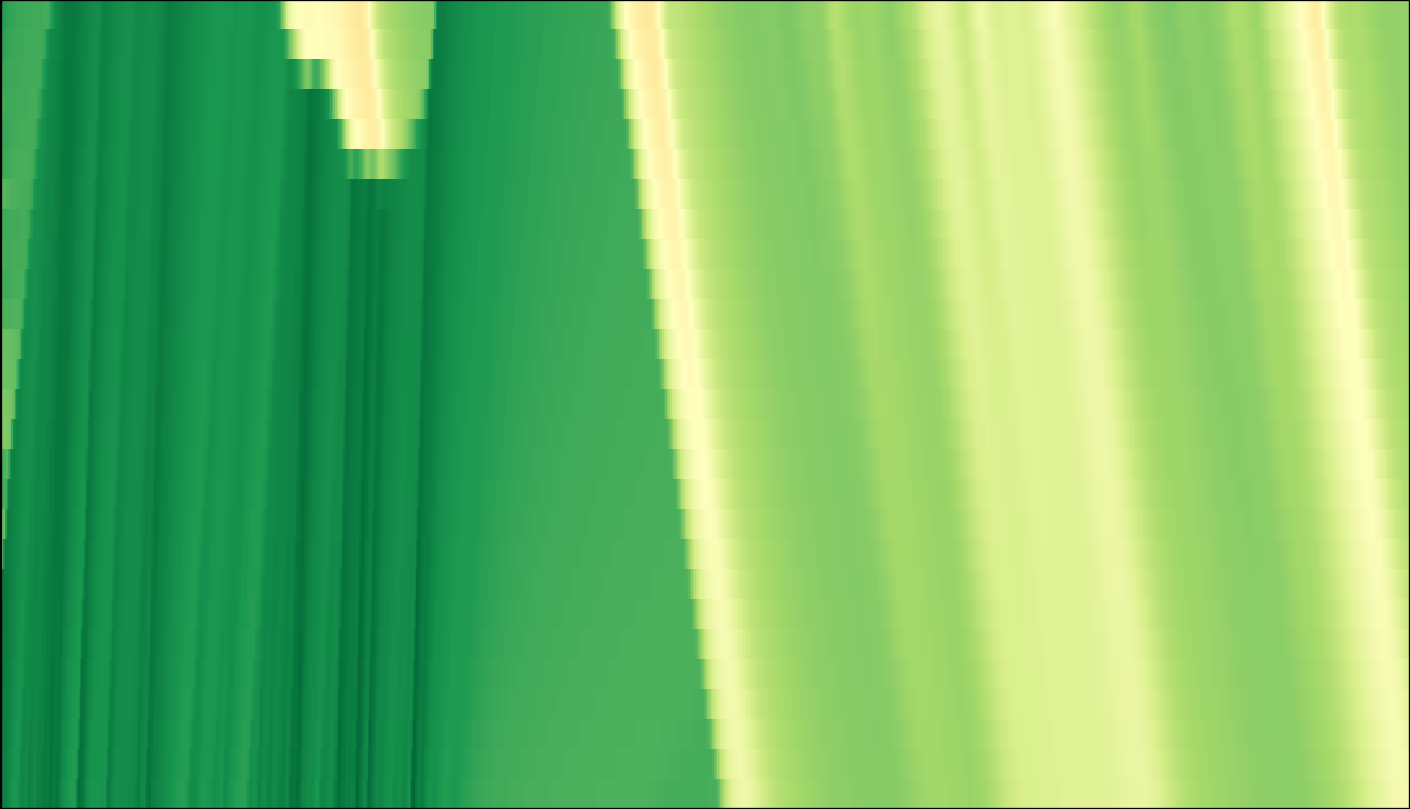}
    \end{subfigure}
  
    \raisebox{0.61cm}{\rotatebox[origin=c]{90}{ENO}}
    \begin{subfigure}[t]{.193\textwidth}
      \centering
      \includegraphics[width=\textwidth]{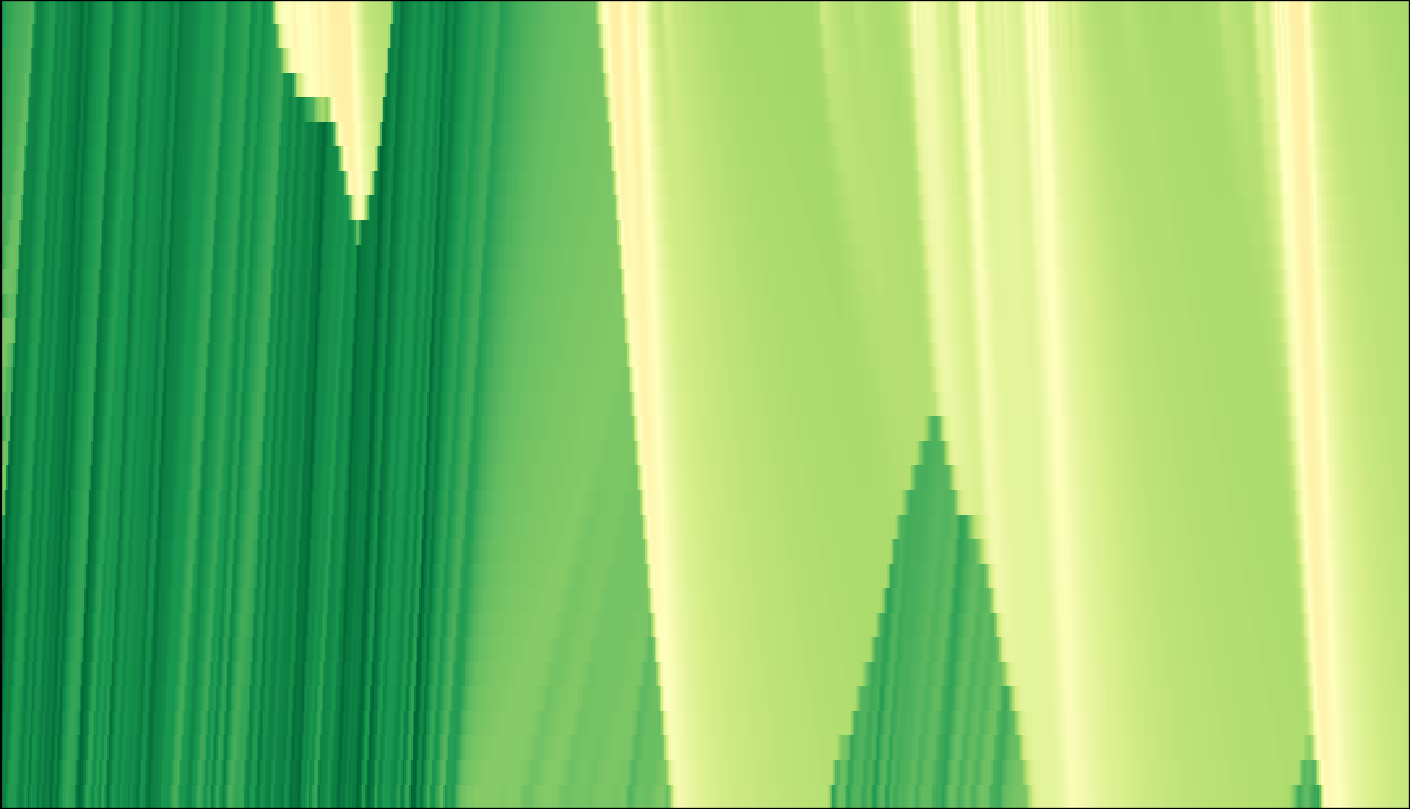}
    \end{subfigure}\hfill
    \begin{subfigure}[t]{.193\textwidth}
      \centering
      \includegraphics[width=\textwidth]{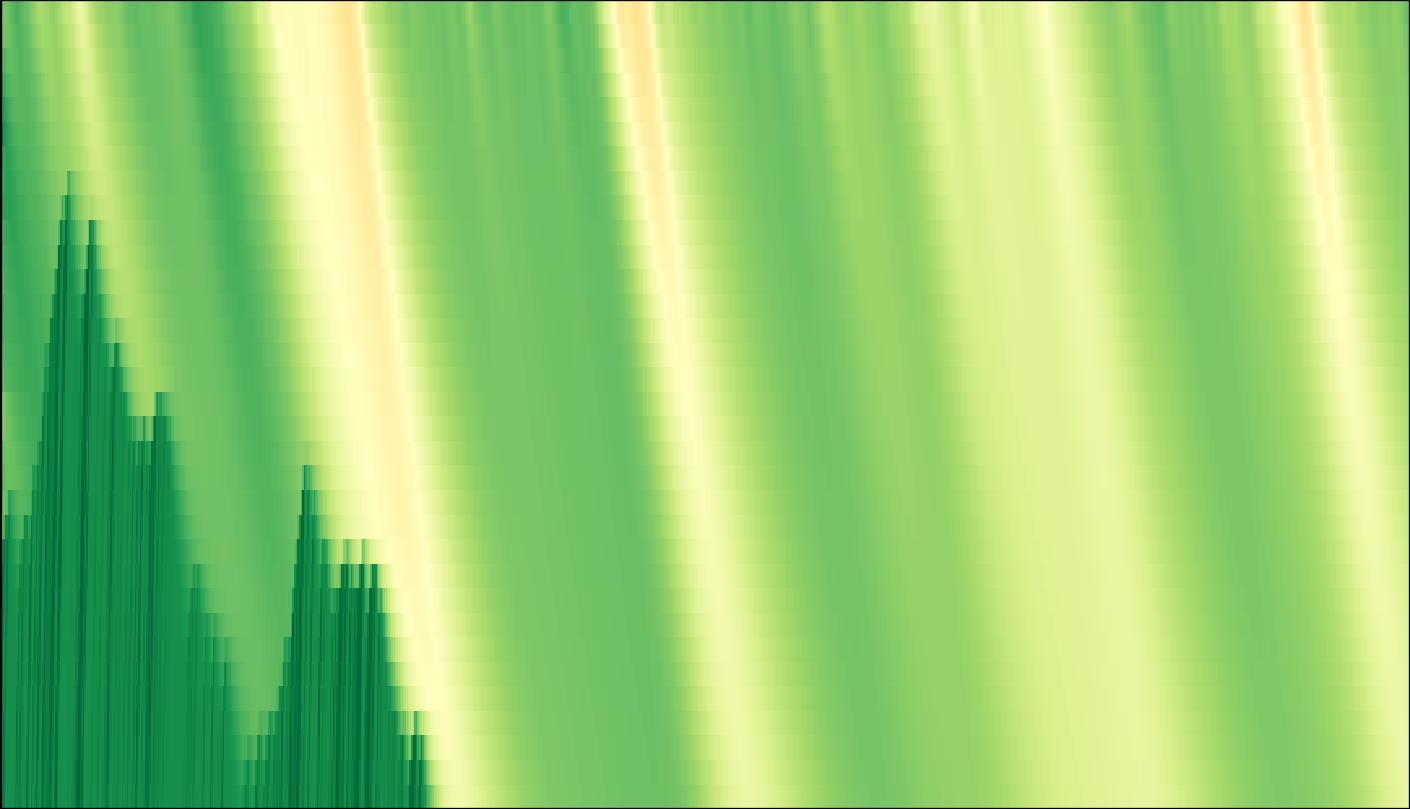}
    \end{subfigure}\hfill
    \begin{subfigure}[t]{.193\textwidth}
      \centering
      \includegraphics[width=\textwidth]{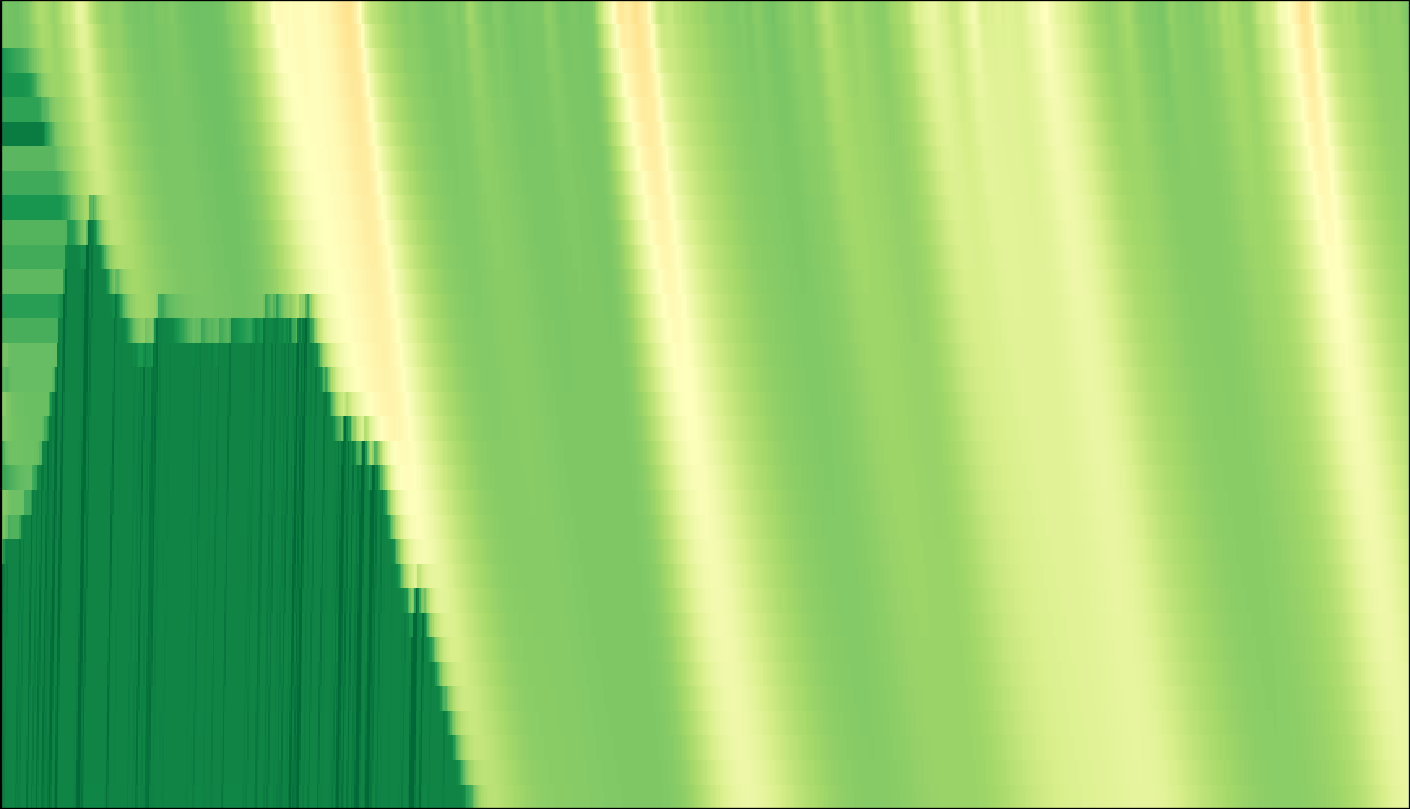}
    \end{subfigure}\hfill
    \begin{subfigure}[t]{.193\textwidth}
      \centering
      \includegraphics[width=\textwidth]{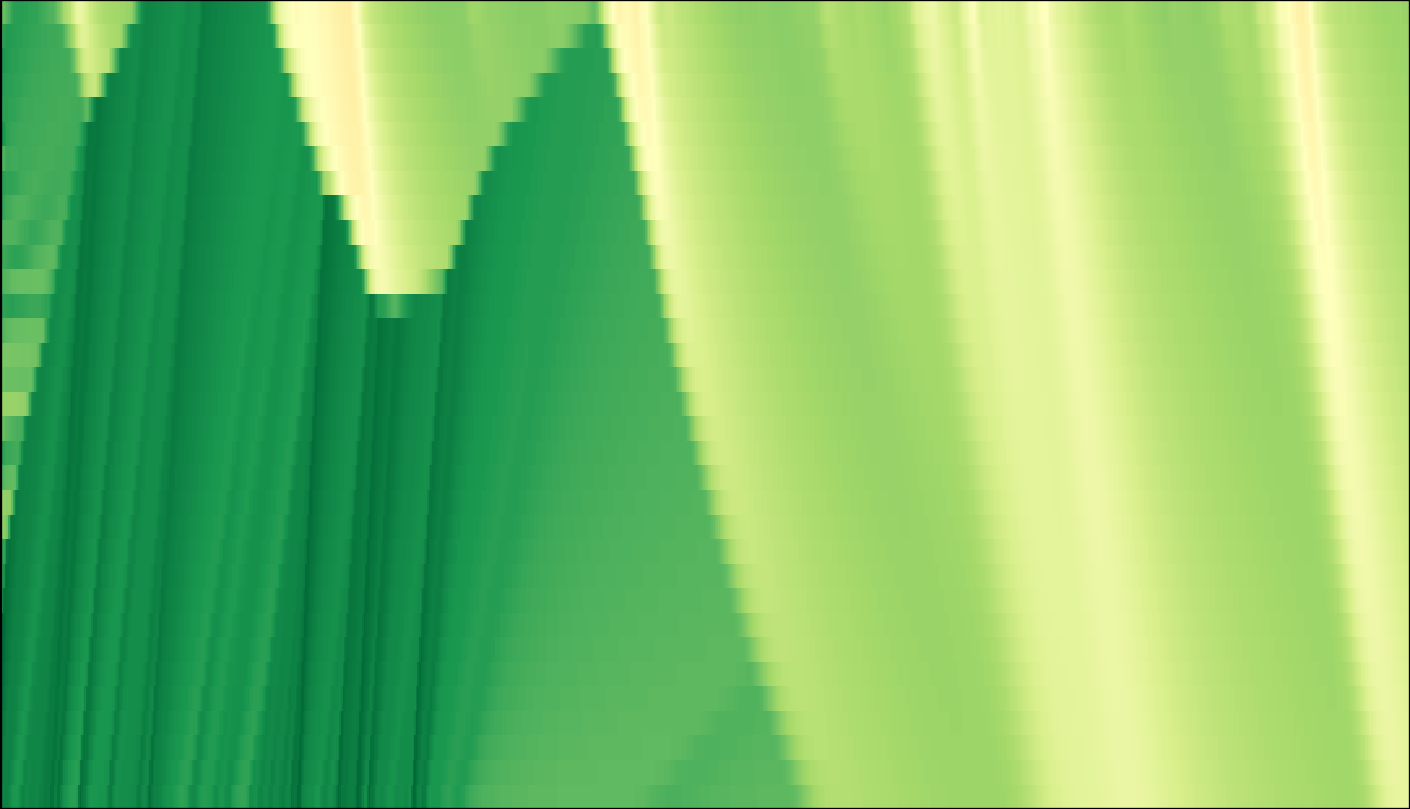}
    \end{subfigure}\hfill
    \begin{subfigure}[t]{.193\textwidth}
      \centering
      \includegraphics[width=\textwidth]{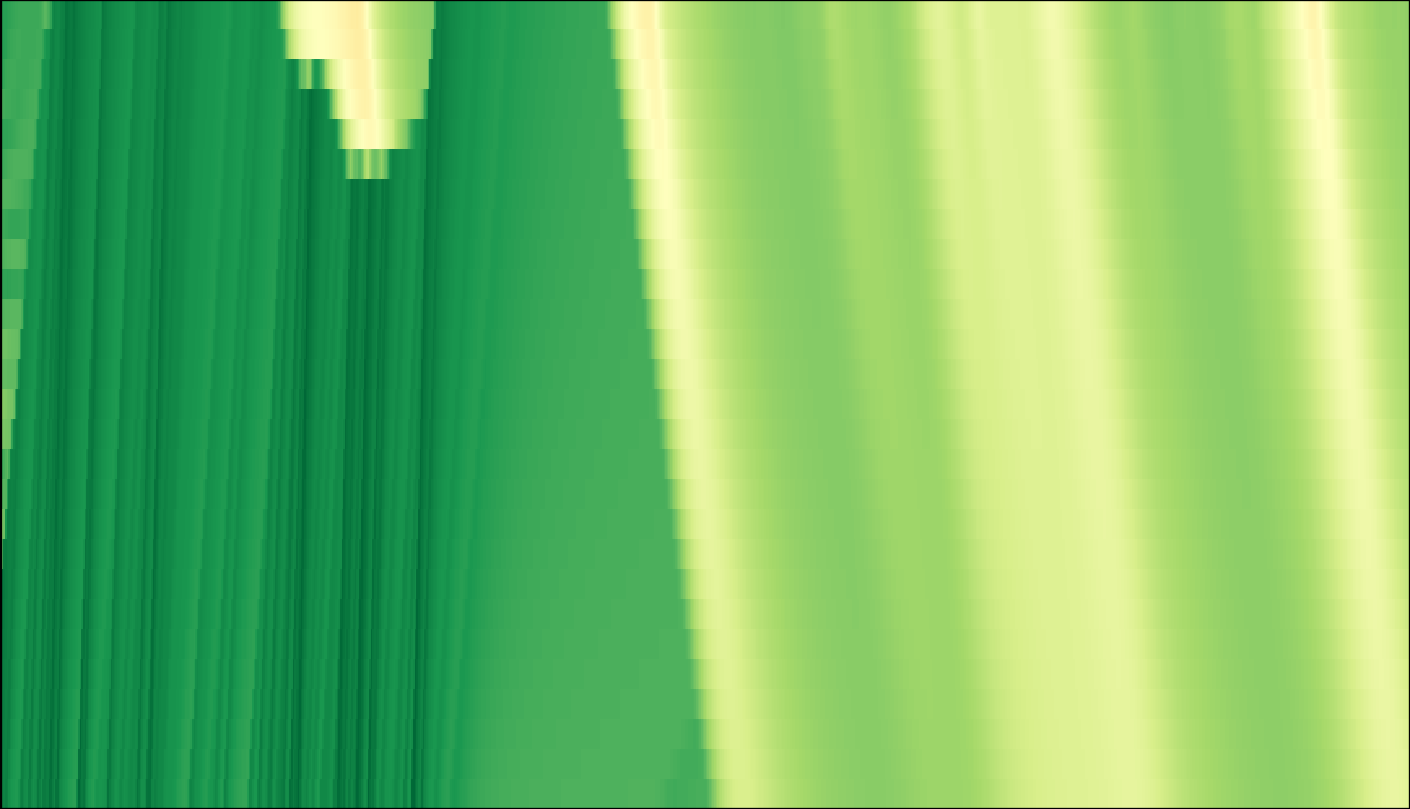}
    \end{subfigure}
  
    \raisebox{0.61cm}{\rotatebox[origin=c]{90}{WENO}}
    \begin{subfigure}[t]{.193\textwidth}
      \centering
      \includegraphics[width=\textwidth]{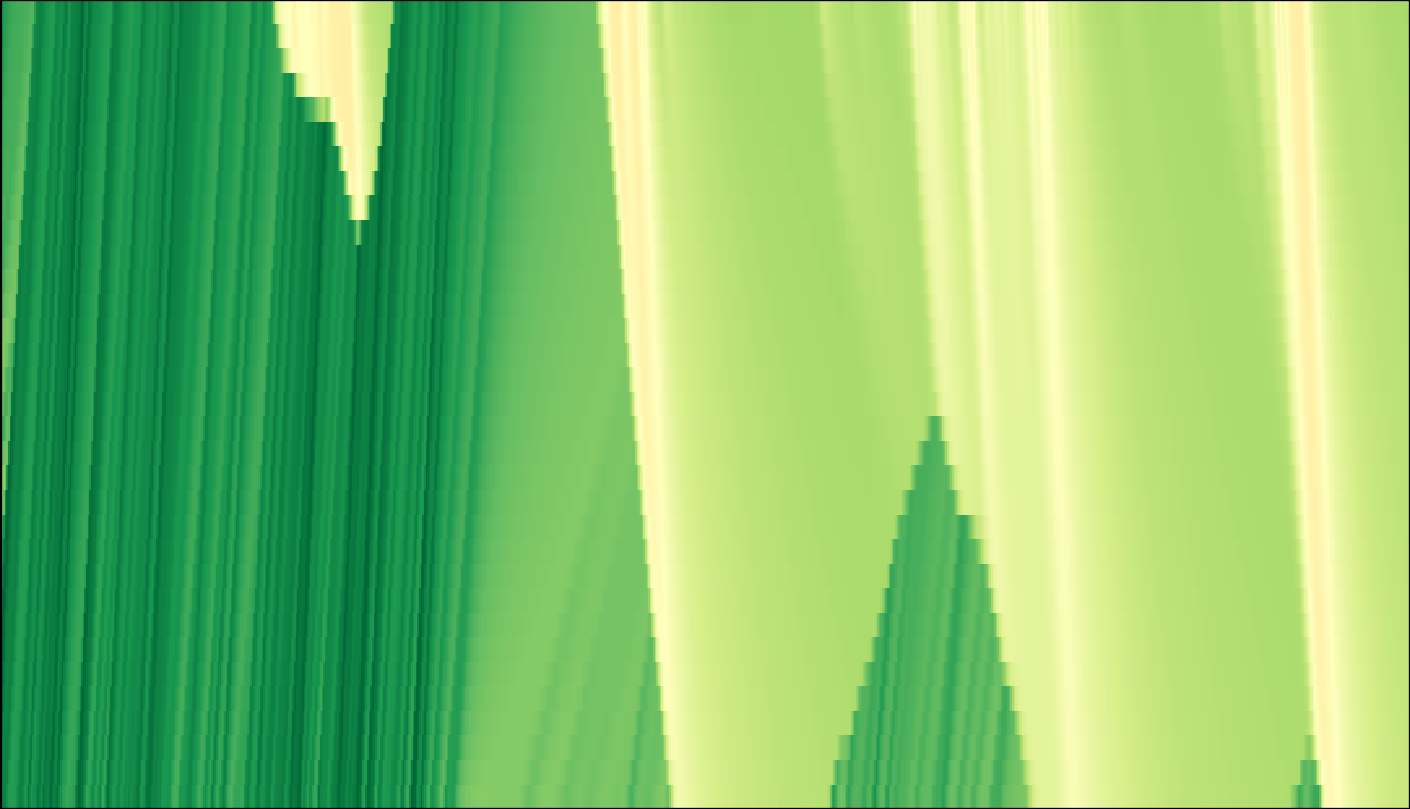}
      \vspace{-8.65cm}\caption*{Greenshields'}
    \end{subfigure}\hfill
    \begin{subfigure}[t]{.193\textwidth}
      \centering
      \includegraphics[width=\textwidth]{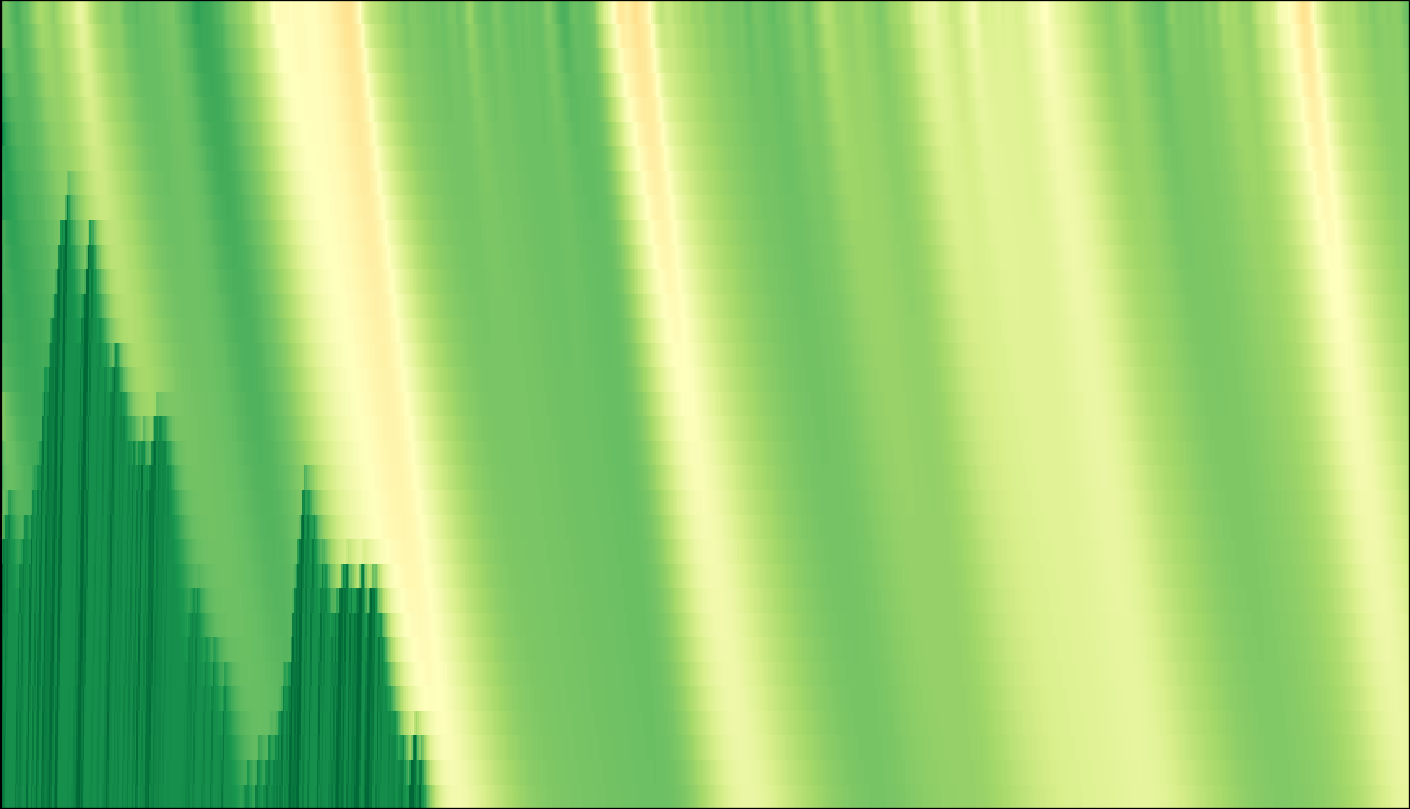}
      \vspace{-8.65cm}\caption*{Triangular}
    \end{subfigure}\hfill
    \begin{subfigure}[t]{.193\textwidth}
      \centering
      \includegraphics[width=\textwidth]{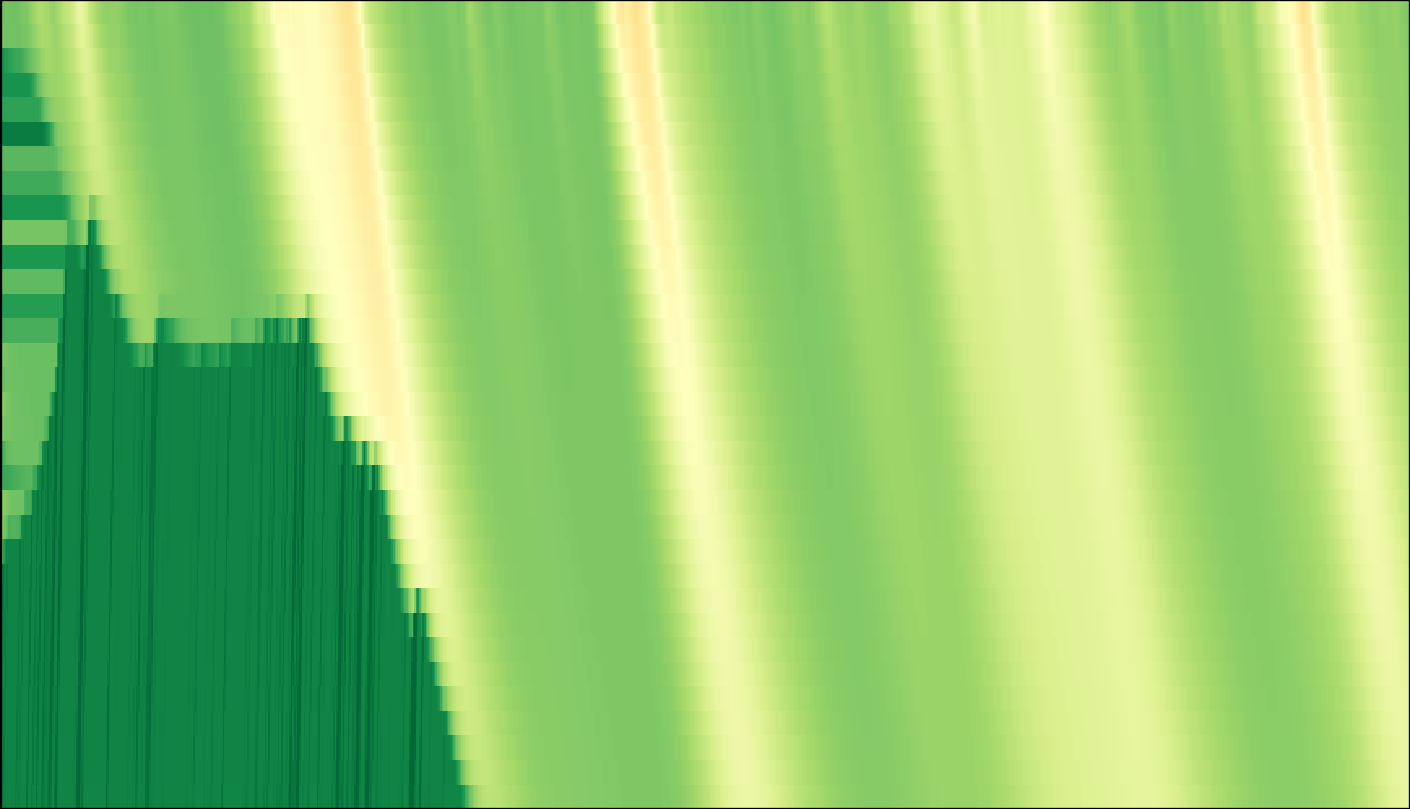}
      \vspace{-8.65cm}\caption*{Trapezoidal}
    \end{subfigure}\hfill
    \begin{subfigure}[t]{.193\textwidth}
      \centering
      \includegraphics[width=\textwidth]{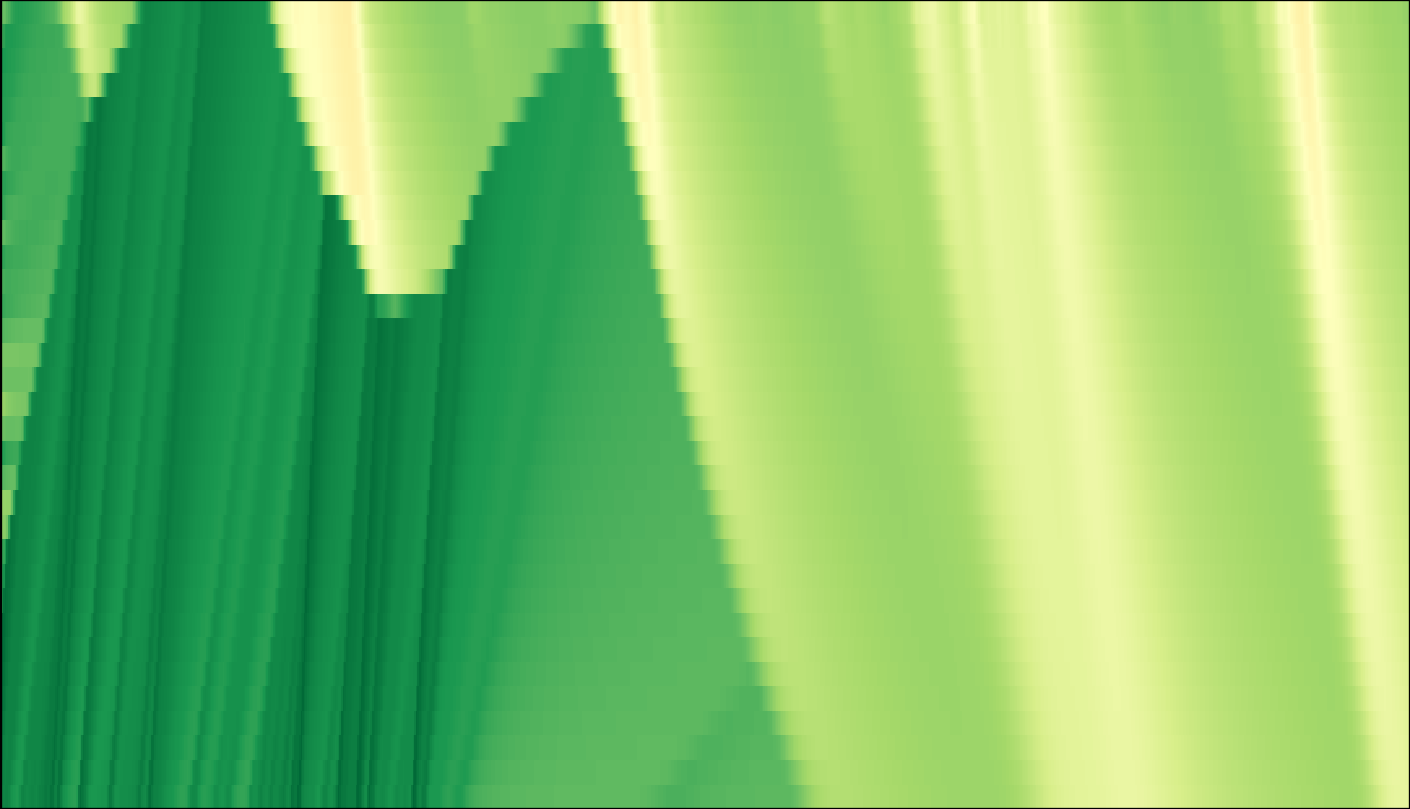}
      \vspace{-8.65cm}\caption*{Greenberg}
    \end{subfigure}\hfill
    \begin{subfigure}[t]{.193\textwidth}
      \centering
      \includegraphics[width=\textwidth]{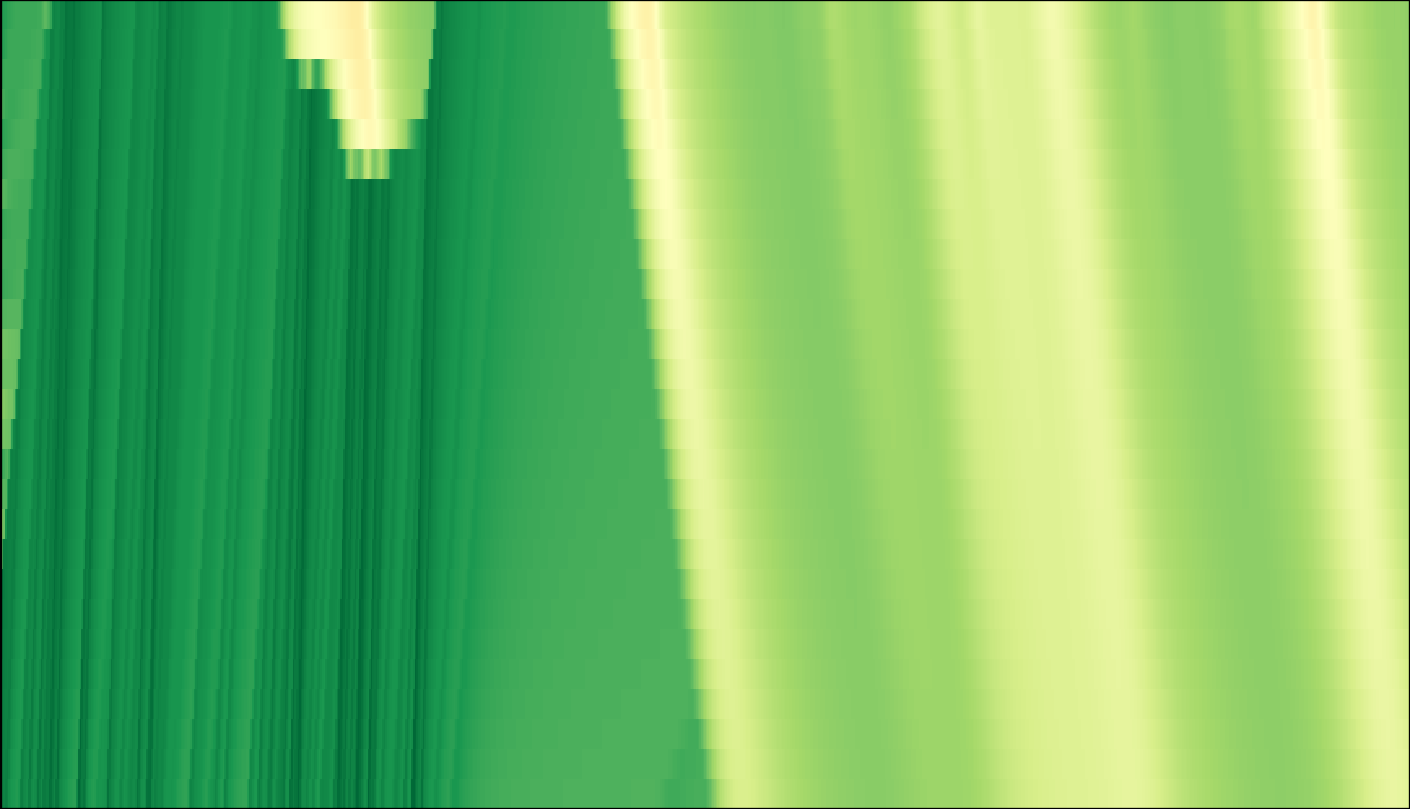}
      \vspace{-8.65cm}\caption*{Underwood}
    \end{subfigure}
  
    \raisebox{1.05cm}{\rotatebox[origin=c]{90}{Flow}}
    \begin{subfigure}[t]{.193\textwidth}
      \centering
      \begin{tikzpicture}
        \begin{axis}[
            grid=both,
            only marks, 
            width=1.625\linewidth,
            height=1.2\linewidth,
            scaled y ticks=false,
            xtick=\empty,
            ytick=\empty,
            xticklabels=\empty,
            yticklabels=\empty,
            xmin=-0.00621504, xmax=0.15537601,
            ymin=-0.05555556, ymax=0.66666667,
          ]
            \addplot[
              mark=*,
              mark size=1pt,
              blue,
              opacity=0.3,
            ] table[
              col sep=comma,
              x expr=\thisrow{density}/4/1609.34,
              y expr=\thisrow{flow}/4/3600
            ] {figs/fundamental_diagram/data_small.csv};
          \end{axis}  
          \begin{axis}[
              grid=both,
              only marks, 
              width=1.625\linewidth,
              height=1.2\linewidth,
              scaled y ticks=false,
              xtick=\empty,
              ytick=\empty,
              xticklabels=\empty,
              yticklabels=\empty,
              xmin=-0.00621504, xmax=0.15537601,
              ymin=-0.05555556, ymax=0.66666667,
            ]
            \addplot[very thick, smooth, color=orange, domain=0:0.10594] {20.009*x*(1.0-x/0.10594)};
            \end{axis}
            \node[anchor=north, yshift=-0.05cm] at (current axis.south) {\footnotesize MSE = 124588};  
      \end{tikzpicture}
    \end{subfigure}\hfill
    \begin{subfigure}[t]{.193\textwidth}
      \centering
      \begin{tikzpicture}
        \begin{axis}[
            grid=both,
            only marks, 
            width=1.625\linewidth,
            height=1.2\linewidth,
            scaled y ticks=false,
            xtick=\empty,
            ytick=\empty,
            xticklabels=\empty,
            yticklabels=\empty,
            xmin=-0.00621504, xmax=0.15537601,
            ymin=-0.05555556, ymax=0.66666667,
          ]
            \addplot[
              mark=*,
              mark size=1pt,
              blue,
              opacity=0.3,
            ] table[
              col sep=comma,
              x expr=\thisrow{density}/4/1609.34,
              y expr=\thisrow{flow}/4/3600
            ] {figs/fundamental_diagram/data_small.csv};
          \end{axis}  
          \begin{axis}[
              grid=both,
              only marks, 
              width=1.625\linewidth,
              height=1.2\linewidth,
              scaled y ticks=false,
              xtick=\empty,
              ytick=\empty,
              xticklabels=\empty,
              yticklabels=\empty,
              xmin=-0.00621504, xmax=0.15537601,
              ymin=-0.05555556, ymax=0.66666667,
            ]
            \addplot[very thick, smooth, color=orange, domain=0:0.01593] {39.42163*x};
            \addplot[very thick, smooth, color=orange, domain=0.01593:0.152867] {-4.5877 * (x-0.152867)};
            \end{axis}   
            \node[anchor=north, yshift=-0.05cm] at (current axis.south) {\footnotesize MSE = 209630};  
      \end{tikzpicture}
    \end{subfigure}\hfill
    \begin{subfigure}[t]{.193\textwidth}
      \centering
      \begin{tikzpicture}
        \begin{axis}[
            grid=both,
            only marks, 
            width=1.625\linewidth,
            height=1.2\linewidth,
            scaled y ticks=false,
            xtick=\empty,
            ytick=\empty,
            xticklabels=\empty,
            yticklabels=\empty,
            xmin=-0.00621504, xmax=0.15537601,
            ymin=-0.05555556, ymax=0.66666667,
          ]
            \addplot[
              mark=*,
              mark size=1pt,
              blue,
              opacity=0.3,
            ] table[
              col sep=comma,
              x expr=\thisrow{density}/4/1609.34,
              y expr=\thisrow{flow}/4/3600
            ] {figs/fundamental_diagram/data_small.csv};
          \end{axis}  
          \begin{axis}[
              grid=both,
              only marks, 
              width=1.625\linewidth,
              height=1.2\linewidth,
              scaled y ticks=false,
              xtick=\empty,
              ytick=\empty,
              xticklabels=\empty,
              yticklabels=\empty,
              xmin=-0.00621504, xmax=0.15537601,
              ymin=-0.05555556, ymax=0.66666667,
            ]            
            \addplot[very thick, smooth, color=orange, domain=0:0.0118300] {39.4904236*x};
            \addplot[very thick, smooth, color=orange, domain=0.0118300:0.04041] {0.467171 + (x - 0.0118300) * (0.467171 - 0.467171)/(0.0118300 - 0.040417)};
            \addplot[very thick, smooth, color=orange, domain=0.04041:0.142490] {-4.576839560 * (x-0.1424900)};
            \end{axis}   
            \node[anchor=north, yshift=-0.05cm] at (current axis.south) {\footnotesize MSE = 119592};  
      \end{tikzpicture}
    \end{subfigure}\hfill
    \begin{subfigure}[t]{.193\textwidth}
      \centering
      \begin{tikzpicture}
        \begin{axis}[
            grid=both,
            only marks, 
            width=1.625\linewidth,
            height=1.2\linewidth,
            scaled y ticks=false,
            xtick=\empty,
            ytick=\empty,
            xticklabels=\empty,
            yticklabels=\empty,
            xmin=-0.00621504, xmax=0.15537601,
            ymin=-0.05555556, ymax=0.66666667,
          ]
            \addplot[
              mark=*,
              mark size=1pt,
              blue,
              opacity=0.3,
            ] table[
              col sep=comma,
              x expr=\thisrow{density}/4/1609.34,
              y expr=\thisrow{flow}/4/3600
            ] {figs/fundamental_diagram/data_small.csv};
          \end{axis}  
          \begin{axis}[
              grid=both,
              only marks, 
              width=1.625\linewidth,
              height=1.2\linewidth,
              scaled y ticks=false,
              xtick=\empty,
              ytick=\empty,
              xticklabels=\empty,
              yticklabels=\empty,
              xmin=-0.00621504, xmax=0.15537601,
              ymin=-0.05555556, ymax=0.66666667,
            ]
            \addplot[very thick, smooth, color=orange, domain=0:0.106515] {7.0944226454 * x * ln(0.10651502/(x + 1e-7))};
            \end{axis}   
            \node[anchor=north, yshift=-0.05cm] at (current axis.south) {\footnotesize MSE = 413663};  
      \end{tikzpicture}
    \end{subfigure}\hfill
    \begin{subfigure}[t]{.193\textwidth}
      \centering
      \begin{tikzpicture}
        \begin{axis}[
            grid=both,
            only marks, 
            width=1.625\linewidth,
            height=1.2\linewidth,
            scaled y ticks=false,
            xtick=\empty,
            ytick=\empty,
            xticklabels=\empty,
            yticklabels=\empty,
            xmin=-0.00621504, xmax=0.15537601,
            ymin=-0.05555556, ymax=0.66666667,
          ]
            \addplot[
              mark=*,
              mark size=1pt,
              blue,
              opacity=0.3,
            ] table[
              col sep=comma,
              x expr=\thisrow{density}/4/1609.34,
              y expr=\thisrow{flow}/4/3600
            ] {figs/fundamental_diagram/data_small.csv};
          \end{axis}  
          \begin{axis}[
              grid=both,
              only marks, 
              width=1.625\linewidth,
              height=1.2\linewidth,
              scaled y ticks=false,
              xtick=\empty,
              ytick=\empty,
              xticklabels=\empty,
              yticklabels=\empty,
              xmin=-0.00621504, xmax=0.15537601,
              ymin=-0.05555556, ymax=0.66666667,
            ]
            \addplot[very thick, smooth, color=orange, domain=0:0.15] {13.874974116187104 * exp(1) * x * exp(-28.990779335400976 * x)};
            \end{axis}   
            \node[anchor=north, yshift=-0.05cm] at (current axis.south) {\footnotesize MSE = 97540};  
      \end{tikzpicture}
    \end{subfigure}

    
    \vspace{-0.1cm}
   \caption*{\textbf{Calibration on fundamental diagram}}
   \vspace{0.1cm}

    \raisebox{0.61cm}{\rotatebox[origin=c]{90}{GD}}
    \begin{subfigure}[t]{.193\textwidth}
      \centering
      \includegraphics[width=\textwidth]{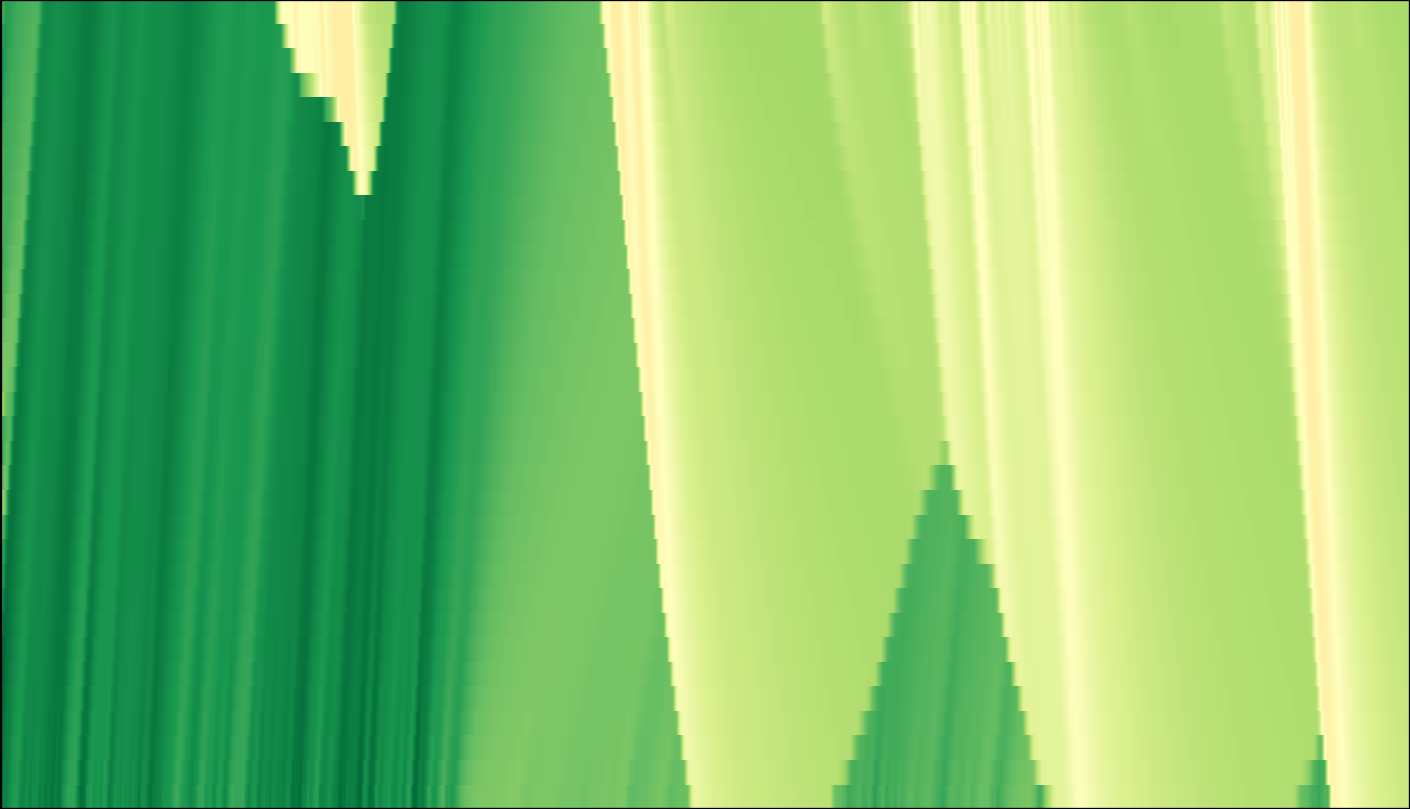}
    \end{subfigure}\hfill
    \begin{subfigure}[t]{.193\textwidth}
      \centering
      \includegraphics[width=\textwidth]{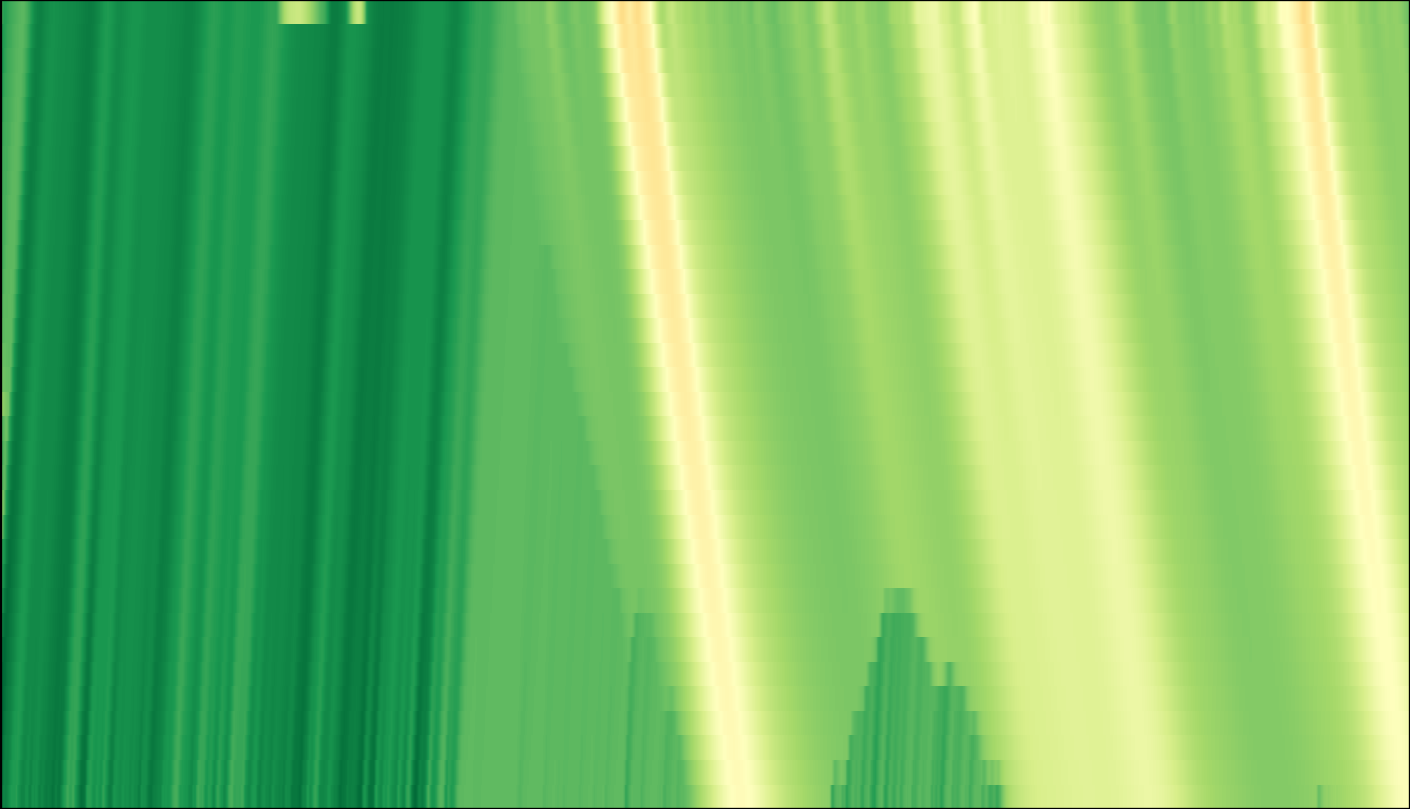}
    \end{subfigure}\hfill
    \begin{subfigure}[t]{.193\textwidth}
      \centering
      \includegraphics[width=\textwidth]{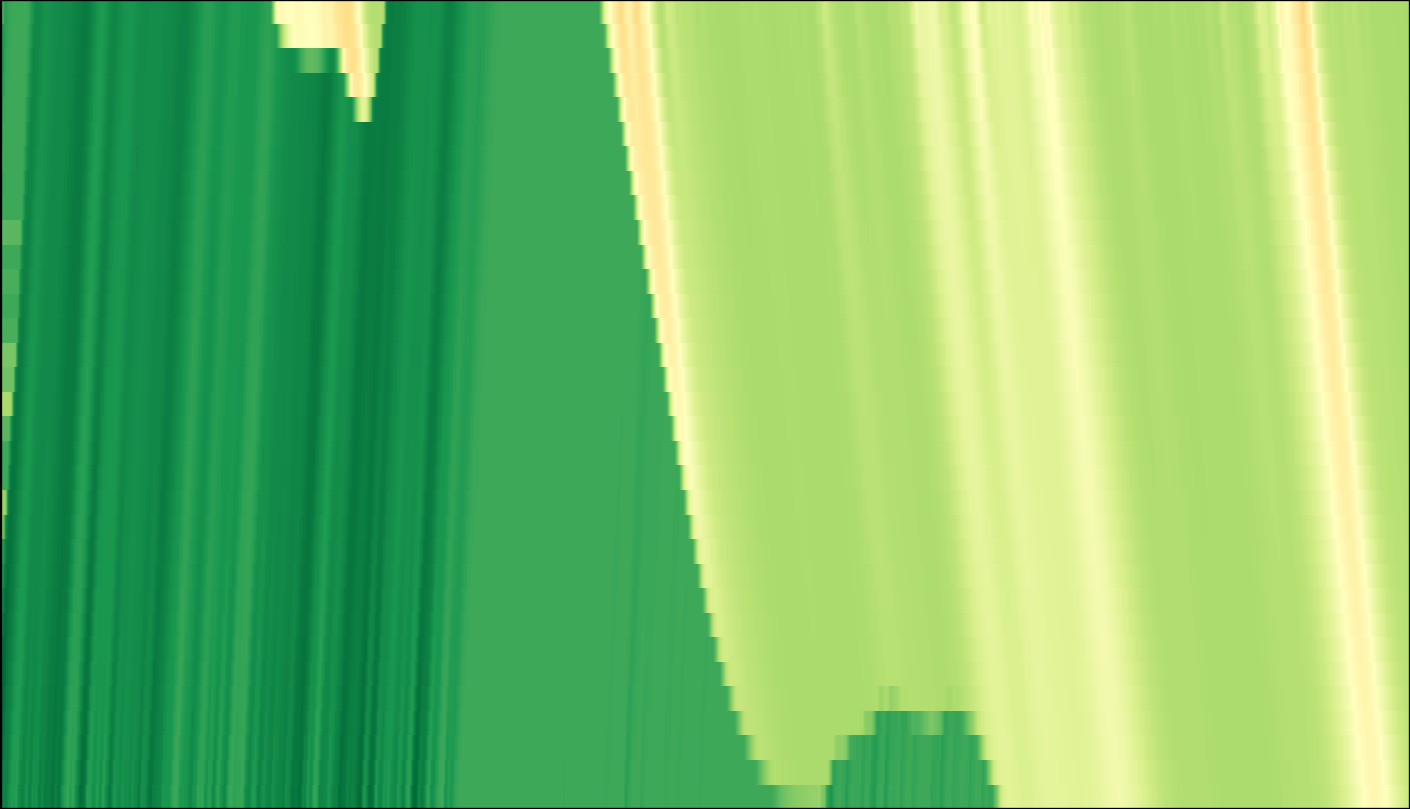}
    \end{subfigure}\hfill
    \begin{subfigure}[t]{.193\textwidth}
      \centering
      \includegraphics[width=\textwidth]{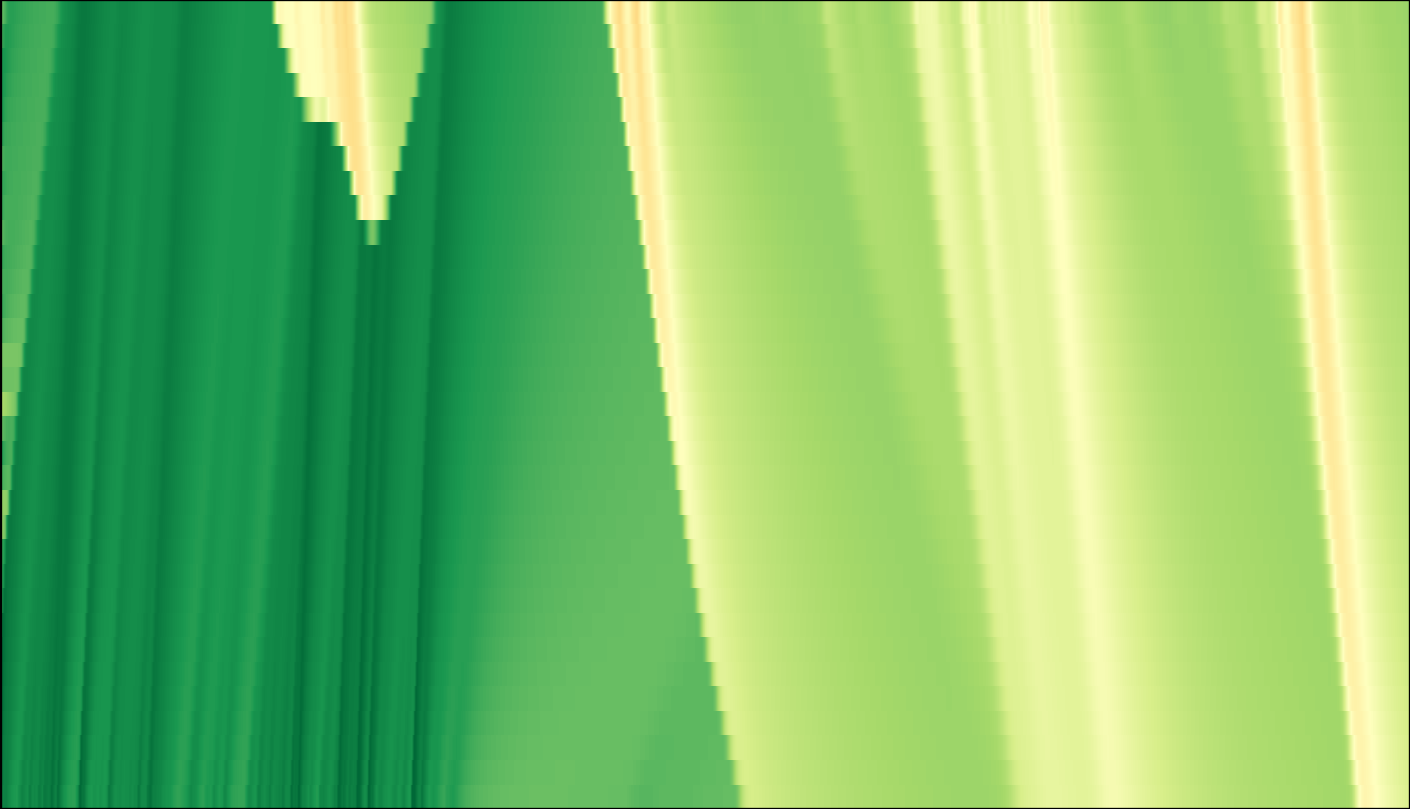}
    \end{subfigure}\hfill
    \begin{subfigure}[t]{.193\textwidth}
      \centering
      \includegraphics[width=\textwidth]{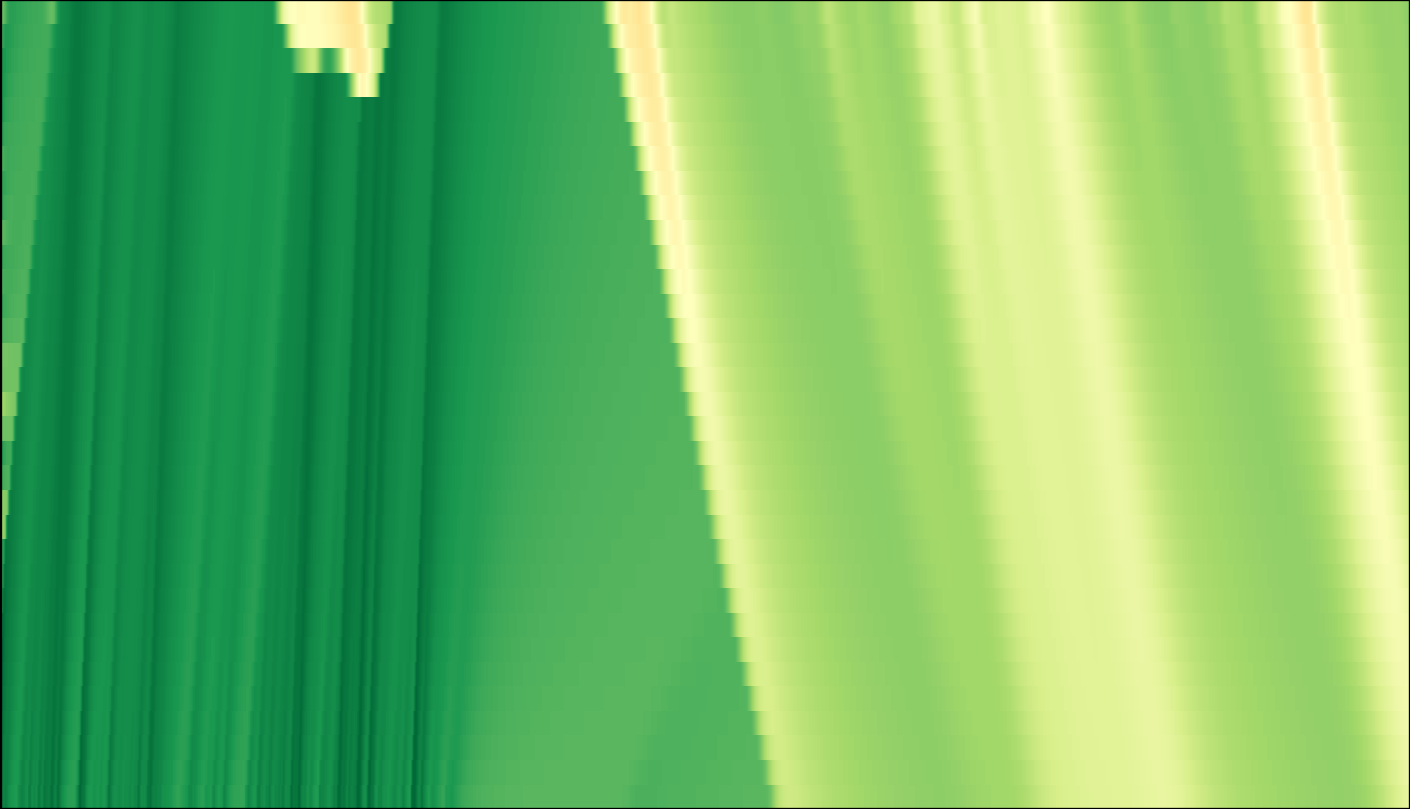}
    \end{subfigure}
  
    \raisebox{0.61cm}{\rotatebox[origin=c]{90}{LxF}}
    \begin{subfigure}[t]{.193\textwidth}
      \centering
      \includegraphics[width=\textwidth]{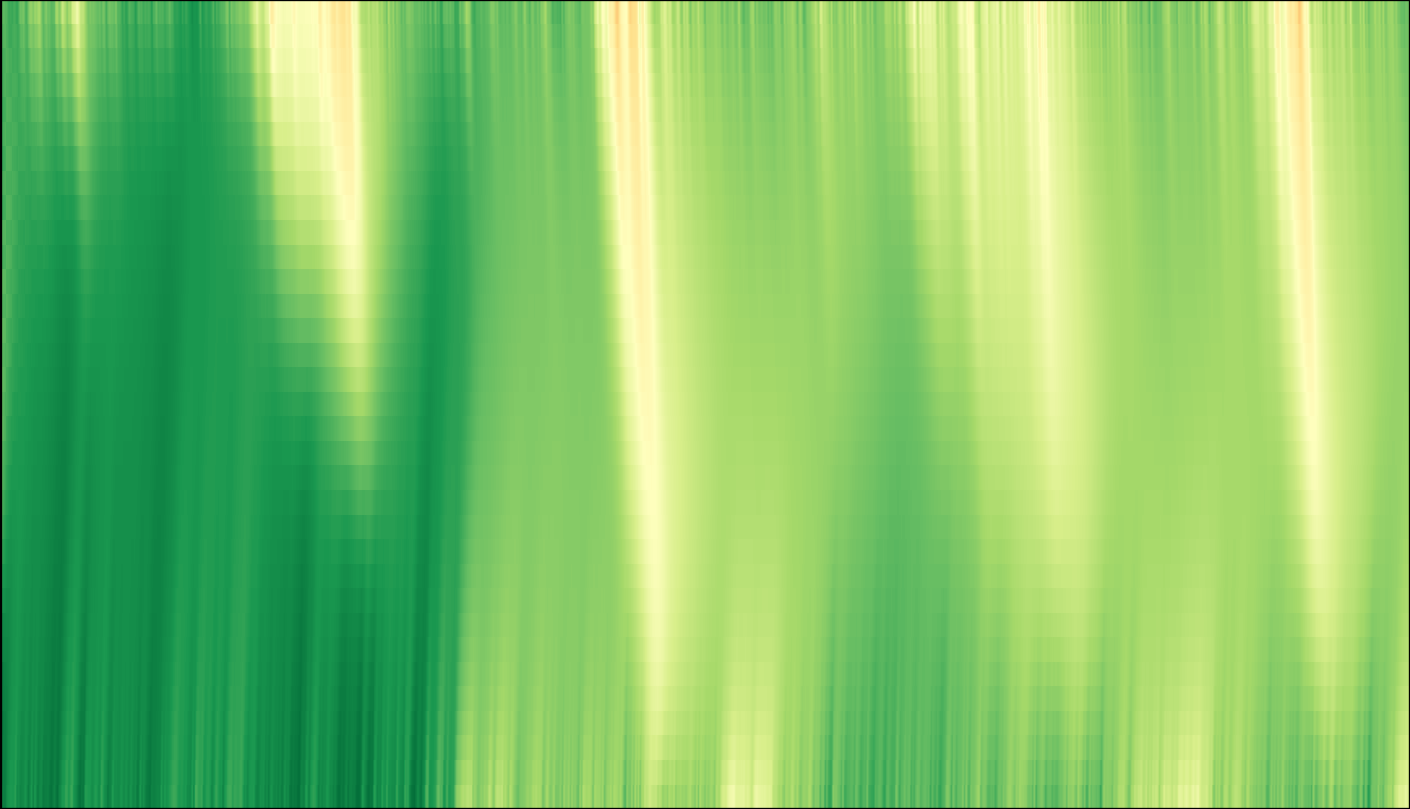}
    \end{subfigure}\hfill
    \begin{subfigure}[t]{.193\textwidth}
      \centering
      \includegraphics[width=\textwidth]{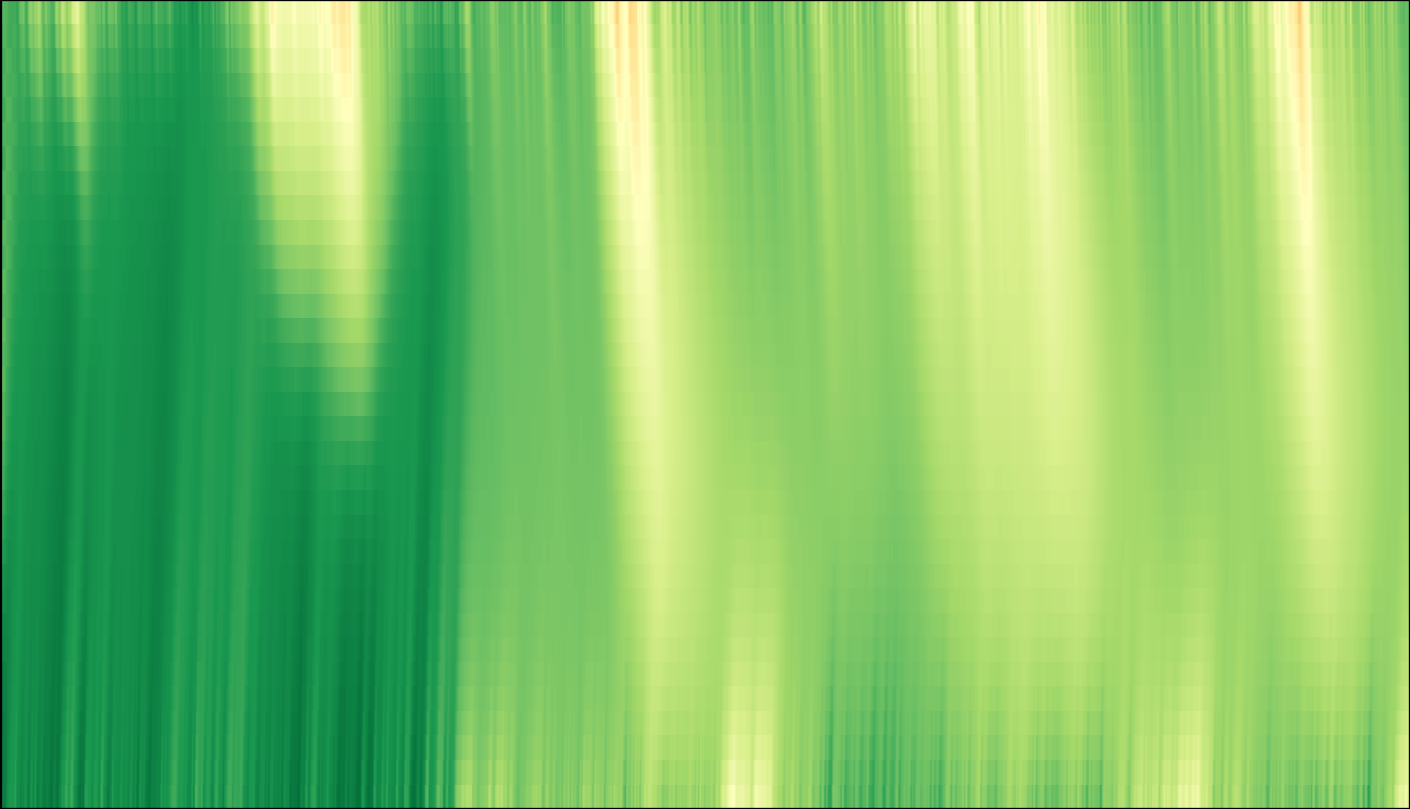}
    \end{subfigure}\hfill
    \begin{subfigure}[t]{.193\textwidth}
      \centering
      \includegraphics[width=\textwidth]{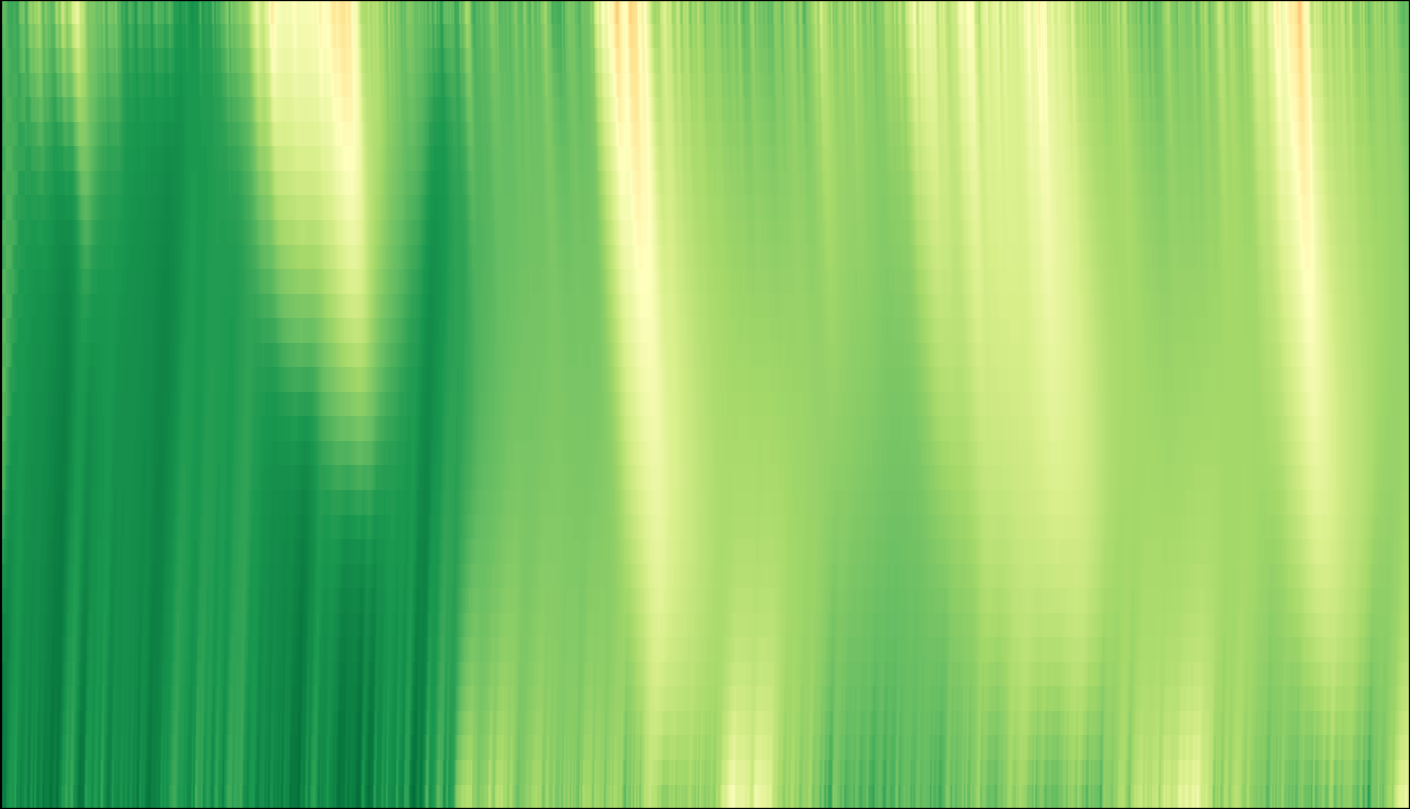}
    \end{subfigure}\hfill
    \begin{subfigure}[t]{.193\textwidth}
      \centering
      \includegraphics[width=\textwidth]{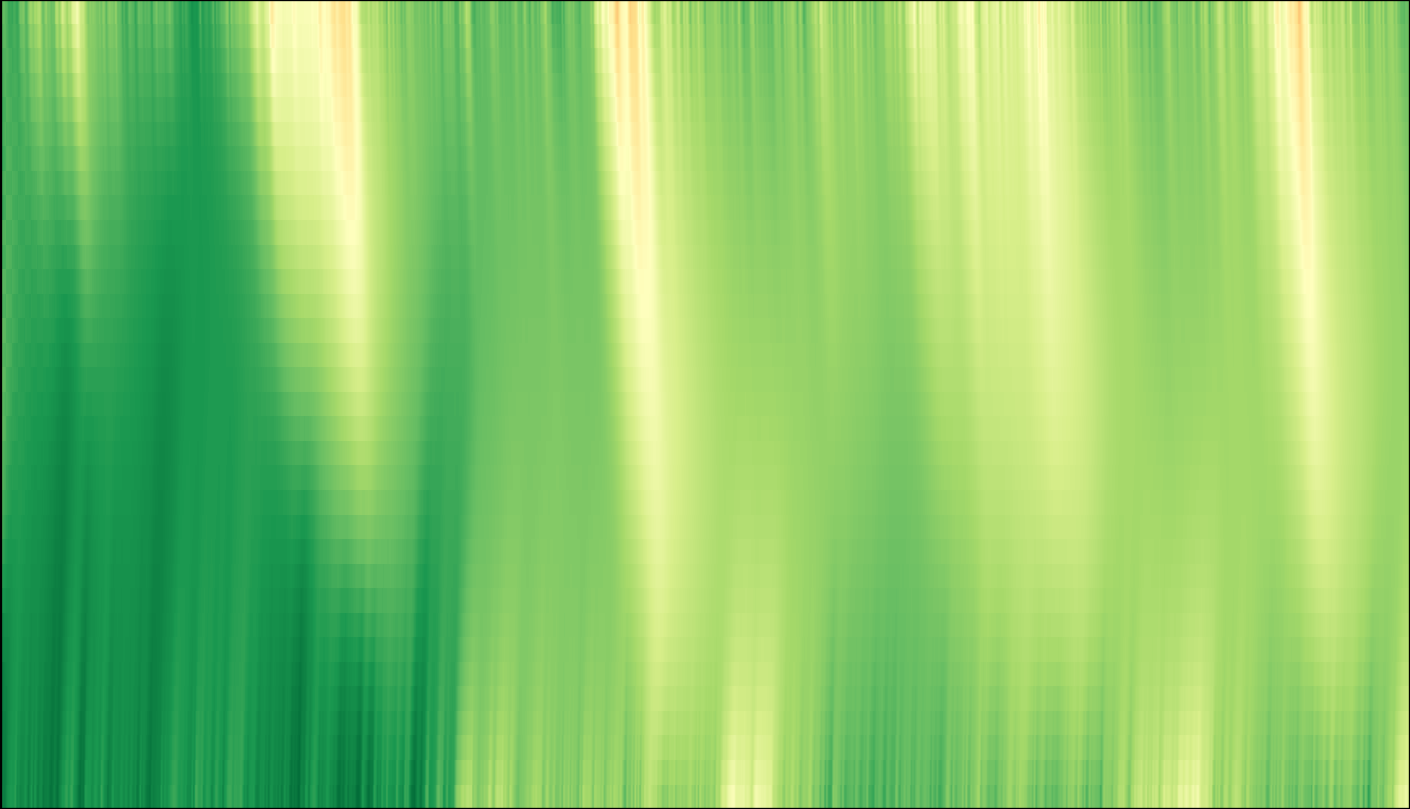}
    \end{subfigure}\hfill
    \begin{subfigure}[t]{.193\textwidth}
      \centering
      \includegraphics[width=\textwidth]{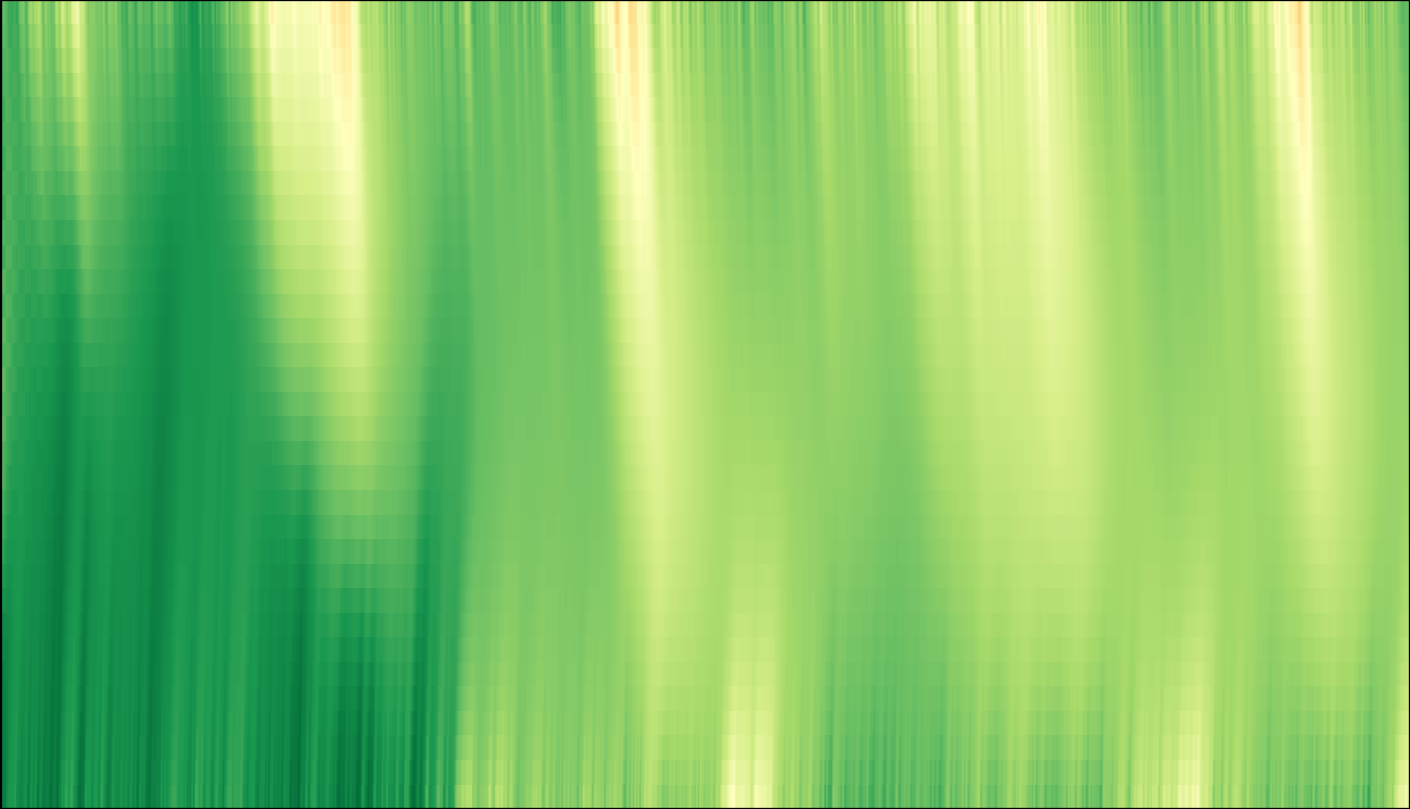}
    \end{subfigure}
  
    \raisebox{0.61cm}{\rotatebox[origin=c]{90}{EO}}
    \begin{subfigure}[t]{.193\textwidth}
      \centering
      \includegraphics[width=\textwidth]{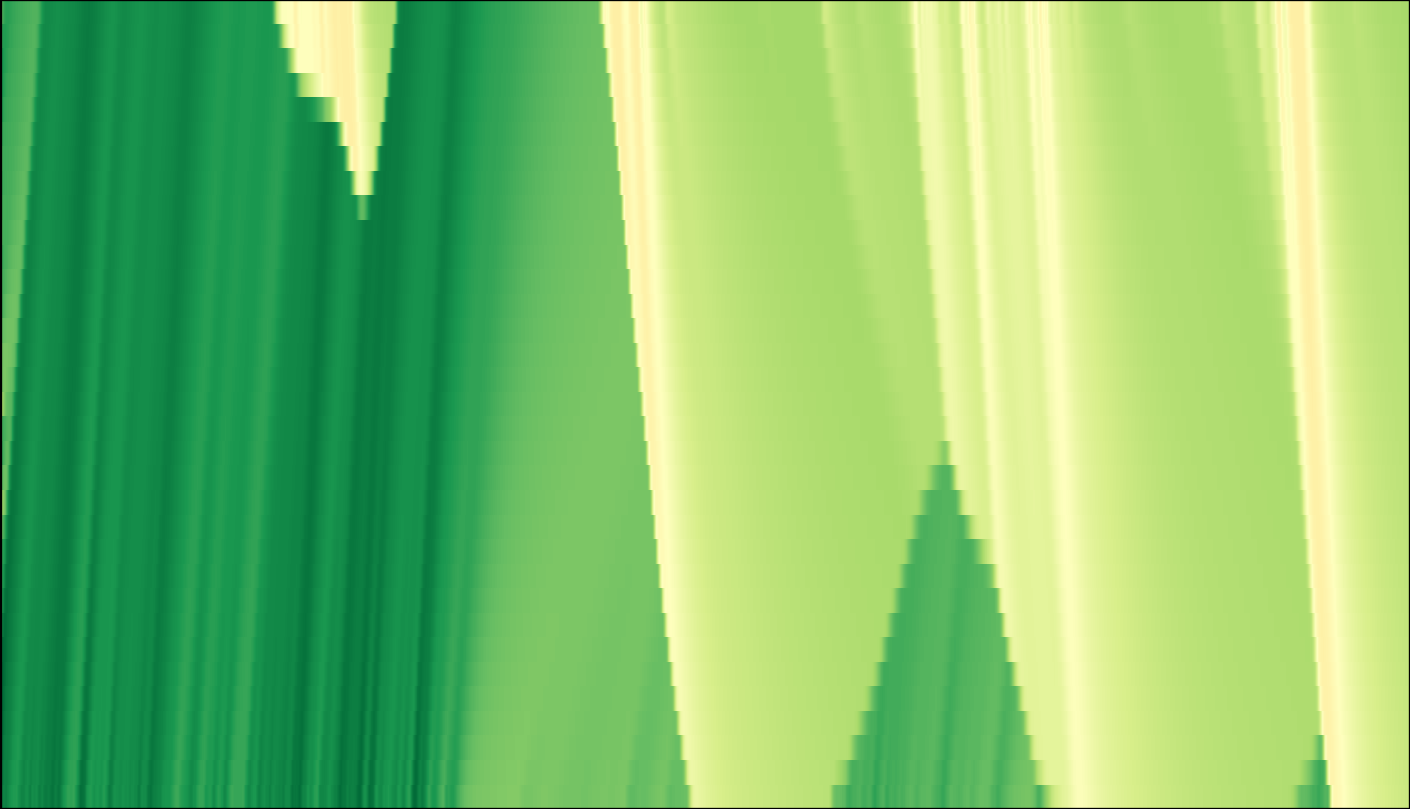}
    \end{subfigure}\hfill
    \begin{subfigure}[t]{.193\textwidth}
      \centering
      \includegraphics[width=\textwidth]{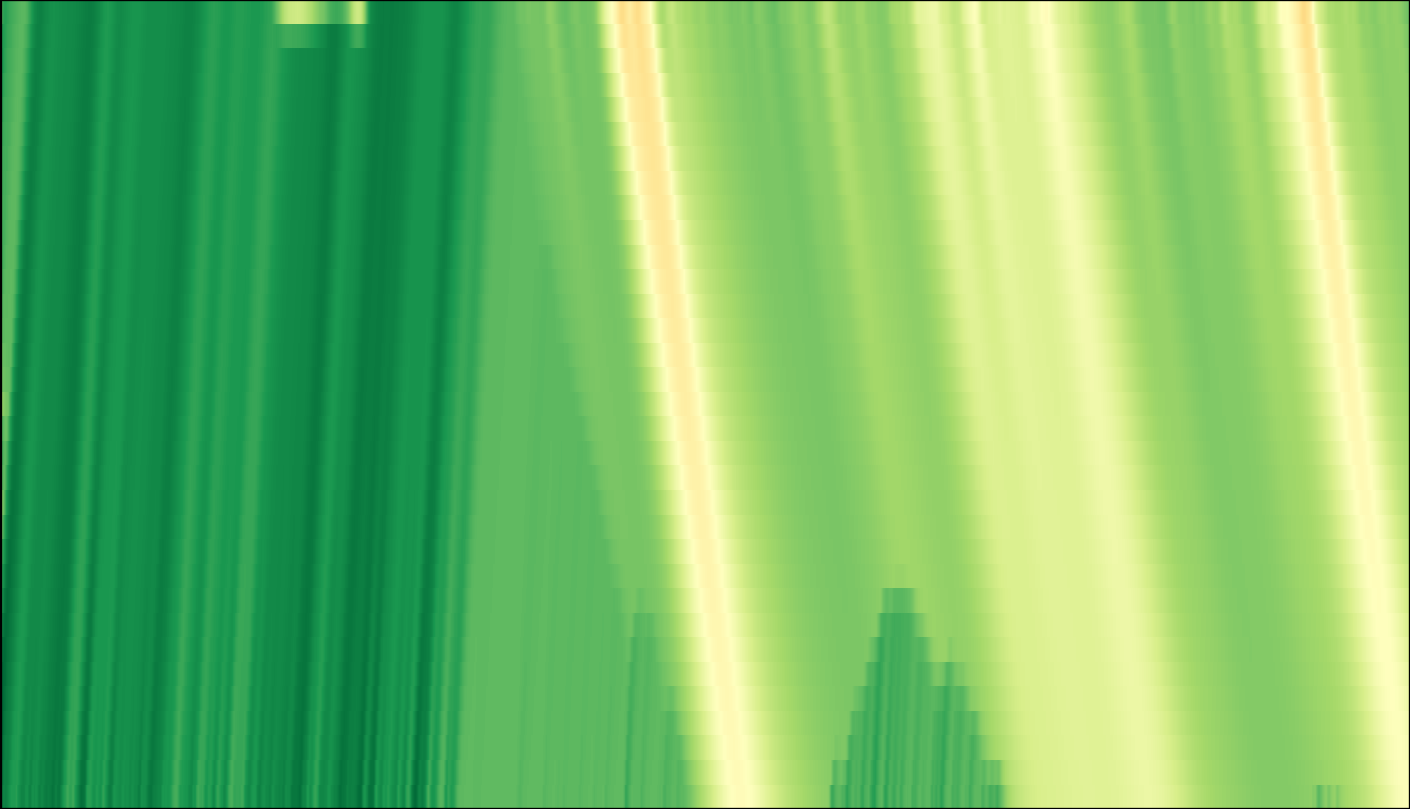}
    \end{subfigure}\hfill
    \begin{subfigure}[t]{.193\textwidth}
      \centering
      \includegraphics[width=\textwidth]{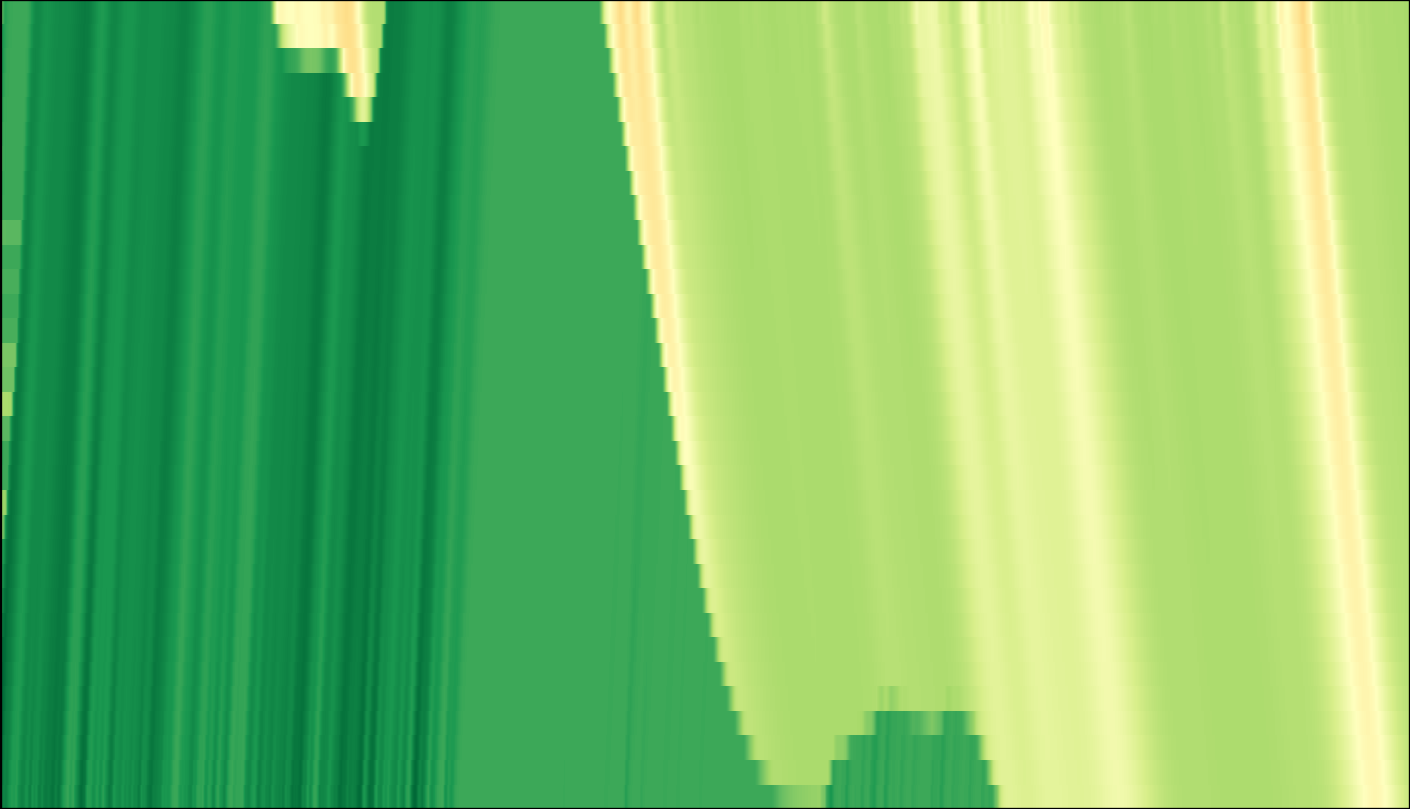}
    \end{subfigure}\hfill
    \begin{subfigure}[t]{.193\textwidth}
      \centering
      \includegraphics[width=\textwidth]{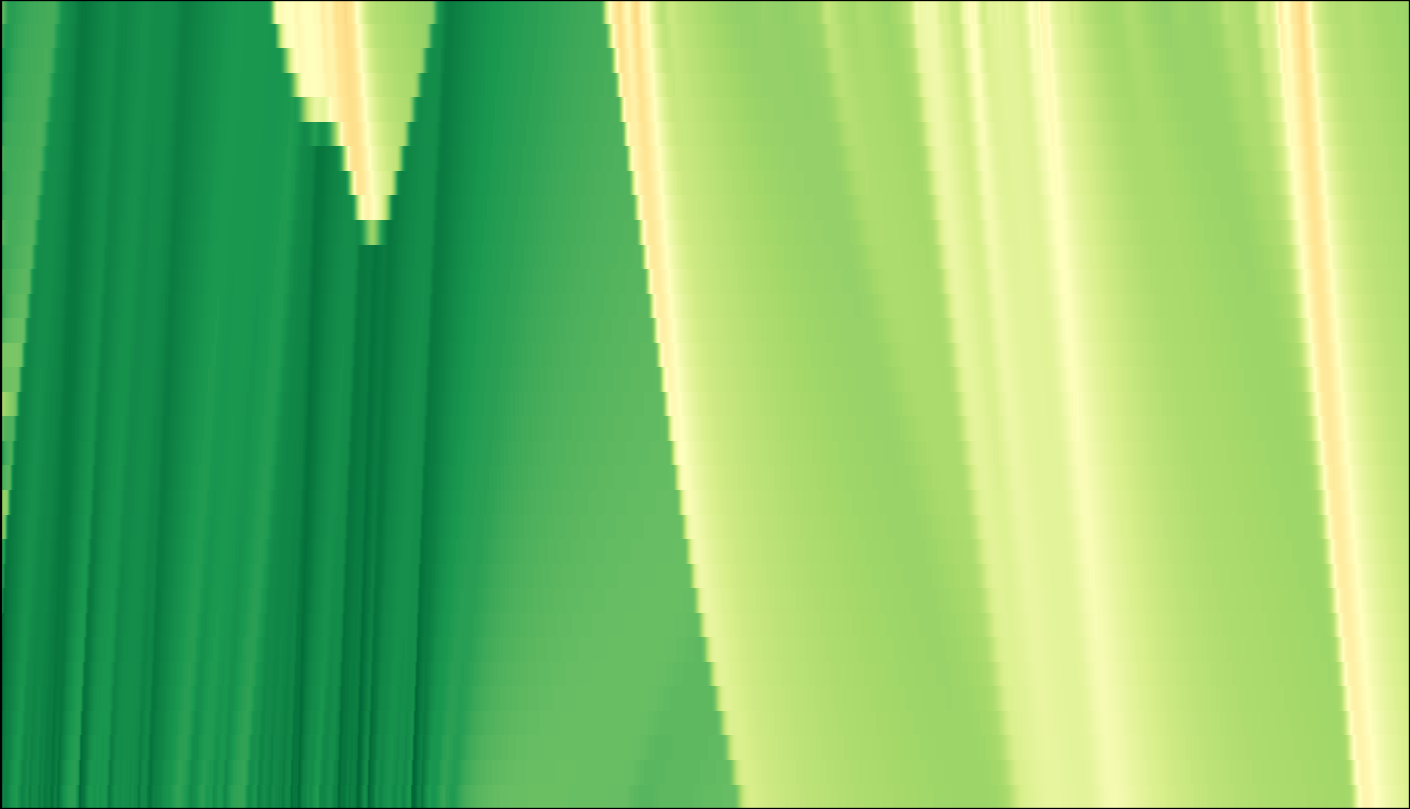}
    \end{subfigure}\hfill
    \begin{subfigure}[t]{.193\textwidth}
      \centering
      \includegraphics[width=\textwidth]{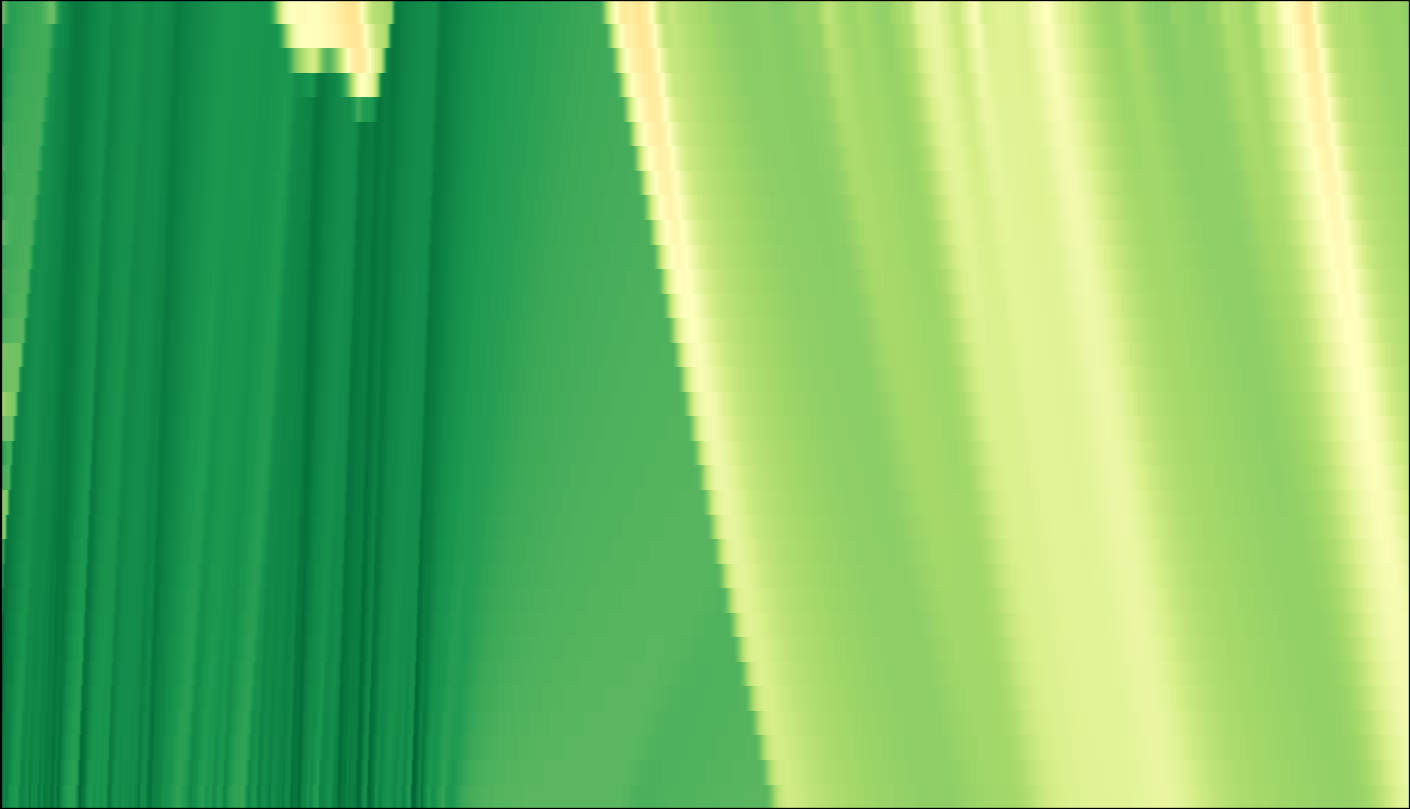}
    \end{subfigure}
  
    \raisebox{0.61cm}{\rotatebox[origin=c]{90}{ENO}}
    \begin{subfigure}[t]{.193\textwidth}
      \centering
      \includegraphics[width=\textwidth]{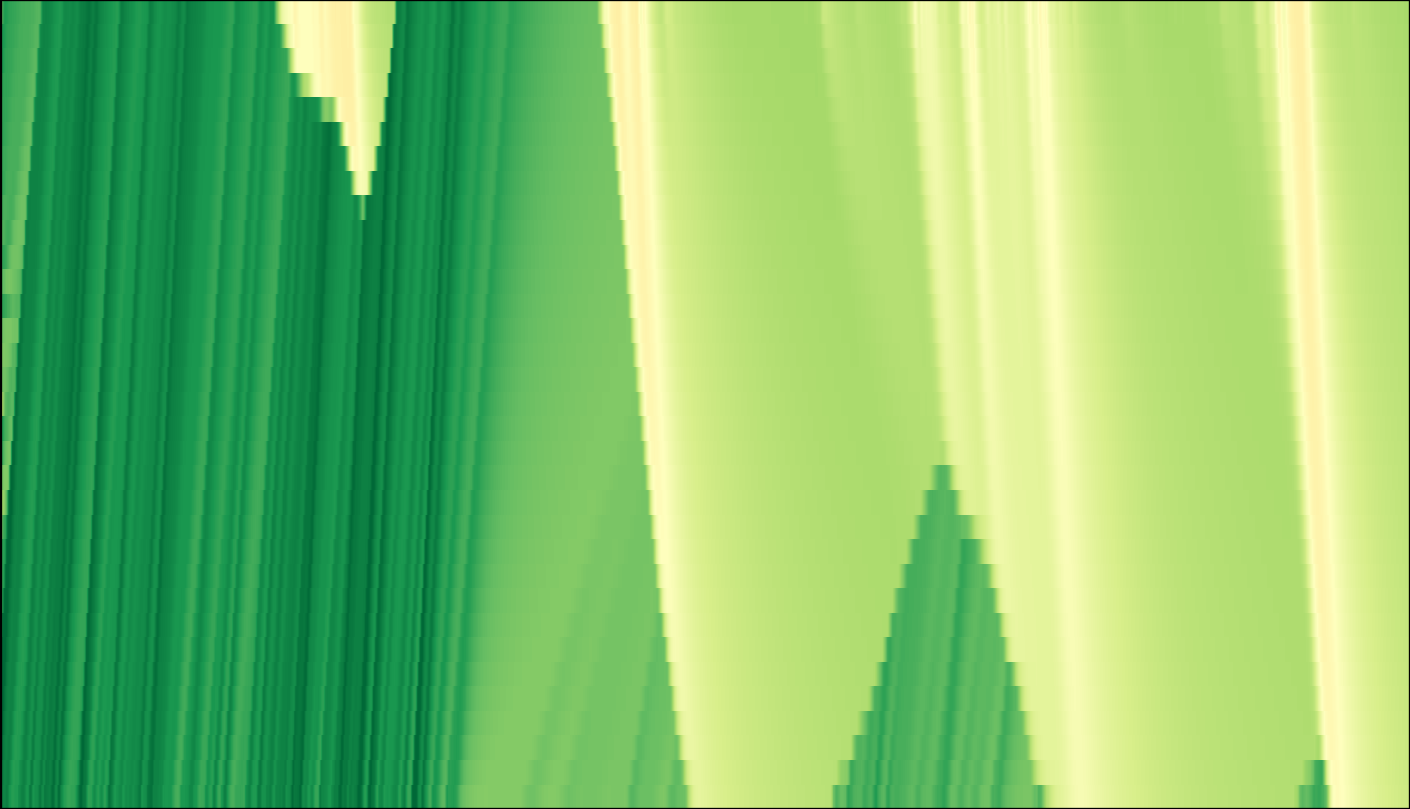}
    \end{subfigure}\hfill
    \begin{subfigure}[t]{.193\textwidth}
      \centering
      \includegraphics[width=\textwidth]{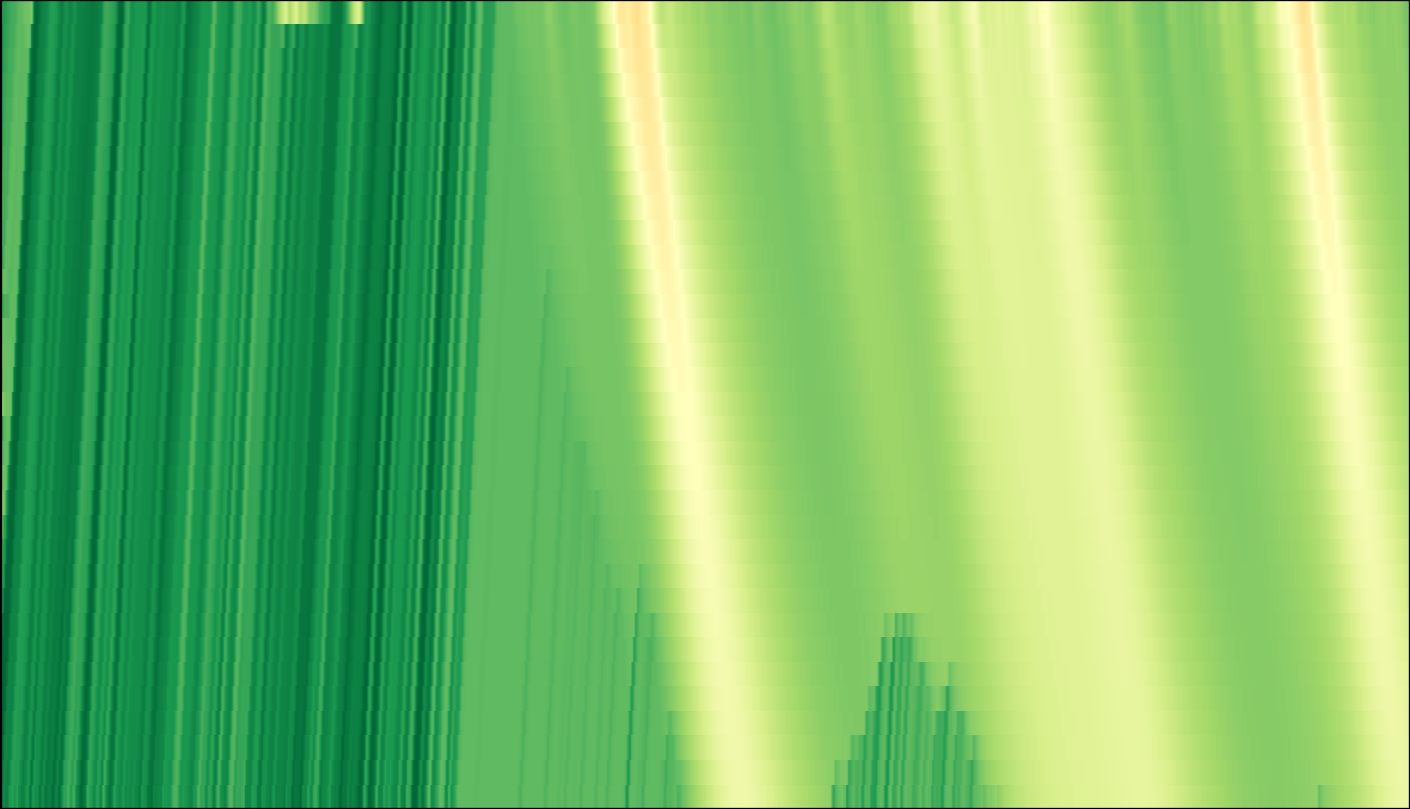}
    \end{subfigure}\hfill
    \begin{subfigure}[t]{.193\textwidth}
      \centering
      \includegraphics[width=\textwidth]{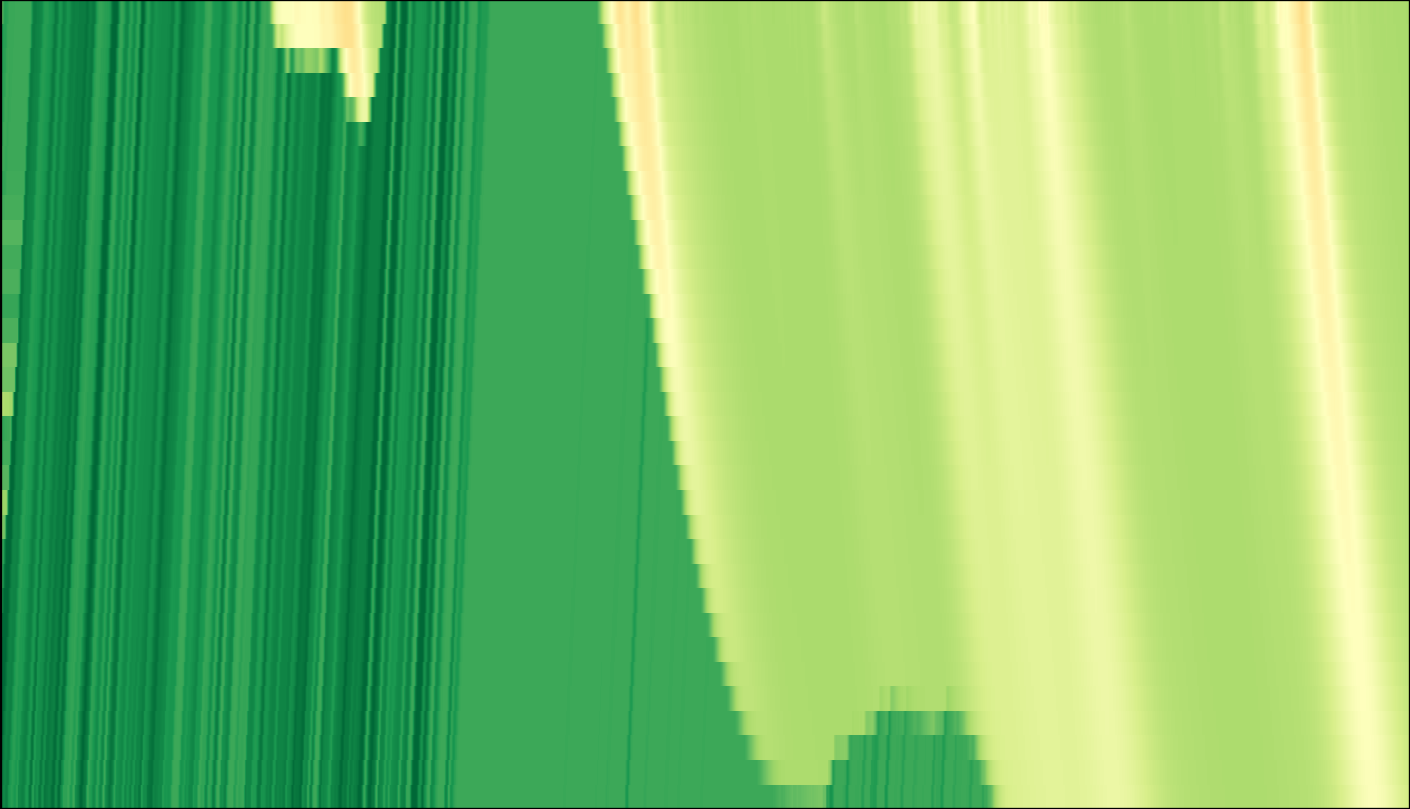}
    \end{subfigure}\hfill
    \begin{subfigure}[t]{.193\textwidth}
      \centering
      \includegraphics[width=\textwidth]{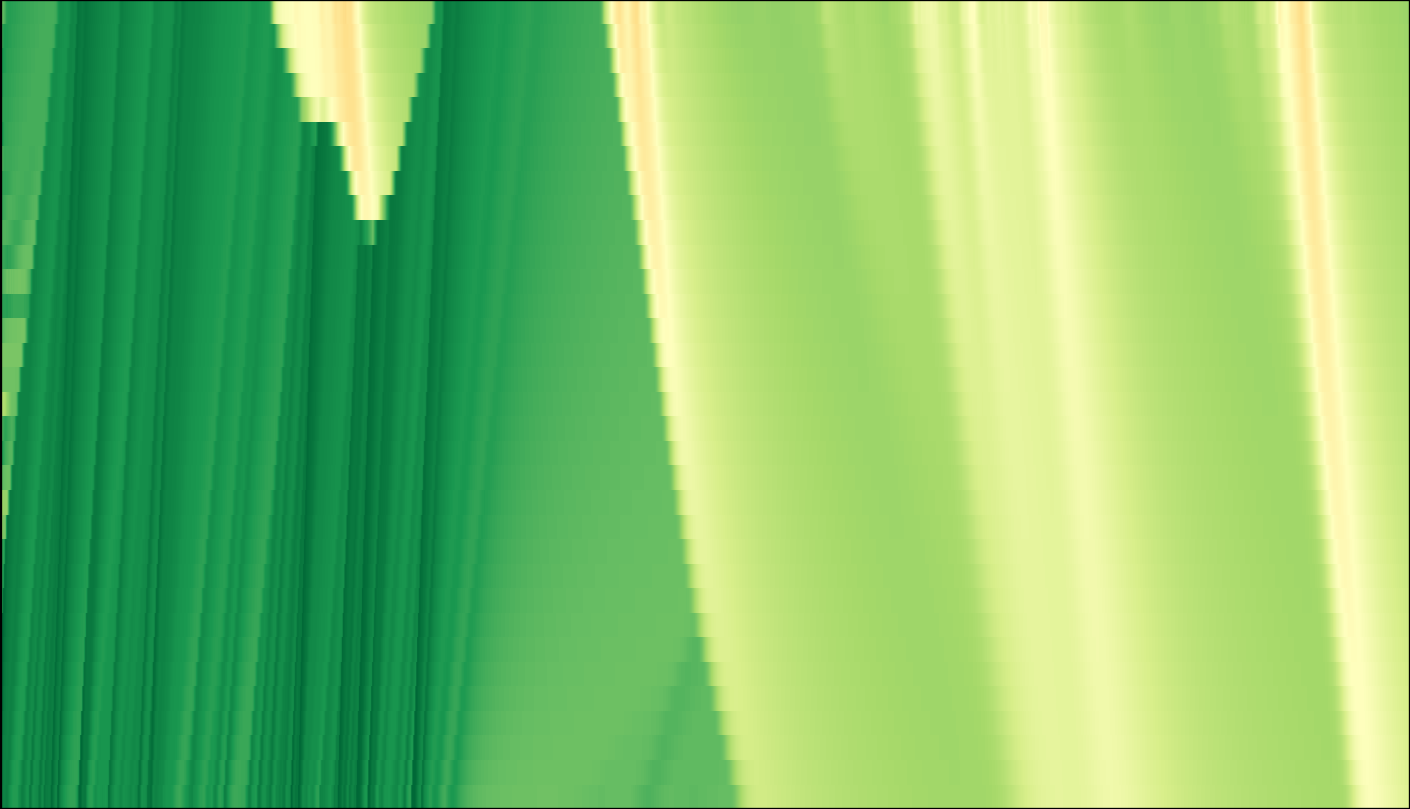}
    \end{subfigure}\hfill
    \begin{subfigure}[t]{.193\textwidth}
      \centering
      \includegraphics[width=\textwidth]{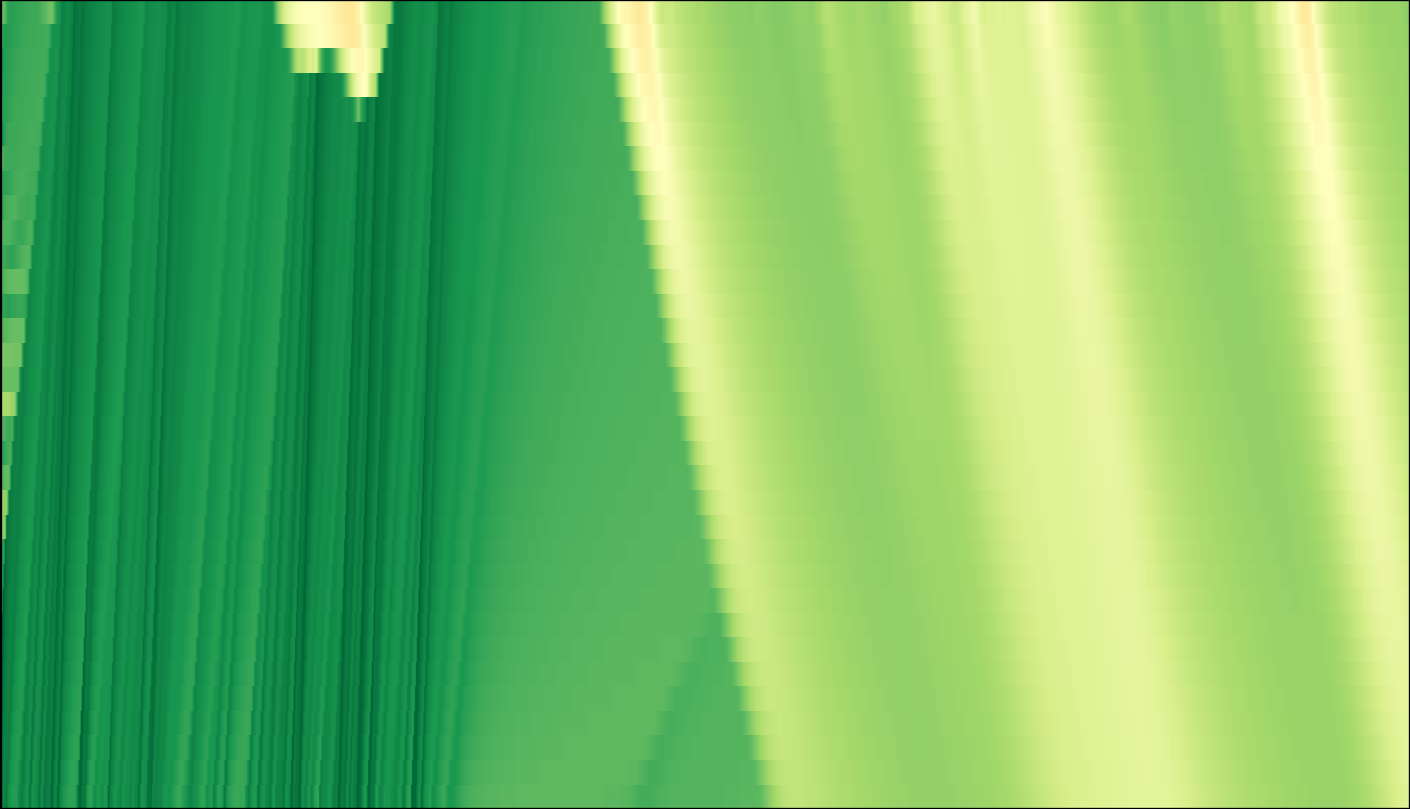}
    \end{subfigure}
  
    \raisebox{0.61cm}{\rotatebox[origin=c]{90}{WENO}}
    \begin{subfigure}[t]{.193\textwidth}
      \centering
      \includegraphics[width=\textwidth]{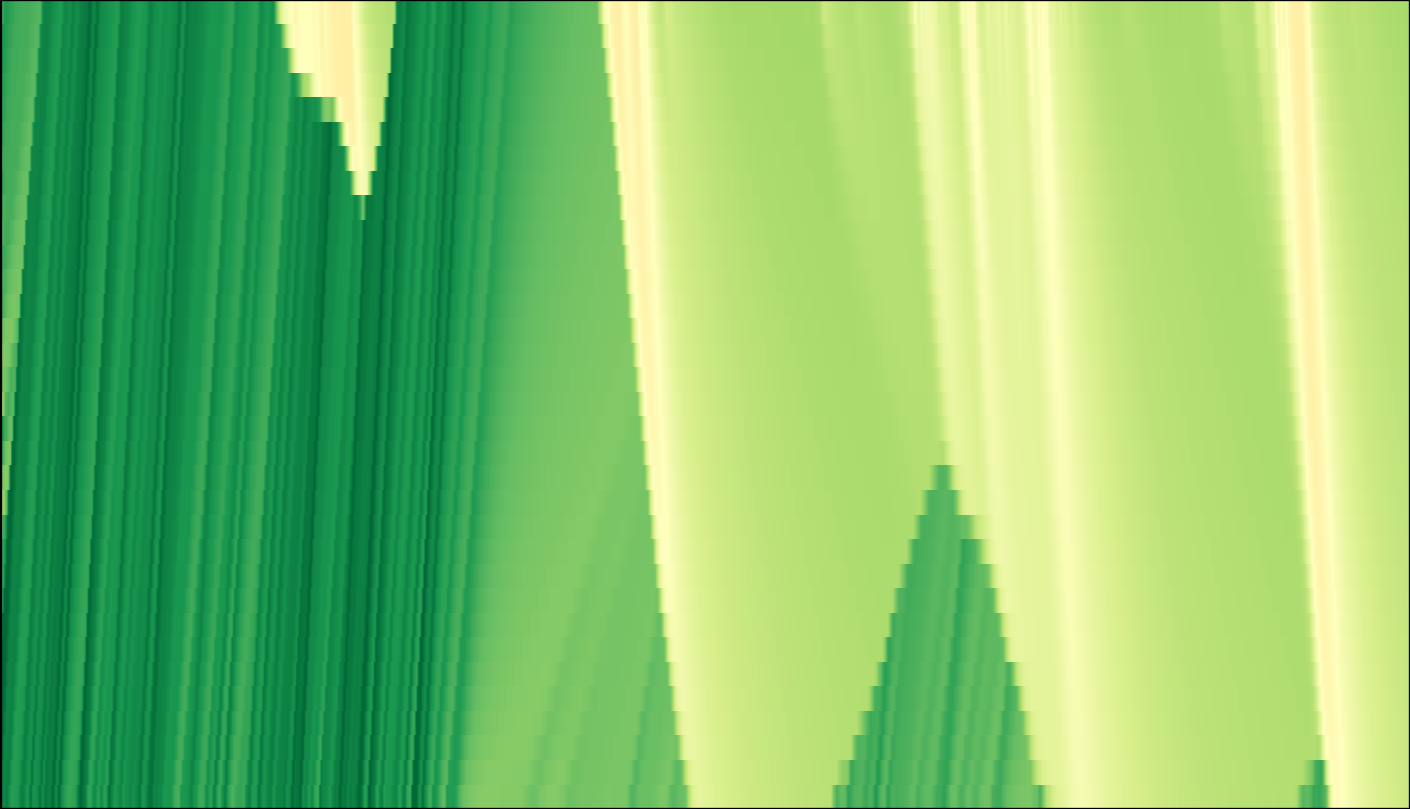}
      \vspace{-8.65cm}\caption*{Greenshields'}
    \end{subfigure}\hfill
    \begin{subfigure}[t]{.193\textwidth}
      \centering
      \includegraphics[width=\textwidth]{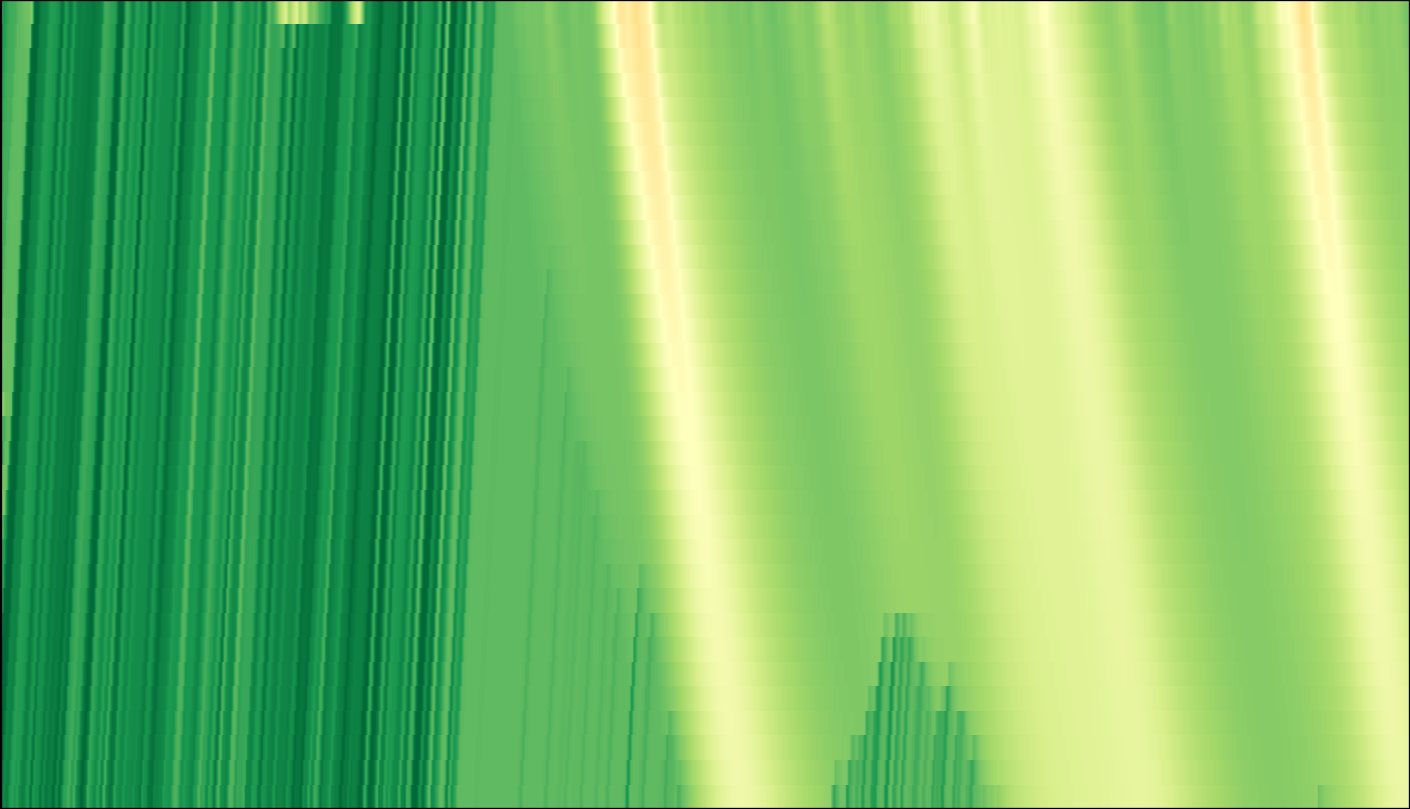}
      \vspace{-8.65cm}\caption*{Triangular}
    \end{subfigure}\hfill
    \begin{subfigure}[t]{.193\textwidth}
      \centering
      \includegraphics[width=\textwidth]{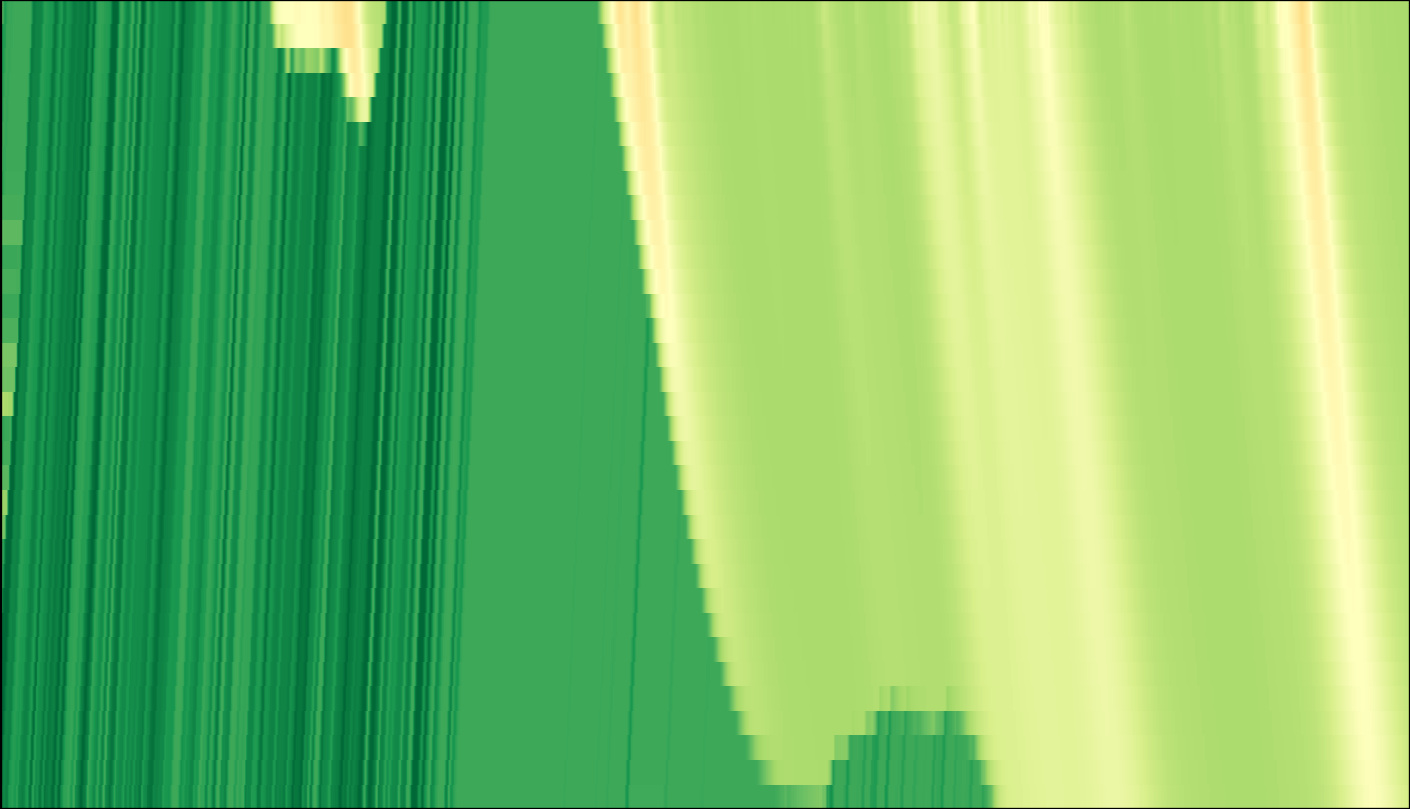}
      \vspace{-8.65cm}\caption*{Trapezoidal}
    \end{subfigure}\hfill
    \begin{subfigure}[t]{.193\textwidth}
      \centering
      \includegraphics[width=\textwidth]{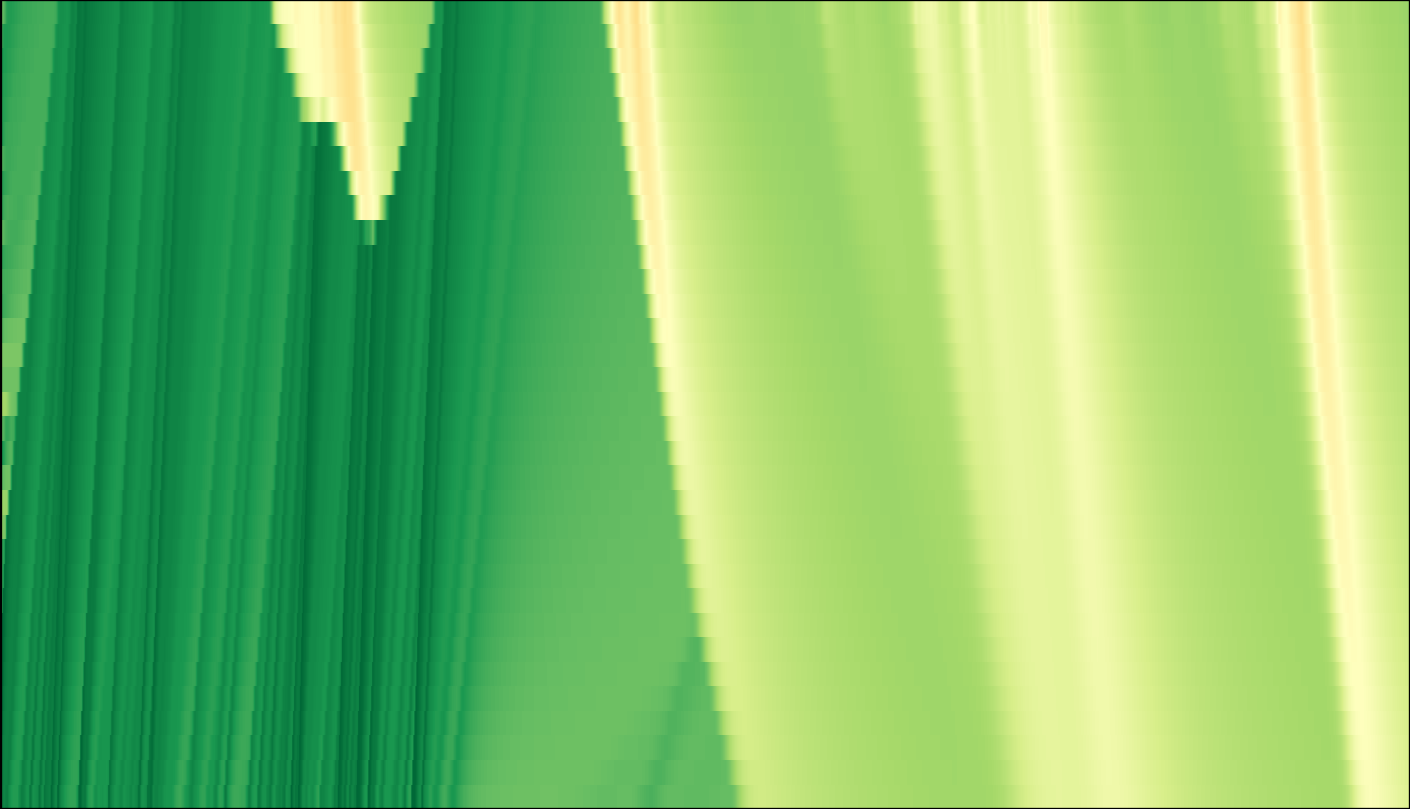}
      \vspace{-8.65cm}\caption*{Greenberg}
    \end{subfigure}\hfill
    \begin{subfigure}[t]{.193\textwidth}
      \centering
      \includegraphics[width=\textwidth]{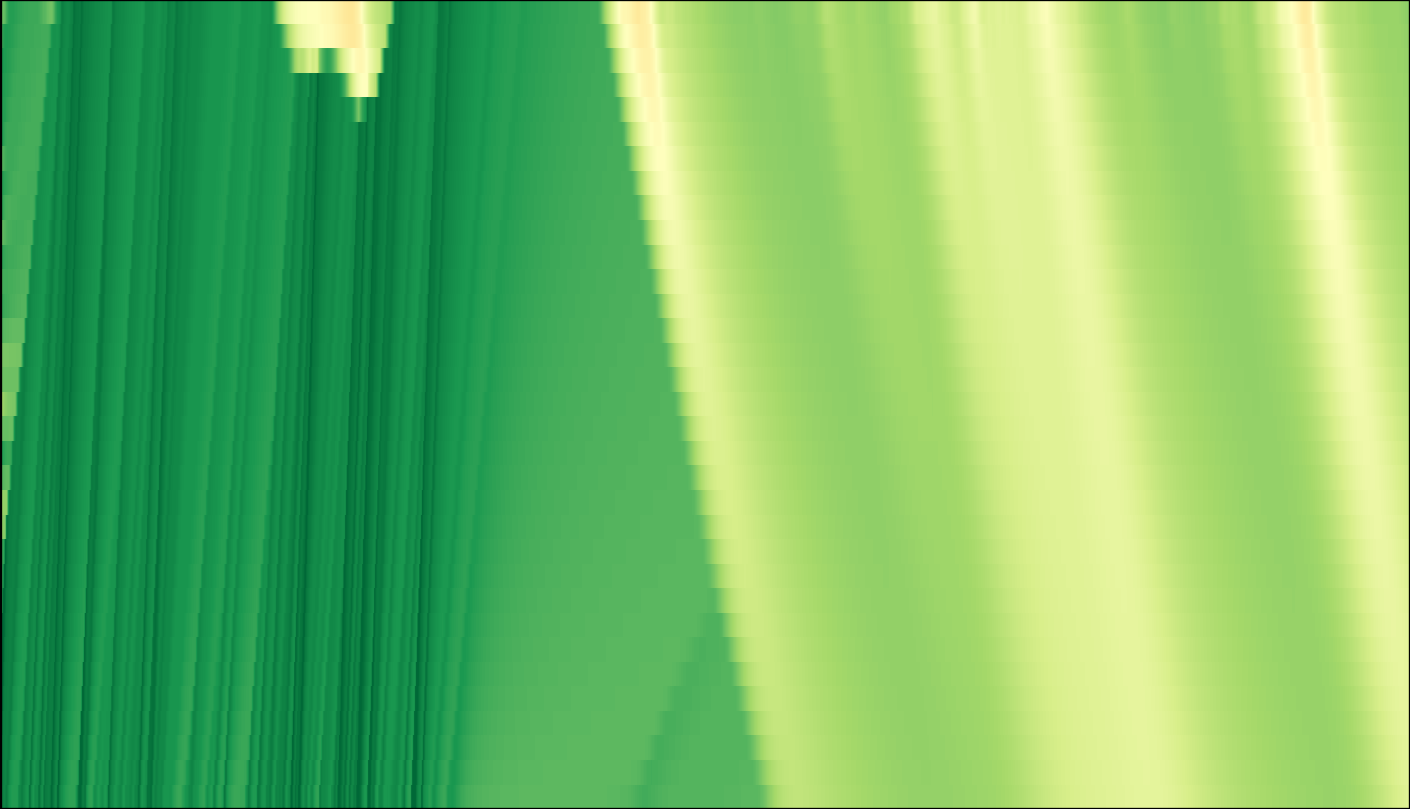}
      \vspace{-8.65cm}\caption*{Underwood}
    \end{subfigure}
  
    \raisebox{1.05cm}{\rotatebox[origin=c]{90}{Flow}}
    \begin{subfigure}[t]{.193\textwidth}
      \centering
      \begin{tikzpicture}
        \begin{axis}[
            grid=both,
            only marks, 
            width=1.625\linewidth,
            height=1.2\linewidth,
            scaled y ticks=false,
            xtick=\empty,
            ytick=\empty,
            xticklabels=\empty,
            yticklabels=\empty,
            xmin=-10, xmax=250,
            ymin=-200, ymax=2400,
          ]
            \addplot[
              mark=*,
              mark size=1pt,
              blue,
              opacity=0.3,
            ] table[
              col sep=comma,
              x expr=\thisrow{density}/4,
              y expr=\thisrow{flow}/4
            ] {figs/fundamental_diagram/data_small.csv};
          \end{axis}  
          \begin{axis}[
              grid=both,
              only marks, 
              width=1.625\linewidth,
              height=1.2\linewidth,
              scaled y ticks=false,
              xtick=\empty,
              ytick=\empty,
              xticklabels=\empty,
              yticklabels=\empty,
              xmin=-10, xmax=250,
              ymin=-200, ymax=2400,
            ]
              \addplot[very thick, smooth, color=orange, domain=0:170.872] {39.0706*x*(1.0-x/170.872)};
            \end{axis}
            \node[anchor=north, yshift=-0.05cm] at (current axis.south) {\footnotesize MSE = 83665};  
      \end{tikzpicture}
    \end{subfigure}\hfill
    \begin{subfigure}[t]{.193\textwidth}
      \centering
      \begin{tikzpicture}
        \begin{axis}[
            grid=both,
            only marks, 
            width=1.625\linewidth,
            height=1.2\linewidth,
            scaled y ticks=false,
            xtick=\empty,
            ytick=\empty,
            xticklabels=\empty,
            yticklabels=\empty,
            xmin=-10, xmax=250,
            ymin=-200, ymax=2400,
          ]
            \addplot[
              mark=*,
              mark size=1pt,
              blue,
              opacity=0.3,
            ] table[
              col sep=comma,
              x expr=\thisrow{density}/4,
              y expr=\thisrow{flow}/4
            ] {figs/fundamental_diagram/data_small.csv};
          \end{axis}  
          \begin{axis}[
              grid=both,
              only marks, 
              width=1.625\linewidth,
              height=1.2\linewidth,
              scaled y ticks=false,
              xtick=\empty,
              ytick=\empty,
              xticklabels=\empty,
              yticklabels=\empty,
              xmin=-10, xmax=250,
              ymin=-200, ymax=2400,
            ]
              \addplot[very thick, smooth, color=orange, domain=0:59.65] {30.67*x};
              \addplot[very thick, smooth, color=orange, domain=59.65:228.26] {-10.85 * (x-228.26)};
            \end{axis}   
            \node[anchor=north, yshift=-0.05cm] at (current axis.south) {\footnotesize MSE = 90131};  
      \end{tikzpicture}
    \end{subfigure}\hfill
    \begin{subfigure}[t]{.193\textwidth}
      \centering
      \begin{tikzpicture}
        \begin{axis}[
            grid=both,
            only marks, 
            width=1.625\linewidth,
            height=1.2\linewidth,
            scaled y ticks=false,
            xtick=\empty,
            ytick=\empty,
            xticklabels=\empty,
            yticklabels=\empty,
            xmin=-10, xmax=250,
            ymin=-200, ymax=2400,
          ]
            \addplot[
              mark=*,
              mark size=1pt,
              blue,
              opacity=0.3,
            ] table[
              col sep=comma,
              x expr=\thisrow{density}/4,
              y expr=\thisrow{flow}/4
            ] {figs/fundamental_diagram/data_small.csv};
          \end{axis}  
          \begin{axis}[
              grid=both,
              only marks, 
              width=1.625\linewidth,
              height=1.2\linewidth,
              scaled y ticks=false,
              xtick=\empty,
              ytick=\empty,
              xticklabels=\empty,
              yticklabels=\empty,
              xmin=-10, xmax=250,
              ymin=-200, ymax=2400,
            ]
              \addplot[very thick, smooth, color=orange, domain=0:45.40] {35.52*x};
              \addplot[very thick, smooth, color=orange, domain=45.40:95.72] {1612.67};
              \addplot[very thick, smooth, color=orange, domain=95.72:203.52] {-14.96 * (x-203.52)};
            \end{axis}   
            \node[anchor=north, yshift=-0.05cm] at (current axis.south) {\footnotesize MSE = 74097};  
      \end{tikzpicture}
    \end{subfigure}\hfill
    \begin{subfigure}[t]{.193\textwidth}
      \centering
      \begin{tikzpicture}
        \begin{axis}[
            grid=both,
            only marks, 
            width=1.625\linewidth,
            height=1.2\linewidth,
            scaled y ticks=false,
            xtick=\empty,
            ytick=\empty,
            xticklabels=\empty,
            yticklabels=\empty,
            xmin=-10, xmax=250,
            ymin=-200, ymax=2400,
          ]
            \addplot[
              mark=*,
              mark size=1pt,
              blue,
              opacity=0.3,
            ] table[
              col sep=comma,
              x expr=\thisrow{density}/4,
              y expr=\thisrow{flow}/4
            ] {figs/fundamental_diagram/data_small.csv};
          \end{axis}  
          \begin{axis}[
              grid=both,
              only marks, 
              width=1.625\linewidth,
              height=1.2\linewidth,
              scaled y ticks=false,
              xtick=\empty,
              ytick=\empty,
              xticklabels=\empty,
              yticklabels=\empty,
              xmin=-10, xmax=250,
              ymin=-200, ymax=2400,
            ]
              \addplot[very thick, smooth, color=orange, domain=0:192.98] {22.71 * x * ln(192.98/(x + 1e-7))};
            \end{axis}   
            \node[anchor=north, yshift=-0.05cm] at (current axis.south) {\footnotesize MSE = 76069};  
      \end{tikzpicture}
    \end{subfigure}\hfill
    \begin{subfigure}[t]{.193\textwidth}
      \centering
      \begin{tikzpicture}
        \begin{axis}[
            grid=both,
            only marks, 
            width=1.625\linewidth,
            height=1.2\linewidth,
            scaled y ticks=false,
            xtick=\empty,
            ytick=\empty,
            xticklabels=\empty,
            yticklabels=\empty,
            xmin=-10, xmax=250,
            ymin=-200, ymax=2400,
          ]
            \addplot[
              mark=*,
              mark size=1pt,
              blue,
              opacity=0.3,
            ] table[
              col sep=comma,
              x expr=\thisrow{density}/4,
              y expr=\thisrow{flow}/4
            ] {figs/fundamental_diagram/data_small.csv};
          \end{axis}  
          \begin{axis}[
              grid=both,
              only marks, 
              width=1.625\linewidth,
              height=1.2\linewidth,
              scaled y ticks=false,
              xtick=\empty,
              ytick=\empty,
              xticklabels=\empty,
              yticklabels=\empty,
              xmin=-10, xmax=250,
              ymin=-200, ymax=2400,
            ]
              \addplot[very thick, smooth, color=orange, domain=0:240] {26.22 * exp(1) * x * exp(-0.0162 * x)};
            \end{axis}   
            \node[anchor=north, yshift=-0.05cm] at (current axis.south) {\footnotesize MSE = 89689};  
      \end{tikzpicture}
    \end{subfigure}
  
    \caption{Forward predictions over the whole dataset from \Cref{fig:drone_dataset_heatmap} using the fitted flow functions with different FV schemes. For each calibration method (see \Cref{sec:calibrated_flux_functions}), the calibrated flows are shown overlapping with the fundamental diagram. The MSE values are computed as in Equation~\eqref{eq:mse_pred}. The heatmaps axes and colorbars are as in in~\Cref{fig:drone_dataset_heatmap}. The flow plots axes are as in~\Cref{fig:fundamental_diagram}. Unit of MSE: $(\text{veh}/\text{hour}/\text{lane})^2$.}
    \label{fig:drone_fitted_all}
  \end{figure}

  \clearpage

  \begin{table}[h]
    \centering
    \begin{tabular}{lccccc}

        \multicolumn{6}{c}{\textbf{Calibration on scheme prediction}} \\ 
        \toprule
        \small \textbf{Scheme} & \small Greenshields' & \small Triangular & \small Trapezoidal & \small Greenberg & \small Underwood \\
        \midrule
        Godunov           & 1018.53 & 404.70 & 428.02 & 927.74 & 887.10 \\
        LxF              & \textbf{719.53} & 620.45 & 607.88 & \textbf{596.30} & \textbf{601.27} \\
        EO               & 1008.98 & \textbf{403.73} & 427.65 & 917.30 & 881.93 \\
        ENO              & 1006.42 & 418.76 & 428.43 & 923.56 & 883.73 \\
        WENO             & 1004.49 & 423.83 & \textbf{425.52} & 922.21 & 880.80 \\
        \bottomrule \\

        \multicolumn{6}{c}{\textbf{Calibration on fundamental diagram}} \\ 
        \toprule
        \small \textbf{Scheme} & \small Greenshields' & \small Triangular & \small Trapezoidal & \small Greenberg & \small Underwood \\
        \midrule
        Godunov           & 1015.57 & 943.43 & 1126.08 & 942.51 & 1083.78 \\
        LxF   & \textbf{698.14} & \textbf{726.66} & \textbf{700.57} & \textbf{616.32} & \textbf{629.95} \\
        EO   & 1005.67 & 937.90 & 1122.44 & 935.62 & 1077.93 \\
        ENO               & 1003.80 & 957.51 & 1114.99 & 932.22 & 1075.46 \\
        WENO              & 1003.38 & 959.16 & 1114.53 & 931.82 & 1073.69 \\
        \bottomrule
    \end{tabular}
    \caption{MSE of the predictions from \Cref{fig:drone_fitted_all} against the ground truth from \Cref{fig:drone_dataset_heatmap}, computed as in Equation~\eqref{eq:mse_pred}. Unit of MSE: $(\text{veh}/\text{km}/\text{lane})^2$.}
    \label{tab:drone_fitted_mse_pred}
  \end{table}

  \begin{table}[H]
    \centering
    \begin{tabular}{ccccc}
        \multicolumn{5}{c}{\textbf{Calibration on scheme prediction}} \\ 
        \toprule
        Greenshields' & Triangular & Trapezoidal & Greenberg & Underwood \\
        \midrule
        $k_{\max}=105.94$ & $k_{\max}=152.87$ & $k_{\max}=142.49$ & $k_{\max}=106.52$ & $k_{\max}=\infty$ \\
        $v_{\max}=72.04$ & $v_{\max}=141.91$ & $v_{\max}=142.16$ & $c_0=25.54$ & $c_1=49.95$ \\
        -- & $w=-16.49$ & $w=-16.45$ & -- & $c_2=0.0290$ \\
        -- & -- & $q_{\max}=1677.54$ & -- & -- \\
        \bottomrule \\

        \multicolumn{5}{c}{\textbf{Calibration on fundamental diagram}} \\ 
        \toprule
        Greenshields' & Triangular & Trapezoidal & Greenberg & Underwood \\
        \midrule
        
        $k_{\max}=106.18$ & $k_{\max}=141.83$ & $k_{\max}=126.46$ & $k_{\max}=119.91$ & $k_{\max}=\infty$ \\
        $v_{\max}=61.27$ & $v_{\max}=49.36$ & $v_{\max}=57.16$ & $c_{0}=36.55$ & $c_{1}=42.20$ \\
        -- & $w=-17.46$ & $w=-24.08$ & -- & $c_{2}=0.0261$ \\
        -- & -- & $q_{\max}=1612.67$ & -- & -- \\
        \bottomrule
    \end{tabular}
    \caption{Fitted flow parameters from \Cref{fig:drone_fitted_all}. Units: Densities ($k_{\max}$) in $\text{veh}/\text{km}/\text{lane}$, speeds ($v_{\max}$) in $\text{km}/\text{h}$, and flows ($q_{\max}$) in $\text{veh}/\text{h}/\text{lane}$. Coefficients $c_0, c_1, c_2$ (and others) are such that if density is in $\text{veh}/\text{km}/\text{lane}$, then flow output is in $\text{veh}/\text{h}/\text{lane}$.}
    \label{tab:drone_fitted_params}
\end{table}

%% file: tex/drone_dtw.tex
\begin{figure}
    \centering
      
    \begin{tikzpicture}
      \begin{axis}[
        width=1.0\linewidth,
        height=0.3\linewidth,
        xlabel={x-slice},
        ylabel={DTW distance},
        grid=major,
        legend pos=north west,
        legend style={draw=none, legend columns=2},
        mark size=1pt,
        xtick=\empty,
        xticklabels={}
      ]
      \addplot[red, mark=square*] table[x=x,y=y2] {figs/drone_dtw/plot_data.dat};
        \addplot[blue, mark=*] table[x=x,y=y1] {figs/drone_dtw/plot_data.dat};
        \legend{$\text{NFVM}_4^5$,FVM}
      \end{axis}
    \end{tikzpicture}
    \caption{Dynamic Time Warping (DTW) distance between the best FVM fit (see~\Cref{fig:drone_heatmap_best}) and $\text{NFVM}_4^5$ solution time series for each fixed $x$ position along the road.}
    \label{fig:drone_dtw}
\end{figure}
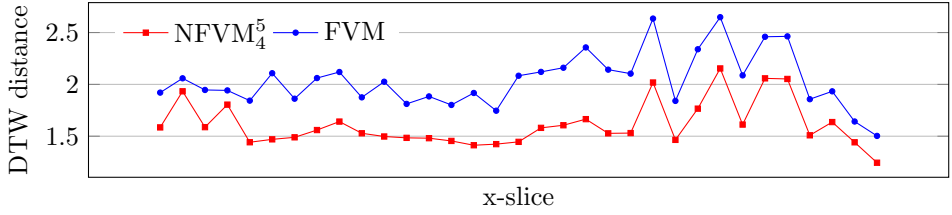 

%% file: tex/drone_heatmap_best.tex
\begin{figure}
    \centering
  \begin{subfigure}[t]{.495\textwidth}
      \begin{tikzpicture}[scale=.95, clip=false]
          \scriptsize
                  \node[anchor=south west, inner sep=0] (img) at (0,0) {\includegraphics[width=.9\textwidth]{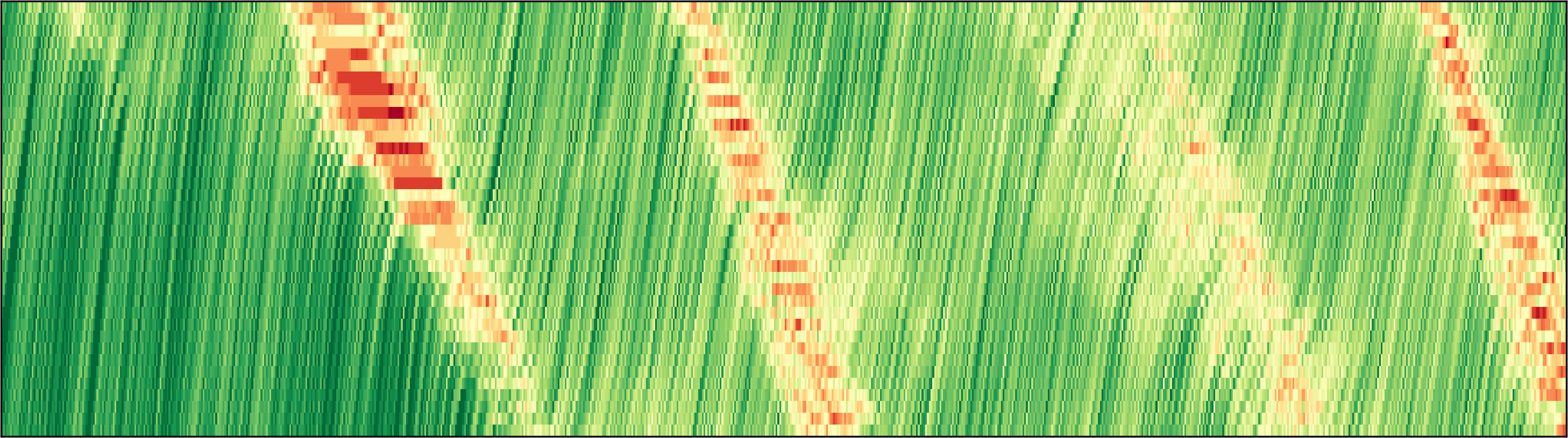}};
                  
                  \draw[black, thick] (img.south west) rectangle (img.north east);
          
                  \path (img.south west) -- (img.south east) node[midway, below, yshift=-0.08cm] {\(t\)};
                  \path (img.south west) -- (img.north west) node[midway, left, xshift=-0.08cm] {\(x\)};
                  
                  \node[below, yshift=-0.08cm] at (img.south west) {\(0\)};
                  \node[below, yshift=-0.08cm, xshift=-.2cm] at (img.south east) {\(870s\)};
                  \node[left, xshift=-0.08cm] at (img.south west) {\(0\)};
                  \node[left, xshift=.4cm, yshift=.15cm] at (img.north west) {\(400m\)};
              \end{tikzpicture}
              \vspace{-3.2cm}
              \caption{Ground truth}
          \end{subfigure}
          \hfill
          \begin{subfigure}[t]{.495\textwidth}
              \begin{tikzpicture}[scale=.95, clip=false]
                  \scriptsize
                          \node[anchor=south west, inner sep=0] (img) at (0,0) {\includegraphics[width=.9\textwidth]{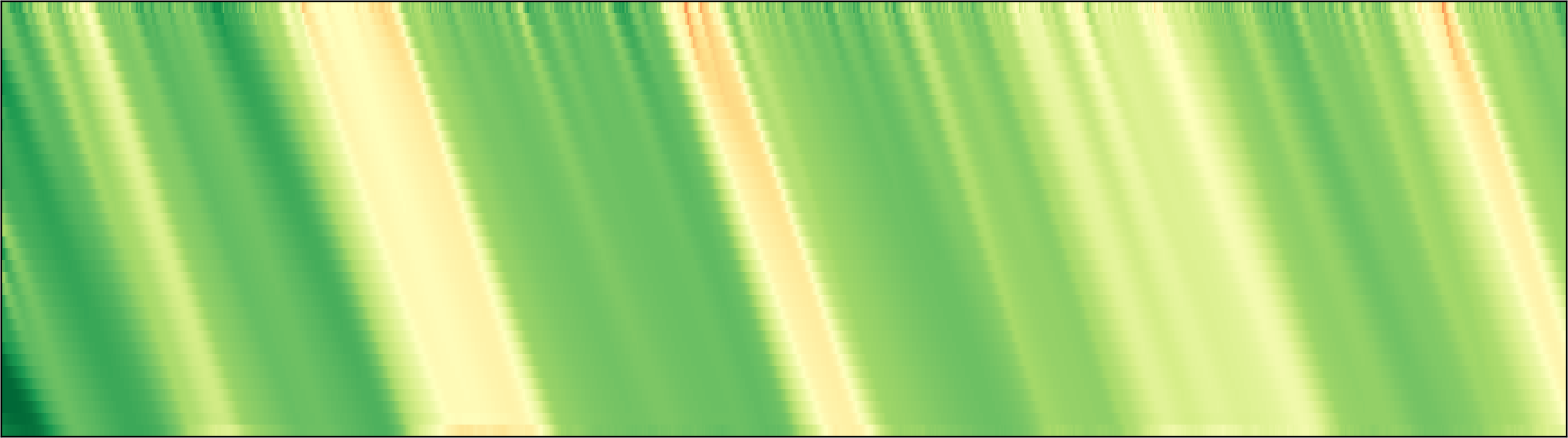}};
                      
                      \draw[black, thick] (img.south west) rectangle (img.north east);
              
                      \path (img.south west) -- (img.south east) node[midway, below, yshift=-0.08cm] {\(t\)};
                      \path (img.south west) -- (img.north west) node[midway, left, xshift=-0.08cm] {\(x\)};
                      
                      \node[below, yshift=-0.08cm] at (img.south west) {\(0\)};
                      \node[below, yshift=-0.08cm, xshift=-.2cm] at (img.south east) {\(870s\)};
                      \node[left, xshift=-0.08cm] at (img.south west) {\(0\)};
                      \node[left, xshift=.4cm, yshift=.15cm] at (img.north west) {\(400m\)};
                  \end{tikzpicture}
                  \vspace{-3.2cm}
                  \caption{$\text{NFVM}_2^1$}
          \end{subfigure}\vspace{0.3cm}
          \begin{subfigure}[t]{.495\textwidth} 
              \begin{tikzpicture}[scale=.95, clip=false]
                  \scriptsize
                          \node[anchor=south west, inner sep=0] (img) at (0,0) {\includegraphics[width=.9\textwidth]{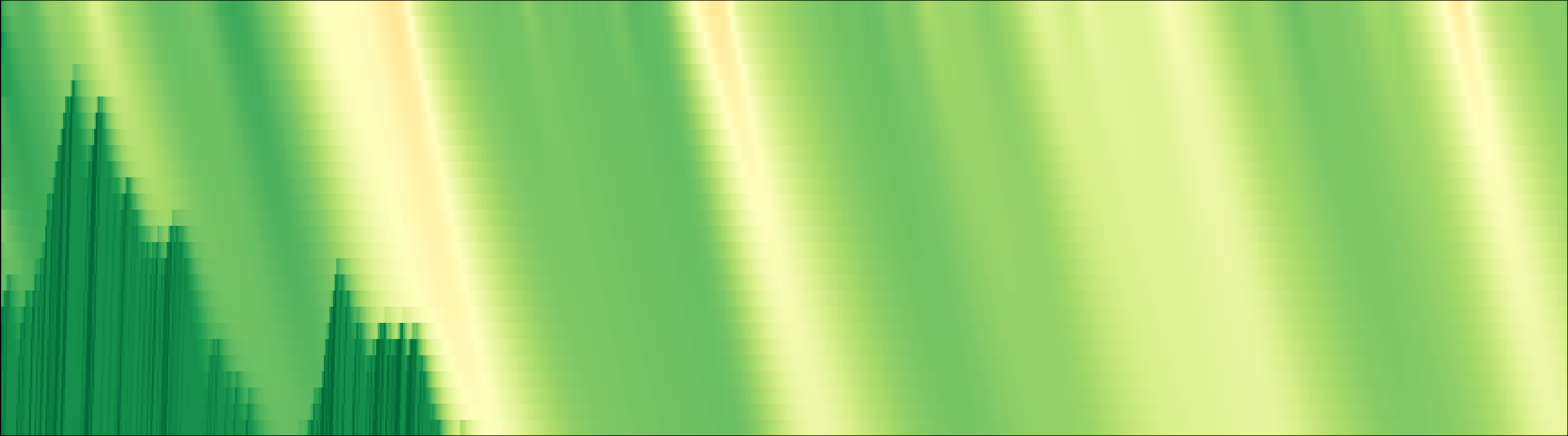}}; 
              
              \draw[black, thick] (img.south west) rectangle (img.north east);
      
              \path (img.south west) -- (img.south east) node[midway, below, yshift=-0.08cm] {\(t\)};
              \path (img.south west) -- (img.north west) node[midway, left, xshift=-0.08cm] {\(x\)};
              
              \node[below, yshift=-0.08cm] at (img.south west) {\(0\)};
              \node[below, yshift=-0.08cm, xshift=-.2cm] at (img.south east) {\(870s\)};
              \node[left, xshift=-0.08cm] at (img.south west) {\(0\)};
              \node[left, xshift=.4cm, yshift=.15cm] at (img.north west) {\(400m\)};
          \end{tikzpicture}
          \vspace{-3.2cm}
          \caption{Representative best FVM fit}
      \end{subfigure}
      \hfill
      \begin{subfigure}[t]{.495\textwidth}
          \begin{tikzpicture}[scale=.95, clip=false]
              \scriptsize
                  \node[anchor=south west, inner sep=0] (img) at (0,0) {\includegraphics[width=.9\textwidth]{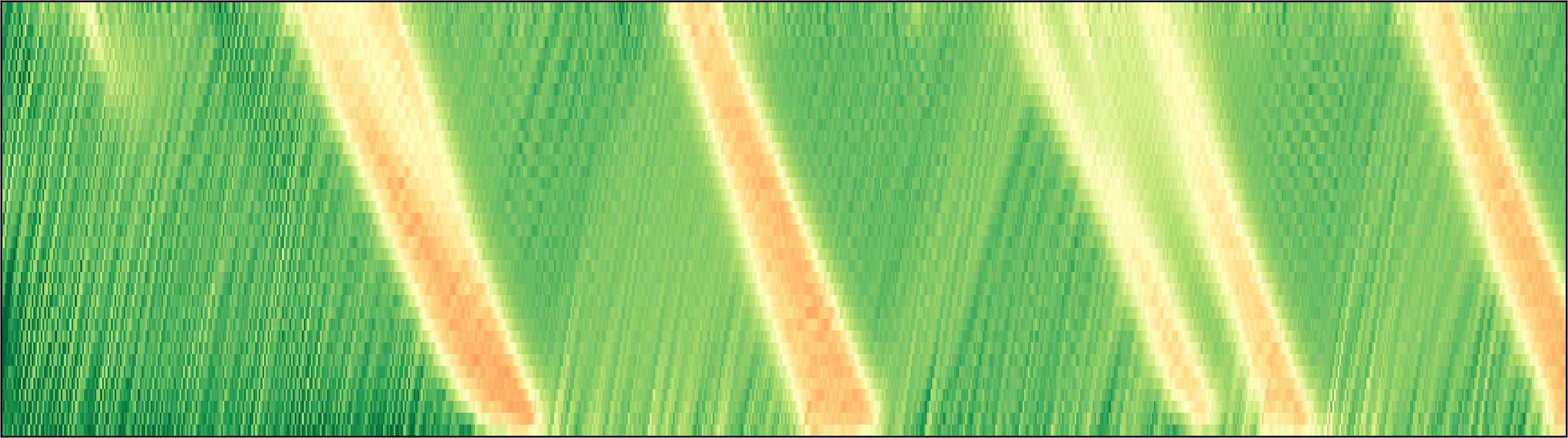}};
                  
                  \draw[black, thick] (img.south west) rectangle (img.north east);
          
                  \path (img.south west) -- (img.south east) node[midway, below, yshift=-0.08cm] {\(t\)};
                  \path (img.south west) -- (img.north west) node[midway, left, xshift=-0.08cm] {\(x\)};
                  
                  \node[below, yshift=-0.08cm] at (img.south west) {\(0\)};
                  \node[below, yshift=-0.08cm, xshift=-.2cm] at (img.south east) {\(870s\)};
                  \node[left, xshift=-0.08cm] at (img.south west) {\(0\)};
                  \node[left, xshift=.4cm, yshift=.15cm] at (img.north west) {\(400m\)};
              \end{tikzpicture}
              \vspace{-3.2cm}
              \caption{$\text{NFVM}_4^5$}
            \end{subfigure}
      \caption{Comparison of the performance of NFVM$_2^1$ and NFVM$_4^5$ against a representative FVM fit. The density is color-coded as in~\Cref{fig:drone_dataset_heatmap}.}
    \label{fig:drone_heatmap_best}
  \end{figure}

%% file: tex/flows.tex
\begin{figure}
    \centering
    \begin{subfigure}[t]{.35\textwidth}
        \centering
        \begin{tikzpicture}
            \begin{axis}[
                width=4.3cm, height=3cm,
                xlabel={Density \(u\)},
                ylabel={Flux \(f(u)\)},
                grid=major,
                samples=200,
                enlargelimits=.1
            ]
                    \addplot[ultra thick, smooth, color=blue, domain=0:1] {x*(1-x)};
            \end{axis}
        \end{tikzpicture}
        \caption{Greenshields'}
    \end{subfigure}
    \hfill
    \begin{subfigure}[t]{.3\textwidth}
        \centering
        \begin{tikzpicture}
            \begin{axis}[
                width=4.3cm, height=3cm,
                xlabel={Density \(u\)},
                grid=major,
                samples=200,
                enlargelimits=.1
            ]
                \addplot[ultra thick, smooth, color=blue, domain=0.:.5] {x};
                \addplot[ultra thick, smooth, color=blue, domain=.5:1., forget plot] {1-x};
            \end{axis}
        \end{tikzpicture}
        \caption{Triangular Sym}
    \end{subfigure}
    \hfill
    \begin{subfigure}[t]{.3\textwidth}
        \centering
        \begin{tikzpicture}
            \begin{axis}[
                width=4.3cm, height=3cm,
                xlabel={Density \(u\)},
                grid=major,
                samples=200,
                enlargelimits=.1
            ]
                \addplot[ultra thick, smooth, color=blue, domain=0.:1/3] {2*x};
                \addplot[ultra thick, smooth, color=blue, domain=1/3:1.] {1-x};
            \end{axis}
        \end{tikzpicture}
        \caption{Triangular Skw}
    \end{subfigure}
    \begin{subfigure}[t]{.35\textwidth}
        \centering
        \begin{tikzpicture}
            \begin{axis}[
                width=4.3cm, height=3cm,
                xlabel={Density \(u\)},
                ylabel={Flux \(f(u)\)},
                grid=major,
                samples=200,
                enlargelimits=.1
            ]
                \addplot[ultra thick, smooth, color=blue, domain=0:.2] {x};
                \addplot[ultra thick, smooth, color=blue, domain=.2:.8, forget plot] {.1*(x-.2) + .2};
                \addplot[ultra thick, smooth, color=blue, domain=.8:1, forget plot] {1.3 - 1.3*x};
            \end{axis}
        \end{tikzpicture}
        \caption{Trapezoidal}
    \end{subfigure}
    \hfill
    \begin{subfigure}[t]{.3\textwidth}
        \centering
        \begin{tikzpicture}
            \begin{axis}[
                width=4.3cm, height=3cm,
                xlabel={Density \(u\)},
                grid=major,
                samples=200,
                enlargelimits=.1
            ]
                \addplot[ultra thick, smooth, color=blue, domain=0:1] {-2*x*ln(x)};
            \end{axis}
        \end{tikzpicture}
        \caption{Greenberg}
    \end{subfigure}
    \hfill
    \begin{subfigure}[t]{.3\textwidth}
        \centering
        \begin{tikzpicture}
            \begin{axis}[
                width=4.3cm, height=3cm,
                xlabel={Density \(u\)},
                grid=major,
                samples=200,
                enlargelimits=.1
            ]
                \addplot[ultra thick, smooth, color=blue, domain=0.:1.] {.25*exp(1 - x)*x};
            \end{axis}
        \end{tikzpicture}
        \caption{Underwood}
    \end{subfigure}
    \caption{LWR flow functions with default parameters from~\Cref{app:lwr_flows}.}
	\label{fig:flows}
\end{figure}
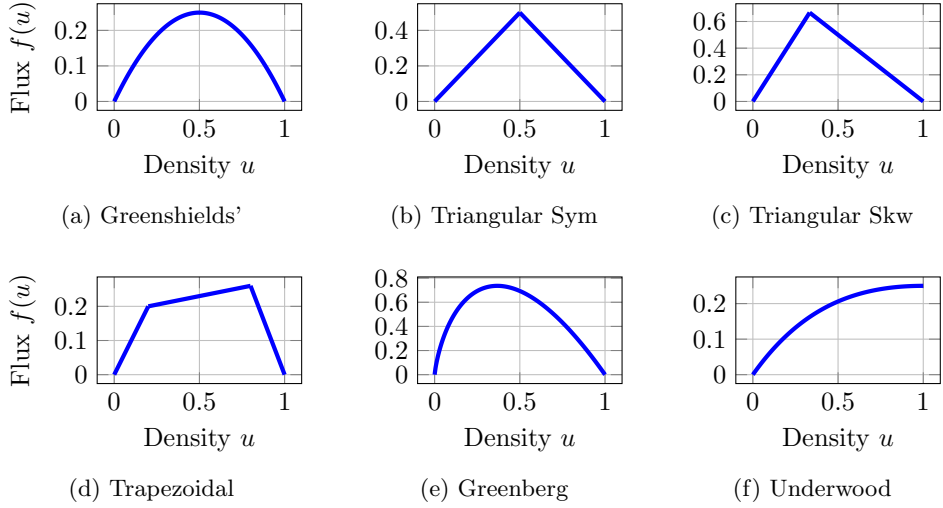

%% file: references.bib
@article{storn1997differential,
  title={Differential evolution--a simple and efficient heuristic for global optimization over continuous spaces},
  author={Storn, Rainer and Price, Kenneth},
  journal={Journal of global optimization},
  volume={11},
  pages={341--359},
  year={1997},
  publisher={Springer}
}

@misc{shin2020error,
  author       = {Y. Shin and Z. Zhang and G. E. Karniadakis},
  title        = {Error Estimates of Residual Minimization Using Neural Networks for Linear PDEs},
  year         = {2020},
  note         = {Preprint},
  eprint       = {2010.08019},
  archivePrefix= {arXiv},
  primaryClass = {math.NA}
}

@article{mishra2023generalization,
  author       = {S. Mishra and R. Molinaro},
  title        = {Estimates on the generalization error of physics-informed neural networks for approximating PDEs},
  journal      = {IMA Journal of Numerical Analysis},
  volume       = {43},
  year         = {2023},
  pages        = {1--43}
}

@article{mishra2022inverse,
  author       = {S. Mishra and R. Molinaro},
  title        = {Estimates on the generalization error of physics-informed neural networks for approximating a class of inverse problems for PDEs},
  journal      = {IMA Journal of Numerical Analysis},
  volume       = {42},
  year         = {2022},
  pages        = {981--1022}
}

@misc{deryck2021kolmogorov,
  author       = {T. De Ryck and S. Mishra},
  title        = {Error Analysis for Physics Informed Neural Networks (PINNs) Approximating Kolmogorov PDEs},
  year         = {2021},
  note         = {Preprint},
  eprint       = {2106.14473},
  archivePrefix= {arXiv},
  primaryClass = {math.NA}
}

@article{nelsen2021random,
  title={The random feature model for input-output maps between Banach spaces},
  author={Nelsen, Nicholas H and Stuart, Andrew M},
  journal={SIAM Journal on Scientific Computing},
  volume={43},
  number={5},
  pages={A3212--A3243},
  year={2021},
  publisher={SIAM}
}

@article{patel2018operator,
  title={Nonlinear integro-differential operator regression with neural networks},
  author={Patel, R.G. and Desjardins, O.},
  journal={arXiv preprint arXiv:1810.08552},
  year={2018}
}

@article{tripura2023wavelet,
  title={Wavelet neural operator for solving parametric partial differential equations in computational mechanics problems},
  author={Tripura, Tapabrata and Chakraborty, Souvik},
  journal={Computer Methods in Applied Mechanics and Engineering},
  volume={404},
  pages={115783},
  year={2023},
  publisher={Elsevier}
}

@article{chen1995universal,
  author={T. Chen and H. Chen},
  title={Universal approximation to nonlinear operators by neural networks with arbitrary activation functions and its application to dynamical systems},
  journal={IEEE Transactions on Neural Networks},
  volume={6},
  number={4},
  pages={911--917},
  year={1995}
}

@article{lanthaler2022error,
  author={S. Lanthaler and S. Mishra and G. E. Karniadakis},
  title={Error estimates for DeepONets: A deep learning framework in infinite dimensions},
  journal={Transactions of Mathematics and Its Applications},
  volume={6},
  number={1},
  year={2022}
}

@article{kovachki2021fourier,
  author={N. B. Kovachki and S. Lanthaler and S. Mishra},
  title={On universal approximation and error bounds for Fourier neural operators},
  journal={Journal of Machine Learning Research},
  volume={22},
  number={1},
  year={2021}
}

@article{bhattacharya2021model,
  author={K. Bhattacharya and B. Hosseini and N. B. Kovachki and A. M. Stuart},
  title={Model reduction and neural networks for parametric PDEs},
  journal={The SMAI Journal of Computational Mathematics},
  volume={7},
  pages={121--157},
  year={2021}
}

@article{herrmann2022neural,
  title={Neural and GPC operator surrogates: Construction and expression rate bounds},
  author={Herrmann, Lukas and Schwab, Christoph and Zech, Jonas},
  journal={arXiv preprint arXiv:2207.04950},
  year={2022}
}

@article{dehoop2023convergence,
  title={Convergence rates for learning linear operators from noisy data},
  author={de Hoop, Maarten V. and Kovachki, Nikola B. and Nelsen, Nicholas H. and Stuart, Andrew M.},
  journal={SIAM/ASA Journal on Uncertainty Quantification},
  volume={11},
  number={2},
  pages={480--513},
  year={2023},
  publisher={SIAM}
}

@book{vershynin2018high,
  title     = {High-Dimensional Probability: An Introduction with Applications in Data Science},
  author    = {Vershynin, Roman},
  year      = {2018},
  publisher = {Cambridge University Press},
  series    = {Cambridge Series in Statistical and Probabilistic Mathematics},
  doi       = {10.1017/9781108231596}
}

@article{Sanders1983,
  author    = {Sanders, Richard},
  title     = {On convergence of monotone finite difference schemes with variable spatial differencing},
  journal   = {Mathematics of Computation},
  volume    = {40},
  number    = {161},
  pages     = {91--106},
  year      = {1983},
  publisher = {American Mathematical Society},
  doi       = {10.2307/2007364},
  url       = {https://www.jstor.org/stable/2007364}
}

@incollection{barth2004finite,
  author    = {Timothy Barth and Raphaèle Herbin and Mario Ohlberger},
  title     = {Finite Volume Methods: Foundation and Analysis},
  booktitle = {Encyclopedia of Computational Mechanics},
  editor    = {E. Stein and R. de Borst and T. J. R. Hughes},
  publisher = {John Wiley \& Sons},
  volume    = {1},
  pages     = {439--466},
  year      = {2004},
  address   = {Chichester, UK}
}

@article{lax1955nonlinear,
  title     = {Integration of Nonlinear Equations of Evolution and Nonlinear Theory of Shock Waves},
  author    = {Lax, Peter D. and Friedrichs, Kurt O.},
  journal   = {Communications on Pure and Applied Mathematics},
  volume    = {8},
  number    = {2},
  pages     = {345--392},
  year      = {1955},
  publisher = {Wiley},
  doi       = {10.1002/cpa.3160080207}
}

@article{boulle2023learning,
  title={Learning elliptic partial differential equations with randomized linear algebra},
  author={Boull{\'e}, Nicolas and Townsend, Alex},
  journal={Foundations of Computational Mathematics},
  volume={23},
  number={2},
  pages={709--739},
  year={2023},
  publisher={Springer}
}

@misc{wang2023operatorlearninghyperbolicpartial,
      title={Operator learning for hyperbolic partial differential equations}, 
      author={Christopher Wang and Alex Townsend},
      year={2023},
      eprint={2312.17489},
      archivePrefix={arXiv},
      primaryClass={math.NA},
      url={https://arxiv.org/abs/2312.17489}, 
}

@article{castro2023kolmogorov,
  author={J. Castro},
  title={The Kolmogorov infinite dimensional equation in a Hilbert space via deep learning methods},
  journal={Journal of Mathematical Analysis and Applications},
  volume={527},
  number={2},
  pages={127413},
  year={2023}
}

@inproceedings{raonic2023iclr,
  title={Convolutional neural operators},
  author={Raonic, B. and Molinaro, R. and Rohner, T. and Mishra, S. and de Bezenac, E.},
  booktitle={ICLR 2023 Workshop on Physics for Machine Learning},
  year={2023}
}

@inproceedings{raonic2023neurips,
  title={Convolutional neural operators for robust and accurate learning of PDEs},
  author={Raonic, B. and Molinaro, R. and De Ryck, T.D. and Rohner, T. and Bartolucci, F. and Alaifari, R. and Mishra, S. and de Bezenac, E.},
  booktitle={Thirty-seventh Conference on Neural Information Processing Systems (NeurIPS)},
  year={2023}
}

@misc{deryck2022wpinnsweakphysicsinformed,
      title={wPINNs: Weak Physics informed neural networks for approximating entropy solutions of hyperbolic conservation laws}, 
      author={Tim De Ryck and Siddhartha Mishra and Roberto Molinaro},
      year={2022},
      eprint={2207.08483},
      archivePrefix={arXiv},
      primaryClass={math.NA},
      url={https://arxiv.org/abs/2207.08483}, 
}

@article{kharazmi2019vpinn,
  title={Variational physics-informed neural networks for solving partial differential equations},
  author={Kharazmi, Ehsan and Zhang, Zongyi and Karniadakis, George Em},
  journal={arXiv preprint arXiv:1912.00873},
  year={2019},
  url={https://arxiv.org/abs/1912.00873}
}

@misc{deryck2022navier,
  author       = {T. De Ryck and A. D. Jagtap and S. Mishra},
  title        = {Error Analysis for Physics Informed Neural Networks Approximating the Navier-Stokes Equations},
  year         = {2022},
  note         = {Preprint},
  eprint       = {2203.09346},
  archivePrefix= {arXiv},
  primaryClass = {math.NA}
}

@misc{deryck2022generic,
  author       = {T. De Ryck and S. Mishra},
  title        = {Generic Bounds on the Approximation Error for Physics-Informed (and) Operator Learning},
  year         = {2022},
  note         = {Preprint},
  eprint       = {2205.11393},
  archivePrefix= {arXiv},
  primaryClass = {math.NA}
}

@article{ardekani2011macroscopic,
  title={Macroscopic speed-flow models for characterization of freeway and managed lanes},
  author={Ardekani, S and Ghandehari, Mostafa and Nepal, S},
  journal={Institutul Politehnic din Iasi. Buletinul. Sectia Constructii. Arhitectura},
  volume={57},
  number={1},
  pages={149},
  year={2011},
  publisher={Institutul Politehnic" Gheorghe Asachi" din Iasi}
}

@article{mazare2011analytical,
  title={Analytical and grid-free solutions to the Lighthill--Whitham--Richards traffic flow model},
  author={Mazar{\'e}, Pierre-Emmanuel and Dehwah, Ahmad H and Claudel, Christian G and Bayen, Alexandre M},
  journal={Transportation Research Part B: Methodological},
  volume={45},
  number={10},
  pages={1727--1748},
  year={2011},
  publisher={Elsevier}
}

@article{claudel2010lax_a,
  title={Lax--Hopf based incorporation of internal boundary conditions into Hamilton--Jacobi equation. Part I: Theory},
  author={Claudel, Christian G and Bayen, Alexandre M},
  journal={IEEE Transactions on Automatic Control},
  volume={55},
  number={5},
  pages={1142--1157},
  year={2010},
  publisher={IEEE}
}

@article{claudel2010lax_b,
  title={Lax--hopf based incorporation of internal boundary conditions into hamilton-jacobi equation. part ii: Computational methods},
  author={Claudel, Christian G and Bayen, Alexandre M},
  journal={IEEE Transactions on Automatic Control},
  volume={55},
  number={5},
  pages={1158--1174},
  year={2010},
  publisher={IEEE}
}

@Misc{amsmath,
  author =	 {{American Mathematical Society}},
  title =	 {User's Guide for the \texttt{amsmath} Package
                  (Version 2.0)},
  url =		 {ftp://ftp.ams.org/pub/tex/doc/amsmath/amsldoc.pdf},
  urldate =	 {2015-07-30},
  year =	 2002}

@Misc{pgfplots,
  author =	 {Christian Feuers\"anger},
  title =	 {Manual for Package \texttt{PGFPLOTS}},
  month =	 may,
  year =	 2015,
  url =		 {http://sourceforge.net/projects/pgfplots}
}

@article{R56,
	author = {Richards, P. I.},
	journal = {Operations Research},
	pages = {42-51},
	title = {Shock waves on the highway},
	volume = {4},
	year = {1956}}

@article{Wilson01061927,
author = {Edwin B. Wilson and},
title = {Probable Inference, the Law of Succession, and Statistical Inference},
journal = {Journal of the American Statistical Association},
volume = {22},
number = {158},
pages = {209--212},
year = {1927},
publisher = {ASA Website},
doi = {10.1080/01621459.1927.10502953},
URL = {       https://www.tandfonline.com/doi/abs/10.1080/01621459.1927.10502953
},
eprint = {       https://www.tandfonline.com/doi/pdf/10.1080/01621459.1927.10502953
}
}

@article{lighthill1955kinematic,
  title={On kinematic waves I. Flood movement in long rivers},
  author={Lighthill, Michael James and Whitham, G Be},
  journal={Proceedings of the Royal Society of London. Series A. Mathematical and Physical Sciences},
  volume={229},
  number={1178},
  pages={281--316},
  year={1955},
  publisher={The Royal Society London}
}

@book{holden2015front,
  title={Front tracking for hyperbolic conservation laws},
  author={Holden, Helge and Risebro, Nils Henrik},
  volume={152},
  year={2015},
  publisher={Springer}
}

@article{Greenberg,
	author = {Harold Greenberg},
	journal = {Operations Research},
	number = {1},
	pages = {79--85},
	title = {An Analysis of Traffic Flow},
	volume = {7},
	year = {1959}}

@article{godunov1959finite,
	author = {Godunov, SK},
	journal = {Sbornik: Mathematics},
	number = {8-9},
	pages = {357--393},
	title = {A finite difference method for the computation of discontinuous solutions of the equations of fluid dynamics.},
	volume = {47},
	year = {1959}}

@article{Underwood1961,
	author = {Underwood, Robin T},
	journal = {Quality and theory of traffic flow},
	title = {Speed, volume and density relationships},
	year = {1961}}

@article{ClaudelBayen2010,
	author = {Claudel, Christian G. and Bayen, Alexandre M.},
	journal = {IEEE Transactions on Automatic Control},
	number = {5},
	pages = {1158-1174},
	title = {{Lax--Hopf Based Incorporation of Internal Boundary Conditions Into Hamilton-Jacobi Equation. Part II: Computational Methods}},
	volume = {55},
	year = {2010}}

@article{SimoniClaudel2017,
	author = {Michele D. Simoni and Christian G. Claudel},
	journal = {Transportation Research Part B: Methodological},
	pages = {238-255},
	title = {A fast simulation algorithm for multiple moving bottlenecks and applications in urban freight traffic management},
	volume = {104},
	year = {2017}}

@article{greenberg1959analysis,
  title={An analysis of traffic flow},
  author={Greenberg, Harold},
  journal={Operations research},
  volume={7},
  number={1},
  pages={79--85},
  year={1959},
  publisher={INFORMS}
}

@book{leveque2002finite,
  title={Finite volume methods for hyperbolic problems},
  author={LeVeque, Randall J},
  volume={31},
  year={2002},
  publisher={Cambridge university press}
}

@inproceedings{greenshields1935study,
  title={A study of traffic capacity},
  author={Greenshields, Bruce D and Bibbins, J Rowland and Channing, WS and Miller, Harvey H},
  booktitle={Highway research board proceedings},
  volume={14},
  number={1},
  pages={448--477},
  year={1935},
  organization={Washington, DC}
}

@article{godounov1959difference,
  title={A difference method for numerical calculation of discontinuous solutions of the equation of hydrodynamics},
  author={Godounov, SK},
  journal={Matematicheskii Sbornik},
  volume={47},
  number={89-3},
  pages={271--306},
  year={1959}
}

@article{wu2022b3d,
  title={Decentralized Vehicle Coordination: The Berkeley DeepDrive Drone Dataset},
  author={Fangyu Wu and Dequan Wang and Minjune Hwang and Chenhui Hao and Jiawei Lu and Jiamu Zhang and Christopher Chou and Trevor Darrell and Alexandre Byen},
  journal={arXiv},
  year={2022}
}

@article{morand2024deep,
  title={Deep learning of first-order nonlinear hyperbolic conservation law solvers},
  author={Morand, Victor and M{\"u}ller, Nils and Weightman, Ryan and Piccoli, Benedetto and Keimer, Alexander and Bayen, Alexandre M},
  journal={Journal of Computational Physics},
  volume={511},
  pages={113114},
  year={2024},
  publisher={Elsevier}
}

@book{evans2022partial,
  title={Partial differential equations},
  author={Evans, Lawrence C},
  volume={19},
  year={2022},
  publisher={American Mathematical Society}
}

@article{de2024wpinns,
  title={wPINNs: Weak physics informed neural networks for approximating entropy solutions of hyperbolic conservation laws},
  author={De Ryck, Tim and Mishra, Siddhartha and Molinaro, Roberto},
  journal={SIAM Journal on Numerical Analysis},
  volume={62},
  number={2},
  pages={811--841},
  year={2024},
  publisher={SIAM}
}

@article{raissi2017physics,
  title={Physics informed deep learning (part i): Data-driven solutions of nonlinear partial differential equations},
  author={Raissi, Maziar and Perdikaris, Paris and Karniadakis, George Em},
  journal={arXiv preprint arXiv:1711.10561},
  year={2017}
}

@article{engquistOsher1981,
  title={One-sided difference approximations for nonlinear conservation laws},
  author={Engquist, Bj{\"o}rn and Osher, Stanley},
  journal={Mathematics of Computation},
  volume={36},
  number={154},
  pages={321--351},
  year={1981}
}

@incollection{shu1999high,
  title={High order ENO and WENO schemes for computational fluid dynamics},
  author={Shu, Chi-Wang},
  booktitle={High-order methods for computational physics},
  pages={439--582},
  year={1999},
  publisher={Springer}
}

@Inbook{Bertoluzza2009,
  title={Numerical solutions of partial differential equations},
  author={Bertoluzza, Silvia and Falletta, Silvia and Russo, Giovanni and Shu, Chi-Wang},
  year={2009},
  publisher={Springer Science \& Business Media},
  pages="89--95",
isbn="978-3-7643-8940-6",
doi="10.1007/978-3-7643-8940-6_7",
url="https://doi.org/10.1007/978-3-7643-8940-6_7"
}

@article{lax1954initial,
  title={The initial value problem for nonlinear hyperbolic equations in two independent variables},
  author={Lax, PD},
  journal={Ann. Math. Studies},
  volume={33},
  number={21},
  pages={1--229},
  year={1954},
  publisher={Princeton U. Press}
}

@article{shin2020convergence,
  title={On the convergence of physics informed neural networks for linear second-order elliptic and parabolic type PDEs},
  author={Shin, Yeonjong and Darbon, Jerome and Karniadakis, George Em},
  journal={arXiv preprint arXiv:2004.01806},
  year={2020}
}

@article{van1989sonic,
  title={Sonic-point capturing},
  author={Van Leer, Bram and Lee, Wen-Tzong and Powell, Kenneth G},
  year={1989},
  publisher={American Institute of Aeronautics and Astronautics}
}

@article{cameron2011notes,
  title={Notes on the Burgers Equation},
  author={Cameron, MARIA},
  journal={University of Marylan},
  year={2011}
}

@article{kruvzkov1970first,
  title={First order quasilinear equations in several independent variables},
  author={Kru{\v{z}}kov, Stanislav N},
  journal={Mathematics of the USSR-Sbornik},
  volume={10},
  number={2},
  pages={217},
  year={1970},
  publisher={IOP Publishing}
}

@article{geroliminis2008existence,
  title={Existence of urban-scale macroscopic fundamental diagrams: Some experimental findings},
  author={Geroliminis, Nikolas and Daganzo, Carlos F},
  journal={Transportation Research Part B: Methodological},
  volume={42},
  number={9},
  pages={759--770},
  year={2008},
  publisher={Elsevier}
}

@article{geroliminis2011properties,
  title={Properties of a well-defined macroscopic fundamental diagram for urban traffic},
  author={Geroliminis, Nikolas and Sun, Jie},
  journal={Transportation Research Part B: Methodological},
  volume={45},
  number={3},
  pages={605--617},
  year={2011},
  publisher={Elsevier}
}

@article{LWR,
  title={Traffic current fluctuation and the Burgers equation},
  author={Musha, Toshimitsu and Higuchi, Hideyo},
  journal={Japanese journal of applied physics},
  volume={17},
  number={5},
  pages={811},
  year={1978},
  publisher={IOP Publishing}
}

@article{ENO,
  author    = {A. Harten and B. Engquist and S. Osher and S. Chakravarthy},
  title     = {Uniformly High Order Accurate Essentially Non-oscillatory Schemes, III},
  journal   = {Journal of Computational Physics},
  volume    = {71},
  number    = {2},
  pages     = {231--303},
  year      = {1987},
  doi       = {10.1016/0021-9991(87)90045-1}
}

@article{WENO,
  author    = {Guan-Shan Jiang and Chi-Wang Shu},
  title     = {Efficient Implementation of Weighted ENO Schemes},
  journal   = {Journal of Computational Physics},
  volume    = {126},
  number    = {1},
  pages     = {202--228},
  year      = {1996},
  doi       = {10.1006/jcph.1996.0130}
}

@article{first_paper_DG,
  author    = {Cockburn, Bernardo and Shu, Chi-Wang},
  title     = {The Runge--Kutta Discontinuous Galerkin Method for Conservation Laws V: Multidimensional Systems},
  journal   = {Journal of Computational Physics},
  volume    = {141},
  number    = {2},
  pages     = {199--224},
  year      = {1998},
  doi       = {10.1006/jcph.1998.5892}
}

@article{kuznetsov1976accuracy,
  author    = {N. N. Kuznetsov},
  title     = {The Accuracy of Some Approximate Methods for Computing Weak Solutions of a First-Order Quasi-Linear Equation},
  journal   = {Zhurnal Vychislitel{\textasciiacute}noi Matematiki i Matematicheskoi Fiziki},
  volume    = {16},
  number    = {6},
  pages     = {1489--1502},
  year      = {1976},
  note      = {English translation: {\em USSR Computational Mathematics and Mathematical Physics}, 16(6):105–119, 1976},
  doi       = {10.1016/0041-5553(76)90046-X}
}

@inproceedings{first_FNO,
  title={Fourier neural operator for parametric partial differential equations},
  author={Li, Zongyi and Kovachki, Nikola and Azizzadenesheli, Kamyar and Liu, Burigede and Stuart, Andrew and Anandkumar, Anima},
  booktitle={International Conference on Learning Representations (ICLR)},
  year={2021}
}

@article{first_Deep0Net,
  title={Learning nonlinear operators via DeepONet based on the universal approximation theorem of operators},
  author={Lu, Lu and Jin, Pengzhan and Karniadakis, George Em},
  journal={Nature Machine Intelligence},
  volume={3},
  number={3},
  pages={218--229},
  year={2021},
  publisher={Nature Publishing Group}
}

@misc{lanthaler2023parametric,
  title        = {The Parametric Complexity of Operator Learning},
  author       = {Lanthaler, Samuel and Stuart, Andrew M.},
  year         = {2023},
  eprint       = {2306.15924},
  archivePrefix= {arXiv},
  primaryClass = {cs.LG},
  note         = {Preprint},
  url          = {https://arxiv.org/abs/2306.15924}
}

@inproceedings{lanthaler2022nonlinear,
  title     = {Nonlinear Reconstruction for Operator Learning of PDEs with Discontinuities},
  author    = {Lanthaler, Samuel and Molinaro, Roberto and Hadorn, Patrik and Mishra, Siddhartha},
  booktitle = {International Conference on Learning Representations (ICLR)},
  year      = {2023},
  arxivId   = {2210.01074},
  url       = {https://arxiv.org/abs/2210.01074},
  note      = {Presented as oral (top 25\%) at ICLR 2023},
}

@article{thodi2023learning,
  title        = {Learning-based solutions to nonlinear hyperbolic PDEs: Empirical insights on generalization errors},
  author       = {Thodi, Bilal Thonnam and Ambadipudi, Sai Venkata Ramana and Jabari, Saif Eddin},
  journal      = {arXiv preprint arXiv:2302.08144},
  year         = {2023},
  month        = feb,
  eprint       = {2302.08144},
  archivePrefix= {arXiv},
  primaryClass = {cs.LG},
  note         = {Presented at NeurIPS Workshop on Machine Learning for the Physical Sciences},
  url          = {https://arxiv.org/abs/2302.08144}
}

@article{kim2025approximating,
  title        = {Approximating Numerical Fluxes Using Fourier Neural Operators for Hyperbolic Conservation Laws},
  author       = {Kim, Taeyoung and Kang, Myungjoo},
  journal      = {Communications in Computational Physics},
  volume       = {37},
  number       = {2},
  pages        = {420--456},
  year         = {2025},
  doi          = {10.4208/cicp.OA-2024-0123},
  eprint       = {2401.01783},
  archivePrefix= {arXiv},
  primaryClass = {math.NA},
  note         = {Published online Jan 2025},
  url          = {https://arxiv.org/abs/2401.01783}
}

@article{lu2019deeponet,
  title        = {DeepONet: Learning nonlinear operators for identifying differential equations based on the universal approximation theorem of operators},
  author       = {Lu, Lu and Jin, Pengzhan and Karniadakis, George Em},
  journal      = {arXiv preprint arXiv:1910.03193},
  year         = {2019},
  month        = oct,
  eprint       = {1910.03193},
  archivePrefix= {arXiv},
  primaryClass = {cs.LG},
  url          = {https://arxiv.org/abs/1910.03193}
}

@article{li2021physics,
  title        = {Physics-Informed Neural Operator for Learning Partial Differential Equations},
  author       = {Li, Zongyi and Zheng, Hongkai and Kovachki, Nikola and Jin, David and Chen, Haoxuan and Liu, Burigede and Azizzadenesheli, Kamyar and Anandkumar, Anima},
  journal      = {arXiv preprint arXiv:2111.03794},
  year         = {2021},
  month        = nov,
  eprint       = {2111.03794},
  archivePrefix= {arXiv},
  primaryClass = {cs.LG},
  url          = {https://arxiv.org/abs/2111.03794}
}

@article{NwaigweMungkasi2021,
  author       = {Chinedu Nwaigwe and Sudi Mungkasi},
  title        = {Comparison of Different Numerical Schemes for 1D Conservation Laws},
  journal      = {Journal of Interdisciplinary Mathematics},
  volume       = {24},
  number       = {3},
  pages        = {537--552},
  year         = {2021},
  month        = {aug},
  doi          = {10.1080/09720502.2020.1792665},
  url          = {https://doi.org/10.1080/09720502.2020.1792665},
  publisher    = {Taylor \& Francis},
  note         = {Published online 30 Aug 2020},
  affiliation  = {Department of Mathematics, Rivers State University, Port Harcourt, Nigeria; Sanata Dharma University}
}

@article{ChenEtAl2024,
  author       = {Zhen Chen and Anne Gelb and Yoonsang Lee},
  title        = {Learning the Dynamics for Unknown Hyperbolic Conservation Laws Using Deep Neural Networks},
  journal      = {SIAM Journal on Scientific Computing},
  volume       = {46},
  number       = {2},
  pages        = {A825--A850},
  year         = {2024},
  doi          = {10.1137/22M1537333},
  url          = {https://doi.org/10.1137/22M1537333},
  publisher    = {Society for Industrial and Applied Mathematics},
  msc2020      = {65M12, 35L65, 62M45}
}

@article{amat2025convergence,
  title={On the convergence of WENO schemes for scalar hyperbolic conservation laws},
  author={Amat, Sergio and Busquier, Sonia and Ruiz, Juan},
  journal={The Journal of Analysis},
  pages={1--20},
  year={2025},
  publisher={Springer}
}

@article{harten1987uniformly,
  title     = {Uniformly High Order Accurate Essentially Non‑oscillatory Schemes, {III}},
  author    = {Harten, Ami and Engquist, Bj{\"o}rn and Osher, Stanley and Chakravarthy, Sukumar R.},
  journal   = {Journal of Computational Physics},
  volume    = {71},
  number    = {2},
  pages     = {231--303},
  year      = {1987},
  doi       = {10.1016/0021-9991(87)90031-3},
  note      = {Dedicated to Peter Lax},
}

@techreport{shu1997essentially,
  author       = {Shu, Chi‑Wang},
  title        = {Essentially Non‑Oscillatory and Weighted Essentially Non‑Oscillatory Schemes for Hyperbolic Conservation Laws},
  institution  = {Institute for Computer Applications in Science and Engineering, NASA Langley Research Center},
  type         = {NASA/CR-97-206253; ICASE Report No. 97-65},
  address      = {Hampton, VA},
  number       = {ICASE-97-65},
  month        = nov,
  year         = {1997},
  note         = {Prepared under NASA Contract NAS1-19480; Langley Technical Monitor: Dennis M. Bushnell},
  doi          = {10.1007/BFb0096355},
  url          = {https://ntrs.nasa.gov/citations/19980007543}
}

@article{Yarotsky16,
  author       = {Dmitry Yarotsky},
  title        = {Error bounds for approximations with deep ReLU networks},
  journal      = {CoRR},
  volume       = {abs/1610.01145},
  year         = {2016},
  url          = {http://arxiv.org/abs/1610.01145},
  eprinttype    = {arXiv},
  eprint       = {1610.01145},
  bibsource    = {dblp computer science bibliography, https://dblp.org}
}

@article{hornik1991approximation,
  author  = {Hornik, Kurt},
  title   = {Approximation Capabilities of Multilayer Feedforward Networks},
  journal = {Neural Networks},
  volume  = {4},
  number  = {2},
  pages   = {251--257},
  year    = {1991},
  month   = mar,
  doi     = {10.1016/0893-6080(91)90009-T},
}

@misc{kovachki2021universalapproximationerrorbounds,
      title={On universal approximation and error bounds for Fourier Neural Operators}, 
      author={Nikola Kovachki and Samuel Lanthaler and Siddhartha Mishra},
      year={2021},
      eprint={2107.07562},
      archivePrefix={arXiv},
      primaryClass={math.NA},
      url={https://arxiv.org/abs/2107.07562}, 
}

@article{kovachki2021neural,
  title={Neural Operator: Learning Maps Between Function Spaces With Applications to PDEs},
  author={Kovachki, Nikola and Li, Zongyi and Liu, Burigede and Azizzadenesheli, Kamyar and Bhattacharya, Kaushik and Stuart, Andrew and Anandkumar, Anima},
  journal={arXiv preprint arXiv:2108.08481},
  year={2021}
}

@book{muller2007dtw,
  title={Information Retrieval for Music and Motion},
  author={M{\"u}ller, Meinard},
  year={2007},
  publisher={Springer},
  doi={10.1007/978-3-540-74048-3},
}

@article{zhang2011maximum,
  title={Maximum-principle-satisfying and positivity-preserving high-order schemes for conservation laws: survey and new developments},
  author={Zhang, Xiangxiong and Shu, Chi-Wang},
  journal={Proceedings of the Royal Society A: Mathematical, Physical and Engineering Sciences},
  volume={467},
  number={2134},
  pages={2752--2776},
  year={2011},
  publisher={The Royal Society Publishing}
}
